\documentclass[12pt,titlepage,a4paper,twoside]{report}
\usepackage[american]{babel}
\usepackage{moreverb}
\usepackage{a4}
\usepackage{graphicx}
\usepackage{latexsym}
\usepackage{paralist}
\usepackage[compact]{titlesec}
\usepackage[plainpages=false]{hyperref}
\usepackage{amsmath, amssymb, amsthm}
\usepackage{stmaryrd}
\usepackage{dsfont}
\usepackage{fancyhdr}
\usepackage{tikz}
\usetikzlibrary{shapes}

\fancypagestyle{plain}{
	\fancyhf{}
	\fancyhead[EL]{\thepage}
	\fancyhead[OR]{\thepage}
	
}

\theoremstyle{plain}
\newtheorem{theorem}{Theorem}
\newtheorem{proposition}{Proposition}[chapter]
\newtheorem{lemma}[proposition]{Lemma}
\newtheorem{corollary}[proposition]{Corollary}

\theoremstyle{definition}

\newtheorem*{definition}{Definition}
\newtheorem*{example}{Example}
\newtheorem*{remark}{Remark}

\begin{document}
\pagenumbering{alph}
\pagestyle{empty}
\begin{titlepage}
\begin{center}
{\huge \bfseries Isometry groups of Lorentzian manifolds of finite volume}\\
\vspace{0.1 cm}
{\Large \bfseries and}\\
\vspace{0.2 cm}
{\huge \bfseries The local geometry of compact homogeneous Lorentz spaces$^1$}\\
\vspace{0.6 cm}
\Large{\textsc{Felix G\"unther$^2$}}\\
\vspace{2.0 cm}

\small{
{\bfseries Abstract}\\
\vspace{0.2 cm}}

\end{center}
\small{

\noindent Based on the work of Adams and Stuck as well as on the work of Zeghib, we classify the Lie groups which can act isometrically and locally effectively on Lorentzian manifolds of finite volume. In the case that the corresponding Lie algebra contains a direct summand isomorphic to the two-dimensional special linear algebra or to a twisted Heisenberg algebra, we also describe the geometric structure of the manifolds if they are compact.

\noindent Using these results, we investigate the local geometry of compact homogeneous Lorentz spaces whose isometry groups have non-compact connected components. It turns out that they all are reductive. We investigate the isotropy representation, curvatures and holonomy. Especially, we obtain that any Ricci-flat compact homogeneous Lorentz space is flat or has compact isometry group.\\
\vspace{0.2 cm}

\noindent Keywords: Lorentz geometry, isometry groups, twisted Heisenberg group, oscillator group, homogeneous Lorentz spaces, reductive, Ricci curvature
}

\footnotetext[1]{This is a corrected version of my diploma thesis which Humboldt-Universit\"at zu Berlin received in May 2011.}
\footnotetext[2]{Institut f\"ur Mathematik, Humboldt-Universit\"at zu Berlin, Unter den Linden 6, 10099 Berlin, Germany.
E-mail: guenthef@math.hu-berlin.de}

\end{titlepage}
\newpage
\mbox{}
\newpage
%%%%%%%%%%%%%%%%%%%%%%%%%%%%%%%%%%%%%%%%%%%%%%%%%%%%%%%%%%%%%%%%%%%%%%%%%%%%%%%%%%%%%%%%%%%%
\setcounter{tocdepth}{2}
\clearpage\pagenumbering{roman}
\pagestyle{plain}
\fancyhead[ER]{Contents}
\fancyhead[OL]{Contents}
\tableofcontents
\newpage

\allowdisplaybreaks[4]
\raggedbottom
\setlength{\parindent}{0pt}
\setlength{\parskip}{1ex}

%%%%%%%%%%%%%%%%%%%%%%%%%%%%%%%%%%%%%%%%%%%%%%%%%%%%%%%%%%%%%%%%%%%%%%%%%%%%%%%%%%%%%%%%%%%%

\pagestyle{plain}
\fancyhead[ER]{Introduction}
\fancyhead[OL]{Introduction}
\chapter*{Introduction}
\addcontentsline{toc}{chapter}{Introduction}

The aim of our work is a classification of Lie groups acting isometrically and locally effectively on Lorentzian manifolds of finite volume as well as providing a detailed investigation of compact homogeneous Lorentzian manifolds whose isometry groups have non-compact connected components.

Symmetries play an important role in geometry and other disciplines. In 1872, Klein proposed in his Erlangen program a way how to characterize the classical geometries (for example, Euclidean geometry and Hyperbolic geometry) by using the underlying group of symmetries. In addition, symmetries are fundamental in theoretical physics, where they are related to conserved quantities. Symmetries are also relevant in architecture and biology.

In semi-Riemannian geometry, symmetries correspond to isometries. The study of isometry groups is of great interest, regarding the structure of these groups as well as the way they act on the manifolds. Semi-Riemannian manifolds may allow only one isometry, but they are easier to investigate if the isometry groups are large, in the sense that they are acting transitively on the manifolds. In this case, we speak about homogeneous manifolds. A homogeneous manifold is very symmetric, that is, the manifold looks the same regardless from which point on the manifold it is viewed. Classical examples for homogeneous spaces are the Euclidean and Hyperbolic space.

The Riemannian case is well understood. The isometry group of each compact Riemannian manifold is itself compact, and conversely, any compact Lie group is acting isometrically and effectively on a compact Riemannian manifold. In the Lorentzian case, it turns out that this is not true anymore: There are compact Lorentzian manifolds whose isometry groups are non-compact.

D'Ambra showed that the isometry group of a simply-connected compact real analytic Lorentzian manifold is compact (\cite{DA88}). It is not yet clear whether the analyticity is necessary. Later, Adams and Stuck (\cite{AdSt97a}) as well as Zeghib (\cite{Ze98a}) independently provided an algebraic classification of Lie groups that act isometrically and (locally) effectively on Lorentzian manifolds that are compact. More generally, Zeghib showed the same result for manifolds of finite volume. We will mainly follow the approach of Zeghib to prove Theorem~\ref{th:algebraic_classification}, where we recover the results of \cite{AdSt97a} and \cite{Ze98a}. Also, we use ideas from \cite{AdSt97a} and own ideas.

Theorem~\ref{th:algebraic_classification} states that if a Lie group $G$ is acting isometrically and locally effectively on a Lorentzian manifold $M=(M,g)$ of finite volume, then there exist Lie algebras $\mathfrak{k}$, $\mathfrak{a}$ and $\mathfrak{s}$, such that the Lie algebra $\mathfrak{g}$ of $G$ is equal to the direct sum $\mathfrak{k}\oplus\mathfrak{a}\oplus\mathfrak{s}$. Here $\mathfrak{k}$ is compact semisimple, $\mathfrak{a}$ is abelian and $\mathfrak{s}$ is isomorphic to one of the following:
\begin{itemize}
\item the trivial algebra,
\item the two-dimensional affine algebra $\mathfrak{aff}(\mathds{R})$,
\item the $(2d+1)$-dimensional Heisenberg algebra $\mathfrak{he}_d$,
\item a certain $(2d+2)$-dimensional twisted Heisenberg algebra $\mathfrak{he}_d^\lambda$ ($\lambda \in \mathds{Z}_+^d$),
\item the two-dimensional special linear algebra $\mathfrak{sl}_2(\mathds{R})$.
\end{itemize}
Furthermore, in the latter two cases, if $G$ is contained in the isometry group of the manifold, the subgroup generated by $\mathfrak{s}$ has compact center if and only if the subgroup is closed in the isometry group.

Zeghib's approach has the advantage, that in contrast to \cite{AdSt97a}, the classification is more descriptive and also works in the case of manifolds of finite volume. In Chapter~\ref{ch:basics}, we introduce the Lie algebras of Lie groups occurring in Theorem~\ref{th:algebraic_classification} and give examples of isometric actions. Moreover, we introduce a certain symmetric bilinear form $\kappa$ on $\mathfrak{isom}(M)$, the Lie algebra of  the isometry group of $(M,g)$. For $X,Y \in \mathfrak{isom}(M)$, \[\kappa(X,Y):=\int\limits_U{g(\widetilde{X},\widetilde{Y})(x) d\mu(x)}.\] Here $\widetilde{X},\widetilde{Y}$ are complete Killing vector fields corresponding to $X,Y$ and $U$ is a $\text{Isom}(M)$-invariant non-empty open subset of $M$, such that $|g(\widetilde{X},\widetilde{Y})|$ is bounded on $U$ by a constant depending only on $X$ and $Y$. We give an example in Section~\ref{sec:induced_form}, that the restriction to an open subset of the manifold in the definition of $\kappa$ is necessary. We will give a proof of the important Proposition~\ref{prop:invariant_open_set}, which is also due to Zeghib, but was only indicated in \cite{Ze98a}. It shows that we can find such an open set. Thus, the symmetric bilinear form $\kappa$ exists and is well-defined. Note that $\kappa$ induces in a canonical way a symmetric bilinear form on the Lie algebra $\mathfrak{g}$ of a Lie group $G$ acting isometrically and locally effectively on $M$.

$\kappa$ is ad-invariant and fulfills the so called condition~\hyperlink{star}{$(\star)$}. A symmetric bilinear form $b$ on the Lie algebra $\mathfrak{g}$ of a Lie group $G$ fulfills this condition, if for any subspace $V$ of $\mathfrak{g}$, such that the set of $X \in V$ generating a non-precompact one-parameter group in $G$ is dense in $V$, the restriction of $b$ to $V \times V$ is positive semidefinite and its kernel has dimension at most one.

In the case of $\kappa$, this property is finally shown in Corollary~\ref{cor:condition_star} using dynamical methods like the Poincar\'e recurrence theorem or the F\"urstenberg lemma.

Condition~\hyperlink{star}{$(\star)$} is the main tool for proving Theorem~\ref{th:algebraic}. Theorem~\ref{th:algebraic} describes the algebraic structure of connected non-compact Lie groups $G$, whose Lie algebras $\mathfrak{g}$ possess an ad-invariant symmetric bilinear form $\kappa$ fulfilling condition~\hyperlink{star}{$(\star)$}. The structure of the Lie algebras is exactly the same as the one in Theorem~\ref{th:algebraic_classification}, as well as the statement about the center of the subgroup generated by $\mathfrak{s}$. We correct the claim formulated in \cite{Ze98a}, that the subgroup generated by $\mathfrak{s}$ has compact center. Counterexamples of this are given in Propositions~\ref{prop:counterexample1} and~\ref{prop:counterexample2}.

Theorem~\ref{th:algebraic_classification} follows in the case that $G$ is a closed subgroup of the isometry group of $M$ from Theorem~\ref{th:algebraic}. For non-closed subgroups, we encounter the problem that one-parameter groups may be non-precompact in $G$, however precompact in the isometry group. With the consideration of these groups in Section~\ref{sec:general_subgroups}, we fill a gap in \cite{Ze98a}.

It turns out that $\kappa$ is either Lorentzian (which leads to the cases that $\mathfrak{s}$ is isomorphic to the two-dimensional special linear algebra or a twisted Heisenberg algebra), positive definite (which yields the case that $\mathfrak{s}$ is trivial) or positive semidefinite with one-dimensional kernel. In the latter case, if $\mathfrak{s}$ is not trivial, we obtain that $\mathfrak{s}$ is isomorphic to either a Heisenberg algebra or the two-dimensional affine algebra, depending on whether the radical of $\mathfrak{g}$ is nilpotent or not.

Our main results are stated in Chapter~\ref{ch:theorems}. Although we do not go into the details of the proof of Theorem~\ref{th:isometry_groups}, we will mention it, because it gives a characterization of Lie groups which can be the entire connected component of the isometry group of a compact Lorentzian manifold. On the level of Lie algebras, all algebras of Theorem~\ref{th:algebraic_classification} but the two-dimensional affine algebra appear. The missing steps of the proof can be found in \cite{AdSt97b} and \cite{Ze98b}.

Chapter~\ref{ch:algebra} is devoted to the proof of Theorem~\ref{th:algebraic}, which is the main ingredient for the proof of Theorem~\ref{th:algebraic_classification}. The algebraic Lemma~\ref{lem:precompact_ad} (that describes elements of a Lie algebra generating a precompact one-parameter group in terms of the adjoint action) is quite powerful and allows us to obtain several results in a different way than it was done in \cite{Ze98a}. Examples are Lemma~\ref{lem:nilpotent_density} stating that almost all one-parameter groups in a non-compact nilpotent group are not precompact, and  Propositions~\ref{prop:kappa_s_Lorentz} and~\ref{prop:sl2R}, which together show, that in the case that the radical of $G$ is compact, $\mathfrak{s}\cong\mathfrak{sl}_2(\mathds{R})$.

Since the proof of Theorem~\ref{th:algebraic_classification} in \cite{Ze98a} relies on Theorem~\ref{th:algebraic}, the statements about the subgroup generated by $\mathfrak{s}$ in Theorems~\ref{th:algebraic_classification} and~\ref{th:isometry_groups} are also different here. Because our counterexamples do not apply for these theorems, it would be interesting to see whether the original statements are still true or whether one can also find counterexamples for them.

Theorem~\ref{th:locally_free}, which was only stated and shown in its general form in \cite{AdSt97a}, says that if the manifold $M$ is compact, the subgroup generated by $\mathfrak{s}$ in Theorem~\ref{th:algebraic_classification} acts locally freely on $M$. This result is important for the geometric characterization of compact Lorentzian manifolds in Theorem~\ref{th:geometric_characterization} (which is due to Zeghib). This theorem considers the case that the Lorentzian manifold $M$ is compact and $\mathfrak{s}$ is isomorphic to the two-dimensional special linear algebra or to a twisted Heisenberg algebra. If $\mathfrak{s}\cong\mathfrak{sl}_2(\mathds{R})$, $M$ is covered isometrically by a warped product of the universal cover of the two-dimensional special linear group and a Riemannian manifold $N$. Else, if $\mathfrak{s}\cong\mathfrak{he}_d^\lambda$, $M$ is covered isometrically by a twisted product $S\times_{Z(S)}N$ of a twisted Heisenberg group $S$ and a Riemannian manifold $N$. Similarly to Theorem~\ref{th:algebraic}, we state and prove Theorem~\ref{th:geometric_characterization} in an adapted version that takes the case that the center of $S$ might be non-compact into account. Furthermore, we show in Proposition~\ref{prop:lorentz_heisenberg}~(ii) that the invariance of a Lorentzian scalar product on a twisted Heisenberg algebra under the adjoint action of the nilradical (what Zeghib called \textit{essential ad-invariance}) is equivalent to ad-invariance, in contrast to the result of \cite{Ze98a}.

The proof of the last two theorems is done in Chapter~\ref{ch:geometry}. For the proof of Theorem~\ref{th:locally_free}, we follow the algebraic approach of \cite{AdSt97a} and differ only in details. The proof of Theorem~\ref{th:geometric_characterization} relies on the work of Zeghib. In our work, we give a more detailed description and correct some steps of the proof. For example, the horizontal space in Proposition~\ref{prop:submersion}~(ii) (showing that $S\times N \to M$ is a submersion) is different than in \cite{Ze98a}, as well as the proof of Proposition~\ref{prop:m_invariant} (showing that the corresponding metrics on $S$ are bi-invariant).

The idea behind the proof of Theorem~\ref{th:geometric_characterization} is considering the orbits $\mathcal{S}$ of the action of $S$, the subgroup in $\text{Isom}(M)$ generated by $\mathfrak{s}$. The orbits have the same dimension as $S$, because $S$ acts locally freely by Theorem~\ref{th:locally_free}. It turns out that if $\mathfrak{s}\cong\mathfrak{sl}_2(\mathds{R})$ or $\mathfrak{s}\cong\mathfrak{he}_d^\lambda$, the orbits have Lorentzian character everywhere on $M$. This allows us to consider the distribution $\mathcal{O}$ orthogonal to $\mathcal{S}$. In the case that $\mathfrak{s}\cong\mathfrak{sl}_2(\mathds{R})$, $\mathcal{O}$ is involutive. If $\mathfrak{s}\cong\mathfrak{he}_d^\lambda$, $\mathcal{O}+\mathcal{Z}$ is involutive, where $\mathcal{Z}$ denotes the distribution defined by the orbits of the center of $S$. In both cases, the involutive distribution defines a foliation by the Frobenius theorem. We choose $N$ to be a leaf of this foliation and furnish it with the metric induced by the metric of the ambient space $(M,g)$ and an arbitrary Riemannian metric on $\mathcal{Z}$ orthogonal to $\mathcal{O}$.

Regarding homogeneous Riemannian manifolds, it is well known that they are geodesically complete. Marsden showed that the same is true for compact homogeneous semi-Riemannian manifolds (\cite{Ma73}). Note that there exists a homogeneous Lorentzian space which is not geodesically complete (cf.~\cite{ON83}, remark~9.37). Additionally, each homogeneous Riemannian manifold that is Ricci-flat is also flat, but certain Cahen-Wallach spaces, which all are symmetric Lorentzian spaces, are Ricci-flat and non-flat. It is not known whether at least all Ricci-flat compact homogeneous Lorentzian manifolds are flat.

In Chapter~\ref{ch:homogeneous}, we provide a detailed analysis of compact homogeneous Lorentz spaces $M$, especially of those whose isometry groups have non-compact connected components. We start by presenting a topological and geometric description of them in Theorem~\ref{th:homogeneous_characterization}, which is also slightly different than the corresponding result of Zeghib. Essentially, it shows that if the isometry group of $M$ has non-compact connected components, $M$ is covered isometrically either by a metric product of the universal cover of the two-dimensional special linear group and a compact homogeneous Riemannian manifold $N$, or by a twisted product $S\times_{Z(S)}N$ of a twisted Heisenberg group $S$ and a compact homogeneous Riemannian manifold $N$. $N$ is constructed in the same way as in the proof of Theorem~\ref{th:geometric_characterization}.

Theorem~\ref{th:homogeneous_reductive} shows what was originally stated in \cite{Ze98a}, namely that the connected components of the isotropy group are compact. We give an elegant proof using the ideas of Adams and Stuck in the proof of Theorem~\ref{th:locally_free}. Moreover, it turns out that every compact homogeneous Lorentzian manifold has a reductive representation defined in a natural way. For this representation, the induced bilinear form $\kappa$ plays an essential role.

Homogeneous pseudo-Riemannian manifolds of constant curvature were studied extensively in \cite{Wo61}. Using a slightly different reductive representation than in Theorem~\ref{th:homogeneous_reductive}, we are able to describe in Section~\ref{sec:homogeneous_geometry} in detail the local geometry of compact homogeneous Lorentzian manifolds whose isometry groups have non-compact connected components, in terms of the curvature and holonomy of the manifold. We also investigate the isotropy representation of the manifold, our result concerning a decomposition into (weakly) irreducible summands being summarized in Theorem~\ref{th:isotropy_representation}.

Our results of Section~\ref{sec:homogeneous_geometry} directly yield the proof of Theorem~\ref{th:homogeneous_not_Ricci_flat} which states that the isometry group of any Ricci-flat compact homogeneous Lorentzian manifold has compact connected components. Together with two results in \cite{PiZe10} and \cite{RoSa96}, it follows that the isometry group of any Ricci-flat compact homogeneous Lorentzian manifold that is non-flat, is in fact compact.

In the appendix, we give a short presentation of the complete Jordan decomposition. The main result is important for our proof of Proposition~\ref{prop:sl2R}.

\section*{Acknowledgment}

I am very grateful to my supervisor Helga Baum for her dedicated guidance, useful comments and strong support. She gave me the opportunity to work on such a fascinating and interesting topic, the treatment of which combines several fields of mathematics. Furthermore, in her demanding and captivating lectures she provided me with the broad knowledge of differential geometry I needed during my work and I will benefit from in the future.

I especially thank Abdelghani Zeghib for his very helpful answers to my questions.

\newpage

%%%%%%%%%%%%%%%%%%%%%%%%%%%%%%%%%%%%%%%%%%%%%%%%%%%%%%%%%%%%%%%%%%%%%%%%%%%%%%%%%%%%%%%%%%%%

\pagestyle{plain}
\renewcommand{\chaptermark}[1]{\markboth{\chaptername\ \thechapter.\ #1}{}}
\renewcommand{\sectionmark}[1]{\markright{\thesection.\ #1}}
\fancyhead[ER]{\leftmark}
\fancyhead[OL]{\rightmark}

\setcounter{secnumdepth}{3}

\clearpage\pagenumbering{arabic}

\chapter{Lie groups acting isometrically on Lorentzian manifolds}\label{ch:basics}

This chapter is devoted to introduce notations and prove elementary theorems which we will use for the classification of Lie algebras of Lie groups acting isometrically and locally effectively on a Lorentzian manifold of finite volume. In Section~\ref{sec:basics}, we describe isometric actions of Lie groups and give several Lie algebraic results focussing especially on compactness properties. We continue in Section~\ref{sec:examples} by describing all the Lie algebras which appear as a direct summand in the Lie algebras we classify as well as giving examples of isometric actions of the corresponding Lie groups. Finally, we define in Section~\ref{sec:induced_form} a certain symmetric bilinear form $\kappa$ on Lie algebras, which will be essential for our proof of Theorem~\ref{th:algebraic_classification} that is stated in Chapter~\ref{ch:theorems}.

\section{Definitions and basic properties}\label{sec:basics}

Let $M=(M,g)$ be a semi-Riemannian manifold of dimension $n$. We will consider only connected smooth real manifolds without boundary. The metric $g$ defines canonically a Lebesgue measure $\mu$ on $M$. Then the volume of $M$ is defined as \[\textnormal{vol}(M):=\mu(M)=\int\limits_M 1 d\mu(x).\] The assertion of finite volume allows us to integrate bounded functions, which will be crucial for defining the ad-invariant symmetric bilinear form $\kappa$ on the Lie algebra of the isometry group $\textnormal{Isom}(M)$. We explain this in more detail in Section~\ref{sec:induced_form}.

All Lie groups and Lie algebras will be real and finite-dimensional. Unless otherwise stated, all actions of Lie groups are smooth.

\begin{definition}
Let $G$ be a connected Lie group, which acts isometrically on a semi-Riemannian manifold $M$. In this case, there exists a Lie group homomorphism $\rho : G \to \textnormal{Isom}(M)$. The group action will be denoted by $\cdot$. The action is \textit{locally effective}, if the kernel of $\rho$ is discrete in $G$. If the kernel is trivial, the action of $G$ on $M$ is \textit{effective}. In this case, we can consider $G$ as a subgroup of $\textnormal{Isom}(M)$ and $\mathfrak{g}$ as a subalgebra of $\mathfrak{isom}(M)$, the Lie algebra of the isometry group.

The action is \textit{locally free}, if there is a neighborhood of the identity element in $G$ acting without fixed points on $M$. If all elements of $G$ but the identity are acting without fixed points on $M$, the action is \textit{free}.
\end{definition}

\begin{remark}
We will usually suppose that the action of $G$ is locally effective. Note that in the situation above, the kernel of $\rho$ is then a discrete normal subgroup and therefore central in $G$. It follows that $\rho:G \to \rho(G)$ is a covering map, so when investigating the Lie algebra $\mathfrak{g}$, we may assume the action to be effective. In either case, we identify a group element $f \in G$ with the isometry $\rho(f)$, which allows us to speak about the differential $df:=d\rho(f)$.
\end{remark}

The following proposition is a classical result which allows us to identify $\mathfrak{isom}(M)$ with $\mathfrak{kill}_c(M)$, the Lie algebra of complete Killing vector fields, as vector spaces. One can find a proof in \cite{ON83}, Proposition~9.33.

\begin{proposition}\label{prop:isometry_Killing}
Let $M$ be a semi-Riemannian manifold. Then the mapping $\mathfrak{isom}(M) \to \mathfrak{kill}_c(M)$ defined by $X \mapsto \widetilde{X}$, $\widetilde{X}(x):=\frac{\partial}{\partial t}(\exp(tX) \cdot x) |_{t=0}$, is an anti-isomorphism of Lie algebras, that is, $[\widetilde{X},\widetilde{Y}]=-\widetilde{[X,Y]}$.
\end{proposition}

Let $O(M)$ denote the bundle of orthonormal frames in $M$. It is well known that for any $x \in M$, the mapping $\phi_x: \textnormal{Isom}(M) \to O(M)$, which is defined by $f \mapsto \left(df_{x}(s_1), \ldots, df_{x}(s_n)\right)$ for a fixed orthonormal basis $\left\{s_1, \ldots, s_n\right\}$ in the tangent space $(T_xM,g_x)$, is an embedding of the isometry group into the bundle of orthonormal frames. If $M$ is a compact Riemannian manifold, $O(M)$ is compact (since the fiber is the orthogonal group $O(n)$, which is compact) and therefore, $\textnormal{Isom}(M)$ is a compact Lie group. Conversely, any connected compact Lie group $G$ acts on a compact Riemannian manifold isometrically and effectively. We can simply take $M=G$ as a manifold and provide it with any left-invariant Riemannian metric.

For the following, one can find a proof in \cite{Se92}, Part~I, Chapter~VI, Paragraph~4, Corollary~1.

\begin{lemma}\label{lem:Levi_decomposition}
Let $\mathfrak{g}$ be a Lie algebra and $\mathfrak{r}$ its radical, that is, the largest solvable ideal with respect to inclusion. Then there exists a subalgebra $\mathfrak{l}$, which is trivial or semisimple, such that $\mathfrak{g}=\mathfrak{l}\inplus\mathfrak{r}$ is a semidirect sum. This decomposition is called \textit{Levi decomposition}.
\end{lemma}

The first part of the following proposition is due to \cite{Got49}, Lemma~4. A proof of the second part one can find in \cite{Vir93}, Proposition~4.2. Note that the second part implies the first one.

\begin{proposition}\label{prop:compact_center_closed_subgroup}
The following is true:
\begin{compactenum}
\item Let $G$ be a Lie group and $H$ a connected semisimple Lie subgroup. If $H$ has finite center, $H$ is closed in $G$.

\item Let $G$ be a Lie group and $H$ a connected semisimple or nilpotent Lie subgroup. If $H$ is not closed in $G$, the center of $H$ is not compact.
\end{compactenum}
\end{proposition}

\begin{definition}
A Lie algebra $\mathfrak{g}$ is \textit{reductive}, if its radical is equal to its center $\mathfrak{z(g)}$ (cf.~\cite{OnVin94}, Chapter~1, Paragraph~3.5).
\end{definition}

\begin{proposition}\label{prop:reductive_algebra}
A Lie algebra $\mathfrak{g}$ is reductive if and only if $\mathfrak{g}=[\mathfrak{g},\mathfrak{g}]\oplus\mathfrak{z(g)}$ is a direct sum with $[\mathfrak{g},\mathfrak{g}]$ trivial or semisimple.
\end{proposition}
\vspace{-1.4em}\begin{proof}
Assume that $\mathfrak{g}$ is reductive. By Lemma~\ref{lem:Levi_decomposition}, the Levi decomposition yields the direct decomposition $\mathfrak{g}=\mathfrak{l}\oplus\mathfrak{z(g)}$ with $\mathfrak{l}$ trivial or semisimple. Clearly, $[\mathfrak{g},\mathfrak{g}]=[\mathfrak{l},\mathfrak{l}]=\mathfrak{l}$.

Now let $\mathfrak{g}=[\mathfrak{g},\mathfrak{g}]\oplus\mathfrak{z(g)}$ be direct. The radical intersects $[\mathfrak{g},\mathfrak{g}]$ in a solvable ideal. But $[\mathfrak{g},\mathfrak{g}]$ is trivial or semisimple, so this intersection is trivial. Since the radical is the largest solvable ideal and the center of $\mathfrak{g}$ is solvable, it follows that the radical is equal to the center.
\end{proof}\vspace{0pt}

\begin{definition}
A Lie algebra $\mathfrak{g}$ is \textit{compact}, if it is isomorphic to the Lie algebra of a compact Lie group $G$.
\end{definition}

The following proposition is a standard Lie algebraic result. One can find a proof in \cite{Bo05}, Chapter~IX, Paragraph~1.3, Proposition~1 (for various characterizations of reductive algebras, see \cite{Bo89}, Chapter~I, Paragraph~6.4, Proposition~5).

\begin{proposition}\label{prop:compact_algebra}
For a Lie algebra $\mathfrak{g}$, the following are equivalent:
\begin{compactenum}
\item $\mathfrak{g}$ is compact.

\item $\mathfrak{g}$ possesses a positive definite symmetric bilinear form $b$ which is ad-invari\-ant, that is, $b([X,Y],Z)=b(X,[Y,Z])$ for all $X,Y,Z \in \mathfrak{g}$.

\item $\mathfrak{g}$ is reductive and the Killing form of $\mathfrak{g}$ is negative semidefinite.
\end{compactenum}
\end{proposition}

\begin{corollary}\label{cor:compact_semisimple}
A Lie algebra $\mathfrak{g}$ is compact semisimple if and only if its Killing form is negative definite.
\end{corollary}
\vspace{-1.4em}\begin{proof}
It is well-known that $\mathfrak{g}$ is semisimple if and only if its Killing form is non-degenerate. Also, a semisimple algebra is reductive. The claim now follows from Proposition~\ref{prop:compact_algebra}.
\end{proof}\vspace{0pt}

\begin{corollary}\label{cor:subalgebra_compact}
Subalgebras of compact algebras are itself compact.
\end{corollary}
\vspace{-1.4em}\begin{proof}
This follows directly from Proposition~\ref{prop:compact_algebra} using (ii).
\end{proof}\vspace{0pt}

\begin{corollary}\label{cor:compact_algebra}
If $\mathfrak{g}$ is compact, $\mathfrak{g}=[\mathfrak{g},\mathfrak{g}]\oplus\mathfrak{z(g)}$ is a direct sum with $[\mathfrak{g},\mathfrak{g}]$ compact semisimple.
\end{corollary}
\vspace{-1.4em}\begin{proof}
Use Proposition~\ref{prop:compact_algebra}~(iii) to see that $\mathfrak{g}$ is reductive and Proposition~\ref{prop:reductive_algebra} and Corollary~\ref{cor:subalgebra_compact} for the direct decomposition.
\end{proof}\vspace{0pt}

The following lemma one can find for example in \cite{On93}, Part~I, Theorem~2.12.

\begin{lemma}\label{lem:abelian_group}
If $G$ is a connected abelian Lie group, then it is isomorphic to the additive group $\mathds{T}^m \times \mathds{R}^{m^\prime}$ of a $m$-dimensional torus $\mathds{T}^m$ and a $m^\prime$-dimensional vector space $\mathds{R}^{m^\prime}$.
\end{lemma}

The next result is a classical one due to \cite{Iw49}, Theorem~4.
\begin{lemma}\label{lem:compact_abelian}
A compact normal abelian subgroup in a connected topological group is contained in its center.
\end{lemma}

\section{Examples}\label{sec:examples}

Theorem~\ref{th:algebraic_classification} that is stated in Chapter~\ref{ch:theorems} describes the structure of the Lie algebra of a Lie group acting isometrically and locally effectively on a Lorentzian manifold of finite volume. It is isomorphic to a direct sum $\mathfrak{k}\oplus\mathfrak{a}\oplus\mathfrak{s}$, where $\mathfrak{k}$ is compact semisimple, $\mathfrak{a}$ abelian and $\mathfrak{s}$ is either trivial, the two-dimensional affine algebra, the two-dimensional special linear algebra, a Heisenberg algebra or a twisted Heisenberg algebra. In Section~\ref{sec:compact_abelian}, we begin the description and investigation of these Lie algebras with the summand $\mathfrak{k}\oplus\mathfrak{a}$ and continue with the other algebras in the order we listed them in the last sentence.

\subsection{Product with a compact Riemannian manifold}\label{sec:compact_abelian}

Let $(M,g)$ be a Lorentzian manifold of finite volume and $(N,h)$ be a compact Riemannian manifold. Then $M \times N$ with the product metric $g \times h$ is a Lorentzian manifold of finite volume and $\textnormal{Isom}(M) \times \textnormal{Isom}(N) \subseteq \textnormal{Isom}(M \times N)$.

Thus, if $H$ is any connected compact Lie group and $G$ a Lie group acting isometrically and (locally) effectively on some Lorentzian manifold of finite volume, so $H \times G$ acts isometrically and (locally) effectively on a (in general different) Lorentzian manifold of finite volume. Due to Corollary~\ref{cor:compact_algebra}, the Lie algebra of $H \times G$ is given by the direct sum $\mathfrak{k}\oplus\mathfrak{a}\oplus\mathfrak{g}$, where $\mathfrak{g}$ is the Lie algebra of G, $\mathfrak{k}$ is compact semisimple and $\mathfrak{a}$ abelian.

\subsection{Two-dimensional affine algebra}\label{sec:affine}

The two-dimensional affine group $\textnormal{Aff}(\mathds{R})$ is the group of orientation preserving affine transformations of the real line. It is centerless, simply-connected and isomorphic to the identity component of the group of upper triangular matrices in the two-dimensional special linear group $\textnormal{SL}_2(\mathds{R})$. Thus, we can identify its Lie algebra $\mathfrak{aff}(\mathds{R})$ as the subalgebra of $\mathfrak{sl}_2(\mathds{R})$ spanned by \[X:= \begin{pmatrix} \frac{1}{2} & 0 \\ 0 & -\frac{1}{2}\end{pmatrix} \text{ and } Y:= \begin{pmatrix} 0 & 1 \\ 0 & 0\end{pmatrix}.\] $X$ corresponds to the infinitesimal generator of the scaling transformations and $Y$ to the infinitesimal generator of the translations. They satisfy the relation $[X,Y]=Y$. Note that the Killing form of $\mathfrak{aff}(\mathds{R})$ is positive semidefinite and the kernel is exactly the span of $Y$.

Considering $\textnormal{Aff}(\mathds{R})$ as a subgroup of $\textnormal{SL}_2(\mathds{R})$, an isometric and locally effective action on a compact Lorentzian manifold is given by restricting the action of $\textnormal{SL}_2(\mathds{R})$ given in the next paragraph.

\subsection{Special linear algebra}\label{sec:sl_2}

The two-dimensional special linear group $\textnormal{SL}_2(\mathds{R})$ is the matrix group consisting of all real $2 \times 2$-matrices of determinant 1.

We can consider its Lie algebra $\mathfrak{sl}_2(\mathds{R})$ as the subalgebra of the general linear algebra $\mathfrak{gl}_2(\mathds{R})$ that is spanned by \[e:= \begin{pmatrix} 0 & 1 \\ 0 & 0\end{pmatrix}, \ f:=\begin{pmatrix} 0 & 0 \\ 1 & 0\end{pmatrix} \text{ and } h:=\begin{pmatrix} 1 & 0 \\ 0 & -1\end{pmatrix}.\] They satisfy the relations $[e,f]=h$, $[h,e]=2e$ and $[h,f]=-2f$ and therefore form an $\mathfrak{sl}_2$-triple (cf.~\cite{OnVin94}, Chapter~6, Paragraph~2.1).

\begin{definition}
A triple $\left\{e,f,h\right\}$ of elements in a Lie algebra $\mathfrak{g}$ satisfying $[e,f]=h$, $[h,e]=2e$ and $[h,f]=-2f$, is called an \textit{$\mathfrak{sl}_2$-triple}.
\end{definition}

$\mathfrak{sl}_2(\mathds{R})$ is a simple Lie algebra. Its Killing form $k$ in the ordered basis $(e,f,h)$ is determined by the matrix \[\begin{pmatrix} 0 & 4 & 0 \\ 4 & 0 & 0 \\ 0 & 0 & 8\end{pmatrix}\] and defines a Lorentzian scalar product.

Since the Killing form is ad-invariant, $k$ defines a bi-invariant Lorentzian metric on the universal cover $\widetilde{\textnormal{SL}_2(\mathds{R})}$. Thus, it passes to a Lorentzian metric on the centerless quotient group $\textnormal{PSL}(2,\mathds{R})$, on which $\widetilde{\textnormal{SL}_2(\mathds{R})}$ and any central quotient group acts isometrically and locally effectively by left multiplication. The same is true for a quotient $\textnormal{PSL}(2,\mathds{R})/\Lambda$, where $\Lambda$ is a uniform lattice, that is, a cocompact discrete subgroup of $\textnormal{PSL}(2,\mathds{R})$. Such a lattice exists (cf.~\cite{EiWa11}, Lemma~11.12).

\subsection{Heisenberg algebra}\label{sec:Heisenberg}

The Heisenberg algebra $\mathfrak{he}_d$ of dimension $2d+1$, $d>0$, is spanned by the elements $Z, X_1, Y_1, \ldots, X_d, Y_d$. The non-vanishing Lie brackets are given by \[[X_k,Y_k]=Z,\] $k=1,\ldots,d$. A Lie algebra $\mathfrak{g}$ is isomorphic to a Heisenberg algebra if and only if its center $\mathfrak{z(g)}$ is one-dimensional (spanned by $Z$) and there is a complementary vector space $V \subset \mathfrak{g}$ and a non-degenerate alternating bilinear form $\omega$, such that $[X,Y]=\omega(X,Y)Z$ for all $X,Y \in V$. To see this, simply notice that there is a symplectic basis $\left\{X_1, Y_1, \ldots, X_d, Y_d\right\}$ of $V$, that is, \[\omega(X_k,Y_j)=\delta_{jk} \text{ and } \omega(X_k,X_j)=0=\omega(Y_k,Y_j)\] for all $j$ and $k$.

\begin{definition}
A set $\left\{Z, X_1, Y_1, \ldots, X_d, Y_d\right\}$ of elements in a Heisenberg algebra fulfilling the same relations as above, is called a \textit{canonical basis}. A \textit{Heisenberg group} is a Lie group with Lie algebra $\mathfrak{he}_d$.
\end{definition}

The Heisenberg algebra is two-step nilpotent. Let $\widetilde{\textnormal{He}_d}$ denote the simply-connect\-ed Lie group with Lie algebra $\mathfrak{he}_d$. As a matrix group, $\widetilde{\textnormal{He}_d}$ can be represented as the set of $(d+2)\times (d+2)$-matrices \[\begin{pmatrix}1 & x & z \\ 0 & I_{d} & y \\ 0 & 0 & 1\end{pmatrix},\] where $x$ and $y$ are real row and column vectors of length $d$, respectively, $z \in \mathds{R}$ and $I_{d}$ is the $d \times d$-identity matrix (cf.~\cite{GorWi86}, Definitions and Notation~2.1~(a)). The center is isomorphic to $\mathds{R}$. Let $\Lambda$ be a (uniform) lattice of the center (isomorphic to $\mathds{Z}$) and consider the quotient group $\widetilde{\textnormal{He}_d}/\Lambda$. This quotient is up to isomorphism independent of the choice of $\Lambda$.

\begin{definition} The quotient $\widetilde{\textnormal{He}_d}/\Lambda$ is denoted by $\textnormal{He}_d$.
\end{definition}

The Heisenberg groups $\textnormal{He}_d$ are subgroups of some twisted Heisenberg groups explained in the next paragraph. The universal cover $\widetilde{\textnormal{He}_d}$ is a subgroup of some twisted Heisenberg group as well. Therefore, an isometric and locally effective action on a compact Lorentzian manifold is given by the corresponding restriction of the action given in the next paragraph.

\subsection{Twisted Heisenberg algebras}\label{sec:twisted_Heisenberg}

Let $\lambda=(\lambda_1,\ldots,\lambda_d) \in \mathds{R}_+^d$, $d> 0$. The corresponding twisted Heisenberg algebra $\mathfrak{he}_d^\lambda$ of dimension $2d+2$ is spanned by the elements $T,Z, X_1, Y_1, \ldots, X_d, Y_d$ and the non-vanishing Lie brackets are given by \[[X_k,Y_k]=\lambda_k Z, \ [T,X_k]=\lambda_k Y_k \text{ and } [T,Y_k]=-\lambda_k X_k,\] $k=1,\ldots,d$. Thus, $\mathfrak{he}_d^\lambda=\mathds{R} T \inplus \mathfrak{he}_d$ is a semidirect sum, where $\mathfrak{he}_d$ can be identified with the subalgebra spanned by $Z, X_1, Y_1, \ldots, X_d, Y_d$.

\begin{definition}
A set $\left\{T, Z, X_1, Y_1, \ldots, X_d, Y_d\right\}$ of elements in a twisted Heisenberg algebra fulfilling the same relations as above, is called a \textit{canonical basis}. A \textit{twisted Heisenberg group} is a Lie group with Lie algebra $\mathfrak{he}_d^\lambda$.
\end{definition}

Now let $\widetilde{\textnormal{He}_d^\lambda}$ be the simply-connected Lie group with Lie algebra $\mathfrak{he}_d^\lambda$ and let $\exp: \mathfrak{he}_d^\lambda \to \widetilde{\textnormal{He}_d^\lambda}$ be the exponential. Then $\widetilde{\textnormal{He}_d^\lambda}=\exp(\mathds{R} T) \ltimes \widetilde{\textnormal{He}_d}$, where $\widetilde{\textnormal{He}_d}$ can be identified with $\exp(\mathfrak{he}_d)$.

For all $k$, the subspace $V_k:=\textnormal{span}\left\{X_k,Y_k\right\}$ is $\textnormal{ad}_T$-invariant and the action of $\textnormal{ad}_{tT}$ on $V_k$ defined by the ordered basis $(X_k,Y_k)$ is described by the matrix \[\begin{pmatrix} 0 & -t\lambda_k \\ t\lambda_k & 0\end{pmatrix},\] $t \in \mathds{R}$. Therefore, the action of $\textnormal{Ad}_{\exp(tT)}=\exp(\textnormal{ad}_{tT})$ on the space $V_k$ corresponds to the matrix \[\begin{pmatrix} \cos(t\lambda_k) & -\sin(t\lambda_k) \\ \sin(t\lambda_k) & \cos(t\lambda_k)\end{pmatrix}.\] Also, $\textnormal{ad}_T$ acts trivially on the center of $\mathfrak{he}_d$.

If $\lambda \in \mathds{Q}^d_+$, then $\Lambda^\prime:=\textnormal{ker}(\exp \circ \textnormal{ad}: \mathds{R}T \to \textnormal{Aut}(\mathfrak{he}_d^\lambda))$ is a (uniform) lattice in $\mathds{R}T$. As in \ref{sec:Heisenberg}, let $\Lambda$ be a lattice of the center of $\widetilde{\textnormal{He}_d}$.

\begin{definition}
Let $\lambda \in \mathds{Q}^d_+$. Then $\textnormal{He}_d^\lambda:=\widetilde{\textnormal{He}_d^\lambda}/(\Lambda^\prime \times \Lambda)$. In an analogous way, we define $\overline{\textnormal{He}_d^\lambda}:=\widetilde{\textnormal{He}_d^\lambda}/\Lambda^\prime $.
\end{definition}
\begin{remark}
We have $\textnormal{He}_d^\lambda\cong \mathds{S}^1 \ltimes \textnormal{He}_d$ and $\overline{\textnormal{He}_d^\lambda} \cong \mathds{S}^1 \ltimes \widetilde{\textnormal{He}_d}$.
\end{remark}

By construction, $\mathds{S}^1$ acts trivially on the center of $\textnormal{He}_d$ and as a rotation on the subgroups generated by $V_k$. We will see in the sequel that the universal cover $\widetilde{\textnormal{He}_d^\lambda}$ cannot be a closed subgroup of the isometry group of a Lorentzian manifold of finite volume, but $\textnormal{He}_d^\lambda$ can.

\begin{lemma}\label{lem:isom_heisenberg}
Let $\lambda,\eta \in \mathds{R}^d_+$. Then $\mathfrak{he}_d^\lambda \cong \mathfrak{he}_d^\eta$ if and only if there is an $a\in \mathds{R}_+$, such that $\left\{\lambda_1,\ldots,\lambda_d\right\}=\left\{a\eta_1,\ldots,a\eta_d\right\}$ as sets.
\end{lemma}
\vspace{-1.4em}\begin{proof}
We first prove the backward direction. Let $\left\{T, Z, X_1, Y_1, \ldots, X_d, Y_d\right\}$ and $\left\{T^\prime, Z^\prime, X_1^\prime, Y_1^\prime, \ldots, X_d^\prime, Y_d^\prime\right\}$ be canonical bases of $\mathfrak{he}_d^\lambda$ and $\mathfrak{he}_d^\eta$, respectively. $\sigma$ denotes the permutation of $\left\{1,\ldots,d\right\}$ such that $\lambda_k=a\eta_{\sigma(k)}$. Then the linear map $f:\mathfrak{he}_d^\lambda \to \mathfrak{he}_d^\eta$ defined by \[f(T)=a T^\prime, \ f(Z)=\frac{1}{a}Z^\prime, \ f(X_k)=X_{\sigma(k)}^\prime \text{ and }f(Y_k)=Y_{\sigma(k)}^\prime\] is an isomorphism of Lie algebras.

For the forward implication, let $\left\{T,Z,X_1,Y_1,\ldots,X_d,Y_d\right\}$ be a canonical basis of $\mathfrak{g}:=\mathfrak{he}_d^\lambda$. Then \[\mathfrak{h}:=[\mathfrak{g},\mathfrak{g}] =\textnormal{span}\left\{Z, X_1, Y_1, \ldots, X_d, Y_d\right\} \cong \mathfrak{he}_d \text{ and } \mathfrak{z(g)}=\mathds{R}Z=[\mathfrak{h},\mathfrak{h}]\] do not depend on the choice of the basis. The same is true for \[\mathfrak{a}:=\left\{X \in \mathfrak{g} | [X,\mathfrak{g}] \cap \mathfrak{z(g)}=\left\{0\right\}\right\}=\mathfrak{z(g)}+\mathds{R}T.\]

Now fix $T^\prime \in \mathfrak{a} \backslash \mathfrak{z(g)}$. Then $T^\prime=aT+bZ$ for real numbers $a,b$; $a\neq0$. Thus, $\textnormal{ad}_{T^\prime}=a \textnormal{ad}_T$. It is easy to see that $\textnormal{ad}_T$ is semisimple. Hence, there is an $\textnormal{ad}_T$-invariant vector space $V$ complementary to the $\textnormal{ad}_T$-invariant subspace $\mathfrak{a}$. It follows that $V=\textnormal{span}\left\{X_1, Y_1, \ldots, X_d, Y_d\right\}$.

$\textnormal{ad}_T$ is an automorphism of $V$ and has the eigenvalues $\lambda_k \textnormal{i}$ with respect to the invariant eigenspaces $V_k:=\textnormal{span}\left\{X_k,Y_k\right\}$, $k=1,\ldots,d$. Thus, $\textnormal{ad}_{T^\prime}$ is also an automorphism of $V$ and has the eigenvalues $a\lambda_k \textnormal{i}$ with respect to the invariant eigenspaces $V_k:=\textnormal{span}\left\{X_k,Y_k\right\}$, $k=1,\ldots,d$. So we have shown that the set $\left\{\lambda_1,\ldots,\lambda_d\right\}$ is determined up to multiplication with a positive real number $a$.
\end{proof}\vspace{0pt}

\begin{corollary}\label{cor:isom_heisenberg}
The isomorphism classes of $\widetilde{\textnormal{He}_d^\lambda}$, $\lambda \in \mathds{Q}^d_+$, are in one-to-one-correspondence with the set $\mathds{Z}_+^d/\sim$, where $\lambda \sim \eta$ if and only if there exists $a \in \mathds{\mathds{Q}}_+$, such that \[ \left\{\lambda_1,\ldots,\lambda_d\right\}=\left\{a\eta_1,\ldots,a\eta_d\right\}.\]
\end{corollary}

In the following, we will describe compact Lorentzian manifolds, on which $\widetilde{\textnormal{He}_d^\lambda}$ acts isometrically and locally effectively.

\begin{proposition}\label{prop:lorentz_heisenberg}
The following is true:
\begin{compactenum}
\item A twisted Heisenberg algebra $\mathfrak{he}_d^\lambda$ admits an ad-invariant Lorentz form.

\item Any ad-invariant Lorentzian scalar product on $\mathfrak{he}_d^\lambda$ is determined by two real parameters $\alpha, \beta$, where $\alpha>0$. Moreover, $\textnormal{ad}(\mathfrak{he}_d)$- and $\textnormal{ad}(\mathfrak{he}_d^\lambda)$-invariance are equivalent.

\item Conversely, if the Lie algebra of a semidirect product $\mathds{S}^1 \ltimes \textnormal{He}_d$ or $\mathds{S}^1 \ltimes \widetilde{\textnormal{He}_d}$, respectively, admits an ad-invariant Lorentz form, then $\mathds{S}^1 \ltimes \textnormal{He}_d$ or $\mathds{S}^1 \ltimes \widetilde{\textnormal{He}_d}$, respectively, is a twisted Heisenberg group.

\item Up to finite index of the lattices, there is a bijective correspondence between lattices in a twisted Heisenberg group $\mathds{S}^1 \ltimes \textnormal{He}_d$ and lattices in the subgroup $\textnormal{He}_d$, which are equivalent to lattices in $\widetilde{\textnormal{He}_d}$. Also, up to finite index of the lattices, there is a bijective correspondence between lattices in a twisted Heisenberg group $\mathds{S}^1 \ltimes \widetilde{\textnormal{He}_d}$ and lattices in the subgroup $\widetilde{\textnormal{He}_d}$.
\end{compactenum}
\end{proposition}
\vspace{-1.4em}\begin{proof}
(i) Let $\left\{T, Z, X_1, Y_1, \ldots, X_d, Y_d\right\}$ be a canonical basis of $\mathfrak{he}_d^\lambda$ and define $V:=\textnormal{span}\left\{X_1, Y_1, \ldots, X_d, Y_d\right\}$.

We define the symmetric bilinear form $\langle \cdot,\cdot \rangle$ by
\begin{align*}
\langle Z,Z \rangle&=0, \\
\langle X_j,X_k\rangle=\langle Y_j,Y_k\rangle&=\delta_{jk} \alpha,\\
\langle X_j,Y_k\rangle&=0,\\
\langle Z,X_k \rangle=\langle Z,Y_k \rangle =\langle T,X_k\rangle=\langle T,Y_k\rangle&=0.
\end{align*}
for all $j,k=1,\ldots,d$, and \[\langle T,Z \rangle=\alpha \text{ and } \langle T,T \rangle=\beta\] with two real parameters $\alpha, \beta$, where $\alpha>0$. Then $\langle \cdot,\cdot \rangle$ is a Lorentzian scalar product. We have to show, that
\begin{equation}\label{eq:ad-invariant}
\langle [W,X],Y \rangle+\langle X,[W,Y] \rangle =0\tag{ad}
\end{equation}
holds for all $W,X,Y \in \mathfrak{he}_d^\lambda$. By linearity, it suffices to check Equation~(\ref{eq:ad-invariant}) for the canonical basis elements.

Consider $W=T$. $\textnormal{ad}_T$ is trivial on $\textnormal{span}\left\{T,Z\right\}$ and by construction, $\textnormal{ad}_T$ restricted to $V:=\textnormal{span}\left\{X_1, Y_1, \ldots, X_d, Y_d\right\}$ is a skew-symmetric matrix with respect to a canonical basis. Thus, (\ref{eq:ad-invariant}) is satisfied.

If $W=Z$, then (\ref{eq:ad-invariant}) is satisfied since $Z$ lies in the center.

Now let $W \in V$. If $X,Y \in \mathfrak{he}_d$, then $[W,X]$ and $[W,Y]$ lie in the center, which is orthogonal to $V$. If $X=Y=T$, then (\ref{eq:ad-invariant}) follows from the fact that $T$ is orthogonal to $V \ni [W,X], [W,Y]$. For symmetry reasons, it therefore suffices to check the equation for $X\neq T$, $Y = T$.

In the case $X=Z$, we can use that $Z$ lies in the center and is orthogonal to $V$.

If $W=X_k$ and $X\neq Y_k$ or $W=Y_k$ and $X\neq X_k$, then $[W,X]=0$ and $[W,Y] \in \mathds{R}W^\prime$ is orthogonal to $X$, where $W^\prime=Y_k$ or $W^\prime=X_k$, respectively.

Finally, we consider $W=X_k$ and $X= Y_k$ or $W=Y_k$ and $X= X_k$. Then
\begin{equation*}
\langle [W,X],Y\rangle+\langle X,[W,Y]\rangle=\langle\pm \lambda_k Z,T\rangle +\langle X,\mp \lambda_k X\rangle=\pm \lambda_k \alpha \mp \lambda_k \alpha=0
\end{equation*}
and we are done.

(ii) As above, we choose a canonical basis $\left\{T, Z, X_1, Y_1, \ldots, X_d, Y_d\right\}$ of $\mathfrak{he}_d^\lambda$. Let $V:=\textnormal{span}\left\{X_1, Y_1, \ldots, X_d, Y_d\right\}$ and $\langle \cdot,\cdot \rangle$ be an $\textnormal{ad}(\mathfrak{he}_d)$-invariant Lorentzian scalar product on $\mathfrak{he}_d^\lambda$. Then Equation~(\ref{eq:ad-invariant}) holds for all $X,Y \in \mathfrak{he}_d^\lambda$ and $W \in \mathfrak{he}_d$.

Let $X=X_k, Y=Y_k$. Equation~(\ref{eq:ad-invariant}) yields \[\langle X_k,Z\rangle=0 \text{ for } W=X_k \quad \text{and} \quad \langle Y_k,Z\rangle=0 \text{ for } W=Y_k.\]
This is true for any $k$. If we choose $W=X_k,X=Y_k,Y=Z$, we obtain \[\langle Z,Z\rangle=0.\]
For $X=T,Y=X_j$ and $W=X_k$, we get for all $j,k$ that \[\langle Y_k, X_j \rangle=0.\]
$X=T,Y=X_j$ and $W=Y_k$, $k \neq j$, yields for all $j\neq k$ that \[\langle X_k, X_j \rangle=0.\]
In the same way, for all $j\neq k$, \[\langle Y_k, Y_j \rangle=0.\]
Choosing $X=T,Y=X_k$ and $W=Y_k$ in Equation~(\ref{eq:ad-invariant}), we obtain for all $k$ that \[\langle X_k, X_k \rangle=\langle T,Z \rangle.\]
Analogously, for all $k$ we have \[\langle Y_k, Y_k \rangle=\langle T,Z \rangle.\]
Thus, $\langle \cdot,\cdot \rangle$ is determined by the two real parameters $\alpha:=\langle T,Z \rangle$ and $\beta:=\langle T,T \rangle$. Since $\langle \cdot,\cdot \rangle$ is Lorentzian, $\alpha>0$.

(iii) Now let $\left\{T,Z,X_1,\ldots,Y_d\right\}$ be a canonical basis of the Lie algebra of the semidirect product $\mathds{S}^1 \ltimes \textnormal{He}_d$ or $\mathds{S}^1 \ltimes \widetilde{\textnormal{He}_d}$, respectively, meaning that $T$ is a generator of the $\mathds{S}^1$-factor and $\left\{Z,X_1,\ldots,Y_d\right\}$ is a canonical basis of the Heisenberg subalgebra. Denote by $\langle \cdot,\cdot\rangle$ the ad-invariant Lorentz form.

Let $a_k Z$ and $b_k Z$ be the $Z$-components (with respect to the canonical basis) of $[T,X_k]$ and $[T,Y_k]$, respectively. Define \[T^\prime := T + \sum\limits_{k=1}^d (a_k Y_k-b_k X_k).\] Then $\textnormal{ad}_{T^\prime}$ is an endomorphism of the subspace $V$, which is defined as above.

Since \[\langle [X_k,Y_k], X \rangle=\langle X_k, [Y_k,X] \rangle = 0\] for any canonical basis element $X$ of $V$ and all $k$, \[\langle Z, X \rangle =0\] for all $X\in \mathfrak{he}_d$. Because $\langle \cdot ,\cdot\rangle$ is a Lorentz form, $\langle \cdot ,\cdot\rangle$ restricted to $V \times V$ is positive definite and \[\alpha := \langle T^\prime ,Z \rangle \neq 0.\] Passing to $-T^\prime$ instead of $T^\prime$ if necessary, we may suppose $\alpha >0$.

We have \[\langle [T^\prime,Z] ,X\rangle =\langle T^\prime ,[Z,X]\rangle=0\] for all $X \in \mathfrak{he}_d$ and using the $\textnormal{ad}_{T^\prime}$-invariance, also \[\langle [T^\prime,Z] ,T^\prime \rangle =-\langle Z,[T^\prime ,T^\prime] \rangle=0.\] Thus, \[[T^\prime,Z]=0.\]

Let $\omega$ be the non-degenerate alternating bilinear form on $V$ defined by the relation \[[X,Y]=\omega(X,Y)Z.\] With respect to an $\langle \cdot ,\cdot\rangle$-orthonormal basis of $V$, $\omega$ corresponds to a skew-symmet\-ric matrix $\Omega$. $\Omega$ has no kernel. Choose an eigenvector $X \in V$ of the symmetric matrix $\Omega^2$. Then $\langle \Omega X,X \rangle =0$ and \[U:=\textnormal{span}\left\{X,\Omega X\right\}\] is $\Omega$-invariant. Furthermore, for any $A \in U^\perp$, \[\langle \Omega A,B \rangle =-\langle A,\Omega B \rangle=0\] for all $B \in U$, so $U^\perp$ is $\Omega$-invariant. Clearly, the operator $\Omega |_{U^\perp}$ is still skew-self-adjoint. Proceeding by induction, we find an $\langle \cdot ,\cdot\rangle$-orthonormal basis $b_1,\ldots,b_{2d}$ such that for $j<k$, $\omega(b_{j},b_{k})=0$, unless $j=2l-1,k=2l$. Thus, without loss of generality, we can choose $\left\{X_1,Y_1,\ldots,X_d,Y_d\right\}$ to be an $\langle \cdot ,\cdot\rangle$-orthogonal basis.

Using $\langle [T^\prime,X],Y \rangle= \langle T^\prime ,[X,Y] \rangle = \omega(X,Y) \alpha$, \[[T^\prime,X_k]=\frac{\alpha}{\langle Y_k, Y_k \rangle} Y_k \text{ and } [T^\prime,Y_k]=-\frac{\alpha}{\langle X_k, X_k \rangle} X_k.\] For any parameters $\eta_1,\ldots,\eta_d \in \mathds{R}_+$, the basis $\left\{Z,\eta_1 X_1, \frac{1}{\eta_1} Y_1, \ldots, \eta_d X_d, \frac{1}{\eta_d} Y_d \right\}$ of $\mathfrak{he}_d$ fulfills the same Lie bracket relations as a canonical basis. Let us choose \[\eta_k:=\sqrt[4]{\frac{{\langle Y_k, Y_k \rangle}}{\langle X_k, X_k \rangle}}\] for $k=1,\ldots,d$. Then for all $k$, \[\langle \eta_k X_k, \eta_k X_k \rangle=\langle \frac{1}{\eta_k} Y_k, \frac{1}{\eta_k} Y_k \rangle.\] Without loss of generality, we can choose the basis $\left\{Z,X_1,Y_1,\ldots,X_d,Y_d\right\}$ in such a way that $\langle X_k,  X_k \rangle=\langle  Y_k, Y_k \rangle$ holds for all $k$. With \[\lambda_k:=\frac{\alpha}{\langle X_k, X_k \rangle},\] $[T^\prime,X_k]=\lambda_k Y_k \text{ and } [T^\prime,Y_k]=-\lambda_k X_k$. $\lambda_k>0$ since $\alpha>0$.

Because of the $\mathds{S}^1$-factor in the semidirect product, all quotients $\frac{\lambda_j}{\lambda_k}$ have to be rational. The claim now follows from Lemma~\ref{lem:isom_heisenberg}.

(iv) The nilradical of the twisted Heisenberg group $\textnormal{He}_d^\lambda=\mathds{S}^1 \ltimes \textnormal{He}_d$ is the subgroup $\textnormal{He}_d$. A theorem of Mostow states, that if $\Gamma$ is a lattice in a connected solvable Lie group $R$ with nilradical $N$, so is $\Gamma \cap N$ a lattice in $N$ (cf.~\cite{OnVin00}, Part~I, Theorem~3.6). If we apply this to our situation, we obtain that given a lattice $\Gamma$ in $\textnormal{He}_d^\lambda$, $\Gamma \cap \textnormal{He}_d$ is a lattice in $\textnormal{He}_d$. Since $\textnormal{He}_d$ is cocompact in $\textnormal{He}_d^\lambda$, any lattice in $\textnormal{He}_d$ is also one in $\textnormal{He}_d^\lambda$. Obviously, $\Gamma \cap \textnormal{He}_d$ has finite index in $\Gamma$. The same argumentation remains true if we replace $\textnormal{He}_d$ by $\widetilde{\textnormal{He}_d}$.

If we consider the universal cover $\pi: \widetilde{\textnormal{He}_d}\to\textnormal{He}_d$, then $\pi^{-1}(\Gamma)$ is a lattice in $\widetilde{\textnormal{He}_d}$ if $\Gamma$ is a lattice in $\textnormal{He}_d$. Conversely, let $\widetilde{\Gamma}$ be a lattice in $\widetilde{\textnormal{He}_d}$. By \cite{Ra72}, Proposition~2.17, $C \cap \widetilde{\Gamma}$ is a lattice in $C$, where $C$ is an element of the ascending central series of $\widetilde{\textnormal{He}_d}$. Choosing $C=Z(\widetilde{\textnormal{He}_d})$ to be the center of $\widetilde{\textnormal{He}_d}$, we see that $\Gamma^\prime:=Z(\widetilde{\textnormal{He}_d}) \cap \widetilde{\Gamma}$ is a lattice in the center. Thus, $\widetilde{\Gamma}$ projects to a lattice of $\textnormal{He}_d \cong \widetilde{\textnormal{He}_d} / \Gamma^\prime$.
\end{proof}\vspace{0pt}

\begin{remark}
In \cite{GorWi86}, Theorem~2.4, a classification of all uniform lattices in $\widetilde{\textnormal{He}_d}$ up to automorphisms of $\widetilde{\textnormal{He}_d}$ is given. Note that any lattice in the simply-connected nilpotent Lie group $\widetilde{\textnormal{He}_d}$ is uniform (cf.~\cite{OnVin00}, Part~I, Chapter~2, Theorem~2.4).

The uniform lattices are classified up to automorphism of $\widetilde{\textnormal{He}_d}$ by the lattices $\Gamma_r$, constructed in the following way: Let $r=(r_1,\ldots,r_d) \in \mathds{Z}_+^d$, such that $r_j$ divides $r_{j+1}$ for all $j=1,\ldots,d-1$. In the matrix model of $\widetilde{\textnormal{He}_d}$ given in Section~\ref{sec:Heisenberg}, $\Gamma_r$ consists of all $(d+2) \times (d+2)$-matrices \[\begin{pmatrix}1 & x & z \\ 0 & I_{d} & y \\ 0 & 0 & 1\end{pmatrix},\] such that $z \in \mathds{Z}$, the row vector $x=(x_1,\ldots, x_d)$ has the property, that $x_j \in r_j\mathds{Z}$ for all $j=1,\ldots,d$, and the column vector $y$ has integer entries.
\end{remark}

From the proof of the third part of Proposition~\ref{prop:lorentz_heisenberg}, we see the following:

\begin{corollary}\label{cor:orthogonal_basis}
Let $\langle \cdot,\cdot \rangle$ be a positive definite scalar product on $V$, where $V$ is a vector space complement to the center $\mathds{R}Z$ in $\mathfrak{he}_d$. Then there exists a canonical basis $\left\{Z,X_1,Y_1,\ldots,X_d,Y_d\right\}$, such that $\left\{X_1,Y_1,\ldots,X_d,Y_d\right\}\subset V$ is $\langle \cdot,\cdot \rangle$-orthogonal.
\end{corollary}

As a consequence of the fourth part of Proposition~\ref{prop:lorentz_heisenberg}, a twisted Heisenberg group $\textnormal{He}_d^\lambda$ or $\overline{\textnormal{He}_d^\lambda}$ admits a uniform lattice $\Lambda$, since $\widetilde{\textnormal{He}_d}$ has a uniform lattice. For example, in the matrix representation given in Section~\ref{sec:Heisenberg}, a uniform lattice is given by $x$ and $y$ having integer entries and $z \in \mathds{Z}$.

Note that by a theorem of Mostow (cf.~\cite{Mo62}, Theorem~6.2), for a solvable Lie group $G$ and a closed subgroup $H \subseteq G$, $G/H$ is compact if it has a finite invariant measure. Thus, any lattice $\Lambda$ in $\textnormal{He}_d^\lambda$ or $\overline{\textnormal{He}_d^\lambda}$, respectively, is uniform.

Any ad-invariant Lorentz form on the Lie algebra given in Proposition~\ref{prop:lorentz_heisenberg}~(i) gives a Lorentzian metric on $\textnormal{He}_d^\lambda/\Lambda$ or $\overline{\textnormal{He}_d^\lambda}/\Lambda$, respectively, such that $\textnormal{He}_d^\lambda$ or $\overline{\textnormal{He}_d^\lambda}$, respectively, acts isometrically and locally effectively by left multiplication.

\section{Induced bilinear form on the Lie algebra}\label{sec:induced_form}

Let $M=(M,g)$ be a Lorentzian manifold of finite volume and $\mu$ the corresponding Lebesgue measure on the manifold. We consider a connected Lie group $G$ with Lie algebra $\mathfrak{g}$ acting isometrically and locally effectively on $M$.

According to Proposition~\ref{prop:isometry_Killing} and the preceding remark, $X,Y \in \mathfrak{g}$ correspond to Killing vector fields $\widetilde{X},\widetilde{Y}$ in $M$.

\begin{definition}
In the situation above, let $U$ be a $G$-invariant non-empty open subset of $M$, such that for any $X,Y \in \mathfrak{g}$, $\left|g(\widetilde{X},\widetilde{Y})(x)\right| \leq C$ for all $x \in U$, where $C$ is a constant depending only on $X,Y$. Then $\kappa$ is the \textit{induced bilinear form} on $\mathfrak{g}$ defined by
\begin{equation*}
\kappa(X,Y):=\int\limits_U{g(\widetilde{X},\widetilde{Y})(x) d\mu(x)}.
\end{equation*}
\end{definition}

\begin{lemma}\label{lem:kappa_invariant}
In the situation above, the following is true:
\begin{compactenum}
\item For any $f \in G$ and $X \in \mathfrak{g}$, $df(\widetilde{X}(x))=\widetilde{\textnormal{Ad}_f(X)}(f \cdot x)$ for all $x \in M$.

\item $\kappa$ is Ad-invariant, that is, $\kappa(\textnormal{Ad}_f(X),\textnormal{Ad}_f(Y))=\kappa(X,Y)$ for all $f \in G$ and $X,Y \in \mathfrak{g}$.

\item $\kappa$ is ad-invariant.
\end{compactenum}
\end{lemma}
\vspace{-1.4em}\begin{proof}
(i) By definition of $\widetilde{X}$, we have:
\begin{align*}
df\left(\widetilde{X}\left(x\right)\right)&=\frac{\partial}{\partial t} \left(f \cdot \left(\exp\left(t X\right) \cdot x\right)\right)|_{t=0}\\
&=\frac{\partial}{\partial t} \left( \left(f \exp\left(t X\right) f^{-1}\right)f  \cdot x\right)|_{t=0}\\
&=\widetilde{\textnormal{Ad}_f\left(X\right)}\left(f \cdot x\right).\end{align*}
(ii) By the first part and because $G$ acts isometrically, 
\begin{align*}
\kappa\left(\textnormal{Ad}_f\left(X\right),\textnormal{Ad}_f\left(Y\right)\right)&=\int\limits_U{g\left(\widetilde{\textnormal{Ad}_f\left(X\right)},\widetilde{\textnormal{Ad}_f\left(Y\right)}\right)\left(x\right) d\mu\left(x\right)}\\
&=\int\limits_U{g\left(df(\widetilde{X}),df(\widetilde{Y})\right)\left(f^{-1} \cdot x\right) d\mu\left(x\right)}\\
&=\int\limits_U{g(\widetilde{X},\widetilde{Y})\left(y\right) d\mu\left(y\right)}\\
&=\kappa\left(X,Y\right).\end{align*}
(iii) By the second part, for any $X,Y,Z \in \mathfrak{g}$ we have:
\begin{align*}
\kappa([X,Y],Z)&=\kappa\left(\frac{\partial}{\partial t}\left( \textnormal{Ad}_{\exp(tX)}(Y)\right)|_{t=0},Z\right)\\
&=\frac{\partial}{\partial t} \left(\kappa(\textnormal{Ad}_{\exp(tX)}(Y),Z)\right)|_{t=0}\\
&=\frac{\partial}{\partial t} \kappa(Y,\textnormal{Ad}_{\exp(-tX)}(Z))|_{t=0}\\
&=\kappa\left(Y,\frac{\partial}{\partial t} \left(\textnormal{Ad}_{\exp(-tX)}(Z)\right)|_{t=0}\right)\\
&=-\kappa(Y,[X,Z]).\qedhere
\end{align*}
\end{proof}\vspace{0pt}

In the case that $M$ is compact, we can simply take $U=M$. But in general, this is not possible for a semi-Riemannian manifold of finite volume.

\begin{example}
Consider $M=\mathds{R}^3$ with cylindrical coordinates $(r,\varphi,z)$ on $M \backslash \left\{0\right\}$ and the Lorentzian metric \[g_{(r,\varphi,z)}=f(r)^2dr^2+r^2d\varphi^2-\exp(-2z^2)dz^2.\] Here $f:\mathds{R}_+ \to \mathds{R}_+$ is a smooth function with $f(r)=1$ if $r<1$ and $f(r)=r^{-3}$ if $r>2$. It follows that $g$ extends to a Lorentzian metric on the whole of $M$.

The volume form of $M=(M,g)$ is determined by \[d\textnormal{vol}=r f(r) \exp(-z^2) dr d\varphi dz.\] Hence, the volume of $M$ is
\begin{align*}
\textnormal{vol}(M)&=\int\limits_M{r f(r) \exp(-z^2) dr d\varphi dz}\\
&=2\pi \sqrt{\pi} \int\limits_{0}^\infty{r f(r)dr}\\
&=2\pi \sqrt{\pi} \left(\frac{1}{2}+\int\limits_{1}^2{r f(r)dr}+\frac{1}{2}\right)<+\infty.\end{align*}
$\mathds{S}^1$ is acting isometrically and effectively on $M$: $\theta \in \mathds{S}^1$ corresponds to the isometry $(r,\varphi,z) \mapsto (r,\varphi+\theta,z)$. Thus, $\frac{\partial}{\partial \varphi}$ is a Killing vector field corresponding to an element of the Lie algebra of $\mathds{S}^1$. $g_{(r,\varphi,z)}\left(\frac{\partial}{\partial \varphi},\frac{\partial}{\partial \varphi}\right)=r^2$. Therefore,
\begin{align*}
\int\limits_M{g\left(\frac{\partial}{\partial \varphi},\frac{\partial}{\partial \varphi}\right)d\textnormal{vol}}&=\int\limits_M{r^3 f(r) \exp(-z^2) dr d\varphi dz}\\
&=2\pi \sqrt{\pi} \int\limits_{0}^\infty{r^3 f(r)dr}\\
&=2\pi \sqrt{\pi} \left(\frac{1}{4}+\int\limits_{1}^2{r^3 f(r)dr}+\int\limits_{2}^\infty{1dr}\right)=+\infty.\end{align*}
\end{example}

The given example can be easily generalized to higher dimensions and different signatures of the metric.

We will now show that an open set $U$ as required in the definition of $\kappa$ always exists in our situation. The idea of the proof is due to Zeghib, but was not mentioned in \cite{Ze98a}.

\begin{proposition}\label{prop:invariant_open_set}
Let $M=(M,g)$ be a semi-Riemannian manifold of finite volume and $G$ a connected Lie group with Lie algebra $\mathfrak{g}$ acting isometrically and locally effectively on $M$. Then there is a $G$-invariant non-empty open subset $U \subseteq M$, such that for any $X,Y \in \mathfrak{g}$, $\left|g(\widetilde{X},\widetilde{Y})(x)\right| \leq C$ for all $x \in U$, where $C$ is a constant depending only on $X,Y$.
\end{proposition}
\vspace{-1.4em}\begin{proof}
We consider the Gau{\ss} map $\textnormal{Ga}: M \to S^2\mathfrak{g}$, which maps $x \in M$ to the symmetric bilinear form $b_x$ defined by \[b_x(X,Y)=g_x(\widetilde{X}(x),\widetilde{Y}(x)).\] $G$ acts from the left on $M$ through isometries and on $S^2\mathfrak{g}$ through the adjoint representation: \[(f \cdot b)(X,Y):=b(\textnormal{Ad}_{f^{-1}}(X),\textnormal{Ad}_{f^{-1}}(Y))\] for any $f \in G, b \in S^2\mathfrak{g}$. Using the first part of Lemma~\ref{lem:kappa_invariant}, we obtain:
\begin{align*}
b_{f \cdot x}(X,Y)&=g_{f\cdot x}(\widetilde{X}(f \cdot x),\widetilde{Y}(f \cdot x))\\
&=g_{x}(df^{-1}(\widetilde{X}(f \cdot x)),df^{-1}(\widetilde{Y}(f \cdot x)))\\
&=g_{x}(\widetilde{\textnormal{Ad}_{f^{-1}}(X)}(x),\widetilde{\textnormal{Ad}_{f^{-1}}(Y)}(x))\\
&=(f\cdot b)_{x}(X,Y).\end{align*}
Thus, Ga is $G$-equivariant. If the image of Ga consists of only one point, we are done. So suppose the image of Ga consists of more than one point. Let $\textnormal{Ga}(\mu)$ denote the push-forward of the Lebesgue measure $\mu$ on $M$ under the Gau{\ss} map. $\textnormal{Ga}(\mu)$ is $G$-invariant.

$S^2\mathfrak{g}$ is a real vector space and $G$ acts linearly on $S^2\mathfrak{g}$, that is, the action of $G$ is associated to a representation $\varrho: G \to \textnormal{GL}(S^2\mathfrak{g})$. The representation \[\mathds{P}(\varrho): G \to \textnormal{GL}(S^2\mathfrak{g}) \to \textnormal{PGL}(S^2\mathfrak{g})\] results in an action of $G$ on the projective space $\mathds{P}(S^2\mathfrak{g})$.

The Gau{\ss} map induces a $G$-equivariant map $M \to \mathds{P}(S^2\mathfrak{g})$. Let $\lambda$ denote the push-forward of $\mu$. $\lambda$ is a $G$-invariant finite measure on $\mathds{P}(S^2\mathfrak{g})$. Without loss of generality, we may assume that $\textnormal{vol}(M)=1$, so $\lambda$ is a probability measure.

It turns out that the F\"urstenberg lemma is quite powerful in situations like our. One can find a proof in \cite{Zi84}, Lemma~3.2.1. For convenience, we will formulate only a special case of it, which is sufficient for our purposes.

\begin{lemma}\label{lem:Furstenberg}
Let $\lambda$ be a probability measure on $\mathds{RP}^{d}$. Consider a sequence $\left\{T_m\right\}_{m=0}^\infty$ of projective transformations leaving $\lambda$ invariant. Then at least one of the following statements is true:
\begin{compactenum}
\item $\left\{T_m\right\}_{m=0}^\infty$ is a precompact sequence in $\textnormal{PGL}(n,\mathds{R})$.

\item There exist linear subspaces $V,W \subset \mathds{R}^d$, $1 \leq \dim V, \dim W \leq d-1$, such that the support of $\lambda$ is contained in $\mathds{P}(V) \cup \mathds{P}(W)$.
\end{compactenum}
\end{lemma}

The lemma applies to our situation considering a sequence $\left\{f_m\right\}_{m=0}^\infty$ of elements of $G$ acting on $\mathds{P}(S^2\mathfrak{g})$. Suppose we are in the case~(ii) of lemma \ref{lem:Furstenberg}. If one of the subspaces $V,W$ is contained in the other, we may ignore this subspace. Then the support of $\lambda$ is contained in say $\mathds{P}(V)$. Because $\lambda$ is $G$-invariant and $G$ is connected and acts by projective transformations on $\mathds{P}(S^2\mathfrak{g})$, the subspace $V$ is $G$-invariant if $V$ is chosen minimal with respect to inclusion. Analogously, if the support of $\lambda$ is contained in $\mathds{P}(V) \cup \mathds{P}(W)$ and neither $V \subseteq W$ nor $W \subseteq V$, $V$ and $W$ are $G$-invariant if they are chosen minimal with respect to inclusion. Thus, we can apply the F\"urstenberg lemma also to $V$ and $W$ instead of $S^2\mathfrak{g}$.

By induction, we may suppose that the support of $\lambda$ is contained in the union $\mathds{P}(V_1) \cup \ldots \cup \mathds{P}(V_k)$ with $G$-invariant subspaces $V_j \subset S^2\mathfrak{g}$, $1 \leq \dim V_j < \dim S^2\mathfrak{g}$, no subspace is contained in another and applying lemma \ref{lem:Furstenberg} to one of the $V_j$ results in the case~(i). Additionally, we suppose that all the $V_j$ are chosen minimal with respect to inclusion. Note that we assume the image of Ga to consist of more than one point, so $k>0$.

The action of $G$ on $V_j$ is associated to the map $\varrho_j:G \to \textnormal{GL}(V_j)$ induced by $\varrho$. Let $\pi_j:\textnormal{GL}(V_j) \to \textnormal{PGL}(V_j)$ denote the canonical projection. By assumption, $\pi_j(\varrho_j(G))$ is precompact in $\textnormal{PGL}(V_j)$.

The push-forward of the Lebesgue measure on $M$ under $\pi^j \circ \textnormal{Ga}: M \to S^2\mathfrak{g} \to V_j$, $\pi^j: S^2\mathfrak{g} \to V_j$ being the canonical projection, is denoted by $\lambda_j$ and is a $G$-invariant probability measure. Assume that $t \textnormal{id}_{V_j} \in \varrho_j(G)$, for some real $t$ with $|t| \neq 1$. Since $\varrho_j(G)$ is a subgroup of $\textnormal{GL}(V_j)$, for any $m \in \mathds{Z}_+$, there are $t_m,t_m^\prime \in \mathds{R}$ such that $t_m \textnormal{id}_{V_j}, t^\prime_m \textnormal{id}_{V_j} \in \varrho_j(G)$ and $0<|t|<\frac{1}{m}$, $|t^\prime|>m$.

But $\lambda_j$ is $G$-invariant, so one easily follows that the measure is concentrated on $0 \in V_j$, contradicting the minimality of $V_j$ with respect to inclusion. Thus, $t \textnormal{id}_{V_j} \in \varrho_j(G)$ only if $|t|=1$. Hence, the kernel of $\pi_j|_{\varrho_j(G)}$ contains at most two elements. Therefore, since $\pi_j(\varrho_j(G))$ is precompact in $\textnormal{PGL}(V_j)$, $\varrho_j(G)$ is precompact in $\textnormal{GL}(V_j)$.

In summary, we have shown that $G$ acts precompactly on $V_1 \cup \ldots \cup V_k$, which contains the support of $\textnormal{Ga}(\mu)$, that is, there is a compact group $K \subset \textnormal{GL}(S^2\mathfrak{g})$ and a group homomorphism $\widetilde{\varrho}: G \to K$, such that for all $f \in G$, \[\varrho(f)|_{\textnormal{supp}(\textnormal{Ga}(\mu))}=\widetilde{\varrho}(f)|_{\textnormal{supp}(\textnormal{Ga}(\mu))}.\]

Let $C$ be an $K$-invariant bounded open subset of $S^2\mathfrak{g}$ intersecting $\textnormal{Ga}(M)$ non-trivially (for example, $C$ is the $K$-orbit of an open ball containing an element of the image of Ga). Let $U_C:=\textnormal{Ga}^{-1}(C)$ and $U$ be the $G$-orbit of $U_C$. Since $C$ is open, $U_C$ and hence $U$ are open as well. Also, $U$ is $G$-invariant and non-empty. Because $C$ is $K$-invariant, \[\textnormal{Ga}(U) \cap \textnormal{supp}(\textnormal{Ga}(\mu)) = \textnormal{Ga}(U_C) \cap \textnormal{supp}(\textnormal{Ga}(\mu)).\] Since $\textnormal{Ga}(\mu)$ is the push-forward of $\mu$, it follows that $U \backslash U_C$ is a set of $\mu$-measure zero. Thus, $U_C$ is dense in $U$.

By construction, for any $X,Y \in \mathfrak{g}$, $\left|g(\widetilde{X},\widetilde{Y})(x)\right|$ is bounded uniformly on $U_C$. Since $U_C$ is dense in $U$, it is also bounded uniformly on $U$.
\end{proof}\vspace{0pt}

The proof of Proposition~\ref{prop:invariant_open_set} yields the following corollary of the F\"urstenberg lemma:

\begin{corollary}\label{cor:Furstenberg}
Suppose that a connected Lie group $G$ is acting continuously and linearly on a real vector space $V$ and is preserving a finite measure $\mu$ on $V$. Then $G$ acts precompactly on the support of $\mu$, that is, there is a Lie group homomorphism $\varrho: G \to K$, $K\subset \textnormal{GL}(V)$ a compact subgroup, such that for all $x$ in the support of $\mu$, $f \cdot x = \varrho(f)(x)$.
\end{corollary}

Our next aim is to show a certain non-degeneracy condition of $\kappa$, which will allow us to classify the possible Lie algebras $\mathfrak{g}$ of Lie groups acting isometrically and (locally) effectively on Lorentzian manifolds of finite volume. The key point is that $\kappa(X,X) \geq 0$ if $X\in\mathfrak{isom}(M)$ generates a non-precompact one-parameter group in $\textnormal{Isom}(G)$.

\begin{proposition}\label{prop:kappa_positive}
Let $M=(M,g)$ be a Lorentzian manifold of finite volume. If $X \in \mathfrak{isom}(M)$ generates a non-precompact one-parameter group in $\textnormal{Isom}(M)$, then $\widetilde{X}$ is non-timelike everywhere, that is, $g(\widetilde{X},\widetilde{X})(x)\geq 0$ for all $x \in M$.
\end{proposition}
\vspace{-1.4em}\begin{proof}
It suffices to show, that if for the generator $X \in \mathfrak{isom}(M)$ of a one-parameter group $\left\{\varphi^t\right\}_{t \in \mathds{R}}$ of isometries the corresponding vector field $\tilde{X}$ is somewhere timelike, then $\left\{\varphi^t\right\}_{t \in \mathds{R}}$ is precompact.

$\left\{\varphi^t\right\}_{t \in \mathds{R}}$ preserves the Lebesgue measure $\mu$ on $M$. By assumption, $\mu(M)<\infty$, so without loss of generality, $\mu$ is a probability measure. As a consequence of the Poincar\'e recurrence theorem, $\mu$-almost all $x \in M$ are recurrent with respect to $\left\{\varphi^t\right\}$ (cf.~\cite{BrSt02}, Proposition~4.2.2), that is, for $\mu$-almost all $x \in M$ there is a sequence $\left\{t_k\right\}_{k=0}^\infty \subset \mathds{R}$, such that $t_k \to \infty$ and $\varphi^{t_k}(x) \to x$ as $k \to \infty$.

The set $T$ of points, where $\widetilde{X}$ is timelike, is non-empty and open. Thus, there is a recurrent $x \in M$ such that $g(\widetilde{X},\widetilde{X})(x)<0$. Let $\left\{t_k\right\}_{k=0}^\infty$ be a sequence of real numbers, such that $t_k \to \infty$ and $x_k:=\varphi^{t_k}(x) \to x$ as $k \to \infty$. Without loss of generality, $x_k \in V$ for all $k$, where $V$ is a precompact neighborhood of $x$ such that for all $y$ in the closure of $V$, $\widetilde{X}(y)$ is timelike.

Now transform for all points of $T$ the Lorentzian metric $g$ into a Riemannian metric $h$ by changing the sign along $\widetilde{X}$ (the metric on $\widetilde{X}^\perp$ remains the same). Then $\left\{\varphi^t\right\}_{t \in \mathds{R}}$ also preserves $h$ (at the points where it is defined). Hence, \[\left\| d\varphi^{t_k}_x\right\|:=\max\limits_{v \in T_xM\backslash\left\{0\right\}}{\frac{h_{x_k}\left(d\varphi^{t_k}_x(v),d\varphi^{t_k}_x(v)\right)}{h_x\left(v,v\right)}}=1.\]

Recall the mapping $\phi_x: \textnormal{Isom}(M) \to O(M)$ defined by $f \mapsto \left(df_{x}(s_1), \ldots, df_{x}(s_n)\right)$ for a fixed orthonormal basis $\left\{s_1, \ldots, s_n\right\}$ in the tangent space $(T_xM,g_x)$. Let $K \subset O(M)$ denote the set of orthonormal frames, such that the base point lies in the closure of $V$ and for all $j=1,\ldots,n$, the $h$-norm of the $j$th vector of each frame in $K$ is bounded by the $h$-norm of $s_j$. $K$ is compact and by construction, $\phi_x(\varphi^{t_k}) \in K$ for all $k$. Therefore, $\left\{\varphi^{t_k}\right\}_{k=0}^\infty$ is precompact in $\textnormal{Isom}(M)$.

Let $L$ be the closure of $\left\{\varphi^t\right\}_{t \in \mathds{R}}$ in $\textnormal{Isom}(M)$. As a connected and abelian group, $L$ is isomorphic to $\mathds{T}^m \times \mathds{R}^{m^\prime}$ by Lemma~\ref{lem:abelian_group}. Since the one-parameter group $\left\{\varphi^t\right\}_{t \in \mathds{R}}$ lies dense in $L$, $L$ is either isomorphic to a torus $\mathds{T}^m$ or to the real line. But $\left\{\varphi^{t_k}\right\}_{k=0}^\infty$ is precompact, so $L$ is a torus, meaning that $\left\{\varphi^t\right\}_{t \in \mathds{R}}$ is precompact.
\end{proof}\vspace{0pt}

\begin{lemma}\label{lem:kernel_kappa}
Let $M=(M,g)$ be a Lorentzian manifold of finite volume and $V$ a subspace of $\mathfrak{isom}(M)$, such that for all $X \in V$, $\widetilde{X}$ is non-timelike everywhere. Then the restriction of the induced bilinear form $\kappa$ to $V \times V$ is positive semidefinite and its kernel has dimension at most 1.
\end{lemma}
\vspace{-1.4em}\begin{proof}
By assumption, $\widetilde{X}$ is lightlike everywhere on the open set $U$ used in the definition of $\kappa$, if $X \in V$ is $\kappa$-isotropic. So if $W$ is a totally $\kappa$-isotropic subspace of $V$ and $x \in U$, then $W_x:=\left\{\widetilde{X}(x) | X \in W\right\}$ is a totally isotropic subspace of $(T_xM,g_x)$. Since $g$ is Lorentzian, $\dim W_x \leq 1$ for all $x \in U$.

Thus, it suffices to show, that if two Killing vector fields $\widetilde{X},\widetilde{Y}$ are linearly dependent on each point of some arbitrary open set $U$, then one is a multiple of the other.

If one of the Killing vector fields vanishes on an open set, then it vanishes everywhere. So let $\widetilde{Y}=\sigma\widetilde{X}\neq 0$ on a connected open set $U^\prime$, $\sigma$ a smooth function on $U^\prime$. Denote by $\nabla$ the Levi-Civita connection of $M$. Since $\widetilde{X},\widetilde{Y}$ are Killing vector fields, $\nabla \widetilde{X}$ and $\nabla \widetilde{Y}$ are skew-self-adjoint. Therefore,
\begin{equation*}
0=g(\nabla_Z \widetilde{Y},Z)=d\sigma(Z) g(\widetilde{X},Z)+\sigma g(\nabla_Z \widetilde{X},Z)=d\sigma(Z) g(\widetilde{X},Z)
\end{equation*}
for any vector field $Z$ on $U^\prime$. Since $\widetilde{X}(x) \neq 0$ for all $x \in U^\prime$, $\widetilde{X}(x)^\perp$ is a subspace of codimension 1 in $T_xM$. Thus, the set of vectors in $T_xM$ not orthogonal to $\widetilde{X}(x)$ is dense in $T_xM$. It follows that $d\sigma(Z)=0$ for all vector fields $Z$ on $U^\prime$. Thus, $\sigma=\sigma_0$ is constant on $U^\prime$. Since $\widetilde{X}-\sigma_0\widetilde{Y}$ is a Killing vector field and vanishes on an open set, $\widetilde{X}=\sigma_0\widetilde{Y}$ on $M$.
\end{proof}\vspace{0pt}

\begin{corollary}\label{cor:condition_star}
Let $M=(M,g)$ be a Lorentzian manifold of finite volume and $V$ a subspace of $\mathfrak{isom}(M)$, such that the set of $X \in V$ generating a non-precompact one-parameter group of isometries is dense in $V$. Then the restriction of the induced bilinear form $\kappa$ to $V \times V$ is positive semidefinite and the dimension of its kernel is at most one.
\end{corollary}
\vspace{-1.4em}\begin{proof}
The claim directly follows from Proposition~\ref{prop:kappa_positive} and Lemma~\ref{lem:kernel_kappa}.
\end{proof}\vspace{0pt}

Note that in general, the induced bilinear form $\kappa$ is not a Lorentz form. But it turns out that in the case of a non-compact connected component of the identity in the isometry group, $\kappa$ is either Lorentzian or positive semidefinite and the dimension of its kernel is at most one.

%%%%%%%%%%%%%%%%%%%%%%%%%%%%%%%%%%%%%%%%%%%%%%%%%%%%%%%%%%%%%%%%%%%%%%%%%%%%%%%%%%%%%%%%%%%%

\chapter{Main theorems}\label{ch:theorems}

In this chapter, we describe the main results of our studies. We start with stating the relevant theorems concerning the structure of the Lie algebras of Lie groups acting isometrically and locally effectively on a Lorentzian manifold of finite volume in Section~\ref{sec:alg_th}. In Section~\ref{sec:geo_th}, we continue in giving a geometric characterization of compact Lorentzian manifolds admitting isometric and effective actions of a cover of $\textnormal{PSL}(2,\mathds{R})$ or of a twisted Heisenberg group. This will be relevant for the description of compact homogeneous Lorentzian manifolds with such an action and the investigation of their geometry in Section~\ref{sec:hom_th}.

\section{Algebraic theorems}\label{sec:alg_th}

The key theorem providing the classification of the Lie algebras is the following algebraic result:

\begin{theorem}\label{th:algebraic}
Let $G$ be a connected non-compact Lie group. Assume there is an ad-invariant symmetric bilinear form $\kappa$ on its Lie algebra $\mathfrak{g}$ that fulfills the following non-degeneracy condition:

\par
\begingroup
\leftskip=1em
\noindent \hypertarget{star}{$(\star)$} Let $V$ be a subspace of $\mathfrak{g}$, such that the set of $X \in V$ generating a non-precompact one-parameter group in $G$ is dense in $V$. Then the restriction of $\kappa$ to $V \times V$ is positive semidefinite and its kernel has dimension at most one.
\par
\endgroup

Then $\mathfrak{g}=\mathfrak{k}\oplus\mathfrak{a}\oplus\mathfrak{s}$ is a $\kappa$-orthogonal direct sum of a compact semisimple Lie algebra $\mathfrak{k}$, an abelian algebra $\mathfrak{a}$ and a Lie algebra $\mathfrak{s}$, which is either trivial, isomorphic to $\mathfrak{aff}(\mathds{R})$, to a Heisenberg algebra $\mathfrak{he}_d$, to a twisted Heisenberg algebra $\mathfrak{he}_d^\lambda$ with $\lambda \in \mathds{Z}_+^d$, or to $\mathfrak{sl}_2(\mathds{R})$.

Furthermore, the following is true:
\begin{compactenum}
\item If $\mathfrak{s}$ is trivial, $\kappa$ is positive semidefinite and its kernel has dimension at most one. If $\mathfrak{s}$ is non-trivial, $\kappa$ restricted to $\mathfrak{a} \times \mathfrak{a}$ and $\mathfrak{k} \times \mathfrak{k}$ is positive definite.

\item Let $\mathfrak{s}\cong\mathfrak{aff}(\mathds{R})$. Then the restriction of $\kappa$ to $\mathfrak{s} \times \mathfrak{s}$ is positive semidefinite and its kernel is exactly the span of the generator of the translations in the affine group.

\item Let $\mathfrak{s}\cong\mathfrak{he}_d$. Then the restriction of $\kappa$ to $\mathfrak{s} \times \mathfrak{s}$ is positive semidefinite and its kernel is exactly the center of $\mathfrak{he}_d$.

\item Let $\mathfrak{s}\cong\mathfrak{he}_d^\lambda$. Then $\kappa$ restricted to $\mathfrak{s} \times \mathfrak{s}$ is a Lorentz form. The subgroup in $G$ generated by $\mathfrak{s}$ is isomorphic to $\textnormal{He}_d^\lambda$ or $\overline{\textnormal{He}_d^\lambda}$, if it is closed in $G$ or not, respectively. Moreover, the abelian subgroup generated by $\mathfrak{a}\oplus\mathfrak{z(s)}$ is compact.

\item Let $\mathfrak{s}\cong\mathfrak{sl}_2(\mathds{R})$. Then $\kappa$ restricted to $\mathfrak{s} \times \mathfrak{s}$ is a positive multiple of the Killing form of $\mathfrak{s}$. The subgroup in $G$ generated by $\mathfrak{s}$ is isomorphic to some $\textnormal{PSL}_k(2,\mathds{R})$, the $k$-covering of $\textnormal{PSL}(2,\mathds{R})$, if and only if it is closed in $G$. Moreover, the abelian subgroup generated by $\mathfrak{a}$ is compact.
\end{compactenum}
\end{theorem}

\begin{remark}
For both cases $\mathfrak{s}\cong\mathfrak{sl}_2(\mathds{R})$ and $\mathfrak{s}\cong\mathfrak{he}_d^\lambda$, $\lambda \in \mathds{Z}_+^d$, 
there exists a connected non-compact Lie group $G$ and a symmetric bilinear form $\kappa$ on its Lie algebra $\mathfrak{g}$ fulfilling condition~\hyperlink{star}{$(\star)$}, such that $\mathfrak{g}=\mathfrak{s}\oplus\mathds{R}$ and the subgroup generated by $\mathfrak{s}$ is not closed in $G$. Especially, it is not isomorphic to a finite covering of $\textnormal{PSL}(2,\mathds{R})$ or isomorphic to $\textnormal{He}_d^\lambda$, respectively. We will show this later in Propositions~\ref{prop:counterexample1} and~\ref{prop:counterexample2}.
\end{remark}

According to Corollary~\ref{cor:condition_star}, Theorem~\ref{th:algebraic} applies to a connected non-compact closed Lie subgroup $G$ of $\textnormal{Isom}(M)$, where $M$ is a Lorentzian manifold of finite volume. More generally, we obtain the following classification result:

\begin{theorem}\label{th:algebraic_classification}
Let $M$ be a Lorentzian manifold of finite volume and $G$ a connected Lie group acting isometrically and locally effectively on $M$. Then its Lie algebra $\mathfrak{g}=\mathfrak{k}\oplus\mathfrak{a}\oplus\mathfrak{s}$ is a direct sum of a compact semisimple Lie algebra $\mathfrak{k}$, an abelian algebra $\mathfrak{a}$ and a Lie algebra $\mathfrak{s}$, which is either trivial, isomorphic to $\mathfrak{aff}(\mathds{R})$, to a Heisenberg algebra $\mathfrak{he}_d$, to a twisted Heisenberg algebra $\mathfrak{he}_d^\lambda$ with $\lambda \in \mathds{Z}_+^d$, or to $\mathfrak{sl}_2(\mathds{R})$.

Suppose $G$ is a Lie subgroup of the isometry group $\textnormal{Isom}(M)$. If $\mathfrak{s}\cong\mathfrak{he}_d^\lambda$, the subgroup generated by $\mathfrak{s}$ is isomorphic to $\textnormal{He}_d^\lambda$ or $\overline{\textnormal{He}_d^\lambda}$, if it is closed in $\textnormal{Isom}(M)$ or not, respectively. If $\mathfrak{s}\cong\mathfrak{sl}_2(\mathds{R})$, the subgroup generated by $\mathfrak{s}$ is isomorphic to some $\textnormal{PSL}_k(2,\mathds{R})$, if and only if it is closed in $\textnormal{Isom}(M)$.
\end{theorem}

Together with the examples given in Section~\ref{sec:examples}, Theorem~\ref{th:algebraic_classification} gives a complete characterization of the Lie algebras of Lie groups acting isometrically and locally effectively on Lorentzian manifolds of finite volume.

Adams and Stuck (cf.~\cite{AdSt97b}, Theorem~1.2) and Zeghib (cf.~\cite{Ze98b}, Theorem~1.1) independently showed, that if one is interested in the Lie algebra of a Lie group, which can be exactly the isometry group of some compact Lorentzian manifold, only the case of $\mathfrak{aff}(\mathds{R})$ does not appear. Additionally, a construction of compact Lorentzian manifolds with given Lie algebra of its isometry group was done in \cite{Ze98b} (Section~4, proof of Theorem~1.5). We will not go into details, but summarize the result which follows from their investigation and Theorem~\ref{th:algebraic_classification}.

\begin{theorem}\label{th:isometry_groups}
Let $M$ be a compact Lorentzian manifold. Then the Lie algebra of the isometry group $\textnormal{Isom}(M)$, $\mathfrak{isom}(M)=\mathfrak{k}\oplus\mathfrak{a}\oplus\mathfrak{s}$, is a direct sum of a compact semisimple Lie algebra $\mathfrak{k}$, an abelian algebra $\mathfrak{a}$ and a Lie algebra $\mathfrak{s}$, which is either trivial, isomorphic to a Heisenberg algebra $\mathfrak{he}_d$, to a twisted Heisenberg algebra $\mathfrak{he}_d^\lambda$ with $\lambda \in \mathds{Z}_+^d$, or to $\mathfrak{sl}_2(\mathds{R})$.

Furthermore, the following is true:
\begin{compactenum}
\item Let $\mathfrak{s}\cong\mathfrak{he}_d^\lambda$. Then the subgroup in $\textnormal{Isom}(M)$ generated by $\mathfrak{s}$ is isomorphic to $\textnormal{He}_d^\lambda$ or $\overline{\textnormal{He}_d^\lambda}$, if it is closed in $\textnormal{Isom}(M)$ or not, respectively. The abelian subgroup generated by $\mathfrak{a}\oplus\mathfrak{z(s)}$ is compact.

\item Let $\mathfrak{s}\cong\mathfrak{sl}_2(\mathds{R})$. Then the subgroup in $\textnormal{Isom}(M)$ generated by $\mathfrak{s}$ is isomorphic to some $\textnormal{PSL}_k(2,\mathds{R})$, if and only if it is closed in $\textnormal{Isom}(M)$. The abelian subgroup generated by $\mathfrak{a}$ is compact.
\end{compactenum}
Conversely, for any Lie algebra $\mathfrak{g}=\mathfrak{k}\oplus\mathfrak{a}\oplus\mathfrak{s}$ as above, there is a compact Lorentzian manifold $M$, such that the Lie algebra of its isometry group is isomorphic to $\mathfrak{g}$.
\end{theorem}

\section{Geometric theorems}\label{sec:geo_th}

The following theorem will be important for the geometric characterization of the investigated manifolds.

\begin{theorem}\label{th:locally_free}
Let $M$ be a compact Lorentzian manifold and $G$ a connected Lie group acting isometrically and locally effectively on $M$. Then the action of the subgroup $S$ generated by the direct summand $\mathfrak{s}$ in its Lie algebra $\mathfrak{g}=\mathfrak{k}\oplus\mathfrak{a}\oplus\mathfrak{s}$ (cf.~Theorem~\ref{th:algebraic_classification}) is locally free.
\end{theorem}

The geometric characterization of compact Lorentzian manifolds with a non-compact connected component of the identity in the isometry group is given by the following:

\begin{theorem}\label{th:geometric_characterization}
Let $M$ be a compact Lorentzian manifold and $G$ a connected closed non-compact Lie subgroup of the isometry group $\textnormal{Isom}(M)$. According to Theorem~\ref{th:algebraic_classification}, its Lie algebra is a direct sum $\mathfrak{g}=\mathfrak{k}\oplus\mathfrak{a}\oplus\mathfrak{s}$ as described in the theorem.
\begin{compactenum}
\item If the induced bilinear form $\kappa$ is positive semidefinite, then $\mathfrak{s}$ is neither isomorphic to $\mathfrak{sl}_2(\mathds{R})$ nor $\mathfrak{he}_d^\lambda$. The orbits of $G$ are nowhere timelike and the kernel of $\kappa$ is either trivial or the span of a lightlike Killing vector field with geodesic orbits.

\item If $\mathfrak{s}\cong\mathfrak{sl}_2(\mathds{R})$, $M$ is covered isometrically by a warped product $N \times_\sigma \widetilde{\textnormal{SL}_2(\mathds{R})}$ of a Riemannian manifold $N=(N,h)$ and the universal cover $\widetilde{\textnormal{SL}_2(\mathds{R})}$ of $\textnormal{SL}_2(\mathds{R})$ furnished with the bi-invariant metric given by the Killing form $k$ of $\mathfrak{sl}_2(\mathds{R})$, that is, $M$ is covered by the manifold $N \times \widetilde{\textnormal{SL}_2(\mathds{R})}$ with the metric $g_{(x,\cdot)}=h_x \times \left(\sigma^2\left(x\right)k\right)$, $\sigma:N \to \mathds{R}_+$ smooth.

Moreover, there is a discrete subgroup $\Gamma$ in $\textnormal{Isom}(N)\times\widetilde{\textnormal{SL}_2(\mathds{R})}$ acting freely on $N \times_\sigma \widetilde{\textnormal{SL}_2(\mathds{R})}$, such that $M\cong\Gamma \backslash\mkern-5mu \left(N\times_\sigma \widetilde{\textnormal{SL}_2(\mathds{R})}\right)$.

\item Let $\mathfrak{s}\cong\mathfrak{he}_d^\lambda$ and $S$ be the subgroup generated by $\mathfrak{s}$. $S$ is isomorphic to $\textnormal{He}_d^\lambda$ or $\overline{\textnormal{He}_d^\lambda}$ and the center $Z(S)$ is isomorphic to $\mathds{S}^1$ or $\mathds{R}$, if $S$ is closed in $\textnormal{Isom}(M)$ or not, respectively.

Consider the space $\mathcal{M}$ of ad-invariant Lorentz forms on $\mathfrak{s}$. Then there exist a Riemannian manifold $N=(N,h)$ with a locally free and isometric $Z(S)$-action and a smooth map $m: N \to \mathcal{M}$ invariant under the $Z(S)$-action, such that $M$ is covered isometrically by the Lorentzian manifold $S \times_{Z(S)} N$, constructed in the following way:

Consider the product $S\times N$ furnished with the metric $g_{(\cdot,x)}=m(x) \times h_x$, where the $S$-factor is provided with the bi-invariant metric defined by $m(x)$. Let $\mathcal{O}$ be the distribution of $N$ orthogonal to the $Z(S)$-orbits and $\widetilde{\mathcal{S}}$ be the distribution of $S$ given by the tangent spaces.

The center $Z(S)$ acts isometrically on $S$ by multiplication of the inverse. The action of a central element $z$, which maps $f$ to $fz^{-1}=z^{-1}f$, corresponds to a translation in the center component (remember the matrix model of $\widetilde{\textnormal{He}_d}$ in Section~\ref{sec:Heisenberg}). Thus, we have a locally free and isometric action of $Z(S)$ on $S\times N$ by the diagonal action. Factorizing through this action, we obtain the quotient space $S\times_{Z(S)} N$. This is a manifold and we can provide it with the metric given by projection of the induced metric on $\widetilde{\mathcal{S}}\oplus\mathcal{O}$.

Furthermore, there is a discrete subgroup $\Gamma$ in $S\times_{Z(S)}\textnormal{Isom}_{Z(S)}(N)$, the quotient group defined in the same way as the quotient manifold above, $\textnormal{Isom}_{Z(S)}(N)$ being the group of $Z(S)$-equivariant isometries of $N$, acting freely on the space $S\times_{Z(S)}N$, such that $M\cong\Gamma \backslash \mkern-5mu  \left(S\times_{Z(S)}N\right)$.
\end{compactenum}
\end{theorem}

\begin{remark}
If $S$ is not closed in $\textnormal{Isom}(M)$, $S \cong \overline{\textnormal{He}_d^\lambda}$ and $Z(S)\cong \mathds{R}$. We obtain $S \times_{Z(S)} N\cong \overline{\textnormal{He}_d^\lambda} \times_\mathds{R} N$ in this case.

If $S$ is closed in $\textnormal{Isom}(M)$, $S \cong \textnormal{He}_d^\lambda$ and $Z(S)\cong \mathds{S}^1$. In this case, we obtain $S \times_{Z(S)} N\cong \textnormal{He}_d^\lambda \times_{\mathds{S}^1} N$. Since we can extend the $\mathds{S}^1$-action on $N$ to an $\mathds{R}$-action by considering the covering map $\mathds{R} \to \mathds{S}^1$, we easily obtain $\textnormal{He}_d^\lambda \times_{\mathds{S}^1} N \cong \overline{\textnormal{He}_d^\lambda} \times_\mathds{R} N$ also in this case.
\end{remark}

\section{Theorems in the homogeneous case}\label{sec:hom_th}

\begin{theorem}\label{th:homogeneous_characterization}
Let $M$ be a compact homogeneous Lorentzian manifold and denote $\textnormal{Isom}^0(M)$ the connected component of the identity in the isometry group. Suppose $\textnormal{Isom}^0(M)$ is not compact and let $\mathfrak{isom}(M)=\mathfrak{k}\oplus\mathfrak{a}\oplus\mathfrak{s}$ be the decomposition of its Lie algebra according to Theorem~\ref{th:algebraic_classification}.

Then either $\mathfrak{s}\cong\mathfrak{sl}_2(\mathds{R})$ or $\mathfrak{s}\cong\mathfrak{he}_d^\lambda$.
\begin{compactenum}
\item If $\mathfrak{s}\cong\mathfrak{sl}_2(\mathds{R})$, $M$ is covered isometrically by the metric product $N \times \widetilde{\textnormal{SL}_2(\mathds{R})}$ of a compact homogeneous Riemannian manifold $N=(N,h)$ and the universal cover $\widetilde{\textnormal{SL}_2(\mathds{R})}$ of $\textnormal{SL}_2(\mathds{R})$ furnished with the metric given by a positive multiple of the Killing form $k$ of $\mathfrak{sl}_2(\mathds{R})$.

Furthermore, there is a uniform lattice $\Gamma_0$ in $\widetilde{\textnormal{SL}_2(\mathds{R})}$ and a group homomorphism $\varrho:\Gamma_0 \to \textnormal{Isom}(N)$, such that $M\cong\Gamma \backslash \mkern-5mu \left(N \times \widetilde{\textnormal{SL}_2(\mathds{R})}\right)$, where $\Gamma$ is the graph of $\varrho$. Also, the centralizer of $\Gamma$ in $\textnormal{Isom}(N \times \widetilde{\textnormal{SL}_2(\mathds{R})})$ acts transitively on $N \times \widetilde{\textnormal{SL}_2(\mathds{R})}$.

$\textnormal{Isom}^0(M)$ is isomorphic to a central quotient of $C \times \widetilde{\textnormal{SL}_2(\mathds{R})}$, where $C$ is the connected component of the identity in the centralizer of $\varrho(\Gamma_0)$ in $\textnormal{Isom}(N)$. $C$ acts transitively on $N$.

\item Let $\mathfrak{s}\cong\mathfrak{he}_d^\lambda$ and $S$ be the subgroup generated by $\mathfrak{s}$. $S$ is isomorphic to $\textnormal{He}_d^\lambda$ or $\overline{\textnormal{He}_d^\lambda}$ and the center $Z(S)$ is isomorphic to $\mathds{S}^1$ or $\mathds{R}$, if $S$ is closed in $\textnormal{Isom}^0(M)$ or not, respectively.

Then there exist a compact homogeneous Riemannian manifold $N=(N,h)$ with a locally free and isometric $Z(S)$-action and an ad-invariant Lorentz form $m$ on $\mathfrak{s}$, such that $M$ is covered isometrically by the Lorentzian manifold $S\times_{Z(S)} N$, constructed in the following way:

Consider the metric product $S\times N$, where $S$ is furnished with the bi-invariant metric defined by $m$. Let $\mathcal{O}$ be the distribution of $N$ orthogonal to the $Z(S)$-orbits and $\widetilde{\mathcal{S}}$ be the distribution of $S$ given by the tangent spaces.

The center $Z(S)$ acts isometrically on $S$ by multiplication of the inverse. Thus, we have a locally free and isometric action of $Z(S)$ on $S\times N$ by the diagonal action. Factorizing through this action, we obtain the quotient space $S\times_{Z(S)} N$ provided with the metric given by projection of the induced metric on $\widetilde{\mathcal{S}}\oplus\mathcal{O}$.

Moreover, there is a discrete subgroup $\Gamma$ in $S \times_{Z(S)} \textnormal{Isom}_{Z(S)}(N)$, where the latter group is constructed similarly to the space $S\times_{Z(S)} N$ and $\textnormal{Isom}_{Z(S)}(N)$ denotes the group of $Z(S)$-equivariant isometries of $N$, such that $M$ is isometric to $\Gamma \backslash \mkern-5mu \left(S\times_{Z(S)}N\right)$ and $\Gamma$ projects isomorphically to a uniform lattice $\Gamma_0$ in $S/Z(S)$. Also, the centralizer of $\Gamma$ in $\textnormal{Isom}(S \times_{Z(S)} N)$ acts transitively on $S \times_{Z(S)} N$.

$\textnormal{Isom}^0(M)$ is isomorphic to a central quotient of $S \times_{Z(S)} C$, where $C$ is the connected component of the identity in the centralizer of the projection of $\Gamma$ to $Z(S) \cdot \textnormal{Isom}_{Z(S)}(N)$ in $\textnormal{Isom}_{Z(S)}(N)$. $C$ acts transitively on $N$.
\end{compactenum}
\end{theorem}

\begin{remark}
Note that in Part~(ii) of Theorem~\ref{th:homogeneous_characterization}, $S/Z(S)$ is isomorphic to a semidirect product $\mathds{S}^1 \ltimes \mathds{R}^{2d}$, where the $\mathds{S}^1$-factor acts through non-trivial rotation on the $d$ $\mathds{R}^2$-components. The nilradical of this group is the subgroup $\mathds{R}^{2d}$. Therefore, if $\Gamma$ is a lattice in $S/Z(S)$, so is $\Gamma \cap \mathds{R}^{2d}$ a lattice in $\mathds{R}^{2d}$ (cf.~\cite{OnVin00}, Part~I, Chapter~2, Theorem~3.6). It follows that lattices in $S/Z(S)$ and $\mathds{R}^{2d}$ are equivalent up to finite index. Any lattice in the additive group $\mathds{R}^{2d}$ is uniform and has the form $\mathds{Z}v_1+\ldots+\mathds{Z}v_{2d}$, where $\left\{v_1,\ldots,v_{2d}\right\}$ form a basis (cf.~\cite{OnVin00}, Part~I, Chapter~2, Theorem~1.1 and Corollary~1.2). Especially, $\Gamma$ is almost abelian in the sense that a subgroup of finite index is abelian.
\end{remark}

\begin{remark}
The construction given in Part~(ii) of Theorem~\ref{th:homogeneous_characterization} generalizes the example $\textnormal{He}_d^\lambda/\Lambda$, $\Lambda$ being a uniform lattice in $\textnormal{He}_d^\lambda$, given at the end of Section~\ref{sec:twisted_Heisenberg}. Simply choose $N=\mathds{S}^1 \subset \mathds{R}^2$ with the induced metric ($\mathds{R}^2$ is provided with the Euclidean metric) and the $\mathds{S}^1$-action given by rotation and choose $m$ to be any ad-invariant Lorentz form on $\mathfrak{he}_d^\lambda$. Then $\textnormal{He}_d^\lambda\times_{\mathds{S}^1} N\cong\textnormal{He}_d^\lambda$ and we can choose $\Gamma$ to be a lattice in $\textnormal{He}_d^\lambda$. Note that because of the bi-invariance of the induced metric on $\textnormal{He}_d^\lambda$, we can consider as well left and right quotients.
\end{remark}

\begin{example}
Consider the three-dimensional sphere $\mathds{S}^3 \subset \mathds{C}^2$ with the standard metric. This is our Riemannian manifold $N$. An element $z \in \mathds{S}^1$, $\mathds{S}^1 \subset \mathds{C}$, acts isometrically via $z \cdot (w_0,w_1) \mapsto (zw_0,zw_1)$. This action is also known as the Hopf fibration (cf.~\cite{Ba09}, example~2.7).

Now \[I:=\textnormal{Isom}_{\mathds{S}^1}(\mathds{S}^3)=\textnormal{U}(2)\cong \mathds{S}^1 \times \textnormal{SU}(2) \cong \mathds{S}^1 \times \mathds{S}^3.\]One may identify $\varphi \in \mathds{S}^1$ with the matrix \[\begin{pmatrix} \exp{(\textnormal{i}\frac{\varphi}{2})} & 0 \\ 0 & \exp{(\textnormal{i}\frac{\varphi}{2})} \end{pmatrix},\] which lies in the center.

We choose a uniform lattice $\Gamma^\prime$ in $\textnormal{He}_d \subset \textnormal{He}_d^\lambda$, such that $\Gamma^\prime$ intersects the center isomorphic to $\mathds{S}^1$ trivially. Since any lattice $\Lambda$ of $\textnormal{He}_d$ intersects the compact center in finitely many points, the statement is always true for some sublattice of finite index in $\Lambda$.

Let $\varrho:\Gamma^\prime \to I$ be a homomorphism with values in $\mathds{S}^1 \subset \mathds{S}^1 \times \mathds{S}^3$. Especially, we guarantee that the centralizer of $\varrho(\Gamma^\prime)$ is equal to the whole $\textnormal{U}(2)$ and therefore acts transitively on $\mathds{S}^3$. Denote by $\Gamma$ the projection of $\Gamma_{\varrho}\subset\textnormal{He}_d^\lambda \times I$, the graph of $\varrho$, to $\textnormal{He}_d^\lambda \times_{\mathds{S}^1} I$.

$\Gamma_{\varrho}$ projects isomorphically to $\Gamma^\prime$ in $\textnormal{He}_d^\lambda$, and $\Gamma$ projects isomorphically to a lattice $\Gamma_0$ of $\textnormal{He}_d^\lambda/Z(\textnormal{He}_d^\lambda)$, because $\Gamma^\prime \cap Z(\textnormal{He}_d^\lambda)$ consists only of the identity. By construction, $\Gamma_0$ lies central in $\textnormal{He}_d^\lambda/Z(\textnormal{He}_d^\lambda)\cong\mathds{S}^1 \ltimes \mathds{R}^{2d}$. It easily follows that the centralizer of $\Gamma$ in $\textnormal{He}_d^\lambda \times_{\mathds{S}^1} I$ acts transitively on $\textnormal{He}_d^\lambda \times_{\mathds{S}^1} N$.

Finally, the resulting manifold $M=\Gamma \backslash \mkern-5mu\left( \textnormal{He}_d^\lambda\times_{\mathds{S}^1}N\right)$ is a compact homogeneous Lorentzian manifold with $\textnormal{He}_d^\lambda \subset \textnormal{Isom}(M)$. To see that $M$ is compact, one uses the obvious fact that $\Gamma_{\varrho} \backslash \mkern-5mu \left(\textnormal{He}_d^\lambda \times N\right)$ is compact.
\end{example}

\begin{theorem}\label{th:homogeneous_reductive}
Let $M$ be a compact homogeneous Lorentzian manifold and denote $G:=\textnormal{Isom}^0(M)$ the connected component of the identity in the isometry group. $G$ acts transitively on $M$ and the connected component of the isotropy group $H \subseteq G$ of some point $x \in M$ is compact.

Let $\mathfrak{h}\subseteq \mathfrak{g}$ denote the Lie algebras of both Lie groups. Then $M\cong G/H$ is reductive, that is, there is an $\textnormal{Ad}(H)$-invariant vector space $\mathfrak{m}\subseteq\mathfrak{g}$ (not necessarily a subalgebra) that is complementary to $\mathfrak{h}$.
\end{theorem}

In the case that the isometry group $\textnormal{Isom}(M)$ has non-compact connected components, we will describe the local geometry of $M$ in detail in Section~\ref{sec:homogeneous_geometry}. Two of our results are:

\begin{theorem}\label{th:isotropy_representation}
Let $M=(M,g)$ be a compact homogeneous Lorentzian manifold and denote $G:=\textnormal{Isom}^0(M)$ the connected component of the identity in the isometry group. Let $H$ be the isotropy group in $G$ of some $x \in M$. Clearly, $M \cong G/H$.
\begin{compactenum}
\item Suppose the Lie algebra of $G$ contains a direct summand $\mathfrak{s}$ isomorphic to $\mathfrak{sl}_2(\mathds{R})$. Then the isotropy representation of $G/H$ allows a decomposition into irreducible invariant subspaces, such that $\mathfrak{s}$ appears as an irreducible summand.

\item Suppose the Lie algebra of $G$ contains a direct summand $\mathfrak{s}$ isomorphic to $\mathfrak{he}_d^\lambda$, $\lambda \in \mathds{Z}_+^d$. Then the isotropy representation of $G/H$ allows a decomposition into weakly irreducible invariant subspaces, such that $\mathfrak{s}$ appears as a weakly irreducible summand. $\mathfrak{s}$ is not irreducible and cannot be decomposed into a direct sum of irreducible subspaces.
\end{compactenum}
\end{theorem}

\begin{theorem}\label{th:homogeneous_not_Ricci_flat}
Let $M$ be a compact homogeneous Lorentzian manifold that is Ricci-flat. Then the connected components of its isometry group $\textnormal{Isom}(M)$ are compact.
\end{theorem}

\begin{corollary}\label{cor:Ricci_flat}
Let $M$ be a compact homogeneous Lorentzian manifold that is Ricci-flat, but not flat. Then its isometry group $\textnormal{Isom}(M)$ is compact.
\end{corollary}
\vspace{-1.4em}\begin{proof}
Assume that $\textnormal{Isom}(M)$ is not compact. It follows from Theorem~\ref{th:homogeneous_not_Ricci_flat}, that $\textnormal{Isom}(M)$ has infinitely many connected components. Corollary~3 of \cite{PiZe10} states, that any compact Lorentzian manifold whose isometry group has infinitely many connected components, possesses an everywhere timelike Killing vector field if it possesses a somewhere timelike Killing vector field. Since any homogeneous Lorentzian manifold has a somewhere timelike Killing vector field, it follows that there exists an everywhere timelike Killing vector field on $M$. Theorem~3.2 of \cite{RoSa96} yields that any Ricci-flat compact homogeneous Lorentzian manifold admitting a timelike Killing vector field is isometric to a flat torus (up to a finite covering).
\end{proof}\vspace{0pt}
\begin{remark}
Note that in Paragraph~1 of \cite{PiZe10} examples of flat compact homogeneous Lorentzian manifolds are given, whose isometry groups are not compact. These manifolds are flat tori provided with the metric defined by certain quadratic forms of Minkowski space.
\end{remark}

%%%%%%%%%%%%%%%%%%%%%%%%%%%%%%%%%%%%%%%%%%%%%%%%%%%%%%%%%%%%%%%%%%%%%%%%%%%%%%%%%%%%%%%%%%%%

\chapter{Algebraic classification of the Lie algebras}\label{ch:algebra}

The aim of this chapter is the proof of Theorems~\ref{th:algebraic} and~\ref{th:algebraic_classification}. In Section~\ref{sec:bilinear_form}, we will show some useful elementary results about ad-invariant symmetric bilinear forms on Lie algebras and give a powerful characterization of elements of a Lie algebra generating precompact one-parameter groups in the corresponding Lie group.

In the subsequent sections, we continue with proving Theorem~\ref{th:algebraic}. Apart from Sections~\ref{sec:bilinear_form} and~\ref{sec:general_subgroups}, we always consider a connected non-compact Lie group $G$ and an ad-invariant symmetric bilinear form $\kappa$ on its Lie algebra $\mathfrak{g}$ that fulfills condition~\hyperlink{star}{$(\star)$} as described in the theorem.

The outline of the proof is the following: In Sections~\ref{sec:nilradical} and~\ref{sec:radical}, we first collect some results about the nilradical and the radical of a Lie algebra furnished with a form $\kappa$ as in Theorem~\ref{th:algebraic}. If the radical of the group is compact, what we assume in Section~\ref{sec:compact_radical}, we are in case~(v) of the theorem and $\kappa$ is Lorentzian. In Section~\ref{sec:noncompact_radical} we will show the other cases, which split as follows: Section~\ref{sec:kappa_Lorentz} deals with the case that $\kappa$ is not positive semidefinite. Then $\kappa$ is Lorentzian as well and we are in case~(iv). Otherwise $\kappa$ is positive semidefinite and its kernel has dimension at most one. We deal with this in Section~\ref{sec:kappa_positive}. If the radical is nilpotent as assumed in Section~\ref{sec:nilpotent_radical}, we are in the cases~(i) and~(iii) of the theorem. Else, we show in Section~\ref{sec:nonnilpotent_radical}, that we are in case~(ii).

The following chart summarizes the outline of the proof. We start with a Levi decomposition $\mathfrak{g}=\mathfrak{l} \inplus \mathfrak{r}$ according to Lemma~\ref{lem:Levi_decomposition}, where $\mathfrak{l}$ is semisimple and $\mathfrak{r}$ is the radical of $\mathfrak{g}$. Denote by $R$ the radical of the group $G$ (the largest connected solvable normal subgroup), its Lie algebra is $\mathfrak{r}$.

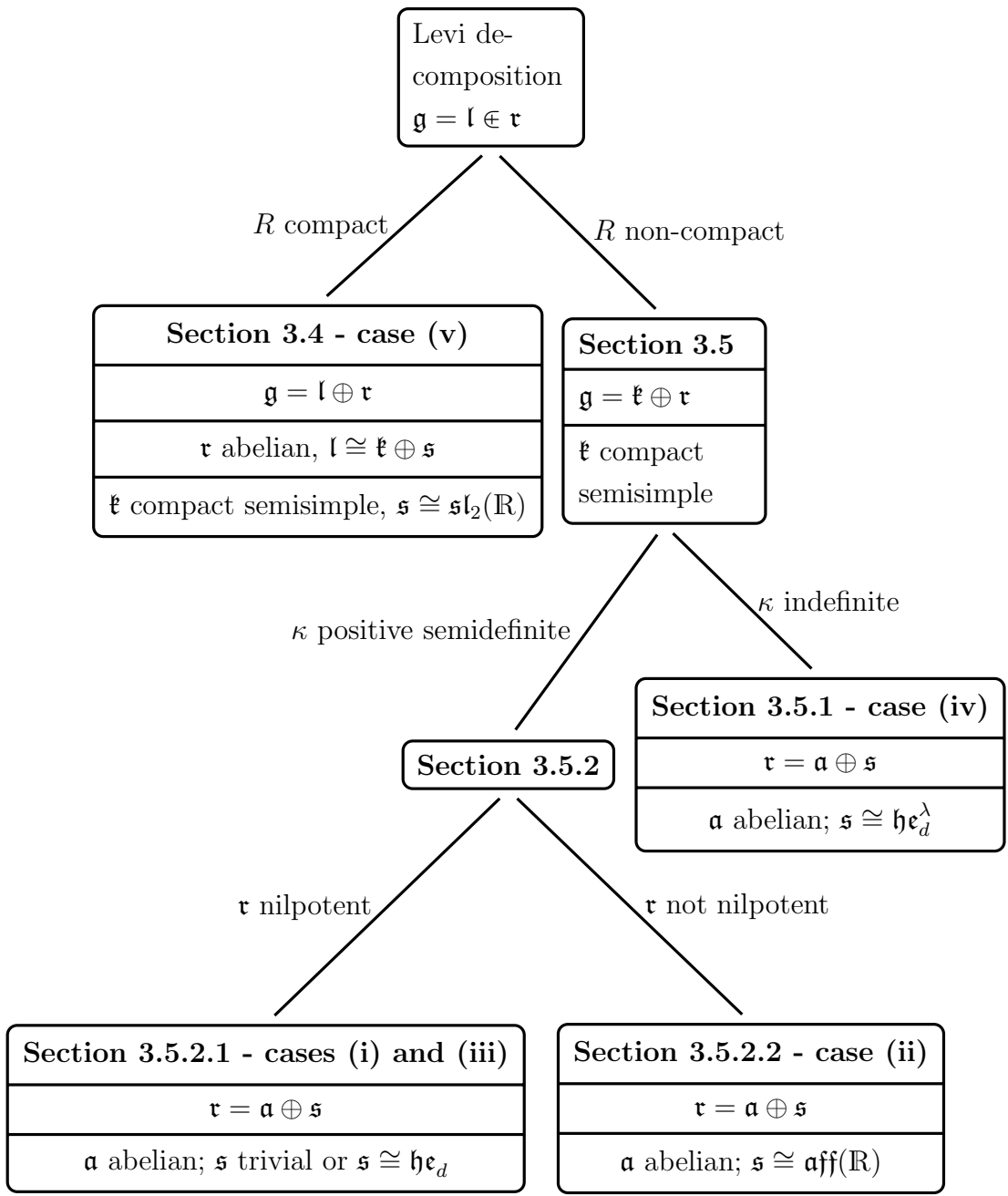
\begin{tikzpicture}[
    level 1/.style={sibling distance=5.0cm,level distance=5.0cm},
    level 2/.style={sibling distance=4.5cm, level distance=5.0cm},
    level 3/.style={sibling distance=7.0cm, level distance=5.0cm},
    edge from parent/.style={very thick,draw=black,
        shorten >=5pt, shorten <=5pt},
    edge from parent path={(\tikzparentnode.south) -- (\tikzchildnode.north)},
    kant/.style={text centered},
    every node/.style={text ragged, inner sep=2mm},
    punkt/.style={rectangle, rounded corners, top color=white,
    draw=black, very thick}
    ]

\node[punkt, text width=5.5em]{
         Levi decomposition $\mathfrak{g}=\mathfrak{l} \inplus \mathfrak{r}$              
         }
    child {
        node[punkt] [rectangle split, rectangle split, rectangle split parts=4,
         text ragged] { \textbf{Section~\ref{sec:compact_radical} - case~(v)}
         		\nodepart{second}
            $\mathfrak{g}=\mathfrak{l}\oplus\mathfrak{r}$
            \nodepart{third}
            $\mathfrak{r} \text{ abelian}$, $\mathfrak{l}\cong\mathfrak{k}\oplus\mathfrak{s}$
            \nodepart{fourth}
            $\mathfrak{k}$ compact semisimple, $ \mathfrak{s}\cong\mathfrak{sl}_2(\mathds{R})$         }
        edge from parent
            node[kant, left] {$R$ compact}
    }
    child {
        node[punkt, text width=6em] [rectangle split, rectangle split, rectangle split parts=3,
         text ragged] { \textbf{Section~\ref{sec:noncompact_radical}}
         \nodepart{second}
         $\mathfrak{g}=\mathfrak{k}\oplus\mathfrak{r}$
         \nodepart{third}
         $\mathfrak{k}$ compact semisimple
         }
        child {
            node [punkt]{\textbf{Section~\ref{sec:kappa_positive}}}
            child {
            node [punkt,rectangle split, rectangle split,
            rectangle split parts=3] {
                \textbf{Section~\ref{sec:nilpotent_radical} - cases~(i) and~(iii)}
                \nodepart{second}
                $\mathfrak{r}=\mathfrak{a}\oplus\mathfrak{s}$
                \nodepart{third}
                $\mathfrak{a}$ abelian; $\mathfrak{s}$ trivial or $\mathfrak{s}\cong\mathfrak{he}_d$
                }
            edge from parent
                node[kant, left] {$\mathfrak{r}$ nilpotent}
            }
            child {
            node [punkt,rectangle split, rectangle split,
            rectangle split parts=3] {
                \textbf{Section~\ref{sec:nonnilpotent_radical} - case~(ii)}
                \nodepart{second}
                $\mathfrak{r}=\mathfrak{a}\oplus\mathfrak{s}$
                \nodepart{third}
                $\mathfrak{a}$ abelian; $\mathfrak{s}\cong\mathfrak{aff}(\mathds{R})$
                }
            edge from parent
                node[kant, right] {$\mathfrak{r}$ not nilpotent}
            }
            edge from parent
                node[kant, left] {$\kappa$ positive semidefinite}
            }
        child {
            node [punkt,rectangle split, rectangle split,
            rectangle split parts=3] {
                \textbf{Section~\ref{sec:kappa_Lorentz} - case~(iv)}
                \nodepart{second}
                $\mathfrak{r}=\mathfrak{a}\oplus\mathfrak{s}$
                \nodepart{third}
                $\mathfrak{a}$ abelian; $\mathfrak{s}\cong\mathfrak{he}_d^\lambda$
            }
            edge from parent
                node[kant, right] {$\kappa$ indefinite}
        }
            edge from parent{
                node[kant, right] {$R$ non-compact}}
     };
\end{tikzpicture}

To close this chapter, we consider general (not necessarily closed) connected Lie subgroups $G$ of the isometry group $\textnormal{Isom}(M)$ of a Lorentzian manifold $M$ of finite volume, and finally prove Theorem~\ref{th:algebraic_classification} in Section~\ref{sec:general_subgroups}. For this proof, we may suppose that $G$ is a subgroup of $\textnormal{Isom}(M)$. The result follows from Corollaries~\ref{cor:subalgebra_compact} and~\ref{cor:compact_algebra}, if $G$ is precompact, and from Corollary~\ref{cor:condition_star} and Theorem~\ref{th:algebraic}, if $G$ is a non-compact closed subgroup of $\textnormal{Isom}(M)$. If $G$ is neither precompact nor closed in $\textnormal{Isom}(M)$, we consider the closure of $G$ in $\textnormal{Isom}(M)$ and use the results of Theorem~\ref{th:algebraic}. One way to show the theorem is to argue analogously as in the proof of Theorem~\ref{th:algebraic}; basically, one replaces the condition~\hyperlink{star}{$(\star)$} by the properties of $\kappa$ given in the theorem. We will give a shorter proof using the Maltsev closure. Also the second part of Theorem~\ref{th:algebraic_classification} can be shown in a similar way as the corresponding statement of Theorem~\ref{th:algebraic}.

Note that we have to consider non-closed subgroups $G$ of the isometry group separately, since in this case, a non-precompact one-parameter group in $G$ can be precompact in the isometry group. Therefore, condition~\hyperlink{star}{$(\star)$} does not necessarily hold for the induced bilinear form $\kappa$ on the Lie algebra of $G$.

\section{Symmetric bilinear forms on Lie algebras}\label{sec:bilinear_form}

\begin{lemma}\label{lem:central_homomorphism}
Let $\mathfrak{l}\subseteq\mathfrak{g}$ be a subalgebra of the Lie algebra $\mathfrak{g}$ and $Y \in \mathfrak{l}$ such that $[Y,\mathfrak{l}]=\left\{0\right\}$. Then for any ad-invariant scalar product $\kappa$, the mapping $\mathfrak{l} \to \mathds{R}$, $X \mapsto \kappa(X,Y)$, is a Lie algebra homomorphism, that is, $\kappa([X,X^\prime],Y)=0$ for all $X,X^\prime \in \mathfrak{l}$.
\end{lemma}
\vspace{-1.4em}\begin{proof}
This follows from $\kappa([X,X^\prime],Y)=\kappa(X,[X^\prime,Y])=0$.
\end{proof}\vspace{0pt}

\begin{lemma}\label{lem:positive_definite}
Let $\mathfrak{g}$ be a Lie algebra with ad-invariant symmetric bilinear form $\kappa$.
\begin{compactenum}
\item The kernel of $\kappa$ is an ideal in $\mathfrak{g}$.

\item If $\mathfrak{g}$ is compact and simple, $\kappa$ is a multiple of the negative definite Killing form of $\mathfrak{g}$.

\item If $\kappa$ is positive definite, $\mathfrak{g}$ is the $\kappa$-orthogonal direct sum $\mathfrak{g}=[\mathfrak{g},\mathfrak{g}]\oplus\mathfrak{z(g)}$ with $[\mathfrak{g},\mathfrak{g}]$ trivial or compact semisimple.
\end{compactenum}
\end{lemma}
\vspace{-1.4em}\begin{proof}
(i) Let $\mathfrak{i}$ denote the kernel of $\kappa$ and let $X \in \mathfrak{i}$ and $Y,Z \in \mathfrak{g}$ be arbitrary. Because of the ad-invariance, $\kappa([X,Y],Z)=\kappa(X,[Y,Z])=0$. Thus, $[X,Y] \in \mathfrak{i}$.

(ii) Let $k$ be the Killing form of $\mathfrak{g}$. Since $\mathfrak{g}$ is compact and simple, $k$ is negative definite by Corollary~\ref{cor:compact_semisimple}. Now choose a $(-k)$-orthonormal basis which is also a $\kappa$-orthogonal basis (principal axis transformation). Obviously, there is $t \in \mathds{R}$ such that $\kappa-t k$ degenerates. The kernel of this ad-invariant form is a non-trivial ideal in $\mathfrak{g}$ and thus equal to $\mathfrak{g}$ since the algebra is simple. Therefore, $\kappa=t k$.

(iii) By Proposition~\ref{prop:compact_algebra}, $\mathfrak{g}$ is reductive. Using Proposition~\ref{prop:reductive_algebra} and Corollary~\ref{cor:subalgebra_compact}, we obtain all but the $\kappa$-orthogonality of the decomposition. Lemma~\ref{lem:central_homomorphism} with $\mathfrak{l}=\mathfrak{g}$ and $Y \in \mathfrak{z(g)}$ yields the orthogonality.
\end{proof}\vspace{0pt}

The next lemma characterizes elements of the Lie algebra generating a precompact one-parameter group in the corresponding Lie group in terms of the adjoint representation of the algebra.

\begin{lemma}\label{lem:precompact_ad}
Let $G$ be a Lie group, $\mathfrak{g}$ its Lie algebra and $X \in \mathfrak{g}$ an element generating a precompact one-parameter group in $G$. Then $\textnormal{ad}_X$ is either trivial or not nilpotent. Furthermore, there is no $0 \neq Y \in \mathfrak{g}$ such that $[X,Y]=\lambda Y$ for some non-zero $\lambda \in \mathds{R}$, or $[X,Y]\neq 0$ and $\textnormal{ad}_X^k(Y)=0$ for some integer $k>0$.

If additionally $\mathfrak{g}=\mathfrak{l}\oplus\mathfrak{m}$ is a direct sum of algebras with $\mathfrak{l}$ semisimple and $X \in \mathfrak{l}$, $\textnormal{ad}_X$ is semisimple (that is, diagonalizable over the complex numbers) and its eigenvalues are all purely imaginary.
\end{lemma}
\vspace{-1.4em}\begin{proof}
We first consider the case that $\mathfrak{g}=\mathfrak{l}\oplus\mathfrak{m}$ with $\mathfrak{l}$ semisimple and $X \in \mathfrak{l}$.

Let $X=S+N$, where $S,N \in \mathfrak{l}$, $[S,N]=0$, $\textnormal{ad}_S$ semisimple and $\textnormal{ad}_N$ nilpotent, be the abstract Jordan decomposition of $X$ (cf.~\cite{Hu72}, Paragraph~5.4). Since $[S,N]=0$, for any real number $t$,
\[\textnormal{Ad}_{\exp(t X)}=\exp(\textnormal{ad}_{t X})=\exp(\textnormal{ad}_{t S}+\textnormal{ad}_{t N})=\exp(t\textnormal{ad}_{S})\exp(t \textnormal{ad}_{N}).\]

Note that $\textnormal{ad}_{S}$ and $\textnormal{ad}_{N}$ commute, since $S$ and $N$ do and ad: $\mathfrak{g} \to \mathfrak{gl(g)}$ is a Lie algebra homomorphism.

Using that Ad: $G \to \textnormal{GL}(\mathfrak{g})$ is a Lie group homomorphism, we obtain that $\left\{\exp(t\textnormal{ad}_{S})\exp(t \textnormal{ad}_{N})\right\}_{t\in\mathds{R}}$ is precompact in $\textnormal{GL}(\mathfrak{g}) \cong \mathds{R}^{k^2} \subset \mathds{C}^{k^2}$, where $k$ is the dimension of $G$.

Let $\lambda_1,\ldots,\lambda_{k^\prime}$ be the eigenvalues of $\textnormal{ad}_X|_{\mathfrak{l}}$, where $k^\prime$ is the dimension of $\mathfrak{l}$. Using $[\mathfrak{l},\mathfrak{m}]=\left\{0\right\}$, $\textnormal{ad}_S+\textnormal{ad}_N$ is the Jordan decomposition of $\textnormal{ad}_X$ in $\mathds{C}^{k^2}$. It follows that $\left\{\exp(t\textnormal{ad}_{S})\exp(t \textnormal{ad}_{N})\right\}_{t\in\mathds{R}}$ is conjugate to the one-parameter group \[\left\{\textnormal{diag}(\exp(t\lambda_1), \ldots, \exp(t\lambda_{k^\prime}),0,\ldots,0) (1+t\textnormal{ad}_N+\ldots+\frac{t^{k^\prime-1}}{(k^\prime-1)!}\textnormal{ad}_N^{k^\prime-1})\right\}_{t\in\mathds{R}}.\]

Since the last group is precompact in $\mathds{C}^{k^2}$, it is especially bounded. It easily follows that $\textnormal{Re}(\lambda_j)=0$ for all $j$ and $\textnormal{ad}_{N}=0$.

In the general case, assume that $\textnormal{ad}_X$ is nilpotent, but not trivial. As above, the one-parameter group $\left\{\exp(t \textnormal{ad}_{X})\right\}_{t\in\mathds{R}}$ is precompact in $\textnormal{GL}(\mathfrak{g}) \cong \mathds{R}^{k^2}$. On the other hand, \[\exp(t \textnormal{ad}_{X})=(1+t\textnormal{ad}_X+\ldots+\frac{t^{k-1}}{(k-1)!}\textnormal{ad}_X^{k-1}),\] contradiction.

If $0 \neq Y \in \mathfrak{g}$ and $\lambda \neq 0$ such that $[X,Y]=\lambda Y$, then \[\textnormal{Ad}_{\exp(tX)}(Y)=\exp(t\lambda)Y,\] contradicting the precompactness of $\left\{\exp(t \textnormal{ad}_{X})\right\}_{t\in\mathds{R}}$.

Let $Y \in \mathfrak{g}$ such that $[X,Y]\neq 0$ and $\textnormal{ad}_X^k(Y)=0$ for some $k>0$. Without loss of generality, $k=2$. But then \[\textnormal{Ad}_{\exp(tX)}(Y)=Y+t[X,Y]\] contradicts the precompactness of $\left\{\exp(t \textnormal{ad}_{X})\right\}_{t\in\mathds{R}}$.
\end{proof}\vspace{0pt}

\section{Nilradical}\label{sec:nilradical}

In this section, we investigate the nilradical of a connected non-compact Lie group $G$, provided with an ad-invariant symmetric bilinear form $\kappa$ on its Lie algebra $\mathfrak{g}$ that fulfills condition~\hyperlink{star}{$(\star)$} as described in Theorem~\ref{th:algebraic}.

\begin{lemma}\label{lem:kappa_abelian}
Let $\mathfrak{p} \subseteq \mathfrak{g}$ be an abelian subalgebra containing an element $X \in \mathfrak{p}$ generating a non-precompact one-parameter group. Then the set of elements in $\mathfrak{p}$ generating a non-precompact one-parameter group is dense in $\mathfrak{p}$. $\kappa|_{\mathfrak{p} \times \mathfrak{p}}$ is positive semidefinite and the dimension of its kernel is at most one.
\end{lemma}
\vspace{-1.4em}\begin{proof}
The closure of the subgroup generated by $\mathfrak{p}$ is a connected abelian subgroup and by Lemma~\ref{lem:abelian_group} isomorphic to $\mathds{T}^m \times \mathds{R}^{m^\prime}$. $m^\prime \neq 0$ since $X$ generates a non-precompact one-parameter group. Let $\mathfrak{p^*}$ denote the Lie algebra of the closure and $\mathfrak{t}\subset\mathfrak{p^*}$ the subalgebra corresponding to the torus factor. All elements of $\mathfrak{p^*}$ that are not in $\mathfrak{t}$ generate a non-precompact one-parameter group. By definition, $\mathfrak{p}$ is not contained in $\mathfrak{t}$, hence, the set of elements in $\mathfrak{p}$ generating a non-precompact one-parameter group is dense in $\mathfrak{p}$. The remainder of the assertion now follows from condition~\hyperlink{star}{$(\star)$}.
\end{proof}\vspace{0pt}

\begin{lemma}\label{lem:nilpotent_density}
Let $N$ be a connected non-compact nilpotent Lie group and $\mathfrak{n}$ its Lie algebra. Then the set of $X \in \mathfrak{n}$ generating a non-precompact one-parameter group is dense in $\mathfrak{n}$.
\end{lemma}
\vspace{-1.4em}\begin{proof}
If $N$ is abelian, $N$ is isomorphic to $\mathds{T}^m \times \mathds{R}^{m^\prime}$ by Lemma~\ref{lem:abelian_group}. $m^\prime \neq 0$ because $N$ is not compact, it follows that the set of elements in $\mathfrak{n}$ generating a non-precompact one-parameter group is dense in $\mathfrak{n}$.

Now assume that $N$ is not abelian and the statement of the lemma is false. Using Lemma~\ref{lem:precompact_ad} and the nilpotency of $\mathfrak{n}$, it follows that $\textnormal{ad}_X$ is trivial for all $X$ of some open ball in $\mathfrak{n}$. Thus, $\mathfrak{z(n)}$ contains elements of an open ball in $\mathfrak{n}$ and therefore has the same dimension as $\mathfrak{n}$. Hence, $\mathfrak{z(n)}=\mathfrak{n}$, that is, $N$ is abelian, contradiction.
\end{proof}\vspace{0pt}

The proof of the last lemma also shows the following:

\begin{corollary}\label{cor:compact_nilpotent}
A connected compact nilpotent Lie group is abelian.
\end{corollary}

In what follows, denote by $N$ the nilradical of $G$, that is, the largest (with respect to inclusion) connected normal Lie subgroup of $G$, which is nilpotent. By maximality, $N$ is closed in $G$. Its Lie algebra $\mathfrak{n}$ is the largest (with respect to inclusion) nilpotent ideal of $\mathfrak{g}$.

According to Corollary~\ref{cor:compact_nilpotent}, $N$ is abelian if it is compact. In the remainder of this section, we assume that $N$ is not compact.

It follows from Lemma~\ref{lem:nilpotent_density} and condition~\hyperlink{star}{$(\star)$}, that the restriction of $\kappa$ to $\mathfrak{n} \times \mathfrak{n}$ is positive semidefinite and its kernel $\mathfrak{i}$ is an ideal of dimension at most one.

\begin{proposition}\label{prop:decomposition_nilradical}
Let $\kappa^\prime$ be any ad-invariant positive semidefinite symmetric bilinear form on $\mathfrak{n}$ with kernel $\mathfrak{i}$ having dimension at most one.

Then $\mathfrak{n}=\mathfrak{a}\oplus\mathfrak{h}$ is the $\kappa^\prime$-orthogonal direct sum of an abelian algebra $\mathfrak{a}$ (which might be trivial) and an subalgebra $\mathfrak{h}$ being either one-dimensional or isomorphic to a Heisenberg algebra $\mathfrak{he}_d$ of dimension $2d+1$. $\mathfrak{h}$ is one-dimensional if and only if $\mathfrak{n}$ is abelian.

The adjoint representation of $G$ on $\mathfrak{n}/\mathfrak{i}$ has precompact image. If $\kappa^\prime$ is not positive definite, then $\mathfrak{i}=\mathds{R}Z$, where $\mathds{R}Z$ is the center of $\mathfrak{h}$.
\end{proposition}
\vspace{-1.4em}\begin{proof}
$\mathfrak{n}/\mathfrak{i}$ is nilpotent and possesses an ad-invariant positive definite symmetric bilinear form. By Propositions~\ref{prop:reductive_algebra} and~\ref{prop:compact_algebra}, $\mathfrak{n}/\mathfrak{i}$ is abelian, so $[\mathfrak{n},\mathfrak{n}]\subseteq \mathfrak{i}$. If $\mathfrak{i}=\left\{0\right\}$, $\mathfrak{n}$ is abelian.

In the abelian case, we choose $\mathds{R}Z$ to be the kernel of $\kappa^\prime$ if it has dimension one and if $\kappa^\prime$ is positive definite, we choose an arbitrary non-zero $Z \in \mathfrak{n}$. Then $\mathfrak{n}=\mathfrak{a}\oplus\mathfrak{h}$, where $\mathfrak{h}=\mathds{R}Z$ and $\mathfrak{a}$ is a vector space complement to $\mathfrak{h}$, which is orthogonal to $Z$ in the case that $\kappa^\prime$ is positive definite.

Suppose $\mathfrak{i}=\mathds{R}Z$ is one-dimensional. Then $Z \in \mathfrak{z(n)}$, since $\mathfrak{n}$ is nilpotent and $[\mathfrak{n},\mathfrak{n}]\subseteq \mathfrak{i}$. It follows that in any case, $\mathfrak{n}$ is at most two-step nilpotent (this actually shows the second part of Theorem~B in \cite{Zi86}). $\mathfrak{z(n)}$ is abelian, so $\mathfrak{z(n)}=\mathfrak{a}\oplus\mathds{R}Z $ with some abelian summand $\mathfrak{a}$. If $\mathfrak{n}=\mathfrak{z(n)}$, $\mathfrak{n}$ is abelian and we are done. So suppose $\mathfrak{n}$ is not abelian.

Let $\mathfrak{h}:=\mathfrak{z(n)}^\perp$ ($\kappa^\prime$-orthogonal complement in $\mathfrak{n}$). Because $\kappa^\prime$ is ad-invariant, $\mathfrak{h}$ is an subalgebra. Moreover, $\mathfrak{n}=\mathfrak{a}\oplus\mathfrak{h}$ is an orthogonal and direct sum by construction. Because $\mathfrak{n}$ is not abelian, $\dim \mathfrak{h} \geq 2$ and so \[[\mathfrak{h},\mathfrak{h}]=[\mathfrak{n},\mathfrak{n}]=\mathds{R}Z \text{ as well as } \mathfrak{z(h)}=\mathfrak{z(n)}\cap\mathfrak{h}=\mathds{R}Z.\] If $V$ is a complementary vector space of $\mathds{R}Z$ in $\mathfrak{h}$, then $\omega:V \times V \to \mathds{R}$ defined by $[X,Y]=\omega(X,Y)Z$ is an alternating bilinear form, which is non-degenerate (an element of the kernel has to lie in the center of $\mathfrak{h}$). Thus, $\mathfrak{h}$ is isomorphic to a Heisenberg algebra.

Finally, the image of the adjoint representation of $G$ on $\mathfrak{n}/\mathfrak{i}$ is precompact, since it preserves a positive definite scalar product and therefore acts by orthogonal transformations.
\end{proof}\vspace{0pt}
\begin{remark}
Note that the abelian summand $\mathfrak{a}$ is in general not canonically defined, but $\mathfrak{a}\oplus\mathds{R}Z=\mathfrak{z(n)}$ and $\mathfrak{h}=\mathfrak{z(n)}^\perp$ are.
\end{remark}

\begin{proposition}\label{prop:nilradical_results}
Suppose $\mathfrak{n}$ is not abelian.
\begin{compactenum}
\item The center $\mathfrak{i}=\mathds{R}Z$ of $\mathfrak{h}\cong\mathfrak{he}_d$ is central in $\mathfrak{g}$.

\item All non-central $X \in \mathfrak{h} \subseteq \mathfrak{n}$ generate non-precompact one-parameter groups.

\item If $Y \in \mathfrak{g}$, $Y \notin \mathds{R}Z$, commutes with a non-central element $X$ of $\mathfrak{n}$, then $\kappa(Y,Y) >0$.
\end{compactenum}
\end{proposition}
\vspace{-1.4em}\begin{proof}
(i) $\mathfrak{h}\cong\mathfrak{he}_d$ by Proposition~\ref{prop:decomposition_nilradical}. We will show that the adjoint action of $G$ on $\mathfrak{i}$ is trivial. Let $X,Y \in \mathfrak{h}$ such that $[X,Y]=Z$. By Proposition~\ref{prop:decomposition_nilradical}, the adjoint action of $G$ on $\mathfrak{n}/\mathfrak{i}$ has precompact image. Therefore, we may choose a compact set $K \subset \mathfrak{n}$ such that $\textnormal{Ad}_f(X) \in K +\mathfrak{i}$ and $\textnormal{Ad}_f(Y) \in K +\mathfrak{i}$ for all $f \in G$. It follows that \[\textnormal{Ad}_f(Z)=\textnormal{Ad}_f([X,Y])=[\textnormal{Ad}_f(X),\textnormal{Ad}_f(Y)]\in [K+\mathfrak{i},K+\mathfrak{i}] \subseteq [K,K].\]

It follows that the adjoint action of $G$ on $\mathfrak{i}$ has precompact image. But $\mathfrak{i}$ is one-dimensional and $G$ is connected, so $\textnormal{Ad}(G)$ acts trivially on ${\mathfrak{i}}$.

(ii) This follows from Lemma~\ref{lem:precompact_ad} and the fact that $\textnormal{ad}_X$ is nilpotent, but not trivial.

(iii) $\mathfrak{p}:=\textnormal{span}\left\{X,Y,Z\right\}$ is an abelian subalgebra of dimension at least 2, and by (ii), $X$ generates a non-precompact one-parameter group. By Lemma~\ref{lem:kappa_abelian}, the restriction of $\kappa$ to $\mathfrak{p} \times \mathfrak{p}$ is positive semidefinite and its kernel has dimension at most one. But the kernel is already given by $\mathds{R}Z$, hence, $\kappa(Y,Y)>0$.
\end{proof}\vspace{0pt}

\begin{proposition}\label{prop:nilradical_semisimple}
Let $\mathfrak{l}\subset\mathfrak{g}$ be a semisimple subalgebra. Then $\mathfrak{g}^\prime=\mathfrak{l}+\mathfrak{n}$ is in fact a $\kappa$-orthogonal direct sum.
\end{proposition}
\vspace{-1.4em}\begin{proof}
$\mathfrak{l} \cap \mathfrak{n}$ is a nilpotent ideal in $\mathfrak{l}$. Since $\mathfrak{l}$ is semisimple, $\mathfrak{l} \cap \mathfrak{n}=\left\{0\right\}$.

As above, let $\mathfrak{i}$ be the kernel of the restriction of $\kappa$ to $\mathfrak{n} \times \mathfrak{n}$. Due to Lemma~\ref{lem:positive_definite}, $\mathfrak{i}$ is an ideal in $\mathfrak{n}$.

Let $X \in \mathfrak{l}, Y \in \mathfrak{n}, Z \in \mathfrak{i}$. Then \[\kappa([Z,X],Y)=\kappa(Z,[X,Y])=0\] since $[\mathfrak{l},\mathfrak{n}]\subseteq\mathfrak{n}$. Thus, $[\mathfrak{l},\mathfrak{i}]\subseteq \mathfrak{i}$, that is, $\mathfrak{i}$ is an ideal in $\mathfrak{g}^\prime$.

The set $\mathfrak{m}$ of $Y \in \mathfrak{l}$ such that $[Y,\mathfrak{i}]=\left\{0\right\}$ is an ideal of $\mathfrak{l}$: Indeed, let $X,Y \in \mathfrak{l}$ such that $[Y,\mathfrak{i}]=\left\{0\right\}$. Since $[\mathfrak{l},\mathfrak{i}]\subseteq\mathfrak{i}$, \[0=[X,[Y,Z]]+[Y,[Z,X]]+[Z,[X,Y]]=[Z,[X,Y]]\] for $Z \in \mathfrak{i}$ by the Jacobi identity.

$\mathfrak{m}$ has codimension at most one, because $\mathfrak{i}$ has dimension at most one. But $\mathfrak{l}$ is semisimple, thus, it contains no ideals of codimension exactly one. Hence, $\mathfrak{l}$ centralizes $\mathfrak{i}$.

Since $\mathfrak{l}$ is semisimple, $\mathfrak{l}=[\mathfrak{l},\mathfrak{l}]$ and by Lemma~\ref{lem:central_homomorphism}, $\mathfrak{l}$ and therefore $\mathfrak{g}^\prime$ are $\kappa$-orthogonal to $\mathfrak{i}$. So we may pass to the quotient $\mathfrak{g}^\prime/\mathfrak{i}=\mathfrak{l}+\mathfrak{n}/\mathfrak{i}$ or equivalently, we can assume without loss of generality that $\kappa$ is positive definite on $\mathfrak{n}$. Note that because of Proposition~\ref{prop:decomposition_nilradical}, $\mathfrak{n}$ is abelian in this case.

We want to show $[\mathfrak{l},\mathfrak{n}]=\left\{0\right\}$. By Proposition~\ref{prop:decomposition_nilradical}, the adjoint action of $G$ on $\mathfrak{n}$ has precompact image. Thus, we may assume that $\mathfrak{l}$ is compact. Treating each simple summand of $\mathfrak{l}$ separately, we also can assume $\mathfrak{l}$ to be simple.

In the following, we construct an ad-invariant positive definite scalar product on $\mathfrak{g}^\prime$. Let $k$ be the Killing form of $\mathfrak{g}^\prime$. If $X \in \mathfrak{n},Y \in \mathfrak{g}^\prime$, then $\textnormal{ad}_X^2(Y)=0$ because $\mathfrak{n}$ is an abelian ideal. Thus, $k$ restricted to $\mathfrak{g}^\prime \times \mathfrak{n}$ is identically zero. Let $X \in \mathfrak{l},Y \in \mathfrak{n}$. If $W \in \mathfrak{n}$, \[\textnormal{ad}_X(\textnormal{ad}_Y(W))=0,\] and if  $W \in \mathfrak{l}$, \[\textnormal{ad}_X(\textnormal{ad}_Y(W))\in \mathfrak{n}.\] Thus, $k(X,Y)=\textnormal{Tr}(\textnormal{ad}_X\circ\textnormal{ad}_Y)=0$.

By Lemma~\ref{lem:positive_definite}~(ii), $k$ restricted to $\mathfrak{l} \times \mathfrak{l}$ is a multiple of the negative definite Killing form of $\mathfrak{l}$. It is a non-zero multiple, because otherwise $k$ would vanish on the whole of $\mathfrak{g}^\prime \times \mathfrak{g}^\prime$ and $\mathfrak{g}^\prime$ would be solvable (cf.~\cite{OnVin94}, Chapter~1, corollary to Theorem~2.1), contradiction.

$-k$ is a positive definite symmetric bilinear form on $\mathfrak{l}$, so we can find a $(-k)$-orthonormal basis which is $\kappa$-orthogonal (principal axis transformation). Thus, there is a $t_0 \in \mathds{R}$, such that $\kappa-t k$ is positive definite on $\mathfrak{l}$ for all $t \geq t_0$.

Obviously, the restriction of $\kappa-t k$ to $\mathfrak{n} \times \mathfrak{n}$ is equal to the restriction of $\kappa$ and therefore positive definite. We want to find a $t\geq t_0$ such that $\kappa-t k$ is positive definite on $\mathfrak{g}^\prime \times \mathfrak{g}^\prime$. It suffices to obtain \[(\kappa-t k)(\lambda X+Y,\lambda X+Y)>0\] for all non-zero $\lambda \in \mathds{R}, X \in \mathfrak{l},Y \in \mathfrak{n}$. By linearity, we may suppose $k(X,X)=-1$ and $\kappa(Y,Y)=1$. Then
\begin{align*}
&(\kappa-t k)(\lambda X+Y,\lambda X+Y)\\
&=\lambda^2 \kappa(X,X)+t\lambda^2+2\lambda\kappa(X,Y)+1\\
&\geq \lambda^2 (t+\min_{X \in \mathfrak{l}, k(X,X)=-1} \kappa(X,X))-2\lambda\max_{X \in \mathfrak{l}, k(X,X)=-1}\max_{Y \in \mathfrak{n}, \kappa(Y,Y)=1} |\kappa(X,Y)|+1.\end{align*}
The last term is a polynomial in $\lambda$ with fixed parameters which discriminant will be less than zero if $t$ is large enough. Thus, $b:=\kappa-t k$ is positive definite if $t$ is large enough.

By Lemma~\ref{lem:positive_definite}~(iii), $\mathfrak{g}^\prime=[\mathfrak{g}^\prime,\mathfrak{g}^\prime]\oplus\mathfrak{z(g^\prime)}$. Thus, the radical of $\mathfrak{g}^\prime$ is equal to its center. Since $\mathfrak{n}$ is a nilpotent ideal, it follows that $\mathfrak{n} \subseteq \mathfrak{z(g^\prime)}$. Also, $\mathfrak{l}\subseteq[\mathfrak{g}^\prime,\mathfrak{g}^\prime]$. But $\mathfrak{g}^\prime=\mathfrak{l}+\mathfrak{n}$, hence, $\mathfrak{g}^\prime=\mathfrak{l}\oplus\mathfrak{n}$ is a direct sum of algebras.

Since $[\mathfrak{l},\mathfrak{n}]=\left\{0\right\}$, we can apply Lemma~\ref{lem:central_homomorphism} to see that $\mathfrak{l}=[\mathfrak{l},\mathfrak{l}]$ is $\kappa$-orthogonal to $\mathfrak{n}$.
\end{proof}\vspace{0pt}

\section{Radical}\label{sec:radical}

In this section, we investigate the radical $R$ of a connected non-compact Lie group $G$, provided with an ad-invariant symmetric bilinear form $\kappa$ on its Lie algebra $\mathfrak{g}$ that fulfills condition~\hyperlink{star}{$(\star)$} as described in Theorem~\ref{th:algebraic}. The radical of $G$ is the largest (with respect to inclusion) connected normal Lie subgroup of $G$, which is solvable. By maximality, $R$ is closed in $G$. Its Lie algebra $\mathfrak{r}$ is the largest (with respect to inclusion) solvable ideal of $\mathfrak{g}$.

First, we show a generalization of Corollary~\ref{cor:compact_nilpotent}:

\begin{proposition}\label{prop:compact_solvable}
A connected compact solvable Lie group $H$ is abelian.
\end{proposition}
\vspace{-1.4em}\begin{proof}
According to Proposition~\ref{prop:compact_algebra}, the Lie algebra $\mathfrak{h}$ of $H$ is reductive. Since $\mathfrak{h}$ is solvable, it follows that $\mathfrak{h}$ itself is abelian.
\end{proof}\vspace{0pt}

A proof of the following can be found in \cite{Ja62}, Chapter~2, Theorem~13:

\begin{lemma}\label{lem:[g,r]<n}
Let $\mathfrak{g}^\prime$ be a Lie algebra, $\mathfrak{r}^\prime$ its radical and $\mathfrak{n}^\prime$ its nilradical. Then $[\mathfrak{g}^\prime,\mathfrak{r}^\prime]\subseteq \mathfrak{n}^\prime$.
\end{lemma}

\begin{proposition}\label{prop:radical_nilradical}
If the radical $R$ is not compact, then also the nilradical $N$ is not compact.
\end{proposition}
\vspace{-1.4em}\begin{proof}
Assume $N$ is compact. By Corollary~\ref{cor:compact_nilpotent}, $N$ is abelian. Using Lemma~\ref{lem:compact_abelian}, $N$ is central in $G$. Hence, $\mathfrak{n} \subseteq \mathfrak{z(g)}$. Because the nilradical of the ideal $\mathfrak{r}$ is the intersection of $\mathfrak{n}$ with $\mathfrak{r}$ (cf.~\cite{Ja62}, Chapter~III, Theorem~7), $\mathfrak{n}$ is also the nilradical of $\mathfrak{r}$. 

Since $\mathfrak{r}$ is a solvable algebra, \[2\cdot\dim \mathfrak{n}\geq\dim\mathfrak{r}+\dim\mathfrak{z(r)}\] (cf.~\cite{OnVin94}, Chapter~2, Theorem~5.2). But $\mathfrak{n} \subseteq \mathfrak{z(g)}$, thus \[2\cdot\dim \mathfrak{n}\geq\dim\mathfrak{r}+\dim\mathfrak{z(r)}\geq\dim\mathfrak{r}+\dim\mathfrak{n}\] and therefore $\dim\mathfrak{n} \geq \dim\mathfrak{r}$. Hence, $\mathfrak{n}=\mathfrak{r}$, contradicting that $N$ is compact, but $R$ not.
\end{proof}\vspace{0pt}

Proposition~\ref{prop:nilradical_semisimple} generalizes to the radical.

\begin{proposition}\label{prop:radical_semisimple}
Let $R$ be non-compact and let $\mathfrak{l}\subset\mathfrak{g}$ be a semisimple subalgebra. Then $\mathfrak{g}^\prime=\mathfrak{l}+\mathfrak{r}$ is in fact a $\kappa$-orthogonal direct sum of algebras.
\end{proposition}
\vspace{-1.4em}\begin{proof}
$\mathfrak{l} \cap \mathfrak{r}$ is a solvable ideal in $\mathfrak{l}$. Since $\mathfrak{l}$ is semisimple, $\mathfrak{l} \cap \mathfrak{r}=\left\{0\right\}$.

We know from Propositions~\ref{prop:nilradical_semisimple} and~\ref{prop:radical_nilradical}, that $[\mathfrak{l},\mathfrak{n}]=\left\{0\right\}$. Additionally, $[\mathfrak{g},\mathfrak{r}]\subseteq\mathfrak{n}$ by Lemma~\ref{lem:[g,r]<n}.

Choose $X,Y \in \mathfrak{l}, W \in \mathfrak{r}$. Then \[0=[[X,Y],W]+[[Y,W],X]+[[W,X],Y]\] by the Jacobi identity. The last two summands are elements of \[[[\mathfrak{l},\mathfrak{r}],\mathfrak{l}]\subseteq[\mathfrak{n},\mathfrak{l}]=\left\{0\right\}.\] We obtain that $[[\mathfrak{l},\mathfrak{l}],\mathfrak{r}]=\left\{0\right\}$. But $[\mathfrak{l},\mathfrak{l}]=\mathfrak{l}$, so $\mathfrak{g}^\prime=\mathfrak{l}\oplus\mathfrak{r}$ is a direct sum of algebras.

The orthogonality of this decomposition follows from Lemma~\ref{lem:central_homomorphism}.
\end{proof}\vspace{0pt}

We conclude the preliminaries with the following decomposition:

\begin{proposition}\label{prop:decomposition_radical_noncompact}
Assume that $R$ is not compact. Then there is a compact semisimple subalgebra $\mathfrak{k}$ such that $\mathfrak{g}=\mathfrak{k}\oplus\mathfrak{r}$ is a $\kappa$-orthogonal direct sum. Furthermore, $\kappa$ restricted to $\mathfrak{k} \times \mathfrak{k}$ is positive definite.
\end{proposition}
\vspace{-1.4em}\begin{proof}
Let $\mathfrak{g}=\mathfrak{k}\inplus\mathfrak{r}$ be a Levi decomposition according to Lemma~\ref{lem:Levi_decomposition}. By Proposition~\ref{prop:radical_semisimple}, $\mathfrak{g}=\mathfrak{k}\oplus\mathfrak{r}$ is a $\kappa$-orthogonal direct sum. It remains to show that $\kappa|_{\mathfrak{k} \times \mathfrak{k}}$ is positive definite.

Since the radical is closed and connected, but not compact, there is an element $X \in \mathfrak{r}$ generating a non-precompact one-parameter group in $G$. Let $Y \in \mathfrak{k}$. Then the subalgebra $\mathfrak{a}:=\textnormal{span}\left\{X,Y\right\}$ is abelian. By Lemma~\ref{lem:kappa_abelian}, the set of elements in $\mathfrak{a}$ generating a non-precompact one-parameter group is dense in $\mathfrak{a}$.

Since this is true for all $Y \in \mathfrak{k}$, we can apply condition~\hyperlink{star}{$(\star)$} to the subalgebra $\mathfrak{p}:=\mathfrak{k}\oplus\mathds{R}X$ and obtain, that $\kappa$ is positive semidefinite on $\mathfrak{p} \times \mathfrak{p}$ and its kernel has dimension at most one. The intersection of the kernel with $\mathfrak{k}$ is an ideal of dimension at most one. But $\mathfrak{k}$ is semisimple, so this ideal is trivial.

Thus, $\kappa$ restricted to $\mathfrak{k} \times \mathfrak{k}$ is positive definite. By Proposition~\ref{prop:compact_algebra}, $\mathfrak{k}$ is compact.
\end{proof}\vspace{0pt}

\section{Compact radical: case of the special linear algebra}\label{sec:compact_radical}

In this section, we assume that the radical $R$ of a connected non-compact Lie group $G$, provided with an ad-invariant symmetric bilinear form $\kappa$ on its Lie algebra $\mathfrak{g}$ that fulfills condition~\hyperlink{star}{$(\star)$} as described in Theorem~\ref{th:algebraic}, is compact. By Proposition~\ref{prop:compact_solvable}, $R$ is abelian.

\begin{proposition}\label{prop:decomposition_case_compact_radical}
$\mathfrak{g}=\mathfrak{k}\oplus\mathfrak{r}\oplus\mathfrak{s}$ is a $\kappa$-orthogonal direct sum, where $\mathfrak{k}$ is compact semisimple and $\mathfrak{s}$ is simple, but non-compact. Furthermore, $\kappa|_{\mathfrak{k} \times \mathfrak{k}}$ is positive definite.
\end{proposition}
\vspace{-1.4em}\begin{proof}
Let $\mathfrak{g}=\mathfrak{l}\inplus\mathfrak{r}$ be a Levi decomposition according to Lemma~\ref{lem:Levi_decomposition}. By Lemma~\ref{lem:compact_abelian}, the radical is central in $G$ and therefore, $\mathfrak{g}=\mathfrak{l}\oplus\mathfrak{r}$ is a direct sum. Lemma~\ref{lem:central_homomorphism} shows that this sum is $\kappa$-orthogonal.

By assumption, $G$ is not compact, but $R$ is. By a theorem of Weyl, any Lie group with a compact semisimple Lie algebra is compact and has finite center (cf.~\cite{Bo05}, Chapter~IX, Paragraph~1.4, Theorem~1). But by Proposition~\ref{prop:compact_center_closed_subgroup}, a semisimple Lie subgroup with finite center is closed in $G$. Thus, $\mathfrak{l}$ contains a simple direct summand $\mathfrak{s}$, which is not compact.

Let $\mathfrak{s}^\prime$ be another simple direct  summand and $X \in \mathfrak{s}$ generating a non-precompact one-parameter group. For any $Y \in \mathfrak{s}^\prime$, the subalgebra $\mathfrak{a}:=\textnormal{span}\left\{X,Y\right\}$ is abelian. Using Lemma~\ref{lem:kappa_abelian}, we obtain that the set of elements in $\mathfrak{a}$ generating a non-precompact one-parameter group is dense in $\mathfrak{a}$. It follows that the same is true for $\mathfrak{p}:=\mathfrak{s}^\prime+\mathds{R}X$. According to condition~\hyperlink{star}{$(\star)$} applied to the subspace $\mathfrak{p}$, $\kappa$ is positive semidefinite on $\mathfrak{p} \times \mathfrak{p}$ and its kernel has dimension at most one. Its intersection with $\mathfrak{s}^\prime$ is an ideal and therefore trivial. Thus, $\kappa$ is positive definite on $\mathfrak{s}^\prime\times\mathfrak{s}^\prime$ and by Proposition~\ref{prop:compact_algebra}, $\mathfrak{s}^\prime$ is compact. With Lemma~\ref{lem:central_homomorphism} we obtain that $\mathfrak{s}$ and $\mathfrak{s}^\prime$ are $\kappa$-orthogonal to each other. In summary, $\mathfrak{g}=\mathfrak{k}\oplus\mathfrak{r}\oplus\mathfrak{s}$ is a $\kappa$-orthogonal direct sum, where $\mathfrak{k}$ is compact semisimple.
\end{proof}\vspace{0pt}

$\mathfrak{s}$ is a non-compact simple Lie algebra. By \cite{Got69}, Proposition~2, $\mathfrak{s}$ contains a subalgebra $\mathfrak{s}_0$ isomorphic to $\mathfrak{sl}_2(\mathds{R})$. Let $\left\{e,f,h\right\} \subset \mathfrak{s}_0$ be an $\mathfrak{sl}_2$-triple, that is, \[[h,e]=2e,[h,f]=-2f,[e,f]=h.\] In the following, we will show that $\mathfrak{s}\cong\mathfrak{sl}_2(\mathds{R})$.

\begin{proposition}\label{prop:kappa_s_Lorentz}
$\kappa|_{\mathfrak{s}_0\times\mathfrak{s}_0}$ is a positive multiple of the Killing form of $\mathfrak{s}_0$, especially a Lorentz form.
\end{proposition}
\vspace{-1.4em}\begin{proof}
Due to the ad-invariance of $\kappa$, $\kappa([e,f],e)=0=\kappa([e,f],f),$ so \[\kappa([h,e)=0=\kappa(h,f).\] Using $\kappa([h,e],e)=0=\kappa([h,f],f),$ we obtain that \[\kappa(e,e)=0=\kappa(f,f).\] Finally, it follows from $\kappa([e,f],h)=\kappa(e,[f,h])$ that \[\kappa(h,h)=2 \kappa(e,f).\] Thus, $\kappa|_{\mathfrak{s}_0\times\mathfrak{s}_0}$ is already determined by the scalar $\kappa(h,h)$. It follows that $\kappa|_{\mathfrak{s}_0\times\mathfrak{s}_0}$ is a multiple of the Killing form of $\mathfrak{s}_0$.

Let $\mathfrak{p}:=\textnormal{span}\left\{h,e\right\}$. Consider $X:= \alpha h + \beta e$, with real numbers $\alpha\neq0$ and $\beta$. The set of such $X$ is dense in $\mathfrak{p}$. Moreover, $\textnormal{ad}_X$ has the eigenvalue $2\alpha$ to the eigenvector $e$. By Lemma~\ref{lem:precompact_ad}, $X$ does not generate a precompact one-parameter group. Thus, we can apply condition~\hyperlink{star}{$(\star)$} to $\mathfrak{p}$ and conclude, that $\kappa(h,h)>0$. Therefore, $\kappa|_{\mathfrak{s}_0\times\mathfrak{s}_0}$ is a positive multiple of the Killing form of $\mathfrak{s}_0$.
\end{proof}\vspace{0pt}

Let $T$ be the $\kappa$-orthogonal complement to $\mathfrak{s}_0$ in $\mathfrak{s}$. Then $T \oplus \mathfrak{s}_0=\mathfrak{s}$ as vector spaces. Since $\kappa$ is ad-invariant, $\kappa([X,Y],Z)=\kappa(X,[Y,Z])=0$ for any $X \in T$ and $Y,Z \in \mathfrak{s}_0$. This means that $T$ is $\textnormal{ad}(\mathfrak{s}_0)$-invariant.

Proposition~\ref{prop:sl2R} will show that \[[X,T]=\left\{0\right\}\] for $X=\frac{h}{2}$ and $X=-\frac{e+f}{2}$. By the Jacobi identity, then also \[[[\frac{h}{2},\frac{e+f}{2}],T]=\left\{0\right\}\] holds. But $[\frac{h}{2},\frac{e+f}{2}]=\frac{e-f}{2}$ and $\left\{h,e+f,e-f\right\}$ form a basis of $\mathfrak{s}_0$. Thus, \[[\mathfrak{s}_0,T]=\left\{0\right\}.\] This implies $[\mathfrak{s}_0,\mathfrak{s}]=\mathfrak{s}_0$, so $\mathfrak{s}_0$ is an ideal in the simple algebra $\mathfrak{s}$, hence \[\mathfrak{s}=\mathfrak{s}_0\cong\mathfrak{sl}_2(\mathds{R}).\]

\begin{proposition}\label{prop:sl2R}
Let $X \in \left\{\frac{h}{2},-\frac{e+f}{2}\right\}$. Then $[X,T]=\left\{0\right\}$.
\end{proposition}
\vspace{-1.4em}\begin{proof}
Let $Y:=e$ if $X=\frac{h}{2}$, $Y:=h+(e-f)$ if $X=-\frac{e+f}{2}$. Then \[[X,Y]=Y\text{ and }\kappa(Y,Y)=0.\] Also, $\textnormal{ad}_Y|_{\mathfrak{s}_0}$ is nilpotent, but not trivial, because of the equations \begin{align*}[e,[e,[e,f]]]=[e,[e,h]]&=-2[e,e]=0\\
\text{and }[h+e-f,[h+e-f,[h+e-f,h]]]&=-2[h+e-f,[h+e-f,e+f]]\\
&=-4[h+e-f,h+e-f]=0.
\end{align*}

Because of \begin{align*}
[h,e]&=2e,\\
[h,f]&=-2f,\\
[h,h]&=0,\\
[e+f,h+(e-f)]&=-2(h+e-f),\\
[e+f,h-(e-f)]&=2(h-(e-f)),\\
[e+f,e+f]&=0,
\end{align*} $\textnormal{ad}_X|_{\mathfrak{s}_0}$ is semisimple with the real eigenvalues $-1,0,1$.

Since $T$ is $\textnormal{ad}(\mathfrak{s}_0)$-invariant, we have a representation $\varrho: \mathfrak{s}_0 \to \mathfrak{gl}(T)$. The kernel of $\varrho$ is an ideal in $\mathfrak{s}_0$. But $\mathfrak{s}_0$ is simple, so the ideal is either $\mathfrak{s}_0$ and we are done, or $\varrho$ is injective, which we suppose from now on. $\varrho(\mathfrak{s}_0)$ is semisimple.

According to Proposition~\ref{prop:Jordan} in the appendix, $X$ acts on $T$ as a semisimple endomorphism with real eigenvalues only and $Y$ acts as a nilpotent endomorphism. It follows that $\textnormal{ad}_Y$ is nilpotent, but not trivial, so due to Lemma~\ref{lem:precompact_ad}, $Y$ does not generate a precompact one-parameter group.  Since $X \neq 0$, there is $0 \neq A \in T$ and $\lambda \in \mathds{R}\backslash\left\{0\right\}$ such that $[X,A]=\lambda A$. \[\lambda \kappa(A,A)=\kappa([X,A],A)=\kappa(X,[A,A])=0,\] because $\kappa$ is ad-invariant.

If $[Y,A]=0$, then $A$ and $Y$ generate an abelian subalgebra. Furthermore, \[\kappa(Y,Y)=0=\kappa(A,A),\ \kappa(A,Y)=0\] since $T$ is $\kappa$-orthogonal to $\mathfrak{s}_0$. This contradicts Lemma~\ref{lem:kappa_abelian}. Therefore, it holds $B:=[Y,A] \neq 0$. Clearly, $B \in T$.

The Jacobi identity yields: \[0=[X,[Y,A]]+[Y,[A,X]]+[A,[X,Y]]=[X,B]-\lambda B-B.\] Therefore, $[X,B]=(\lambda+1)B$. Unless $\lambda \neq -1$, we could have started with $B$ instead of $A$. But $\textnormal{ad}_X$ has only finitely many eigenvalues, so it follows that all eigenvalues of $\textnormal{ad}_X$ are negative integers. It suffices to consider $\lambda=-1$ and $[X,B]=0$.

Applying the Jacobi identity another time, \[0=[X,[Y,B]]+[Y,[B,X]]+[B,[X,Y]]=[X,[Y,B]]+[B,Y].\] Hence, either $[Y,B]=0$ or $[Y,B]$ is an eigenvector of $\textnormal{ad}_X$ with eigenvalue 1. But all eigenvalues of $\textnormal{ad}_X$, so $[Y,B]=0$.

Because of $B \in T$, we have \[\kappa(Y,Y)=0=\kappa(Y,B).\] Since $\kappa$ is ad-invariant, \[\kappa(B,B)=\kappa(B,[Y,A])=\kappa([B,Y],A)=0.\] As above, we obtain a contradiction to Lemma~\ref{lem:kappa_abelian}.
\end{proof}\vspace{0pt}

We have shown that $\mathfrak{g}=\mathfrak{k}\oplus\mathfrak{r}\oplus\mathfrak{s}$, where $\mathfrak{k}$ is compact semisimple, the radical $\mathfrak{r}$ is abelian and $\mathfrak{s}\cong\mathfrak{sl}_2(\mathds{R})$.

\begin{lemma}\label{lem:kappa_radical}
$\kappa$ is positive definite on $\mathfrak{r} \times \mathfrak{r}$.
\end{lemma}
\vspace{-1.4em}\begin{proof}
$\textnormal{ad}_e$ is nilpotent, but not trivial. So for any $X \in \mathfrak{r}$, we can apply condition~\hyperlink{star}{$(\star)$} to the abelian subalgebra $\mathfrak{a}:=\textnormal{span}\left\{e,X\right\}$ by Lemma~\ref{lem:kappa_abelian}. Since $\kappa(e,e)=0=\kappa(e,X)$, it follows that $\kappa(X,X)>0$ for all $0 \neq X \in \mathfrak{r}$.
\end{proof}\vspace{0pt}

To conclude the proof of Theorem~\ref{th:algebraic} in the case that the radical is compact, we have to show that the subgroup generated by $\mathfrak{s}$ is a finite covering of $\textnormal{PSL}(2,\mathds{R})$ if and only if it is closed.

$\textnormal{PSL}(2,\mathds{R})$ is centerless and has fundamental group $\mathds{Z}$ (for the latter, see for example \cite{OnVin94}, Chapter~4, Paragraph~3.2, example~4). Let $\widetilde{\textnormal{SL}_2(\mathds{R})}$ be the universal cover of $\textnormal{PSL}(2,\mathds{R})$. By a theorem of Wolf, $\widetilde{\textnormal{SL}_2(\mathds{R})}$ has no non-trivial compact subgroup (cf.~\cite{Wo63}).

Thus, if the subgroup generated by $\mathfrak{s}$ is closed, it cannot be isomorphic to a finite central quotient of $\widetilde{\textnormal{SL}_2(\mathds{R})}$, since then all non-trivial one-parameter groups would be not precompact, contradicting condition~\hyperlink{star}{$(\star)$} (remember that $\kappa|_{\mathfrak{s}\times\mathfrak{s}}$ is a Lorentz form).

Conversely, if the subgroup generated by $\mathfrak{s}$ is isomorphic to some $\textnormal{PSL}_k(2,\mathds{R})$, it is a semisimple group with finite center and by Proposition~\ref{prop:compact_center_closed_subgroup} closed in $G$.

Since $\textnormal{PSL}(2,\mathds{R})$ is centerless and has fundamental group $\mathds{Z}$, the center of $\widetilde{\textnormal{SL}_2(\mathds{R})}$ is isomorphic to $\mathds{Z}$. The following group is motivated by an example of a non-closed semisimple Levi factor in a Lie group given in \cite{OnVin94}, Chapter~1, Paragraph~4.1. It shows that in Theorem~\ref{th:algebraic}~(v) both cases (some $\textnormal{PSL}_k(2,\mathds{R})$ or some finite quotient of $\widetilde{\textnormal{SL}_2(\mathds{R})}$) are possible.

\begin{proposition}\label{prop:counterexample1}
Let $z$ be a generator of the center of $\widetilde{\textnormal{SL}_2(\mathds{R})}$ and $\exp(\textnormal{i}\varphi)$ in $\mathds{S}^1 \subset \mathds{C}$ with $\varphi \in \mathds{R}$, $\frac{\varphi}{\pi}\notin \mathds{Q}$. Denote by $\Gamma \subset \widetilde{\textnormal{SL}_2(\mathds{R})} \times \mathds{S}^1$ the discrete central subgroup generated by $(z,\exp(\textnormal{i}\varphi))$. Let $G^\prime:=\left(\widetilde{\textnormal{SL}_2(\mathds{R})} \times \mathds{S}^1\right)\mkern-3mu/\Gamma$.

Then $G^\prime$ is a non-compact Lie group with Lie algebra $\mathfrak{g}^\prime=\mathfrak{sl}_2(\mathds{R})\oplus\mathds{R}$ and there is a symmetric bilinear form $\kappa^\prime$ on $\mathfrak{g}^\prime$, such that $\kappa^\prime$ is ad-invariant and fulfills condition~\hyperlink{star}{$(\star)$}. Moreover, the subgroup generated by $\mathfrak{sl}_2(\mathds{R})$ is not closed in $G^\prime$ and is not isomorphic to a finite covering of $\textnormal{PSL}(2,\mathds{R})$.
\end{proposition}
\vspace{-1.4em}\begin{proof}
It is clear that $G^\prime$ is a Lie group with Lie algebra $\mathfrak{g}^\prime=\mathfrak{sl}_2(\mathds{R})\oplus\mathds{R}$. $G^\prime$ is not compact, since its Lie algebra contains the non-compact summand $\mathfrak{sl}_2(\mathds{R})$ (consider Corollary~\ref{cor:subalgebra_compact}).

Remember the isometric action of $\textnormal{PSL}(2,\mathds{R})$ on the compact Lorentzian manifold $M=PSL(2,\mathds{R})/\Lambda$ given in Section~\ref{sec:sl_2}. This action is locally effective. Since the kernel of $\rho: \textnormal{PSL}(2,\mathds{R}) \to \textnormal{Isom}(M)$ is discrete and normal and hence central in $\textnormal{PSL}(2,\mathds{R})$, it has to be trivial, so we can consider $\textnormal{PSL}(2,\mathds{R})$ as a closed subgroup of $\textnormal{Isom}(M)$ (since it is semisimple and centerless; see Proposition~\ref{prop:compact_center_closed_subgroup}). By the way, this also shows that the case of a finite covering of $\textnormal{PSL}(2,\mathds{R})$ can appear in Theorem~\ref{th:algebraic}~(v).

The Killing form of $\mathfrak{sl}_2(\mathds{R})$ defines a Lorentzian metric on $M$. It follows from Corollary~\ref{cor:condition_star} and Proposition~\ref{prop:kappa_s_Lorentz}, that the induced bilinear form $\kappa$ fulfills condition~\hyperlink{star}{$(\star)$} and is Lorentzian. This means, that any timelike $X \in \mathfrak{sl}_2(\mathds{R})$ generates a precompact one-parameter group in $\textnormal{PSL}(2,\mathds{R})$.

Now define $\kappa^\prime$ as follows: $\kappa^\prime|_{\mathfrak{sl}_2(\mathds{R})\times\mathfrak{sl}_2(\mathds{R})}$ is equal to $\kappa$, $\kappa^\prime|_{\mathds{R}\times\mathds{R}}$ is positive definite and $\mathfrak{sl}_2(\mathds{R})$ and $\mathds{R}$ are $\kappa^\prime$-orthogonal. Then $\kappa^\prime$ is an ad-invariant Lorentzian scalar product on $\mathfrak{g}^\prime$.

For showing that $\kappa^\prime$ fulfills condition~\hyperlink{star}{$(\star)$}, it suffices to show that any timelike $X \in \mathfrak{sl}_2(\mathds{R})$ generates a precompact one-parameter group in $G^\prime$.

Let $X \in \mathfrak{sl}_2(\mathds{R})$ be timelike and $\left\{\psi^t\right\}_{t \in \mathds{R}}$ be the one-parameter group generated by $X$ in $\widetilde{\textnormal{SL}_2(\mathds{R})}$. Consider any sequence $\left\{t_k\right\}_{k=0}^\infty \subset \mathds{R}$. We want to show, that there is a convergent subsequence of $\left\{\exp(t_k X)\right\}_{k=0}^\infty$ in $G^\prime$.

From above, we know already that $X$ generates a precompact one-parameter group in $\textnormal{PSL}(2,\mathds{R})=\widetilde{\textnormal{SL}_2(\mathds{R})}/Z(\widetilde{\textnormal{SL}_2(\mathds{R})})$. Without loss of generality, the projection of $\left\{\psi^{t_k}\right\}_{k=0}^\infty$ onto $\textnormal{PSL}(2,\mathds{R})$ converges there. It follows that there is a sequence $\left\{z_k\right\}_{k=0}^\infty \subseteq Z(\widetilde{\textnormal{SL}_2(\mathds{R})})$, such that $\psi^{t_k}z_k^{-1}$ converges to some $\overline{\psi} \in \widetilde{\textnormal{SL}_2(\mathds{R})}$. For this, remember that $\widetilde{\textnormal{SL}_2(\mathds{R})} \to \textnormal{PSL}(2,\mathds{R})$ is a covering map.

By construction, \[\exp(t_k X)=(\psi^{t_k},0)\Gamma=(\psi^{t_k}z_k^{-1},\exp(\textnormal{i}\varphi_k))\Gamma\] in $G^\prime$ for some $\exp(\textnormal{i}\varphi_k) \in \mathds{S}^1$. Since $\mathds{S}^1$ is compact, we can choose a convergent subsequence of $\left\{\exp(\textnormal{i}\varphi_k)\right\}_{k=0}^\infty$. Without loss of generality, $\exp(\textnormal{i}\varphi_k) \to \exp(\textnormal{i}\overline{\varphi})$ as $k \to \infty$. But then $\exp(t_k X)$ converges to $(\overline{\psi},\exp(\textnormal{i}\overline{\varphi}))\Gamma$ as $k \to \infty$.

We have shown that any subsequence of the one-parameter group generated by $X$ contains a convergent subsequence. Therefore, $\left\{\exp(tX)\right\}_{t \in \mathds{R}}$ is precompact.

Finally, the image $S$ of $\widetilde{\textnormal{SL}_2(\mathds(R))}$ in $G^\prime$ is dense in $G^\prime$, because $\left\{\exp(\textnormal{i}k\varphi)\right\}_{k=0}^\infty$ is dense in $\mathds{S}^1$. But $S$ is equal to the subgroup generated by $\mathfrak{sl}_2(\mathds{R})$ in $G^\prime$. Since it is not closed, it follows by Proposition~\ref{prop:compact_center_closed_subgroup}, that the center of $S$ is not finite, so it is a finite central quotient of $\widetilde{\textnormal{SL}_2(\mathds{R})}$.
\end{proof}\vspace{0pt}

\section{Non-compact radical}\label{sec:noncompact_radical}

In this section, we assume that the radical $R$ of a connected non-compact Lie group $G$, provided with an ad-invariant symmetric bilinear form $\kappa$ on its Lie algebra $\mathfrak{g}$ that fulfills condition~\hyperlink{star}{$(\star)$} as described in Theorem~\ref{th:algebraic}, is not compact.

According to Proposition~\ref{prop:decomposition_radical_noncompact}, there is a $\kappa$-orthogonal direct decomposition of the Lie algebra $\mathfrak{g}=\mathfrak{k}\oplus\mathfrak{r}$ with $\mathfrak{k}$ compact semisimple and the restriction of $\kappa$ to $\mathfrak{k} \times \mathfrak{k}$ is positive definite.

For proving Theorem~\ref{th:algebraic}, we only have to show that the radical decomposes as announced. Thus, we may ignore the summand $\mathfrak{k}$ and suppose in the sequel that $G$ is solvable. Note that the radical is closed in $G$.

By Proposition~\ref{prop:radical_nilradical}, the nilradical $N$ is not compact. Hence, we can apply Proposition~\ref{prop:decomposition_nilradical} and obtain the $\kappa$-orthogonal direct sum $\mathfrak{n}=\mathfrak{a}\oplus\mathfrak{h}$ of a possibly trivial abelian algebra $\mathfrak{a}$ and an subalgebra $\mathfrak{h}$ being isomorphic to a Heisenberg algebra $\mathfrak{he}_d$ of dimension $2d+1$. Remember that $\mathfrak{a}\oplus\mathds{R}Z$, $\mathds{R}Z=\mathfrak{z(h)}$, was canonically defined, but $\mathfrak{a}$ in general not.

\subsection{Form not positive semidefinite: case of the twisted Heisenberg algebra}\label{sec:kappa_Lorentz}

Suppose that $\kappa$ is not positive semidefinite. Then there is a $T \in \mathfrak{g}$ such that $\kappa(T,T)<0$. Because of condition~\hyperlink{star}{$(\star)$}, $T$ generates a precompact one-parameter group. Thus, according to Lemma~\ref{lem:abelian_group}, the closure of $\left\{\exp(tT)\right\}_{t \in \mathds{R}}$ is isomorphic to a torus $\mathds{T}^m$ with Lie algebra $\mathfrak{t}$. The set of $X \in \mathfrak{t}$ generating a one-parameter group isomorphic to $\mathds{S}^1$ is dense in $\mathfrak{t}$ (this lies in the fact that the rationals are dense in the real numbers). Especially, we can choose a $T^\prime \in \mathfrak{t}$ generating a compact one-parameter group such that $\kappa(T^\prime,T^\prime)<0$. Without loss of generality, we may assume $T=T^\prime$.

\begin{proposition}\label{prop:T_centralizes}
$T$ centralizes $\mathfrak{a}\oplus\mathds{R}Z$. Additionally, $\mathfrak{a}\oplus\mathds{R}Z$ generates a compact central subgroup of $G$.
\end{proposition}
\vspace{-1.4em}\begin{proof}
We have seen in Proposition~\ref{prop:decomposition_nilradical}, that $\mathfrak{a}\oplus\mathds{R}Z=\mathfrak{z(n)}$. Since $\mathfrak{n}$ is an ideal in $\mathfrak{g}$, for any $X \in \mathfrak{z(n)}$, $Y \in \mathfrak{n}$ and $W \in \mathfrak{g}$, \[0=[[W,X],Y]+[[X,Y],W]+[[Y,W],X]=[[W,X],Y]\]
by the Jacobi identity. Thus, $\mathfrak{z(n)}$ is an ideal in $\mathfrak{g}$ and $\textnormal{ad}_T: \mathfrak{z(n)} \to \mathfrak{z(n)}$.

Due to Proposition~\ref{prop:decomposition_nilradical}, $\kappa|_{\mathfrak{z(n)} \times \mathfrak{z(n)}}$ is positive semidefinite and its kernel either vanishes or is equal to $\mathds{R}Z$.

Since $\kappa$ is ad-invariant, $\kappa([T,X],[T,X])=\kappa(T,[X,[T,X]])=0$ for any $X$ in $\mathfrak{z(n)}$. If the restriction of $\kappa$ to $\mathfrak{z(n)} \times \mathfrak{z(n)}$ is positive definite, it follows that $T$ centralizes $\mathfrak{z(n)}$. Otherwise $\textnormal{ad}_T: \mathfrak{z(n)} \to \mathds{R}Z$. But because $T$ generates a compact one-parameter group, $\textnormal{ad}_T$ has no real non-zero eigenvalue due to Lemma~\ref{lem:precompact_ad}. Especially, $\textnormal{ad}_T^2|_{\mathfrak{z(n)}}=0$. Thus, there exists no element $X\in \mathfrak{z(n)}$ such that $[T,X]\neq 0$, that is, $\textnormal{ad}_T|_{\mathfrak{z(n)}}=0$.

The nilradical $N$ is closed in $G$. Since the center of $N$ is closed in $N$, it follows that $\mathfrak{z(n)}$ generates a closed subgroup of $G$.

$\mathds{R}T \oplus \mathfrak{z(n)}$ is an abelian subalgebra. If any $X \in \mathfrak{z(n)}$ generated a non-precompact one-parameter group, then $\kappa$ restricted to $\mathfrak{b} \times \mathfrak{b}$, where $\mathfrak{b}:=\textnormal{span}\left\{T,X\right\}$, would be positive semidefinite by Lemma~\ref{lem:kappa_abelian}, contradicting $\kappa(T,T)<0$. That $\mathfrak{z(n)}$ generates a central subgroup now follows from Lemma~\ref{lem:compact_abelian}.
\end{proof}\vspace{0pt}

\begin{corollary}\label{cor:d>0}
The summand $\mathfrak{h}$ (being one-dimensional or a Heisenberg algebra) in $\mathfrak{n}=\mathfrak{a}\oplus\mathfrak{h}$ is not abelian, that is, $\mathfrak{h} \cong \mathfrak{he}_d$.
\end{corollary}
\vspace{-1.4em}\begin{proof}
The group generated by $\mathfrak{a}\oplus\mathds{R}Z$ is compact, but the nilradical is not.
\end{proof}\vspace{0pt}

\begin{proposition}\label{prop:group_generated}
Let $\mathfrak{s}:=\mathds{R}T\oplus\mathfrak{h}$. Then $\mathfrak{h}$ is an ideal in $\mathfrak{g}$, $\mathfrak{s}$ is a subalgebra and $\kappa|_{\mathfrak{s}\times\mathfrak{s}}$ is a Lorentz form. Furthermore, the subgroup generated by $\mathfrak{s}$ is isomorphic to a twisted Heisenberg group $\textnormal{He}_d^\lambda$ if the center of the subgroup is closed in $G$, and isomorphic to a twisted Heisenberg group $\overline{\textnormal{He}_d^\lambda}$ otherwise. 
\end{proposition}
\vspace{-1.4em}\begin{proof}
Since $\mathfrak{n}$ is an ideal, $[\mathfrak{g},\mathfrak{h}]\subseteq\mathfrak{n}$. By Proposition~\ref{prop:decomposition_nilradical}, $\mathfrak{h}$ is the $\kappa$-orthogonal complement to $(\mathfrak{a}\oplus\mathds{R}Z)$ in $\mathfrak{n}$. $\kappa$ is ad-invariant and $\mathfrak{g}$ centralizes $(\mathfrak{a}\oplus\mathds{R}Z)$ by Proposition~\ref{prop:T_centralizes}, therefore, $\mathfrak{h}$ is an ideal in $\mathfrak{g}$. Especially, $\mathfrak{s}$ is a subalgebra.

Let $Z^\perp$ be the $\kappa$-orthogonal complement of $\mathds{R}Z$ in $\mathfrak{g}$. Then $\mathfrak{n}\subseteq Z^\perp$. Since $Z \in \mathfrak{z(g)}$ due to Proposition~\ref{prop:nilradical_results}~(i), \[\kappa([X,Y],Z)=\kappa(X,[Y,Z])=0\] for all $X,Y \in \mathfrak{g}$. Especially, $Z^\perp$ is an ideal.

By Lemma~\ref{lem:[g,r]<n}, we have $[\mathfrak{g},\mathfrak{g}] \subseteq \mathfrak{n}$ since $\mathfrak{g}$ is solvable. Therefore, $[Z^\perp,Z^\perp] \subseteq \mathfrak{n}$. For any $X,Y \in \mathfrak{n}$ and $W \in Z^\perp$, \[\kappa([W,X],Y)=\kappa(W,[X,Y])=0\] because $[\mathfrak{n},\mathfrak{n}]=\mathds{R}Z$. But the kernel of $\kappa|_{\mathfrak{n}\times\mathfrak{n}}$ is equal to $\mathds{R}Z$, so \[[Z^\perp,[Z^\perp,Z^\perp]] \subseteq \mathds{R}Z.\] Since $Z$ is central, it follows that $Z^\perp$ is nilpotent, so $Z^\perp=\mathfrak{n}$.

Especially, $\kappa(T,Z)\neq 0$ and $\kappa|_{\mathfrak{s}\times\mathfrak{s}}$ is non-degenerate (since the kernel of $\kappa|_{\mathfrak{h}\times\mathfrak{h}}$ is $\mathds{R}Z$). Without loss of generality, $\kappa(T,Z)>0$ (otherwise we choose $-T$ instead of $T$).

Assume that there would exist two timelike elements in $\mathfrak{s}$ which are orthogonal to each other. Clearly, they cannot lie in $\mathfrak{h}$. Without loss of generality, $T+X$ and $T+Y$ with $X,Y\in \mathfrak{h}$ are these elements. Then
\begin{align*}
0&>\kappa(T+X,T+X)+\kappa(T+Y,T+Y)\\
&=2 \kappa(T,T) + 2\kappa(T,X)+2\kappa(T,Y)+\kappa(X,X)+\kappa(Y,Y)\\
&=2 \kappa(T+X,T+Y) - 2\kappa(X,Y)+\kappa(X,X)+\kappa(Y,Y)\\
&=\kappa(X-Y,X-Y)\geq 0,\end{align*}
contradiction. Hence, $\kappa|_{\mathfrak{s}\times\mathfrak{s}}$ is a Lorentzian scalar product.

Let $V$ be the $\kappa$-orthogonal complement to $\textnormal{span}\left\{T,Z\right\}$ in $\mathfrak{s}$. $\textnormal{span}\left\{T,Z\right\}$ is centralized by $T$, so $V$ is $\textnormal{ad}_T$-invariant. Also, $V \subseteq Z^\perp$, hence $V \subseteq \mathfrak{h}$. But $Z \notin V$, hence $V$ is a vector space complement to $\mathds{R}Z$ in $\mathfrak{h}$. Obviously, $\kappa|_{V \times V}$ is positive definite.

By Corollary~\ref{cor:orthogonal_basis}, we may choose a canonical basis $\left\{Z,X_1,Y_1,\ldots,X_d,Y_d\right\}$ of $\mathfrak{h}\cong\mathfrak{he}_d$, such that $\left\{X_1,Y_1,\ldots,X_d,Y_d\right\}\subset V$ is $\kappa$-orthogonal. For $1 \leq k \leq d$, let \[\widehat{X_k}:=\frac{X_k}{\sqrt{\kappa(X_k,X_k)}}\text{ and }\widehat{Y_k}:=\frac{Y_k}{\sqrt{\kappa(Y_k,Y_k)}}.\] By construction, $\left\{\widehat{X_1},\widehat{Y_1},\ldots,\widehat{X_d},\widehat{Y_d}\right\}$ is $\kappa$-orthonormal. We define \[\lambda^\prime_k:=\frac{1}{\sqrt{\kappa(X_k,X_k)\kappa(Y_k,Y_k)}},\] such that $[\widehat{X_k},\widehat{Y_k}]=\lambda^\prime_k Z$. Let \[\lambda_k:=\lambda^\prime_k \kappa(T,Z)>0.\]

Because of the ad-invariance of $\kappa$,
\begin{align*} \begin{array}{rcl} \kappa([T,\widehat{X_j}],\widehat{X_k})&=\kappa(T,[\widehat{X_j},\widehat{X_k}])&=0\\
\text{and} \quad \kappa([T,\widehat{X_j}],\widehat{Y_k})&=\kappa(T,[\widehat{X_j},\widehat{Y_k}])&=\delta_{jk}\lambda^\prime_k \kappa(T,Z)
\end{array}
\end{align*} for all $j$ and $k$. It follows that \[[T,\widehat{X_k}]=\lambda_k \widehat{Y_k}\] for all $k$. \[[T,\widehat{Y_k}]=-\lambda_k \widehat{X_k}\] in a similar way. Therefore, $\mathfrak{s}\cong\mathfrak{he}_d^\lambda$.

For all $t \in \mathds{R}$ and all $k$, the action of $\textnormal{Ad}_{\exp(tT)}=\exp(\textnormal{ad}_{tT})$ on $\textnormal{span}\left\{\widehat{X_k},\widehat{Y_k}\right\}$ corresponds to the matrix \[\begin{pmatrix} \cos(t\lambda_k) & -\sin(t\lambda_k) \\ \sin(t\lambda_k) & \cos(t\lambda_k)\end{pmatrix}.\] Since $\left\{\textnormal{Ad}_{\exp(tT)}\right\}_{t \in \mathds{R}}$ is isomorphic to a circle, it follows that $\lambda \in r\mathds{Q}_+^d$ for some $r \in \mathds{R}_+$. Due to Lemma~\ref{lem:isom_heisenberg}, we may assume already $\lambda \in \mathds{Z}_+^d$.

The subgroup $S$ generated by $\mathfrak{s}$ is isomorphic to $\textnormal{He}_d^\lambda$ if $Z$ generates a circle and isomorphic to $\overline{\textnormal{He}_d^\lambda}$ otherwise. If $S$ is not closed in $G$, the nilradical of $S$ is neither, since $S\cong\mathds{S}^1 \ltimes N$ and $\mathds{S}^1$ is compact. Thus, Proposition~\ref{prop:compact_center_closed_subgroup} implies that the center of $N$ is not compact, especially not a circle. If $S$ is closed in $G$, the center of $S$ is also closed in $G$. But the center of $S$ is generated by $Z$ and $Z$ generates a precompact one-parameter group by Proposition~\ref{prop:T_centralizes}. Thus, the center of $S$ is a circle.
\end{proof}\vspace{0pt}

Denote by $\mathfrak{s}^\perp$ the $\kappa$-orthogonal complement to $\mathfrak{s}$ in $\mathfrak{g}$. Since $\kappa$ is ad-invariant and $\mathfrak{s}$ is a subalgebra, $\mathfrak{s}^\perp$ is $\textnormal{ad}(\mathfrak{s})$-invariant.

\begin{proposition}\label{prop:radical_noncompact_decomposition}
$\mathfrak{s}^\perp$ centralizes $\mathfrak{n}$. Furthermore, $\mathfrak{s}^\perp \subset \mathfrak{n}$.
\end{proposition}
\vspace{-1.4em}\begin{proof}
Let $X \in \mathfrak{h}$, $Y \in \mathfrak{s}^\perp$. Since $\mathfrak{h}$ is an ideal due to Proposition~\ref{prop:group_generated} and $\mathfrak{h}\subset \mathfrak{s}$, $[X,Y] \in \mathfrak{s} \cap \mathfrak{s}^\perp=\left\{0\right\}$. $\mathfrak{a}$ lies central in $\mathfrak{g}$ by Proposition~\ref{prop:T_centralizes}, so $\mathfrak{n}=\mathfrak{a}\oplus\mathfrak{h}$ is centralized by $\mathfrak{s}^\perp$. So for any $X \in \mathfrak{s}^\perp$, $\mathds{R}X+\mathfrak{n}$ is nilpotent. But $\mathfrak{g}$ is solvable, so $[\mathfrak{g},\mathfrak{g}] \subseteq \mathfrak{n}$ by Lemma~\ref{lem:[g,r]<n}. Hence, $\mathds{R}X+\mathfrak{n}$ is an ideal. It follows that $X \in \mathfrak{n}$ and $\mathfrak{s}^\perp \subset \mathfrak{n}$.
\end{proof}\vspace{0pt}

Clearly, $\mathfrak{s}^\perp \cap \mathfrak{h}=\left\{0\right\}$. So without loss of generality, we may take $\mathfrak{a}=\mathfrak{s}^\perp$ (remember that $\mathfrak{a}$ was not canonically defined yet). We obtain $\mathfrak{g}=\mathfrak{s}^\perp\oplus\mathfrak{s}=\mathfrak{a}\oplus\mathfrak{s}$ as a $\kappa$-orthogonal direct sum of algebras. This finishes the proof of Theorem~\ref{th:algebraic} in the case that the radical is non-compact and $\kappa$ is not positive semidefinite.

To conclude this section, we show an analogous result to Proposition~\ref{prop:counterexample1} for the case of twisted Heisenberg algebras. It shows that in Theorem~\ref{th:algebraic} (iv) both cases ($\textnormal{He}_d^\lambda$ and $\overline{\textnormal{He}_d^\lambda}$) are possible.

\begin{proposition}\label{prop:counterexample2}
Let $z$ be a non-trivial element of the center of $\overline{\textnormal{He}_d^\lambda}$ and $\exp(\textnormal{i}\varphi) \in \mathds{S}^1 \subset \mathds{C}$ with $\varphi \in \mathds{R}$, $\frac{\varphi}{\pi}\notin \mathds{Q}$. Denote by $\Gamma \subset \overline{\textnormal{He}_d^\lambda} \times \mathds{S}^1$ the discrete central subgroup generated by $(z,\exp(\textnormal{i}\varphi))$. Let $G^\prime:=\left(\overline{\textnormal{He}_d^\lambda} \times \mathds{S}^1\right)\mkern-3mu/\Gamma$.

Then $G^\prime$ is a non-compact Lie group with Lie algebra $\mathfrak{g}^\prime=\mathfrak{he}_d^\lambda\oplus\mathds{R}$ and there is a symmetric bilinear form $\kappa^\prime$ on $\mathfrak{g}^\prime$, such that $\kappa^\prime$ is ad-invariant and fulfills condition~\hyperlink{star}{$(\star)$}. Furthermore, the subgroup generated by $\mathfrak{he}_d^\lambda$ is not closed in $G^\prime$ and is not isomorphic to $\textnormal{He}_d^\lambda$.
\end{proposition}
\vspace{-1.4em}\begin{proof}
The proof is done in an analogous manner as the proof of Proposition~\ref{prop:counterexample1}.

It is clear that $G^\prime$ is a Lie group with Lie algebra $\mathfrak{g}^\prime=\mathfrak{he}_d^\lambda\oplus\mathds{R}$. $G^\prime$ is not compact, since it is solvable, but not abelian (consider Proposition~\ref{prop:compact_solvable}).

Remember the isometric action of $\textnormal{He}_d^\lambda$ (which we can consider as the quotient of $\widetilde{\textnormal{He}_d^\lambda}$ by the discrete subgroup generated by $z$) on the compact Lorentzian manifold $M=\textnormal{He}_d^\lambda/\Lambda$ given in Section~\ref{sec:twisted_Heisenberg}. This action is locally effective. Since the kernel of $\rho: \textnormal{He}_d^\lambda \to \textnormal{Isom}(M)$ is discrete and normal, it is central in $\textnormal{He}_d^\lambda$. Therefore, also the image is isomorphic to $\textnormal{He}_d^\lambda$ (remember that $\textnormal{He}_d^\lambda$ was defined as the quotient of $\widetilde{\textnormal{He}_d^\lambda}$ by any lattice of the center). Thus, we may also consider the image instead or equivalently suppose that the kernel of $\rho$ is trivial. In this case, $\textnormal{He}_d^\lambda$ as a closed subgroup of $\textnormal{Isom}(M)$ (since $\textnormal{He}_d^\lambda=\mathds{S}^1\ltimes \textnormal{He}_d$ and $\textnormal{He}_d$ is closed in $\textnormal{Isom}(M)$ because it is nilpotent and has compact center; see Proposition~\ref{prop:compact_center_closed_subgroup}). By the way, this also shows that the case of $\textnormal{He}_d^\lambda$ can occur in Theorem~\ref{th:algebraic}~(iv).

Any ad-invariant Lorentz form given in Proposition~\ref{prop:lorentz_heisenberg} defines a Lorentzian metric on $M$. It follows from Corollary~\ref{cor:condition_star} and Proposition~\ref{prop:group_generated}, that the induced bilinear form $\kappa$ fulfills condition~\hyperlink{star}{$(\star)$} and is Lorentzian. This means, that any timelike $X \in \mathfrak{he}_d^\lambda$ generates a precompact one-parameter group in $\textnormal{He}_d^\lambda$.

Now define $\kappa^\prime$ as follows: $\kappa^\prime|_{\mathfrak{he}_d^\lambda\times\mathfrak{he}_d^\lambda}$ is equal to $\kappa$, $\kappa^\prime|_{\mathds{R}\times\mathds{R}}$ is positive definite and $\mathds{R}$ and $\mathfrak{he}_d^\lambda$ are $\kappa^\prime$-orthogonal. Then $\kappa^\prime$ is an ad-invariant Lorentzian scalar product on $\mathfrak{g}^\prime$.

For showing that $\kappa^\prime$ fulfills condition~\hyperlink{star}{$(\star)$}, it suffices to show that any timelike $X \in \mathfrak{he}_d^\lambda$ generates a precompact one-parameter group in $G^\prime$.

Let $X \in \mathfrak{he}_d^\lambda$ be timelike and $\left\{\psi^t\right\}_{t \in \mathds{R}}$ be the one-parameter group generated in $\overline{\textnormal{He}_d^\lambda}$. Consider any sequence $\left\{t_k\right\}_{k=0}^\infty \subset \mathds{R}$. We want to show, that there is a convergent subsequence of $\left\{\exp(t_k X)\right\}_{k=0}^\infty$ in $G^\prime$.

From above, we know already that $X$ generates a precompact one-parameter group in $\textnormal{He}_d^\lambda=\overline{\textnormal{He}_d^\lambda}/\left\{z^l | l \in \mathds{Z}\right\}$. Without loss of generality, the projection of $\left\{\psi^{t_k}\right\}_{k=0}^\infty$ onto $\textnormal{He}_d^\lambda$ converges there. So there is a sequence $\left\{z_k\right\}_{k=0}^\infty \subseteq \left\{z^l | l \in \mathds{Z}\right\}$, such that $\psi^{t_k}z_k^{-1}$ converges to some $\overline{\psi} \in \overline{\textnormal{He}_d^\lambda}$. For this, simply notice that the map $\overline{\textnormal{He}_d^\lambda} \to \textnormal{He}_d^\lambda$ is a covering map.

By construction, \[\exp(t_k X)=(\psi^{t_k},0)\Gamma=(\psi^{t_k}z_k^{-1},\exp(\textnormal{i}\varphi_k))\Gamma\] in $G^\prime$ for some $\exp(\textnormal{i}\varphi_k) \in \mathds{S}^1$. Since $\mathds{S}^1$ is compact, we can choose a convergent subsequence of $\left\{\exp(\textnormal{i}\varphi_k)\right\}_{k=0}^\infty$. Without loss of generality, $\exp(\textnormal{i}\varphi_k) \to \exp(\textnormal{i}\overline{\varphi})$ as $k \to \infty$. But then $\exp(t_k X)$ converges to $(\overline{\psi},\exp(\textnormal{i}\overline{\varphi}))\Gamma$ as $k \to \infty$.

We have shown that any subsequence of the one-parameter group generated by $X$ contains a convergent subsequence. Therefore, $\left\{\exp(tX)\right\}_{t \in \mathds{R}}$ is precompact.

Finally, the image $S$ of $\overline{\textnormal{He}_d^\lambda}$ in $G^\prime$ is dense in $G^\prime$, because $\left\{\exp(\textnormal{i}k\varphi)\right\}_{k=0}^\infty$ is dense in $\mathds{S}^1$. But $S$ is equal to the subgroup generated by $\mathfrak{he}_d^\lambda$ in $G^\prime$. Since it is not closed, the nilradical $N$ is neither (note that $S \cong \mathds{S}^1\ltimes N$ and $\mathds{S}^1$ is compact), so it follows by Proposition~\ref{prop:compact_center_closed_subgroup}, that the center of $N$, which is equal to the center of $S$, is not compact, so $S$ is not isomorphic to $\textnormal{He}_d^\lambda$.
\end{proof}\vspace{0pt}

\subsection{Form positive semidefinite}\label{sec:kappa_positive}

In the following, we suppose that $\kappa$ is positive semidefinite. Note that we still assume $\mathfrak{g}$ to be solvable.

\subsubsection{Nilpotent radical: trivial case and case of the Heisenberg algebra}\label{sec:nilpotent_radical}

If $\mathfrak{g}$ is nilpotent, Proposition~\ref{prop:decomposition_nilradical} yields the proof of Theorem~\ref{th:algebraic}. Note that we obtain the cases~(i) and~(iii).

\subsubsection{Non-nilpotent radical: case of the affine algebra}\label{sec:nonnilpotent_radical}

It remains the case when $\mathfrak{g}$ is not nilpotent. Especially, $\mathfrak{g}$ is not abelian, so by Propositions~\ref{prop:reductive_algebra} and~\ref{prop:compact_algebra}, $\kappa$ is not positive definite. Let $\mathfrak{i}$ be the kernel of $\kappa$.

\begin{proposition}\label{prop:case_affine_algebra_kernel}
The following is true:
\begin{compactenum}
\item $[\mathfrak{g},\mathfrak{g}] \subseteq \mathfrak{i}$.

\item $\mathfrak{n}$ is abelian.

\item $\mathfrak{i}\subseteq\mathfrak{n}$. Furthermore, $\mathfrak{i}$ is one-dimensional and not central in $\mathfrak{g}$.
\end{compactenum}
\end{proposition}
\vspace{-1.4em}\begin{proof}
(i) Clearly, $\mathfrak{g}/\mathfrak{i}$ is solvable and $\kappa$ induces an ad-invariant positive definite symmetric bilinear form. By Propositions~\ref{prop:reductive_algebra} and~\ref{prop:compact_algebra}, $\mathfrak{g}/\mathfrak{i}$ is abelian, that is, $[\mathfrak{g},\mathfrak{g}] \subseteq \mathfrak{i}$.

(ii) By Lemma~\ref{lem:[g,r]<n}, $[\mathfrak{g},\mathfrak{g}] \subseteq \mathfrak{n}$. Since $\mathfrak{g}$ is not abelian, it follows that $\mathfrak{i} \cap \mathfrak{n} =\mathds{R}Z$ is one-dimensional by Proposition~\ref{prop:decomposition_nilradical}. Because $\mathfrak{g}$ is not nilpotent, $Z$ is not central in $\mathfrak{g}$. So by Proposition~\ref{prop:nilradical_results}~(i), $\mathfrak{n}$ is abelian.

(iii) Assume that $\dim \mathfrak{i} >1$. Then $\mathfrak{i}+\mathfrak{n}$ is an ideal, which is not nilpotent. It follows from \[[\mathfrak{i}+\mathfrak{n},\mathfrak{i}+\mathfrak{n}]\subseteq[\mathfrak{g},\mathfrak{g}]\subseteq \mathfrak{i} \cap \mathfrak{n}=\mathds{R}Z,\] that there exists $X \in \mathfrak{i}+\mathfrak{n}$, such that $[X,Z]\neq 0$. Clearly, we may take $X \in \mathfrak{i} \backslash \mathfrak{n}$ such that $[X,Z]=Z$. Consider the two-dimensional subalgebra $\mathfrak{p}$ generated by $X$ and $Z$. Then $[\alpha X+ \beta Z,Z]=\alpha Z$ for all $\alpha,\beta \in \mathds{R}$, so by Lemma~\ref{lem:precompact_ad}, $\alpha X+ \beta Z$ generates a non-precompact one-parameter group if $\alpha \neq 0$. Especially, we can apply condition~\hyperlink{star}{$(\star)$} to $\mathfrak{p}$ and obtain a contradiction to $\mathfrak{p} \subseteq \mathfrak{i}$. Thus, $\mathfrak{i}=\mathds{R}Z$. We have seen in~(ii), that $Z$ is not central in $\mathfrak{g}$.
\end{proof}\vspace{0pt}

We are now able to finish the proof of Theorem~\ref{th:algebraic}.

\begin{proposition}\label{prop:case_affine_algebra}
$\mathfrak{g}=\mathfrak{a}\oplus\mathfrak{s}$ with $\mathfrak{a}$ abelian and $\mathfrak{s} \cong \mathfrak{aff}(\mathds{R})$ is a $\kappa$-orthogonal direct sum. Moreover, the kernel of $\kappa|_{\mathfrak{s} \times \mathfrak{s}}$ corresponds exactly to the span of the generator of the translations in the affine group.
\end{proposition}
\vspace{-1.4em}\begin{proof}
We know from Proposition~\ref{prop:case_affine_algebra_kernel}, that $\mathfrak{i}=\mathds{R}Z$ is one-dimensional. Since $[\mathfrak{g},\mathfrak{g}]\subseteq \mathfrak{i}$ and $\mathfrak{g}$ is not abelian, the linear map $\mathfrak{g} \to \mathfrak{i}$ defined by $W \mapsto [Z,W]$ is well defined and surjective. So its kernel $P$ has codimension 1 in $\mathfrak{g}$ and $\mathfrak{i}\subseteq P$. We choose $X \in \mathfrak{g}$ $\kappa$-orthogonal to $P$ such that $[X,Z]=Z$ and define the subalgebra $\mathfrak{s}:=\textnormal{span}\left\{X,Z\right\} \cong \mathfrak{aff}(\mathds{R})$. The kernel of $\kappa|_{\mathfrak{s} \times \mathfrak{s}}$ is equal to $\mathds{R}Z$.

Consider now the linear map $P \to \mathfrak{i}$, $W \mapsto [X,W]$. Since $Z \in P$, this map is surjective and its kernel $\mathfrak{a}$ has codimension 1 in $P$. We obtain \[\mathfrak{g}=P\oplus\mathds{R}X=\mathfrak{a}\oplus\mathfrak{s}\] as a $\kappa$-orthogonal direct sum of vector spaces. It remains to show that $\mathfrak{a}$ is a central subalgebra in $\mathfrak{g}$.

Let $W,W^\prime \in \mathfrak{a}$. Then \[[X,[W,W^\prime]]=-[W,[W^\prime,X]]-[W^\prime,[X,W]]=0\] by the Jacobi identity. Since \[[\mathfrak{g},\mathfrak{g}]=\mathds{R}Z\text{ and }[X,[W,W^\prime]]=0\neq[X,Z],\] $[W,W^\prime]=0$. Thus, $\mathfrak{a}$ is an abelian algebra. By construction, \[[X,\mathfrak{a}]=\left\{0\right\}=[Z,\mathfrak{a}],\] therefore, $\mathfrak{a}$ is central in $\mathfrak{g}$.
\end{proof}\vspace{0pt}

\section{General subgroups of the isometry group}\label{sec:general_subgroups}

This section is devoted to prove Theorem~\ref{th:algebraic_classification}. It suffices to consider connected Lie subgroups $G \subseteq \textnormal{Isom}(M)$.

If $G$ is precompact in $\textnormal{Isom}(M)$, the result follows from Corollaries~\ref{cor:subalgebra_compact} and~\ref{cor:compact_algebra}. So suppose $G$ is not precompact.

Let $\overline{G}$ be the closure of $G$ in $\textnormal{Isom}(M)$. Denote by $\overline{\mathfrak{g}}$ the Lie algebra of $\overline{G}$, by $\overline{\mathfrak{r}}$ the radical and by $\overline{\mathfrak{n}}$ the nilradical of $\overline{\mathfrak{g}}$. Since $G$ normalizes $\mathfrak{r}$ and $\mathfrak{n}$, also $\overline{G}$ does, that is, $\mathfrak{r}$ and $\mathfrak{n}$ are ideals in $\overline{\mathfrak{g}}$. Thus, $\mathfrak{r}\subseteq\overline{\mathfrak{r}}$ and $\mathfrak{n}\subseteq\overline{\mathfrak{n}}$.

The following proposition is due to Maltsev and a proof can be found in \cite{On93}, Part~I, Chapter~1, Theorem~5.3.

\begin{proposition}\label{prop:Maltsev}
$[\overline{\mathfrak{g}},\overline{\mathfrak{g}}]=[\mathfrak{g},\mathfrak{g}]$.
\end{proposition}

$\overline{G}$ is non-compact and by Corollary~\ref{cor:condition_star}, the induced bilinear form $\kappa$ fulfills condition~\hyperlink{star}{$(\star)$} in Theorem~\ref{th:algebraic}. We now distinguish between the cases~(i) to~(v) of Theorem~\ref{th:algebraic} applied to $\overline{G}$.

\subsection{Trivial case}

We have $\overline{\mathfrak{g}}=\mathfrak{k}\oplus\mathfrak{a}$ with $\mathfrak{k}$ compact semisimple and $\mathfrak{a}$ abelian. Thus, $\overline{\mathfrak{g}}$ is reductive, and by Proposition~\ref{prop:compact_algebra}, $\overline{\mathfrak{g}}$ is compact. Due to Corollary~\ref{cor:subalgebra_compact}, $\mathfrak{g}$ is compact as well and the result follows from Propositions~\ref{prop:reductive_algebra} and~\ref{prop:compact_algebra}.

\subsection{Case of the affine algebra}

We have $\overline{\mathfrak{g}}=\mathfrak{k}\oplus\mathfrak{a}\oplus\mathfrak{s}$ with $\mathfrak{k}$ compact semisimple, $\mathfrak{a}$ abelian and $\mathfrak{s}\cong\mathfrak{aff}(\mathds{R})$. Especially, $\overline{\mathfrak{g}}=\mathfrak{k}\oplus\overline{\mathfrak{r}}$ with $\overline{\mathfrak{r}}=\mathfrak{a}\oplus\mathfrak{s}$. By Proposition~\ref{prop:Maltsev}, $\mathfrak{k}\subseteq \mathfrak{g}$. Since $\mathfrak{r} \subseteq \overline{\mathfrak{r}}$, it follows that $\mathfrak{k}\oplus\mathfrak{r} \subseteq \mathfrak{g}$ and hence, $\mathfrak{k}\oplus\mathfrak{r}=\mathfrak{g}$.

$\mathfrak{r}$ is an ideal in $\overline{\mathfrak{r}}=\mathfrak{a}\oplus\mathfrak{s}$. Using $\mathfrak{s}\cong\mathfrak{aff}(\mathds{R})$, we obtain that $\mathfrak{r}$ is either abelian or isomorphic to the direct sum of an abelian algebra (possibly trivial) and $\mathfrak{aff}(\mathds{R})$.

\subsection{Case of the Heisenberg algebra}

We have $\overline{\mathfrak{g}}=\mathfrak{k}\oplus\mathfrak{a}\oplus\mathfrak{s}$ with $\mathfrak{k}$ compact semisimple, $\mathfrak{a}$ abelian and $\mathfrak{s}\cong\mathfrak{he}_d$. Moreover, $\kappa|_{\mathfrak{a}\times\mathfrak{a}}$ is positive definite and $\kappa|_{\mathfrak{s}\times\mathfrak{s}}$ positive semidefinite with one-dimensional kernel.

Clearly, $\overline{\mathfrak{g}}=\mathfrak{k}\oplus\overline{\mathfrak{r}}$ and $\overline{\mathfrak{r}}=\overline{\mathfrak{n}}=\mathfrak{a}\oplus\mathfrak{s}$. By Proposition~\ref{prop:Maltsev}, $\mathfrak{k}\subseteq \mathfrak{g}$. Since $\mathfrak{r} \subseteq \overline{\mathfrak{r}}$, it follows that $\mathfrak{k}\oplus\mathfrak{r} \subseteq \mathfrak{g}$ and hence, $\mathfrak{k}\oplus\mathfrak{r}=\mathfrak{g}$.

$\mathfrak{r}$ is an ideal in $\overline{\mathfrak{n}}=\mathfrak{a}\oplus\mathfrak{s}$. Therefore, $\mathfrak{r}=\mathfrak{n}$. $\kappa|_{\mathfrak{n}\times\mathfrak{n}}$ is positive semidefinite with kernel having dimension at most one, so the result follows from Proposition~\ref{prop:decomposition_nilradical}.

\subsection{Case of the twisted Heisenberg algebra}

We have $\overline{\mathfrak{g}}=\mathfrak{k}\oplus\mathfrak{a}\oplus\mathfrak{s}$ with $\mathfrak{k}$ compact semisimple, $\mathfrak{a}$ abelian and $\mathfrak{s}\cong\mathfrak{he}_d^\lambda$, $\lambda \in \mathds{Z}_+^d$. Additionally, $\kappa|_{\mathfrak{a}\times\mathfrak{a}}$ is positive definite and $\kappa|_{\mathfrak{s}\times\mathfrak{s}}$ is a Lorentz form. Also, the subgroup generated by $\mathfrak{a}\oplus\mathfrak{z(s)}$ is compact.

Especially, $\overline{\mathfrak{g}}=\mathfrak{k}\oplus\overline{\mathfrak{r}}$. By Proposition~\ref{prop:Maltsev}, $\mathfrak{k} \oplus \mathfrak{h}\subseteq \mathfrak{g}$, where $\mathfrak{h}\subset\mathfrak{s}$ with $\mathfrak{h}\cong\mathfrak{he}_d$. Since $\mathfrak{r} \subseteq \overline{\mathfrak{r}}$, it follows that $\mathfrak{k}\oplus\mathfrak{r} \subseteq \mathfrak{g}$ and hence, $\mathfrak{k}\oplus\mathfrak{r}=\mathfrak{g}$.

$\mathfrak{r}$ is an ideal in $\overline{\mathfrak{r}}=\mathfrak{a}\oplus\mathfrak{s}$. We first suppose that $\mathfrak{r}=\mathfrak{n}$. It follows that $\mathfrak{r} \subseteq \mathfrak{a}\oplus\mathfrak{h}$. Then $\kappa|_{\mathfrak{r}\times\mathfrak{r}}$ is positive semidefinite and its kernel has dimension at most one, so the result follows from Proposition~\ref{prop:decomposition_nilradical}.

Now suppose $\mathfrak{r}$ is not nilpotent. Let $\left\{T,Z,X_1,Y_1,\ldots,X_d,Y_d\right\}$ be a canonical basis of $\mathfrak{s}\cong\mathfrak{he}_d^\lambda$. Since $\mathfrak{r}$ is not nilpotent, $T+X \in \mathfrak{r}$ for some $X \in \mathfrak{a}\oplus\mathfrak{h}$. If $V:=\textnormal{span}\left\{X_1,Y_1,\ldots,X_d,Y_d\right\}$, \[[T+X,V]=V,\] hence $V \subset \mathfrak{r}$ because the latter is an ideal. $[V,V]=Z$, such that $\mathfrak{h}\subset\mathfrak{r}$. It follows that $\mathfrak{r}$ is isomorphic to the direct sum of $\mathfrak{he}_d^\lambda$ and an abelian algebra.

Denote by $\mathfrak{s}^\prime$ the corresponding subalgebra in $\mathfrak{r}$ isomorphic to $\mathfrak{he}_d^\lambda$. The subgroup $S^\prime$ generated by $\mathfrak{s}^\prime$ is isomorphic to $\textnormal{He}_d^\lambda$ if $\mathfrak{z(s^\prime)}=\mathfrak{z(s)}=\mathds{R}Z$ generates a circle and isomorphic to $\overline{\textnormal{He}_d^\lambda}$ otherwise. If $S^\prime$ is not closed in $\textnormal{Isom}(M)$, the nilradical $N^\prime$ of $S^\prime\cong \mathds{S}^1 \ltimes N^\prime$ is neither, so Proposition~\ref{prop:compact_center_closed_subgroup} yields that the center of $N^\prime$, which is equal to the center of $S^\prime$, is not compact, especially not a circle. If $S^\prime$ is closed in $\textnormal{Isom}(M)$, the center of $S^\prime$ is also closed in $\textnormal{Isom}(M)$. But the center of $S^\prime$ is generated by $Z$ and $Z$ generates a precompact one-parameter group since the subgroup generated by $\mathfrak{a}\oplus\mathfrak{z(s)}$ is compact. Thus, the center of $S^\prime$ is a circle.

\subsection{Case of the special linear algebra}

We have $\overline{\mathfrak{g}}=\mathfrak{k}\oplus\mathfrak{a}\oplus\mathfrak{s}$ with $\mathfrak{k}$ compact semisimple, $\mathfrak{a}$ abelian and $\mathfrak{s}\cong\mathfrak{sl}_2(\mathds{R})$. $\kappa|_{\mathfrak{s}\times\mathfrak{s}}$ is a Lorentzian scalar product.

Thus, $\overline{\mathfrak{g}}$ is reductive. By Proposition~\ref{prop:Maltsev}, $[\mathfrak{g},\mathfrak{g}]=\mathfrak{k}\oplus\mathfrak{s}$. Since $\mathfrak{r} \subseteq \overline{\mathfrak{r}}=\mathfrak{a}$, it follows that $\mathfrak{r}$ is abelian and $\mathfrak{g}$ is reductive. Due to Proposition~\ref{prop:reductive_algebra}, $\mathfrak{k}\oplus\mathfrak{s}\oplus\mathfrak{z(g)}=\mathfrak{g}$.

By a theorem of Wolf, $\widetilde{\textnormal{SL}_2(\mathds{R})}$ has no non-trivial compact subgroup (cf.~\cite{Wo63}). Hence, if the subgroup generated by $\mathfrak{s}$ is closed in $\textnormal{Isom}(M)$, it cannot be isomorphic to a finite central quotient of $\widetilde{\textnormal{SL}_2(\mathds{R})}$, since then all non-trivial one-parameter groups would be not precompact, contradicting condition~\hyperlink{star}{$(\star)$} ($\kappa|_{\mathfrak{s}\times\mathfrak{s}}$ is a Lorentz form).

Conversely, if the subgroup generated by $\mathfrak{s}$ is isomorphic to some $\textnormal{PSL}_k(2,\mathds{R})$, it is a semisimple group with finite center. By Proposition~\ref{prop:compact_center_closed_subgroup}, it is closed in $\textnormal{Isom}(M)$.

%%%%%%%%%%%%%%%%%%%%%%%%%%%%%%%%%%%%%%%%%%%%%%%%%%%%%%%%%%%%%%%%%%%%%%%%%%%%%%%%%%%%%%%%%%%%

\chapter{Geometric characterization of the manifolds}\label{ch:geometry}

In this chapter, we will prove Theorems~\ref{th:locally_free} and~\ref{th:geometric_characterization}. We will first prove the first part of Theorem~\ref{th:geometric_characterization} in Section~\ref{sec:kappa_semidefinite}. It turns out that our investigation will be helpful for the proof of Theorem~\ref{th:locally_free} in Section~\ref{sec:loc_free}. The proof of this theorem mainly relies on the algebraic structure of the involved Lie algebras and a remarkable property concerning the action of one-parameter groups with lightlike orbits (Lemma~\ref{lem:lightlike_no_fixed_points}).

In Section~\ref{sec:kappa_indefinite} we continue with the proof of the second and third part of Theorem~\ref{th:geometric_characterization}. We start with showing that the orbits of a central quotient of $\widetilde{\text{SL}_2(\mathds{R})}$ or a twisted Heisenberg group have Lorentzian character on the compact Lorentzian manifold they act isometrically and effectively on, in Section~\ref{sec:Lorentz_orbit}. Later, in Section~\ref{sec:orthogonal_distribution}, we investigate the distribution orthogonal to the orbits of $S$ and prove that it is (almost) involutive. This allows us to show the topological structure of the manifold as announced in the theorem. Also, the most part of the geometric result is proved in Section~\ref{sec:structure}. Finally, we look in Section~\ref{sec:Lorentz_metric_Heis} at the metric defined on $S$ if $S$ is a twisted Heisenberg group and finish the proof of the theorem.

\section{Induced bilinear form is positive semidefinite}\label{sec:kappa_semidefinite}

Let $M$ be a compact Lorentzian manifold and $G \subseteq \textnormal{Isom}(M)$ a connected closed non-compact Lie subgroup. We also consider its Lie algebra $\mathfrak{g}=\mathfrak{k}\oplus\mathfrak{a}\oplus\mathfrak{s}$ as described in Theorem~\ref{th:algebraic_classification}. As usual, denote by $\kappa$ the induced bilinear form on $\mathfrak{g}$.

Suppose that $\kappa$ is positive semidefinite. It follows from Corollary~\ref{cor:condition_star} and Theorem~\ref{th:algebraic}, that the summand $\mathfrak{s}$ is either trivial, isomorphic to $\mathfrak{aff}(\mathds{R})$ or to $\mathfrak{he}_d$. Additionally, the kernel of $\kappa$ has dimension at most one.

In the case that $\mathfrak{s}$ is trivial, $\mathfrak{a}$ generates the center of $G$. Since $G$ is non-compact and $\mathfrak{k}$ is compact semisimple, the center of $G$ is not compact. Thus, there is some $X \in \mathfrak{a}$ generating a non-precompact one-parameter group. For any $Y \in \mathfrak{g}$, $\mathfrak{p}:=\textnormal{span}\left\{X,Y\right\}$ is an abelian subalgebra and by Lemma~\ref{lem:kappa_abelian}, the set of elements of $\mathfrak{p}$ generating a non-precompact one-parameter group is dense in $\mathfrak{p}$. It follows that the set of elements in $\mathfrak{g}$ generating a non-precompact one-parameter group is dense in $\mathfrak{g}$. Proposition~\ref{prop:kappa_positive} yields that for all $X \in \mathfrak{g}$, the corresponding Killing vector field $\widetilde{X}$ (cf. Proposition~\ref{prop:isometry_Killing}) is everywhere non-timelike. Especially, if $X$ is $\kappa$-isotropic, then $\widetilde{X}$ is lightlike everywhere.

\begin{lemma}\label{lem:aff_heis_lightlike}
Let $G$ be a Lie group acting isometrically and locally effectively on a Lorentzian manifold $M$ of finite volume. Assume that the Lie algebra $\mathfrak{g}=\mathfrak{m}\oplus\mathfrak{s}$ is a direct sum of some arbitrary subalgebra $\mathfrak{m}$ and a subalgebra $\mathfrak{s}$, which is either isomorphic to $\mathfrak{aff}(\mathds{R})$ or $\mathfrak{he}_d$.

Then the induced bilinear form $\kappa$ is positive semidefinite, the orbits of $G$ are non-timelike everywhere and the orbit of any $\kappa$-isotropic $Z \in \mathfrak{g}$ is lightlike everywhere. Moreover, $\kappa(Z,Z)=0$ if $Z \in [\mathfrak{s},\mathfrak{s}]$.
\end{lemma}
\vspace{-1.4em}\begin{proof}
Without loss of generality, we may assume $G \subseteq \textnormal{Isom}(M)$. Especially, $\mathfrak{g}\subseteq\mathfrak{isom}(M)$.

First, suppose $\mathfrak{s} \cong \mathfrak{aff}(\mathds{R})$. Let $X,Y \in \mathfrak{s}$ such that $[X,Y]=Y$. For any $W \in \mathfrak{m}$ and $\lambda\in\mathds{R}$, \[[W+X+\lambda Y,Y]=Y,\] so by Lemma~\ref{lem:precompact_ad}, $W+X+\lambda Y$ generates a non-precompact one-parameter group in $\textnormal{Isom}(M)$. Thus, the set of $X \in \mathfrak{g}$ generating a non-precompact one-parameter group in $\textnormal{Isom}(M)$ is dense in $\mathfrak{g}$.

Corollary~\ref{cor:condition_star} yields that $\kappa$ is positive semidefinite and for all $X \in \mathfrak{g}$, the corresponding Killing vector field $\widetilde{X}$ is everywhere non-timelike by Proposition~\ref{prop:kappa_positive}. Clearly, if $Z$ is $\kappa$-isotropic, then $\widetilde{Z}$ is lightlike everywhere. Since $\kappa$ is ad-invariant, \[\kappa(Y,Y)=\kappa([X,Y],Y)=\kappa(X,[Y,Y])=0.\]

Now let $\mathfrak{s} \cong \mathfrak{he}_d$. Then for any non-central $X \in \mathfrak{he}_d$, there is $X^\prime \in \mathfrak{he}_d$, such that $[X,X^\prime]\neq 0$. But $[X,X^\prime]$ is central, hence, $\textnormal{ad}_X^2(X^\prime)=0$. So by Lemma~\ref{lem:precompact_ad}, $X$ generates a non-precompact one-parameter group in $\textnormal{Isom}(M)$.

If $Y \in \mathfrak{m}\oplus\mathfrak{z(s)}$, the subalgebra $\mathfrak{p}:=\textnormal{span}\left\{X,Y\right\}$ is abelian. Thus, by Lemma~\ref{lem:kappa_abelian}, the set of elements in $\mathfrak{p}$ generating a non-precompact one-parameter group in $\textnormal{Isom}(M)$ is dense in $\mathfrak{p}$. It follows that the set of elements in $\mathfrak{g}$ generating a non-precompact one-parameter group in $\textnormal{Isom}(M)$ is dense in $\mathfrak{g}$, and we can argue exactly as above.

Finally, if $Z\in\mathfrak{z(s)}$, there are $X,Y \in \mathfrak{he}_d$ such that $[X,Y]=Z$ and therefore, \[\kappa(Z,Z)=\kappa([X,Y],Z)=\kappa(X,[Y,Z])=0.\qedhere\]
\end{proof}\vspace{0pt}

The next lemma finishes the proof of Theorem~\ref{th:geometric_characterization}~(i).

\begin{lemma}\label{lem:lightlike_geodesic}
Let $(M,g)$ be a semi-Riemannian manifold and $\widetilde{X}$ a Killing vector field, such that $g(\widetilde{X},\widetilde{X})$ is constant. Then the orbits of the (local) flow of $\widetilde{X}$ are geodesics.
\end{lemma}
\vspace{-1.4em}\begin{proof}
We have to show that $\nabla_{\widetilde{X}} \widetilde{X}$ vanishes on $M$, where $\nabla$ denotes the Levi-Civita connection of $M$. For any vector field $Y$ on $M$, \[0\equiv Y(g(\widetilde{X},\widetilde{X}))\equiv 2g(\nabla_Y \widetilde{X},\widetilde{X}).\]

Since $\widetilde{X}$ is a Killing vector field, \[g(\nabla_{\widetilde{X}} \widetilde{X}, Y) \equiv -g(\widetilde{X}, \nabla_Y \widetilde{X}) \equiv 0.\] Thus, $\nabla_{\widetilde{X}} \widetilde{X} \equiv 0.$
\end{proof}\vspace{0pt}

\section{Locally free action}\label{sec:loc_free}

\begin{lemma}\label{lem:lightlike_no_fixed_points}
Let $M=(M,g)$ be a connected Lorentzian manifold and $\widetilde{X}$ a non-trivial lightlike complete Killing vector field. Then $\widetilde{X}$ vanishes nowhere, that is, $\widetilde{X}(x) \neq 0$ for all $x \in M$. Equivalently, the action of the corresponding one-parameter group $\psi:=\left\{\exp(tX)\right\}_{t \in \mathds{R}}\subseteq \textnormal{Isom}(M)$ (cf. Proposition~\ref{prop:isometry_Killing}) is locally free.
\end{lemma}
\vspace{-1.4em}\begin{proof}
Let $\psi^t:=\exp(tX)$ and assume, there is $x \in M$ such that $\widetilde{X}(x)=0$. Then $\psi^t(x)=x$ for all $t$ in $\mathds{R}$.

$g_x$ defines a quadratic form $Q:T_xM \to \mathds{R}$. Since $\psi$ preserves the Lorentzian metric $g$ on $M$, it follows that $d\psi_x=\left\{d\psi_x^t\right\}_{t \in \mathds{R}}$ preserves the scalar product $g_x$ as well as the quadratic form $Q$. Thus, $d\psi_x$ is a one-parameter group of isometries of $(T_xM,g_x)$. According to Proposition~\ref{prop:isometry_Killing}, it is associated to a complete Killing vector field $X^\prime$ on $T_xM$.

The isometry $\psi^t$ of $M$ is uniquely defined by $\psi^t(x)$ and $d\psi^t_x$. Since the action of $\psi$ is locally free, it follows that $d\psi^t_x \neq \textnormal{id}_{T_xM}$ for all $t \in \mathds{R}\backslash \Gamma$, where $\Gamma$ is a free subgroup in $\mathds{R}$ generated by at most one element. Thus, the flow $d\psi_x$ is non-trivial and $X^\prime$ does not vanish on an open neighborhood of $T_xM$ (because the Killing vector field $X^\prime$ is determined by $X^\prime(v)$ and $\nabla_v X^\prime$ for some $v \in T_xM$, $\nabla$ being the Levi-Civita connection on $M$).

Let $U \subset M$ be a star-shaped open neighborhood of $x$, on which the exponential map $\exp_x$ is invertible, and define $q:=Q \circ \exp_x^{-1}$.

Let $v \in \exp_x^{-1}(U)$. Then \[\gamma:[0,1]\to M, \ \gamma(\tau)=\exp_x(\tau v),\] is the unique geodesic with \[\gamma(0)=x, \ \frac{\partial}{\partial \tau}\gamma(\tau)|_{\tau=0}=v.\] $\psi^t$ is an isometry with fixed point $x$, so $(\psi^t \circ \gamma)$ is a geodesic with \[(\psi^t \circ \gamma)(0)=x \text{ and }\frac{\partial}{\partial \tau}(\psi^t \circ \gamma)(\tau)|_{\tau=0}=d\psi_x^t(v).\] Thus, we have \[\psi^t(\exp_x(v))=\psi^t(\gamma(1))=\exp_x(d\psi_x^t(v))\] for small $t$. Differentiating with respect to $t$ in $t=0$ yields
\begin{equation}\label{eq:Killing}
\widetilde{X}(\exp_x(v))=(d\exp_x)_v(X^\prime(v)).\tag{1}
\end{equation}

For $w \in T_v(T_xM)\cong T_xM$, \[g_x((\textnormal{grad }Q)(v),w)=\frac{\partial}{\partial t}Q(v+tw)|_{t=0}=\frac{\partial}{\partial t}g_x(v+tw,v+tw)|_{t=0}=g_x(2v,w).\] Hence, $\textnormal{grad }Q$ is radial. By the lemma of Gau{\ss},
\begin{align*}
&g_{\exp_x(v)}((d\exp_x)_v((\textnormal{grad }Q)(v)),(d\exp_x)_v(w))\\
=&g_x((\textnormal{grad }Q)(v),w)\\
=&dQ_v(w)\\
=&dq_{\exp_x(v)} \circ (d\exp_x)_v (w)\\
=&g_{\exp_x(v)}(((\textnormal{grad }q)(\exp_x(v))),(d\exp_x)_v(w)).\end{align*}
It follows that \begin{align}
g_x((\textnormal{grad }Q)(v),w)&=g_{\exp_x(v)}(((\textnormal{grad }q)(\exp_x(v))),(d\exp_x)_v(w)),\tag{2}\label{eq:grad1}\\
(d\exp_x)_v((\textnormal{grad }Q)(v))&=(\textnormal{grad }q)(\exp_x(v)).\tag{3}\label{eq:grad2}
\end{align}
Choosing $w=X^\prime(v)$ and later $w=(\textnormal{grad }Q)(v)$ in Equation~(\ref{eq:grad1}), we obtain successively with the help of Equation~(\ref{eq:Killing})
\begin{align}
g_{\exp_x(v)}((\textnormal{grad }q)(\exp_x(v)),\widetilde{X}(\exp_x(v)))&=g_x((\textnormal{grad }Q)(v),X^\prime(v)),\tag{4}\label{eq:grad-Killing1}\\
g_{\exp_x(v)}((\textnormal{grad }q)(\exp_x(v)),(\textnormal{grad }q)(\exp_x(v)))&=g_x((\textnormal{grad }Q)(v),(\textnormal{grad }Q)(v)).\tag{5}\label{eq:grad3}\end{align}
Furthermore,
\begin{align*}
g_x((\textnormal{grad }Q)(v),X^\prime(v))&=dQ_v(X^\prime(v))\\
&=dQ_v\left(\frac{\partial}{\partial t}\left(d\psi_x^t(v)\right)|_{t=0}\right)\\
&=\frac{\partial}{\partial t} \left( (Q \circ d\psi_x^t) (v)\right)|_{t=0}=0,\end{align*}
because $d\psi_x$ preserves $Q$.

It follows together with Equation~(\ref{eq:grad-Killing1}), that
\begin{align}
0=g_{u}((\textnormal{grad }q)(u),\widetilde{X}(u)) \tag{6}\label{eq:grad-Killing2}
\end{align} for all $u \in U$.

By assumption, $0=g_u(\widetilde{X}(u),\widetilde{X}(u))$ since $\widetilde{X}$ is lightlike. Now choose a timelike $v \in \exp_x^{-1}(U)$ such that $X^\prime(v) \neq 0$ (remember that on any open neighborhood, $X^\prime$ does not vanish everywhere). Then because of Equation~(\ref{eq:Killing}), $\widetilde{X}(u) \neq 0$ for $u=\exp_x(v)$, but the radial vector $(\textnormal{grad }Q)(v)$ and so $(\textnormal{grad }q)(u)$ by Equation~(\ref{eq:grad3}) are timelike. Since $(\textnormal{grad }q)(u)$ is in the $g_u$-orthogonal complement of the lightlike vector $\widetilde{X}(u)$ by Equation~(\ref{eq:grad-Killing2}), $(\textnormal{grad }q)(u)$ is not timelike, contradiction.
\end{proof}\vspace{0pt}

\begin{remark}
Lemma~\ref{lem:lightlike_no_fixed_points} is also true if $\widetilde{X}$ is not complete, in this case one consider the local flow $\psi$ defined by $\widetilde{X}$.
\end{remark}

\begin{definition}
Let $G$ be a Lie group with Lie algebra $\mathfrak{g}$ acting continuously on a compact topological space $M$. For $p \in M$ denote by $G_p$ the subgroup of $G$ fixing $p$ and $\mathfrak{g}_p$ its Lie algebra.
\end{definition}
\begin{remark}
$G_p$ is a closed subgroup, especially a Lie subgroup. Thus, its Lie algebra is defined.
\end{remark}

\begin{lemma}\label{lem:stabilizer}
Let $G$ be a connected Lie group with Lie algebra $\mathfrak{g}$ acting continuously on a compact topological space $M$.

Suppose there are $x \in M$, $X \in \mathfrak{g}, Y \in \mathfrak{g}_x$ and $k \in \mathds{Z}_+$, such that $[X,Z]=0$ for $Z=\textnormal{ad}_X^k(Y)$. Then there exists $y \in M$ such that $Z \in \mathfrak{g}_y$.
\end{lemma}
\vspace{-1.4em}\begin{proof}
It suffices to consider $Z \neq 0$. Define \[Y_j:=\textnormal{ad}_X^j(Y)\] for any non-negative integer $i$. By assumption, $Y_k=Z$ and $Y_j=0$ if $j>k$.

Let $t \in \mathds{R}$. Since $Y \in \mathfrak{g}_x$, $\exp(\tau Y) \cdot x=x$, so it holds for all $\tau \in \mathds{R}$:
\[\exp(tX) \exp(\tau Y) \exp(-tX) \cdot (\exp(tX) \cdot x)=\exp (tX) \cdot x.\]

This is equivalent to \[\exp(tX) \exp(\tau Y) \exp(-tX) \in G_{\exp(tX) \cdot x}.\] Differentiating with respect to $\tau$ in $\tau=0$ yields \[\textnormal{Ad}_{\exp(tX)} Y \in \mathfrak{g}_{\exp(tX) \cdot x}.\] Moreover,
\[\textnormal{Ad}_{\exp(tX)} Y=\exp(\textnormal{ad}_{tX}) Y=\sum\limits_{j=0}^{k-1} \left(\frac{t^j}{j!} Y_j \right)+\frac{t^k}{k!} Z.\]

It follows that \[\mathfrak{g}_{\exp(tX) \cdot x} \ni \frac{k!}{t^k} \textnormal{Ad}_{\exp(tX)} Y=Z+\sum\limits_{j=0}^{k-1} \left(\frac{k!}{j!t^{k-j}} Y_j\right)=:Z_t\] for all $t \in \mathds{R}$. Since $M$ is compact, we can choose an increasing sequence $\left\{t_n\right\}_{n=0}^\infty$ of real numbers, such that $t_n \to \infty$ and \[\exp(t_n X) \cdot x \to y \in M\] as $n \to \infty$.

Since
\begin{align*}
Z=\lim\limits_{n\to\infty} \left(Z+\sum\limits_{j=0}^{k-1} \left(\frac{k!}{j!t_n^{k-j}} Y_j\right)\right)&=\lim\limits_{n\to\infty} Z_{t_n},\\
\exp(\tau Z_{t_n}) \cdot (\exp(t_n X) \cdot x) &\to \exp(\tau Z) \cdot y \end{align*} as $n \to \infty$ for all $\tau \in \mathds{R}$. But for all $n$, it holds \[\exp(\tau Z_{t_n}) \cdot (\exp(t_n X) \cdot x) = \exp(t_n X) \cdot x.\] It follows that $\exp(\tau Z) \cdot y = y$ for all $\tau \in \mathds{R}$, so $Z \in \mathfrak{g}_y$.
\end{proof}\vspace{0pt}

\begin{corollary}\label{cor:affR_locally_free}
Let $M$ be a compact Lorentzian manifold and $S$ a Lie group with Lie algebra $\mathfrak{s} \cong \mathfrak{aff}(\mathds{R})$ acting isometrically and locally effectively on $M$. Then this action is locally free. 
\end{corollary}
\vspace{-1.4em}\begin{proof}
Let $X,Y \in \mathfrak{s}$ such that $[X,Y]=Y$. By Lemma~\ref{lem:aff_heis_lightlike}, the subgroup generated by $Y$ has lightlike orbits, and due to Lemma~\ref{lem:lightlike_no_fixed_points}, it acts locally freely. Especially, $Y \notin \mathfrak{s}_x$ for all $x \in M$.

Suppose $X+\lambda Y \in \mathfrak{s}_x$, $\lambda \in \mathds{R}$, for some $x \in M$. Since \[[Y,X+\lambda Y]=-Y \text{ and } [Y,-Y]=0,\] we can apply Lemma~\ref{lem:stabilizer} and obtain $-Y \in \mathfrak{s}_y$ for some $y \in M$, contradiction.

Thus, $\mathfrak{s}_x$ is trivial for all $x \in M$, that is, $S$ acts locally freely on $M$.
\end{proof}\vspace{0pt}

\begin{corollary}\label{cor:sl2R_locally_free}
Let $M$ be a compact Lorentzian manifold and $S$ a Lie group with Lie algebra $\mathfrak{s} \cong \mathfrak{sl}_2(\mathds{R})$ acting isometrically and locally effectively on $M$. Then this action is locally free.
\end{corollary}
\vspace{-1.4em}\begin{proof}
Without loss of generality, $S \subseteq \textnormal{Isom}(M)$. Let $e,f,h \in \mathfrak{s}$ be elements of an $\mathfrak{sl}_2$-triple, that is, $[h,e]=2e$, $[h,f]=-2f$ and $[e,f]=h$.

Let $\kappa$ be the induced bilinear form on $\mathfrak{s}$. \[2\kappa(e,e)=\kappa([h,e],e)=\kappa(h,[e,e])=0\] because $\kappa$ is ad-invariant. Analogously, \[\kappa(f,f)=0.\] Because of $[e,h]=-2e$ and $[e,-2e]=0$, $e$ generates a non-precompact one-parame\-ter group in $\textnormal{Isom}(M)$ by Lemma~\ref{lem:precompact_ad}. The same is true for $f$. It follows from Proposition~\ref{prop:kappa_positive}, that the subgroups generated by $e$ and $f$ have lightlike orbits everywhere. By Lemma~\ref{lem:lightlike_no_fixed_points}, $e,f \notin \mathfrak{s}_x$ for all $x \in M$.

Let \[0\neq Y:=\alpha h + \beta e + \gamma f\] for some real $\alpha,\beta,\gamma$. We have to show that $Y \notin \mathfrak{s}_x$ for all $x \in M$. If $\alpha=\gamma=0$, we have shown this already.

Assume $\gamma \neq 0$ and choose $X:=e$, $Z:=-2 \gamma e$. Then \[[X,Y]=-2\alpha e+\gamma h, \ [X,[X,Y]]=Z \text{ and }[X,Z]=0.\] Lemma~\ref{lem:stabilizer} and $Z \notin \mathfrak{s}_x$ for all $x$ yield the required result.

Finally, assume $\gamma=0$, but $\alpha \neq 0$. Choose $X:=e$ and $Z:=-2\alpha e$. Then \[[X,Y]=Z \text{ and }[X,Z]=0,\] so we can apply Lemma~\ref{lem:stabilizer} to see that $Y \notin \mathfrak{s}_x$ for all $x \in M$, since $Z \notin \mathfrak{s}_x$ for all $x$.
\end{proof}\vspace{0pt}

\begin{lemma}\label{lem:nilpotent_locally_free}
Let $G$ be a connected nilpotent Lie group acting continuously on a compact manifold $M$. If the action of the center is locally free, so is the action of $G$.
\end{lemma}
\vspace{-1.4em}\begin{proof}
We may assume that $G$ is non-trivial. Let $\mathfrak{g}_0:=\mathfrak{g}$, $\mathfrak{g}_j:=[\mathfrak{g},\mathfrak{g}_{j-1}]$ for $j>0$ be the descending central series of $\mathfrak{g}$. Since $\mathfrak{g}$ is nilpotent, there is an integer $k \geq 0$, such that $\mathfrak{g}_k \neq \left\{0\right\}=\mathfrak{g}_{k+1}$.

Since $\left\{0\right\}=[\mathfrak{g},\mathfrak{g}_k]$, $\mathfrak{g}_k \subseteq \mathfrak{z(g)}$. Thus, the subgroup generated by $\mathfrak{g}_k$ acts locally freely by assumption. In the following, we prove by induction, that the action of the subgroup generated by $\mathfrak{g}_{j-1}$ is locally free, if the action of the subgroup generated by $\mathfrak{g}_{j}$ is.

Assume the contrary, that is, there is $x \in M$ and $Y \in \mathfrak{g}_{j-1}$ such that $Y \in \mathfrak{g}_x$. Since $Y \notin \mathfrak{z(g)}$, there is $X \in \mathfrak{g}$ such that $[X,Y]\neq 0$. $\mathfrak{g}$ is nilpotent, so there exists an integer $l>0$ such that \[Z:=\textnormal{ad}_X^l (Y)\neq0=\textnormal{ad}_X^{l+1} (Y).\] By Lemma~\ref{lem:stabilizer}, $Z \in \mathfrak{g}_y$ for some $y \in M$. But $Z \in \mathfrak{g}_{j-1+l} \subseteq \mathfrak{g}_j$, contradiction.
\end{proof}\vspace{0pt}

We can now finish the proof of Theorem~\ref{th:locally_free} by showing that the action of the subgroup generated by $\mathfrak{s}$ is locally free, if $\mathfrak{s} \cong \mathfrak{he}_d$ or $\mathfrak{s} \cong \mathfrak{he}_d^\lambda$.

\begin{corollary}\label{cor:heis_locally_free}
Let $M$ be a compact Lorentzian manifold and $S$ a Lie group with Lie algebra $\mathfrak{s} \cong \mathfrak{he}_d$, acting isometrically and locally effectively on $M$. Then this action is locally free. 
\end{corollary}
\vspace{-1.4em}\begin{proof}
By Lemma~\ref{lem:aff_heis_lightlike}, the center of $\mathfrak{z(s)}$ is $\kappa$-isotropic and the subgroup generated by $0 \neq Z \in \mathfrak{z(s)}$ has lightlike orbits. Due to Lemma~\ref{lem:lightlike_no_fixed_points}, it acts locally free. The result now follows from Lemma~\ref{lem:nilpotent_locally_free}.
\end{proof}\vspace{0pt}

\begin{corollary}\label{cor:twisted_heis_locally_free}
Let $M$ be a compact Lorentzian manifold and $S$ a Lie group with Lie algebra $\mathfrak{s} \cong \mathfrak{he}_d^\lambda$, acting isometrically and locally effectively on $M$. Then this action is locally free.
\end{corollary}
\vspace{-1.4em}\begin{proof}
We identify $\mathfrak{he}_d^\lambda$ with $\mathfrak{s}$. By Corollary~\ref{cor:heis_locally_free}, the subgroup generated by $\mathfrak{he}_d \subset \mathfrak{he}_d^\lambda \cong \mathfrak{s}$ acts locally freely.

Let $Y \in \mathfrak{he}_d^\lambda \backslash \mathfrak{he}_d$. Then there is $X \in \mathfrak{he}_d$ such that $[X,Y] \neq 0$ (otherwise $Y \in \mathfrak{z(s)}\subset\mathfrak{he}_d$, contradiction). But $[X,Y] \in \mathfrak{he}_d$ and $\mathfrak{he}_d$ is nilpotent, hence \[Z:=\textnormal{ad}_X^k(Y)\neq 0=\textnormal{ad}_X^{k+1}(Y)\] for some integer $k>0$. The subgroup generated by $Z \in \mathfrak{he}_d$ acts locally freely by Corollary~\ref{cor:heis_locally_free}, therefore, for all $x \in M$, $Y \notin \mathfrak{s}_x$ due to Lemma~\ref{lem:stabilizer}.
\end{proof}\vspace{0pt}

\section{Induced bilinear form is indefinite}\label{sec:kappa_indefinite}

In this section, let $M=(M,g)$ be a compact Lorentzian manifold and $G$ a connected closed non-compact Lie subgroup of $\textnormal{Isom}(M)$. According to Theorem~\ref{th:algebraic_classification}, its Lie algebra is a direct sum $\mathfrak{g}=\mathfrak{k}\oplus\mathfrak{a}\oplus\mathfrak{s}$ with $\mathfrak{k}$ compact semisimple and $\mathfrak{a}$ abelian. We assume that $\mathfrak{s}$ is either isomorphic to $\mathfrak{sl}_2(\mathds{R})$ or to a twisted Heisenberg algebra $\mathfrak{he}_d^\lambda$ with $\lambda \in \mathds{Z}_+^d$.

We identify $\mathfrak{sl}_2(\mathds{R})$ and $\mathfrak{he}_d^\lambda$ with $\mathfrak{s}$, respectively. In the first case, we choose an $\mathfrak{sl}_2$-triple $\left\{e,f,h\right\}$, and we choose a canonical basis $\left\{T,Z,X_1,Y_1,\ldots,X_d,Y_d\right\}$ of $\mathfrak{s}$ in the latter case.

Denote by $S$ the subgroup generated by $\mathfrak{s}$. By Theorem~\ref{th:algebraic_classification}, $S \cong \textnormal{PSL}_k(2,\mathds{R})$ if and only if $S$ is closed in $\textnormal{Isom}(M)$ in the first case and in the latter case, $S \cong \textnormal{He}_d^\lambda$ if $S$ is closed in $\textnormal{Isom}(M)$ and $S \cong \overline{\textnormal{He}_d^\lambda}$ otherwise.

\subsection{Lorentzian character of orbits}\label{sec:Lorentz_orbit}

\begin{lemma}\label{lem:light_cone_multiple}
Let $V$ be a real vector space, $b_1$ a symmetric bilinear form and $b_2$ a Lorentzian scalar product on $V$. If any element of the light cone of $b_2$ is $b_1$-isotropic, $b_1=\lambda b_2$ for some real number $\lambda$.
\end{lemma}
\vspace{-1.4em}\begin{proof}
Choose a $b_2$-orthonormal basis $\left\{v_0, \ldots, v_m\right\}$ of $V$ such that $b_2(v_0,v_0)=-1$. In what follows, let $j, k \in \left\{1,\ldots,m\right\}$, $j\neq k$.

Then \[b_2(v_0 + v_j,v_0+v_j)=0=b_2(v_0 - v_j, v_0-v_j).\] By assumption, it follows that \[b_1(v_0 + v_j,v_0+v_j)=0=b_1(v_0 - v_j, v_0-v_j).\] Especially, $b_1(v_0,v_j)=0$. Using this, $\lambda:=-b_1(v_0,v_0)=b_1(v_j,v_j)$.
\[b_2(\sqrt{2}v_0 + (v_j+v_k),\sqrt{2}v_0 + (v_j+v_k))=2b_2(v_0,v_0)+b_2(v_j,v_j)+b_2(v_k,v_k)=0\] yields that \[b_1(\sqrt{2}v_0 + (v_j+v_k),\sqrt{2}v_0 + (v_j+v_k))=0.\] Since \[b_1(v_0,v_j)=0=b_1(v_0,v_k),\] we have \[2b_1(v_0,v_0)+b_1(v_j,v_j)+b_1(v_k,v_k)+2b_1(v_j,v_k)=0.\] But we know already \[2b_1(v_0,v_0)+b_1(v_j,v_j)+b_1(v_k,v_k)=-2\lambda+\lambda+\lambda=0.\] Therefore, $b_1(v_j,v_k)=0$.

In summary, $\left\{v_0, \ldots, v_m\right\}$ is a $b_1$-orthogonal basis of $V$. Also, it holds \[-b_1(v_0,v_0)=\lambda= b_1(v_j,v_j).\] Thus, $b_1=\lambda b_2$.
\end{proof}\vspace{0pt}

\begin{proposition}\label{prop:orbit_Lorentz}
In the situation of this section, the following is true:
\begin{compactenum}
\item The orbits of $S$ have Lorentzian character everywhere on $M$ if $\mathfrak{s} \cong \mathfrak{sl}_2(\mathds{R})$.

\item The orbits of $S$ have Lorentzian character everywhere on $M$ if $\mathfrak{s} \cong \mathfrak{he}_d^\lambda$.
\end{compactenum}
\end{proposition}
\vspace{-1.4em}\begin{proof}
Let $x \in M$ and consider the linear map \[\iota : \mathfrak{s} \to T_xM, \ X \mapsto \widetilde{X}(x)\] (cf. Proposition~\ref{prop:isometry_Killing}). Due to Corollaries~\ref{cor:sl2R_locally_free} and~\ref{cor:twisted_heis_locally_free}, the action of $S$ is locally free, hence, $\iota$ is injective. Let $b$ be the symmetric bilinear form on $\mathfrak{s}$ defined by \[b(X,Y):=g_x(\widetilde{X}(x),\widetilde{Y}(x)).\] The character of the orbit of $S$ in $x$ is Lorentzian if and only if $b$ is a Lorentz form.

(i) By Theorem~\ref{th:algebraic}~(v), the induced bilinear form $\kappa$ is a positive multiple of the Killing form $k$ on $\mathfrak{s}$. Let \[0\neq X:=\alpha e+ \beta f+\gamma h\] for real parameters $\alpha,\beta,\gamma$. With respect to the ordered basis $(e,f,h)$, $\textnormal{ad}_X$ corresponds to the $3 \times 3$-matrix \[A:=\begin{pmatrix} 2\gamma & 0 &-2\alpha \\ 0 & -2\gamma & 2\beta\\ -\beta & \alpha & 0 \end{pmatrix}\] An elementary calculation yields \begin{align*} A^2&=\begin{pmatrix} 4\gamma^2+2\alpha\beta & -2\alpha^2 &-4\alpha\gamma \\ -2\beta^2 & 4\gamma^2+2\alpha\beta & -4\beta\gamma\\ -2\beta\gamma & -2\alpha\gamma & 4\alpha\beta \end{pmatrix}\\ \text{and} \quad A^3&=\begin{pmatrix} 8\gamma(\alpha\beta+\gamma^2) & 0 &-8\alpha(\alpha\beta+\gamma^2) \\ 0 & -8\gamma(\alpha\beta+\gamma^2) & 8\beta(\alpha\beta+\gamma^2)\\ -4\beta(\alpha\beta+\gamma^2) & 4\alpha(\alpha\beta+\gamma^2) & 0 \end{pmatrix}.\end{align*}
It follows that \[k(X,X)=\textnormal{Tr} (\textnormal{ad}_X^2)=8 \gamma^2+8\alpha \beta.\] Thus, $X$ is $\kappa$-isotropic if and only if $\gamma^2=-\alpha \beta$. Moreover, $\textnormal{ad}_X^3 =0$ if and only if $\gamma^2=-\alpha \beta$.

In summary, $\textnormal{ad}_X$ is nilpotent (but non-trivial), if and only if $X\neq 0$ is $\kappa$-isotropic. By Lemma~\ref{lem:precompact_ad}, it follows that all non-trivial $\kappa$-isotropic $X$ generate non-precompact one-parameter groups. Due to Proposition~\ref{prop:kappa_positive}, $b(X,X)=0$ if $\kappa(X,X)=0$. The latter is equivalent to $k(X,X)=0$. We can apply Lemma~\ref{lem:light_cone_multiple} to see that $b=\lambda k$. Because $b$ is the restriction of a Lorentzian scalar product, $\lambda >0$. 

(ii) It follows from Lemma~\ref{lem:aff_heis_lightlike}, that $b$ is positive semidefinite on $\mathfrak{he}_d \times \mathfrak{he}_d$ and $b(Z,Z)=0$. $g_x$ is a Lorentzian scalar product and $\iota$ is injective, therefore, $b(X,X)>0$ for all $X \in \mathfrak{he}_d \backslash \mathfrak{z(s)}$. It follows that we do not have a Lorentzian orbit if and only if $b$ is positive semidefinite.

Assume $b$ is positive semidefinite. Then the kernel of $b$ has dimension at most one. But $\mathfrak{z(s)}$ has to lie in the kernel: \[0 \leq b(X+\lambda Z,X+\lambda Z)=b(X,X)+2\lambda(X,Z)\] for all $\lambda \in \mathds{R}$ and $X \in \mathfrak{s}$ requires $b(X,Z)=0$ for all $X \in \mathfrak{s}$.

Now let $D$ be the set of $y \in M$ where the orbit of $S$ is not Lorentzian. $D$ is non-empty, because $x \in D$. By continuity, $D$ is closed in $M$, hence $D$ is a compact metric space. Furthermore, $D$ is $S$-invariant and $S$ acts continuously on $D$.

$S$ is amenable as a solvable group (cf.~\cite{Zi86}, Corollary~4.1.7), that is, any continuous $S$-action on a compact metrizable space has an $S$-invariant probability measure (cf.~\cite{Zi86}, Definition~4.1.1). Let $\mu$ be the $S$-invariant probability measure on $D$. $\mu$ allows us to construct on $\mathfrak{s}$ the symmetric bilinear form \[B(X,Y):=\int\limits_D g_y (\widetilde{X}(y),\widetilde{Y}(y)) d\mu(y).\]

By construction, $B$ is ad-invariant (compare with Lemma~\ref{lem:kappa_invariant}~(iii)) and positive semidefinite. Its kernel is equal to $\mathfrak{z(s)}$. Thus, the solvable algebra $\mathfrak{s}/\mathfrak{z(s)}$ possesses a positive definite scalar product and due to Propositions~\ref{prop:reductive_algebra} and~\ref{prop:compact_algebra}, $\mathfrak{s}/\mathfrak{z(s)}$ is abelian. Therefore, $\mathfrak{he}_d=[\mathfrak{s},\mathfrak{s}]\subseteq\mathfrak{z(s)}$, contradiction.
\end{proof}\vspace{0pt}

\subsection{Orthogonal distribution}\label{sec:orthogonal_distribution}

The distribution defined by the orbits of $S$ will be denoted by $\mathcal{S}$ and by $\mathcal{O}$ we denote the distribution on $M$ orthogonal to $\mathcal{S}$. Note that $\mathcal{S}$ has the same dimension as $S$, because the action of $S$ is locally free due to Corollaries~\ref{cor:sl2R_locally_free} and~\ref{cor:twisted_heis_locally_free}. Since the metric $g$ restricted to $\mathcal{S}\times\mathcal{S}$ is Lorentzian by Proposition~\ref{prop:orbit_Lorentz}, $g$ restricted to $\mathcal{O}\times\mathcal{O}$ is Riemannian. Additionally, denote by $\mathcal{Z}$ the orbits of the center of $S$ in the case that $S$ is a twisted Heisenberg group.

\begin{proposition}\label{prop:integrable_distribution}
The distributions $\mathcal{O}$ for $\mathfrak{s}\cong\mathfrak{sl}_2(\mathds{R})$ and $\mathcal{O}+\mathcal{Z}$ for $\mathfrak{s}\cong\mathfrak{he}_d^\lambda$, respectively, are involutive.
\end{proposition}
\vspace{-1.4em}\begin{proof}
For any point $x$ in the manifold $M$, we define a vector-valued symmetric bilinear form $\omega_x: \mathcal{O}_x \times \mathcal{O}_x \to \mathcal{S}_x$ by \[\omega_x(v,w)=\textnormal{proj}_{\mathcal{S}_x}\left([\widetilde{V},\widetilde{W}]\left(x\right)\right),\] where $\textnormal{proj}_{\mathcal{S}_x}$ denotes the orthogonal projection to $\mathcal{S}_x$ and $\widetilde{V}$ and $\widetilde{W}$ are vector fields in $\mathcal{O}$ extending $v$ and $w$, respectively. Let $\omega:\mathcal{O}\times\mathcal{O}\to\mathcal{S}$ be the form defined by $\left\{\omega_x\right\}_{x \in M}$. Note that $\omega$ is correctly defined, since it is $C^\infty (M)$-linear: For any $\sigma \in C^\infty(M)$ and vector fields $\widetilde{V}$, $\widetilde{W}$ in $\mathcal{O}$, \[\omega(\sigma \widetilde{V}, \widetilde{W})=\textnormal{proj}_{\mathcal{S}}\left([\sigma\widetilde{V},\widetilde{W}]\right)=\textnormal{proj}_{\mathcal{S}}\left(\sigma[\widetilde{V},\widetilde{W}]-\widetilde{W}(\sigma)\widetilde{V}\right)=\sigma\textnormal{proj}_{\mathcal{S}}\left([\widetilde{V},\widetilde{W}]\right).\] In the same way, $\omega( \widetilde{V},\sigma \widetilde{W})=\sigma\omega( \widetilde{V}, \widetilde{W})$.

$\mathcal{O}$ and $\mathcal{O}+\mathcal{Z}$, respectively, are involutive if and only if $\omega$ is trivial or takes only values in $\mathcal{Z}$, respectively.

Using the action of $S$ on $M$, we can identify the vector spaces $\mathcal{S}_x$ and $\mathfrak{s}$ in a canonical way: $X \in \mathfrak{s}$ corresponds to $\frac{\partial}{\partial t}(\exp(tX) \cdot x) |_{t=0} \in T_xM$ (compare with Proposition~\ref{prop:isometry_Killing}). Thus, we can consider $\omega:\mathcal{O}\times\mathcal{O}\to\mathfrak{s}$.

Let $\widetilde{V}$ and $\widetilde{W}$ be vector fields in $\mathcal{O}$ extending $v,w \in \mathcal{O}_x$ and $f \in S$. Then
\begin{align*}
\omega_{f \cdot x}(df_x(v),df_x(w))&=\textnormal{proj}_{\mathcal{S}_{f \cdot x}}\left([df(\widetilde{V}),df(\widetilde{W})](f \cdot x)\right)\\
&=\textnormal{proj}_{\mathcal{S}_{f \cdot x}}\left(df_x\left([\widetilde{V},\widetilde{W}] ( x)\right)\right)\\
&=df_x\left(\textnormal{proj}_{\mathcal{S}_x}\left([\widetilde{V},\widetilde{W}](x)\right)\right)\\
&=df_x(\omega_x(v,w)),\end{align*}
where we have used that $\textnormal{proj}_{\mathcal{S}_{f \cdot x}}\circ df_x=df_x\circ\textnormal{proj}_{\mathcal{S}_x}$, which is true because the isometry $f$ preserves the distributions $\mathcal{S}$ and $\mathcal{O}=\mathcal{S}^\perp$.

It follows by Lemma~\ref{lem:kappa_invariant}~(i) that $df_x(\omega_x(v,w))\mathop{\widehat{=}}\textnormal{Ad}_f(\omega(v,w))$.

Since $g$ restricted to $\mathcal{O} \times \mathcal{O}$ is Riemannian and the action of $S$ preserves $g$ and $\mathcal{O}$ as well, $df_x(v) \in K(v)$, where $K(v):=\left\{u \in \mathcal{O} | g(u,u)=g_x(v,v)\right\}$. Because $M$ is compact, $K(v)$ is compact as well. Thus, the set $\left\{(df_x(v),df_x(w))\right\}_{f \in S}$ contained in $K(v)\times K(w)$ is precompact. $\omega$ is continuous, so $\left\{\textnormal{Ad}_f(\omega(v,w))\right\}_{f \in S} \subset \mathfrak{s}$ has to be precompact as well.

Assume $\mathfrak{s}\cong\mathfrak{sl}_2(\mathds{R})$ and let $X \in \left\{ e,f\right\}$, so $\textnormal{ad}_X$ is three-step nilpotent. Then the set consisting of all
\begin{align*}
\textnormal{Ad}_{\exp(tX)}(\omega(v,w))&=\exp(\textnormal{ad}_{tX})(\omega(v,w))\\
&=\omega(v,w)+t[X,\omega(v,w)]+\frac{t^2}{2}[X,[X,\omega(v,w)]],\end{align*}
$t \in \mathds{R}$, can only be precompact, if $[X,\omega(v,w)]=0$. But \[[e, \alpha e +\beta f +\gamma h]=\beta h -2 \gamma e \text{ and }[f, \alpha e +\beta f +\gamma h]=-\alpha h +2 \gamma f\] for any real $\alpha,\beta,\gamma$. Thus, $\omega(v,w) =0$.

Suppose $\mathfrak{s}\cong\mathfrak{he}_d^\lambda$ and let $X \in \mathfrak{he}_d$. Then $\textnormal{ad}_X^3=0$ and the set consisting of all
\begin{align*}
\textnormal{Ad}_{\exp(tX)}(\omega(v,w))&=\exp(\textnormal{ad}_{tX})(\omega(v,w))\\
&=\omega(v,w)+t[X,\omega(v,w)]+\frac{t^2}{2}[X,[X,\omega(v,w)]],\end{align*}
$t \in \mathds{R}$, can only be precompact, if $[X,\omega(v,w)]=0$. Since this is true for all $X \in \mathfrak{he}_d$, $\omega(v,w)$ has to lie in the center, that is, $\omega(v,w) \in \mathds{R}Z$.
\end{proof}\vspace{0pt}

\subsection{Structure of the manifold}\label{sec:structure}

In this section, we want to show that the compact manifold $M$ is diffeomorphic to $\Gamma\backslash\mkern-5mu\left(N \times S\right)$ if $\mathfrak{s}\cong\mathfrak{sl}_2(\mathds{R})$ and to $\Gamma\backslash\mkern-5mu\left(S \times_{Z(S)} N\right)$ if $\mathfrak{s}\cong\mathfrak{he}_d^\lambda$, respectively. For this, we choose auxiliary Riemannian metrics on $N \times S$ and $S \times_{Z(S)} N$, respectively.  $\Gamma$ will be a certain discrete subgroup in the corresponding isometry group of these spaces.

Let $\beta$ be any positive definite scalar product on $\mathfrak{s}$. We furnish $S$ with the right-invariant metric $\beta_S$ defined by $\beta$.

We transform the Lorentzian metric $g$ on $M$ into a Riemannian metric $g^\prime$ as follows: $g^\prime$ and $g$ coincide on $\mathcal{O}\times\mathcal{O}$ and $g^\prime$ on $\mathcal{S} \times \mathcal{S}$ is given by $\beta$ and the identification of $\mathcal{S}_x$ with $\mathfrak{s}$ as in Section~\ref{sec:orthogonal_distribution}. Additionally, $\mathcal{O}$ is orthogonal to $\mathcal{S}$ in both metrics.

\begin{lemma}\label{lem:S_covering}
For any $x\in M$, the canonical mapping $S \to S \cdot x$, $f \mapsto f\cdot x$, is an isometric covering map. $S \cdot x$ is furnished with the metric induced by the ambient space $(M,g^\prime)$.
\end{lemma}
\vspace{-1.4em}\begin{proof}
Since the action of $S$ on $M$ is locally free by Corollaries~\ref{cor:sl2R_locally_free} and~\ref{cor:twisted_heis_locally_free}, $S \to S \cdot x$ is a covering map.

To show that the mapping is isometric, it suffices to consider right-invariant vector fields on $S$. Let $X,Y \in \mathfrak{s}$ and $\overline{X},\overline{Y}$ the corresponding right-invariant vector fields on $S$. By construction, $\beta_S(\overline{X},\overline{Y})\equiv \beta(X,Y)$. Since $\overline{X}$ is right-invariant, $\overline{X}(f)=\frac{\partial}{\partial t}\left(\exp(tX)f\right)|_{t=0}$ for any $f\in S$. Thus, the vector field $\overline{X}$ on $S$ corresponds to the vector field $\widetilde{X}$ on $S \cdot x$ (cf.~Proposition~\ref{prop:isometry_Killing}). The result follows from $g^\prime(\widetilde{X},\widetilde{Y})\equiv \beta(X,Y)$.
\end{proof}\vspace{0pt}

Since $\mathcal{O}$ or $\mathcal{O}+\mathcal{Z}$, respectively, are involutive by Proposition~\ref{prop:integrable_distribution}, the corresponding distribution induces a foliation on $M$ due to the Frobenius theorem (cf.~\cite{Ba09}, Proposition~A.9). Let $N$ be the leaf of $\mathcal{O}$ or $\mathcal{O}+\mathcal{Z}$, respectively, passing through a point $x_0 \in M$. We furnish $N$ with the metric $h$ induced by $g^\prime$.

\begin{lemma}\label{lem:action_center_N}
Suppose $\mathfrak{s}\cong\mathfrak{he}_d^\lambda$. Then $Z(S)$, the center of $S$, acts isometrically on $N=(N,h)$.
\end{lemma}
\vspace{-1.4em}\begin{proof}
Let $\lambda \in \mathds{R}$ and $x \in N \subset M$. The curve $t \mapsto \exp(t\lambda Z)\cdot x$ has its tangential vector in $\mathcal{Z}$, thus, the action of the center of $S$ on $M$ reduces to an action on $N$. It is isometric by Proposition~\ref{lem:S_covering}.
\end{proof}\vspace{0pt}

\begin{proposition}\label{prop:submersion}
As above, let $M=(M,g^\prime)$. Then the following is true:
\begin{compactenum}
\item Let $\mathfrak{s}\cong\mathfrak{sl}_2(\mathds{R})$. The mapping $p: S \times N \to M$, $p(f,x)=f\cdot x$, is surjective and a local isometry.

\item Let $\mathfrak{s}\cong\mathfrak{he}_d^\lambda$. The mapping $p: S \times N \to M$, $p(f,x)=f\cdot x$, is a Riemannian submersion. The horizontal space can be identified with $\widetilde{\mathcal{S}}+\mathcal{O}_N$ in a natural way, where $\widetilde{\mathcal{S}}$ denotes the distribution on $S$ given by the tangent spaces and $\mathcal{O}_N$ denotes the restriction of the distribution $\mathcal{O}$ to $N$.
\end{compactenum}
\end{proposition}
\vspace{-1.4em}\begin{proof}
Obviously, $p$ is smooth. We know already that $S$ preserves $\mathcal{O}$. Using that $Z$ generates $Z(S)$, we obtain in the case~(ii) that for all $f \in S$,
\begin{align*}
df_x(\widetilde{Z}(x))&=df_x(\frac{\partial}{\partial t} \left(\exp(tZ) \cdot x\right)|_{t=0})\\
&=\frac{\partial}{\partial t} \left(f\exp(tZ) \cdot x\right)|_{t=0}\\
&=\frac{\partial}{\partial t}\left( \exp(tZ) \cdot (f\cdot x)\right)|_{t=0}\\
&=\widetilde{Z}(f \cdot x).\end{align*}
Thus, $S$ preserves $\mathcal{Z}$ as well. It follows that in both cases, $f \cdot N$ is a leaf as well.

Since the action of $S$ is locally free by Corollaries~\ref{cor:sl2R_locally_free} and~\ref{cor:twisted_heis_locally_free} and $\mathcal{O}=\mathcal{S}^\perp$, any leaf of $\mathcal{O}$ or $\mathcal{O}+\mathcal{Z}$, respectively, in a neighborhood of $x \in M$ is given by $f \cdot N^\prime$ for some $f \in S$ and $N^\prime$ being the leaf through $x$. Thus, the image of $p$ is open as well as the complement of the image of $p$ in $M$. But $M$ is connected and the image is non-empty, hence, $p$ is surjective. Also, $dp$ is surjective everywhere. It follows that $p$ is a submersion.

(i) The action of $S$ on $M$ preserves $\mathcal{O}$ and leaves. Hence, $dp|_{\mathcal{O}_N}$ is isometric. It follows from Lemma~\ref{lem:S_covering}, that $dp|_{\widetilde{\mathcal{S}}}$ is isometric as well. Using the orthogonality of $\widetilde{\mathcal{S}}$ and $\mathcal{O}$ on $M$, we obtain that $p$ is a local isometry.

(ii) Locally, only the center $Z(S)$ preserves a leaf. Therefore, the connected component of $(f,x) \in S \times N$ in the preimage of $f \cdot x$ is equal to $\left\{(fz^{-1},z \cdot x)\right\}_{z \in Z(S)}$, which is isomorphic to the center of $S$.

Let $\widetilde{\mathcal{Z}}$ denote the distribution defined by the action of the center on $S$, $f \mapsto z \cdot f$ for $f \in S$ and $z \in Z(S)$. Now consider the distribution $\widetilde{\mathcal{Z}}^{\perp}$ on $S$ orthogonal to $\widetilde{\mathcal{Z}}$ with respect to the metric $\beta_S$. Moreover, let $\mathcal{Z}^\prime$ be the one-dimensional distribution on $S \times N$ defined by $(X,X) \in \widetilde{\mathcal{Z}} \times \mathcal{Z}$ for any $X \in \mathds{R}Z$. Then the horizontal space of the submersion $p$ is $\widetilde{\mathcal{Z}}^{\perp}\oplus\mathcal{Z}^\prime\oplus\mathcal{O}_N$ and can be naturally identified with $\widetilde{\mathcal{S}}\oplus\mathcal{O}_N$ using the canonical identification of $\mathcal{Z}^\prime$ with $\widetilde{\mathcal{Z}}$ by projecting to the first component.

To show that $p$ is a Riemannian submersion, we have to show that $dp|_{\widetilde{\mathcal{Z}}^{\perp}\oplus\mathcal{Z}^\prime\oplus\mathcal{O}_N}$ is isometric. As in the first part, $dp|_{\widetilde{\mathcal{Z}}^{\perp}}$ and $dp|_{\mathcal{O}_N}$ are isometric. Using Lemma~\ref{lem:S_covering} and the fact, that the center of $S$ acts isometrically on $N$, $dp|_{\mathcal{Z}^\prime}$ is also isometric. An orthogonality argument as above yields the result.
\end{proof}\vspace{0pt}

Assume $\mathfrak{s}\cong\mathfrak{he}_d^\lambda$. The diagonal action $(f,x) \mapsto (fz^{-1},z \cdot x)$ of the center $Z(S)$ on $S \times N$ is isometric. If we factorize $S \times N$ through this action, we obtain a set denoted by $S \times_{Z(S)} N$. The action of $Z(S)$ on $S \times N$ induces a foliation with leaves isomorphic to $Z(S)$. Since $Z(S)$ is a closed Lie subgroup of $S$, it follows that the leaves are closed in $S \times N$. The Frobenius theorem (cf.~\cite{Ba09}, Proposition~A.9) gives then a canonical manifold structure on $S \times_{Z(S)} N$. Note that this construction is a special case of constructing manifolds with a certain group action (cf.~\cite{Ad01}).

As in the proof of Proposition~\ref{prop:submersion}~(ii), let $\widetilde{\mathcal{Z}}$ denote the distribution defined by the action of the center on $S$, $\widetilde{\mathcal{Z}}^{\perp}$ the distribution on $S$ $\beta_S$-orthogonal to $\widetilde{\mathcal{Z}}$ and $\mathcal{Z}^\prime$ the one-dimensional distribution on $S \times N$ defined by $(X,X) \in \widetilde{\mathcal{Z}} \times \mathcal{Z}$ for any $X \in \mathds{R}Z$.

We furnish $S \times_{Z(S)} N$ with the metric induced by the metric on $\widetilde{\mathcal{Z}}^{\perp}\oplus\mathcal{Z}^\prime\oplus\mathcal{O}_N$, that is, $S \times N \to S \times_{Z(S)} N$ is a Riemannian submersion with horizontal space $\widetilde{\mathcal{Z}}^{\perp}\oplus\mathcal{Z}^\prime\oplus\mathcal{O}_N$. $p$ induces a smooth map $\pi: S \times_{Z(S)} N \to M$, which is surjective and a local isometry.

\begin{lemma}\label{lem:N_geodesically_complete}
$N=(N,h)$ is geodesically complete.
\end{lemma}
\vspace{-1.4em}\begin{proof}
Assume there is a maximal geodesic $\gamma$ in $N$ which is not defined on the entire real line. Without loss of generality, $\gamma: (a,b) \to N$, $-\infty\leq a<b <+\infty$. Since the metric space $N$ inherits the metric of $(M,g^\prime)$, $\left\{\gamma(t_k)\right\}_{k=0}^\infty$ is a Cauchy sequence in $N$ and $M$ as well, if $t_k \to b$ increases as $k \to \infty$ (note that in a Riemannian space, a geodesic is locally a length-minimizing curve). $M$ is compact, thus, $\gamma(t_k) \to x$ in $M$. It follows that we can extend the curve $\gamma$ continuously on $M$, such that $\gamma(b)$ is defined. Since the tangent vector of $\gamma$ is contained in $\mathcal{O}+\mathcal{Z}$ on $(a,b]$, $x \in N$ by the Frobenius theorem (cf.~\cite{Ba09}, Proposition~A.9). Hence, we can extend $\gamma$ on $N$, contradiction.
\end{proof}\vspace{0pt}

\begin{definition}
If $\mathfrak{s}\cong\mathfrak{he}_d^\lambda$, $\textnormal{Isom}_{Z(S)}(N)$ denotes the subgroup of $\textnormal{Isom}(N)$ consisting of all $Z(S)$-equivariant isometries of $N$.
\end{definition}

\begin{proposition}\label{prop:covering}
Let $M=(M,g^\prime)$ as above. Then the following is true:
\begin{compactenum}
\item $p: S \times N \to M$ for $\mathfrak{s}\cong\mathfrak{sl}_2(\mathds{R})$ and $\pi: S \times_{Z(S)} N \to M$ for $\mathfrak{s}\cong\mathfrak{he}_d^\lambda$, respectively, are covering maps.

\item The group of deck transformations acts transitively on each fiber. Furthermore, any deck transformation comes from an element of $S \times \textnormal{Isom}(N)$ for $\mathfrak{s}\cong\mathfrak{sl}_2(\mathds{R})$ or $S \times \textnormal{Isom}_{Z(S)}(N)$ for $\mathfrak{s}\cong\mathfrak{he}_d^\lambda$, respectively.
\end{compactenum}
\end{proposition}
\vspace{-1.4em}\begin{proof}
(i) All manifolds we consider are connected. Moreover, $S$ is a Riemannian homogeneous space and therefore geodesically complete. $N$ is geodesically complete by Lemma~\ref{lem:N_geodesically_complete}. Therefore, $S \times N$ is geodesically complete. Since $S \times N \to S \times_{Z(S)} N$ is a Riemannian submersion, $S \times_{Z(S)} N$ is geodesically complete as well (cf.~\cite{Her60}, Theorem~1).

$p$ and $\pi$, respectively, are locally isometric. Furthermore, $S \times N$ and $S \times_{Z(S)} N$, respectively, are geodesically complete, so $p$ and $\pi$, respectively, are Riemannian covering maps.

(ii) Assume \[p(f,x)=p(f^\prime,x^\prime) \text{ or } \pi([f,x])=\pi([f^\prime,x^\prime]),\] respectively, for $(f,x) \in S \times N$ and $[f,x]$ being the corresponding equivalence class in $S \times_{Z(S)} N$. In both cases, $f \cdot x=f^\prime \cdot x^\prime$. Equivalently, \[x= f^{-1}f^\prime \cdot x^\prime.\] In the proof of Proposition~\ref{prop:submersion} we have seen that $f^{-1}f^\prime \cdot N$ is a leaf of $\mathcal{O}$ or $\mathcal{O}+\mathcal{Z}$, respectively, as well. It follows that it is equal to $N$. By construction, $f^{-1}f^\prime$ corresponds to an isometry of $N$ (being $Z(S)$-equivariant if $\mathfrak{s}\cong\mathfrak{he}_d^\lambda$, because $Z(S)$ commutes with $f^{-1}f^\prime$ in $S$). Therefore, there is $\psi \in \textnormal{Isom}(N)$ or $\psi \in \textnormal{Isom}_{Z(S)}(N)$, respectively, such that $f^{-1}f^\prime \cdot y=\psi(y)$ for all $y \in N$. This yields for all $y$ \[{f^\prime}^{-1}f \cdot \psi(y)=y.\]

$({f^\prime}^{-1}f,\psi)$ acts as an isometry $\varphi$ of $S \times_{Z(S)} N$ by \[\varphi([\widetilde{f},y])=[\widetilde{f}{f^\prime}^{-1}f,\psi(y)].\] The action on $S \times N$ in the other case is analogous. We know that \[\pi(\varphi([\widetilde{f}f^\prime,y]))=\pi([\widetilde{f}f,\psi(y)])=\pi([\widetilde{f}f^\prime,y])\] and an analogous result for the case of $S \times N$, holds for all $y \in N$ and $\widetilde{f} \in S$, because $f \cdot \psi(y)=f^\prime \cdot y$ implies $\widetilde{f}f \cdot \psi(y)=\widetilde{f}f^\prime \cdot y$.

It follows that $\varphi$ is a deck transformation. Since each deck transformation is determined by the value at one point, we are done.
\end{proof}\vspace{0pt}

Suppose $\mathfrak{s}\cong\mathfrak{sl}_2(\mathds{R})$ and consider the covering map $\widetilde{SL_2(\mathds{R})} \to S$. It follows from our results above that there is a discrete subgroup $\Gamma$ in $\textnormal{Isom}(N)\times\widetilde{\textnormal{SL}_2(\mathds{R})}$ acting freely on $N\times\widetilde{\textnormal{SL}_2(\mathds{R})}$, such that $M$ is diffeomorphic to $\Gamma \backslash \mkern-5mu\left(N\times  \widetilde{\textnormal{SL}_2(\mathds{R})}\right)$. We have seen in the proof of Proposition~\ref{prop:orbit_Lorentz}~(i), that if we pull back the metric $g$ on $M$ to $S$ through the canonical map $S \to S \cdot x$, we obtain a positive multiple of the Killing form $k$. It follows that $M=(M,g)$ is covered isometrically by the manifold $N \times \widetilde{\textnormal{SL}_2(\mathds{R})}$ provided with the metric \[g_{(x,\cdot)}= h_x \times (\sigma^2(x)k),\] $\sigma:N \to \mathds{R}_+$ smooth. Therefore, we have finished the proof of Theorem~\ref{th:geometric_characterization}~(ii).

Assume that $\mathfrak{s}\cong\mathfrak{he}_d^\lambda$ and consider the canonical defined map \[S \times \textnormal{Isom}_{Z(S)}(N) \to \textnormal{Isom}(S \times_{Z(S)} N),\] \[ (f,\psi) \in S \times \textnormal{Isom}_{Z(S)}(N) \text{ maps } [f^\prime,y] \in S \times_{Z(S)} N \text{ to } [f^\prime f,\psi(y)].\] It has kernel exactly $\left\{(z^{-1},z)\right\}_{z \in Z(S)} \cong Z(S)$. Here $z \in \textnormal{Isom}_{Z(S)}(N)$ corresponds to the isometry $x \mapsto z \cdot x$. $\textnormal{Isom}_{Z(S)}(N)$ is a closed subgroup of $\textnormal{Isom}(N)$ and the kernel $\left\{(z^{-1},z)\right\}_{z \in Z(S)}$ is a closed central subgroup of $S \times \textnormal{Isom}_{Z(S)}(N)$. Thus, if we factor out the kernel, we obtain a Lie group $S \times_{Z(S)} \textnormal{Isom}_{Z(S)}(N)$. The induced group homomorphism $S \times_{Z(S)} \textnormal{Isom}_{Z(S)}(N) \to \textnormal{Isom}(S \times_{Z(S)} N)$ is injective. Thus, there is a discrete subgroup $\Gamma \subset S \times_{Z(S)} \textnormal{Isom}_{Z(S)}(N)$ and a group isomorphism between $\Gamma$ and the deck transformations of $\pi$. Since every non-trivial deck transformation has no fixed point, $\Gamma$ acts freely on $S \times_{Z(S)} N$. It follows that $M$ is diffeomorphic to $\Gamma \backslash \mkern-5mu\left( S \times_{Z(S)} N\right)$.

\subsection{Lorentzian metrics on the twisted Heisenberg group}\label{sec:Lorentz_metric_Heis}

To finish the proof of Theorem~\ref{th:geometric_characterization}~(iii), we have to determine the Lorentzian metric on $S \times_{Z(S)} N$ given by the covering map $\pi$ and the metric $g$ on $M$. For this, it suffices to investigate the metric $m=m_x$ on $S$ given by pulling back the metric of $M$ under the canonical map $S \to S\cdot x$ for some $x \in M$. Remember that by Lemma~\ref{lem:S_covering}, the map $S \to S\cdot x$ is a covering map.

By Proposition~\ref{prop:orbit_Lorentz}, the metric $m$ is Lorentzian. Moreover, $m$ is invariant under left translation of $S$, since $S$ acts isometrically on $(M,g)$. Also, $m$ is invariant under right translation of $S$ with an element of $\Gamma_0 \cdot Z(S)$, where $\Gamma_0$ denotes the projection of $\Gamma$ under the continuous map $S \times_{Z(S)} N \to S/Z(S)$ defined by $[f,x] \to [f]$. Let $\overline{\Gamma_0}$ be the topological closure of $\Gamma_0$ in the Lie group $S/Z(S)$. Note that $S/Z(S)\cong Z(S)\backslash S$.

\begin{lemma}\label{lem:heis_gamma_0_cocompact}
$\overline{\Gamma_0}$ is cocompact in $S/Z(S)$.
\end{lemma}
\vspace{-1.4em}\begin{proof}
Consider the covering map $\pi: S \times_{Z(S)} N \to M$ and cover $M$ by evenly covered precompact open sets. Since $M$ is compact, we may choose finitely many of them, such that they still cover $M$. For any such precompact open set, choose one sheet in $S \times_{Z(S)} N$. The collection of these sheets is precompact and we denote the closure of them by $A$. Since $M$ is diffeomorphic to $\Gamma \backslash \mkern-5mu\left( S \times_{Z(S)} N\right)$, the deck transformations are acting transitively on the fibers, so $\Gamma \cdot A=S \times_{Z(S)} N$.

Let $B$ denote the projection of $A$ under the continuous map $S \times_{Z(S)} N \to S/Z(S)$ defined by $[f,x] \to [f]$. $B$ is compact and $B \cdot \Gamma_0 = S/Z(S)$ (remember that $\Gamma_0$ is acting from the right on $S/Z(S)$).

Thus, $B \cdot \overline{\Gamma_0}=S/Z(S)$. To show that $(S/Z(S))/\overline{\Gamma_0}$ is compact, we have to show that any sequence $\left\{f_k\overline{\Gamma_0}\right\}_{k=0}^\infty$, $f_k \in S/Z(S)$, has a convergent subsequence. Choose $b_k \in B$ and $\gamma_k \in \overline{\Gamma_0}$ such that $b_k \gamma_k=f_k$. $B$ is compact, so we may choose a convergent subsequence $\left\{b_{k_j}\overline{\Gamma_0}\right\}_{j=0}^\infty$. Thus, $f_{k_j}\overline{\Gamma_0}=b_{k_j}\overline{\Gamma_0}$ converges as $j \to \infty$.
\end{proof}\vspace{0pt}

\begin{proposition}\label{prop:finite_measure}
$(S/Z(S))/\overline{\Gamma_0}$ admits an $S$-invariant probability measure.
\end{proposition}
\vspace{-1.4em}\begin{proof}
The Lie algebra of $S/Z(S)$ is isomorphic to $\mathfrak{he}_d^\lambda/\mathds{R}Z$. Then $\textnormal{ad}_{[X]}=0$ if $X \in \mathfrak{he}_d$ and $\textnormal{ad}_{[T]}$ is semisimple and has only purely imaginary eigenvalues since \[\textnormal{ad}_{[T]} ([X_k])=\lambda_k[Y_k] \text{ and } \textnormal{ad}_{[T]} ([Y_k])=-\lambda_k[X_k]\] for all $1\leq k\leq d$. Since $S/Z(S)$ is connected, it follows that for any $f \in S/Z(S)$, $\textnormal{Ad}_f$ is semisimple and all its eigenvalues have absolute value $1$. Thus, $S/Z(S)$ as well as $\overline{\Gamma_0}$ are unimodular, so $(S/Z(S))/\overline{\Gamma_0}$ admits a non-trivial $S/Z(S)$-invariant Borel measure, which is unique up to a constant factor (cf.~\cite{Bo04}, Chapter~VII, Paragraph~2, Theorem~3, Corollary~2). Note that this measure is induced by the bi-invariant Haar measure on $S/Z(S)$.

By Lemma~\ref{lem:heis_gamma_0_cocompact}, $\overline{\Gamma_0}$ is cocompact in $S/Z(S)$, hence the measure on $(S/Z(S))/\overline{\Gamma_0}$ is finite. Without loss of generality, we may assume that the measure is a probability measure. Note that the quotient map $S \to S/Z(S)$ yields an $S$-action on the space $(S/Z(S))/\overline{\Gamma_0}$ by multiplication from the left.
\end{proof}\vspace{0pt}

Consider the space $S^2\mathfrak{s}$ of symmetric bilinear forms on $\mathfrak{s}$. $S$ acts on $S^2\mathfrak{s}$ via the adjoint representation: \[(f \cdot b)(X,Y):=b(\textnormal{Ad}_{f^{-1}}(X),\textnormal{Ad}_{f^{-1}}(Y))\] for any $f \in \mathfrak{s}, b \in S^2\mathfrak{s}$. We investigate the $S$-orbit of $m$ in $S^2\mathfrak{s}$.

Since the center $Z(S)$ is equal to the kernel of $\textnormal{Ad}$, the $S$-action factors to an action of $S/Z(S)$.

We know that $m$ is invariant under right translation of $S$ with an element of $\Gamma_0\cdot Z(S)$, hence $m$ is $\textnormal{Ad}(\Gamma_0\cdot Z(S))$-invariant. It follows that $\overline{\Gamma_0}$ acts trivially on the orbit of $m$. Thus, by Proposition~\ref{prop:finite_measure}, $(S/Z(S))/\overline{\Gamma_0}$ induces an $S$-invariant probability measure $\mu$ on $S^2\mathfrak{s}$. Because the Haar measure of $S/Z(S)$ induces this measure, the support of $\mu$ is equal to the whole orbit of $m$.

Let $H$ be the subgroup generated by $\mathfrak{he}_d$ in $S$. Then the $S$-action above yields an $H$-action on $S^2\mathfrak{s}$. Applying the F\"urstenberg lemma in form of Corollary~\ref{cor:Furstenberg}, we obtain a Lie group homomorphism $\varrho: H \to K$ into a compact group $K$ in $\textnormal{GL}(S^2\mathfrak{s})$ such that the $H$-action on the orbit of $m$ coincides with the action of $K$.

\begin{remark}
The center $Z(S)$ acts trivially on $S^2\mathfrak{s}$, so $\varrho$ induces a Lie group homomorphism $\varrho^\prime: \widetilde{\textnormal{He}_d} \to K$. Note that a Lie group homomorphism $\widetilde{\textnormal{He}_d} \to K$ with a compact group $K$, which is trivial on the center, is not necessarily trivial. For example, if $d=1$, consider the homomorphism $\widetilde{\textnormal{He}_1} \to \textnormal{SO}(4)$ defined by \[\begin{pmatrix}1 & x & z \\ 0 & 1 & y \\ 0 & 0 & 1\end{pmatrix} \mapsto \begin{pmatrix}\cos(x) & -\sin(x) & 0 & 0 \\ \sin(x) & \cos(x) & 0 & 0 \\ 0 & 0 & \cos(y) & -\sin(y)\\ 0 & 0 & \sin(y) & \cos(y)\end{pmatrix}.\]
Examples for higher dimensions can be constructed easily in the same way.
\end{remark}

\begin{proposition}\label{prop:m_invariant}
$m$ is an ad-invariant Lorentz form on $\mathfrak{s}$.
\end{proposition}
\vspace{-1.4em}\begin{proof}
Let $V:=\textnormal{span}\left\{X_1,Y_1,\ldots,X_d,Y_d\right\}$.

We know already that $m$ is Lorentzian. From the proof of Proposition~\ref{prop:lorentz_heisenberg}~(ii) follows, that it suffices to show
\[m(T,X_k)=0=m(T,Y_k) \text{ and } m(Z,Z)=0=m(Z,X_k)=m(Z,Y_k)\] for all $1 \leq k \leq d,$ \[m(X_k,Y_j)=0 \text{ and } m(X_k,X_j)=\delta_{jk}m(T,Z)=m(Y_k,Y_j)\] for all $j,k=1,\ldots,d.$
Since $H$ is acting precompactly on the orbit of $m$, the value $(f \cdot m)(T,T)$ is bounded by a constant depending only on $T$, for all $f \in H$. Choosing $f=\exp(-tX_k)$ for an arbitrary $k$ and using \[[X_k,T]=-\lambda_k Y_k, \ [X_k,[X_k,T]]=-\lambda_k^2 Z,\] we obtain that
\[m(\textnormal{Ad}_{\exp(tX_k)}(T),\textnormal{Ad}_{\exp(tX_k)}(T))=m(T-t\lambda_kY_k-\frac{t^2}{2}\lambda_k^2 Z,T-t\lambda_kY_k-\frac{t^2}{2}\lambda_k^2 Z)\]
is bounded for all $t \in \mathds{R}$. But the latter is equal to \begin{align*} m(T,T)&-2t\lambda_k m(T,Y_k)+t^2 \lambda_k^2 \left(m(Y_k,Y_k)-m(Z,T)\right)\\
&+t^3 \lambda_k^3 m(Y_k,Z)+\frac{t^4}{4} \lambda_k^4 m(Z,Z) \end{align*} and a polynomial is bounded if and only if it is constant. Hence, \[m(T,Y_k)=m(Y_k,Z)=m(Z,Z)=0 \text{ and } m(Y_k,Y_k)=m(Z,T).\] Analogously, if we choose $f=\exp(tY_k)$, we obtain \[m(T,X_k)=m(X_k,Z)=0 \text{ and } m(X_k,X_k)=m(Z,T).\]
In the same way as above, for any $j$ and $k$,
\begin{align*}
m(\textnormal{Ad}_{\exp(tX_k)}(T),\textnormal{Ad}_{\exp(tX_k)}(X_j))&=m(T-t\lambda_kY_k-\frac{t^2}{2}\lambda_k^2 Z,X_j)\\
&=m(T,X_j)-t\lambda_k m(Y_k,X_j)-\frac{t^2}{2}\lambda_k^2 m(Z,X_j)\end{align*}
is bounded for all real $t$. Thus, \[m(X_j,Y_k)=0.\] Replacing $X_j$ by $Y_j$, $j\neq k$, we obtain \[m(Y_j,Y_k)=0\] if $j \neq k$. Analogously, \[m(X_j,X_k)=0\] if $j \neq k$.
\end{proof}\vspace{0pt}

To conclude this section, we show that there exists essentially only one ad-invariant Lorentzian scalar product on $\mathfrak{s}$.

\begin{proposition}\label{prop:Lorentz_s_only_one}
Any two ad-invariant Lorentzian scalar products $b_1,b_2$ on $\mathfrak{s}$ are equivalent, that is, there is an automorphism $L:\mathfrak{s} \to \mathfrak{s}$, such that for all $X,Y \in \mathfrak{s}$, $b_1(X,Y)=b_2(L(X),L(Y))$.
\end{proposition}
\vspace{-1.4em}\begin{proof}
Let $b$ be an ad-invariant Lorentzian scalar product on $\mathfrak{s}$. The second part of Proposition~\ref{prop:lorentz_heisenberg} implies that any two $\alpha, \beta \in \mathds{R}$, $\alpha>0$, determine an ad-invariant Lorentzian scalar product. In the case of $b$, $\alpha=b(T,Z)$ and $\beta=b(T,T)$.

The linear map $L:\mathfrak{s}\to\mathfrak{s}$ defined by \[T \mapsto T-\frac{\beta}{2\alpha}Z \text{ and } X \mapsto X \text{ for all }X\in\mathfrak{he}_d,\] is an automorphism. Moreover, \[b(L(T),L(T))=\beta-2\frac{\beta}{2\alpha}b(T,Z)=0.\] Thus, up to equivalence, we may assume $\beta=0$.

If we define $L$ by \[T \mapsto T, \ X_k \mapsto \frac{1}{\sqrt{\alpha}}X_k, \ Y_k \mapsto \frac{1}{\sqrt{\alpha}}Y_k, \ Z \mapsto \frac{1}{\alpha}Z \text{ for all }1\leq k \leq d,\] we also obtain an automorphism of $\mathfrak{s}$. We have \[b(L(T),L(Z))=\frac{1}{\alpha} b(T,Z)=1.\] So up to equivalence, $\alpha=1$, and $b$ is uniquely determined.
\end{proof}\vspace{0pt}

%%%%%%%%%%%%%%%%%%%%%%%%%%%%%%%%%%%%%%%%%%%%%%%%%%%%%%%%%%%%%%%%%%%%%%%%%%%%%%%%%%%%%%%%%%%%

\chapter{Compact homogeneous Lorentzian manifolds}\label{ch:homogeneous}

In the following, we will investigate compact homogeneous Lorentzian manifolds, especially those, whose isometry groups have non-compact connected components. We use the geometric description of Theorem~\ref{th:geometric_characterization} to prove Theorem~\ref{th:homogeneous_characterization} in Section~\ref{sec:proof_homogeneous}. 

We continue in Section~\ref{sec:reductive} with giving a reductive representation for any compact homogeneous Lorentzian manifold. In the case that the connected component of the identity in the isometry group is compact, the statement can be shown in the same way as in the Riemannian case. For the more interesting case of non-compact connected components of the isometry group, the induced bilinear form $\kappa$ yields a reductive representation. Additionally, we show that the isotropy group of a point has compact connected components.

Although the representation given in Section~\ref{sec:reductive} seems to be quite natural, we use a different one in Section~\ref{sec:homogeneous_geometry} that is more convenient with respect to calculational purposes. This representation allows us to determine the local geometry of compact homogeneous Lorentzian manifolds whose isometry groups have non-compact connected components.

\section{Structure of homogeneous manifolds}\label{sec:proof_homogeneous}

Let $M=(M,g)$ be a compact homogeneous Lorentzian manifold and denote by $\textnormal{Isom}^0(M)$ the identity component in the isometry group. We suppose that $\textnormal{Isom}^0(M)$ is not compact. Additionally, let $\mathfrak{isom}(M)=\mathfrak{k}\oplus\mathfrak{a}\oplus\mathfrak{s}$ be the decomposition of its Lie algebra according to Theorem~\ref{th:algebraic_classification}. Due to Corollary~\ref{cor:condition_star} and Theorem~\ref{th:algebraic}, this decomposition is $\kappa$-orthogonal, where $\kappa$ is the induced bilinear form on $\mathfrak{isom}(M)$.

In a homogeneous semi-Riemannian manifold, each tangent vector can be extended to a Killing vector field (cf.~\cite{ON83}, Corollary~9.38). Thus, the orbits of the isometry group on the Lorentzian manifold $M$ have Lorentzian character. It follows from Theorem~\ref{th:geometric_characterization}, that $\mathfrak{s}\cong\mathfrak{sl}_2(\mathds{R})$ or $\mathfrak{s}\cong\mathfrak{he}_d^\lambda$.

Let $S$ be the subgroup generated by $\mathfrak{s}$ and denote by $\mathcal{S}$ the distribution defined by the orbits of $S$. We can apply Theorem~\ref{th:geometric_characterization}~(ii) and~(iii) to see that $M$ is isometric to $\Gamma \backslash \mkern-5mu \left(N \times_\sigma \widetilde{\textnormal{SL}_2(\mathds{R})}\right)$ or $\Gamma \backslash \mkern-5mu \left(S \times_{Z(S)} N\right)$, respectively. $N$ is a Riemannian manifold and $\Gamma$ is a discrete subgroup in $\textnormal{Isom}(N)\times\widetilde{\textnormal{SL}_2(\mathds{R})}$ or $S\times_{Z(S)}\textnormal{Isom}_{Z(S)}(N)$, respectively.

Remember that in Theorem~\ref{th:geometric_characterization}~(ii) and~(iii), $N$ was a leaf of the involutive distribution $\mathcal{O}$ ($\mathfrak{s}\cong\mathfrak{sl}_2(\mathds{R})$) or $\mathcal{O}+\mathcal{Z}$ ($\mathfrak{s}\cong\mathfrak{he}_d^\lambda$), respectively. $\mathcal{O}$ was the distribution orthogonal to $\mathcal{S}$ and $\mathcal{Z}$ denoted the orbit of the center of $S$, if $S$ is a twisted Heisenberg group.

Let $H$ be the isotropy group of a point $x_0 \in N \subset M$ in $\textnormal{Isom}^0(M)$. Denote by $\mathfrak{h}$ its Lie algebra.

\begin{proposition}\label{prop:H_compact}
In the situation of this section, the following is true:
\begin{compactenum}
\item Assume that $\mathfrak{s}\cong\mathfrak{sl}_2(\mathds{R})$. Then $\mathfrak{h} \subseteq \mathfrak{k}\oplus\mathfrak{a}$.

\item Assume that $\mathfrak{s}\cong\mathfrak{he}_d^\lambda$. Then $\mathfrak{h} \subseteq \mathfrak{k}\oplus\mathfrak{a}\oplus\mathfrak{z(s)}$.
\end{compactenum}
In both cases, it follows that the identity component of $H$ is compact.
\end{proposition}
\vspace{-1.4em}\begin{proof}
(i) Take $Y=A+W \in \mathfrak{h}$ with $A \in \mathfrak{k}\oplus\mathfrak{a}$, $W \in \mathfrak{s}$. If $\left\{e,f,h\right\}$ is an $\mathfrak{sl}_2$-triple of $\mathfrak{s}$ (fulfilling $[h,e]=2e$, $[h,f]=-2f$, $[e,f]=h$), then $\textnormal{ad}_e$ and $\textnormal{ad}_f$ are nilpotent.

Due to Theorem~\ref{th:locally_free}, the subgroup generated by $\mathfrak{s}$ acts locally freely on $M$. Using $[\mathfrak{s},Y]\subseteq \mathfrak{s}$ and Lemma~\ref{lem:stabilizer} with $X=e$ or $X=f$, we obtain $[e,Y]=0$ as well as $[f,Y]=0$. By the Jacobi identity, \[[h,Y]=[[e,f],Y]=[[e,Y],f]+[[Y,f],e]=0,\] so $W=0$.

Therefore, $\mathfrak{h}\subseteq\mathfrak{k}\oplus\mathfrak{a}$. Theorem~\ref{th:algebraic} together with Corollary~\ref{cor:condition_star} states that the subgroup generated by $\mathfrak{a}$ is compact. The compact semisimple algebra $\mathfrak{k}$ also generates a compact subgroup (cf.~\cite{Bo05}, Chapter~IX, Paragraph~1.4, Theorem~1). Thus, the subgroup generated by $\mathfrak{h}$ is compact as well, since it is closed.

(ii) Take $Y=A+W \in \mathfrak{h}$ with $A \in \mathfrak{k}\oplus\mathfrak{a}$, $W \in \mathfrak{s}$. Let $X$ be an arbitrary element of the nilradical of $\mathfrak{s}$, which is isomorphic to $\mathfrak{he}_d$.

Due to Theorem~\ref{th:locally_free}, the subgroup generated by $\mathfrak{s}$ acts locally freely on $M$. Using $[\mathfrak{s},Y]\subseteq \mathfrak{s}$ and Lemma~\ref{lem:stabilizer}, we obtain $[X,Y]=0$. Since this is true for any $X$ as above, $W$ centralizes the nilradical of $\mathfrak{s}$. It follows that $W \in \mathfrak{z(s)}$ because $\mathfrak{s}$ is a twisted Heisenberg algebra.

Therefore, $\mathfrak{h}\subseteq\mathfrak{k}\oplus\mathfrak{a}\oplus\mathfrak{z(s)}$. Theorem~\ref{th:algebraic} together with Corollary~\ref{cor:condition_star} states that the subgroup generated by $\mathfrak{a}\oplus\mathfrak{z(s)}$ is compact. We can now conclude exactly as in~(ii), that $\mathfrak{h}$ generates a compact subgroup.
\end{proof}\vspace{0pt}

\begin{proposition}\label{prop:N_compact}
In the situation of this section, $N$ is a compact homogeneous Riemannian manifold.
\end{proposition}
\vspace{-1.4em}\begin{proof}

Consider the decomposition $\mathfrak{isom}(M)=\mathfrak{k}\oplus\mathfrak{a}\oplus\mathfrak{s}$, where $\mathfrak{s}$ is isomorphic to $\mathfrak{sl}_2(\mathds{R})$ or $\mathfrak{he}_d^\lambda$, respectively. Let $C$ be the subgroup generated by \[\mathfrak{c}:=\mathfrak{k}\oplus\mathfrak{a} \text{ or } \mathfrak{c}:=\mathfrak{k}\oplus\mathfrak{a}\oplus\mathfrak{z(s)},\] respectively. But $\mathfrak{a}$ and $\mathfrak{a}\oplus\mathfrak{z(s)}$, respectively, generate compact subgroups by Corollary~\ref{cor:condition_star} and Theorem~\ref{th:algebraic}. $\mathfrak{k}$ is compact semisimple and also generates a compact subgroup (cf.~\cite{Bo05}, Chapter~IX, Paragraph~1.4, Theorem~1). Thus, $C$ is compact.

The centralizer of $S$ in $\textnormal{Isom}^0(M)$ is a Lie subgroup, whose Lie algebra coincides with the centralizer of $\mathfrak{s}$ in $\mathfrak{isom}(M)$ (cf.~\cite{On93}, Part~I, Chapter~2, Theorem~2.8). But the centralizer of $\mathfrak{s}$ is given by $\mathfrak{c}$, therefore, $C$ is contained in the centralizer of $S$.

Let $f \in C$. Using $ff^\prime \cdot x=f^\prime f \cdot x$ for all $x \in M$ and $f^\prime \in S$, we see that $f$ preserves the orbits $\mathcal{S}$ of $S$. But $f$ is an isometry, so $\mathcal{O}=\mathcal{S}^\perp$ is preserved as well. Restricting $f^\prime$ to elements of the center of $S$, we also obtain that $f$ preserves $\mathcal{Z}$ in the case $\mathfrak{s}\cong\mathfrak{he}_d^\lambda$. Thus, for all $f \in C$, $f \cdot N$ is a leaf of the foliation induced by the involutive distribution $\mathcal{O}$ or $\mathcal{O}+\mathcal{Z}$, respectively.

If $f \in S$, then $\mathcal{O}$ is obviously preserved by $f$. Using that $Z(S)$ is central, we see that $f$ preserves $\mathcal{Z}$ in the case $\mathfrak{s}\cong\mathfrak{he}_d^\lambda$ as well (see the proof of Proposition~\ref{prop:submersion}). As above, $f \cdot N$ is a leaf for all $f \in S$.

Because $\mathfrak{isom}(M)=\mathfrak{c}+\mathfrak{s}$, it follows that $f \cdot N$ is a leaf for all $f \in \textnormal{Isom}^0(M)$.

Let $\widehat{C}^0$ be the identity component of the subgroup $\widehat{C}$ in $\textnormal{Isom}^0(M)$, which maps the leaf $N$ to itself. Since $M$ is homogeneous, for any $x,y \in N$ there is $f \in \textnormal{Isom}(M)$, such that $f \cdot x=y$. But $f \cdot N$ is also a leaf, hence $f \cdot N=N$ and $f \in \widehat{C}$. It follows that $N$ is a homogeneous space. Since $N$ is connected, $\widehat{C}^0$ acts transitively on $N$.

Let $\widehat{\mathfrak{c}}$ be the Lie algebra of $\widehat{C}$. By definition, any Killing vector field generated by an element of $\widehat{\mathfrak{c}}$ has to lie in $\mathcal{O}$ or $\mathcal{O}+\mathcal{Z}$, respectively. But $\mathcal{O}=\mathcal{S}^\perp$ and $\mathcal{Z}$ is the orbit of the center of $S$, if it is a twisted Heisenberg group. Thus, $\widehat{\mathfrak{c}}$ is $\kappa$-orthogonal to $\mathfrak{s}$ in the case $\mathfrak{s}\cong\mathfrak{sl}_2(\mathds{R})$ and $\widehat{\mathfrak{c}}\subseteq\mathfrak{o}+\mathfrak{z(s)}$ if $\mathfrak{s}\cong\mathfrak{he}_d^\lambda$, where $\mathfrak{o}$ is the $\kappa$-orthogonal complement of $\mathfrak{s}$ in $\mathfrak{isom}(M)$. It follows that \[\widehat{\mathfrak{c}}\subseteq \mathfrak{c} \text{ and }\widehat{C}^0\subseteq C.\]

We want to show that $\widehat{C}^0= C$. Suppose the contrary. Then for any neighborhood $U$ of the identity in $C$, there is $f \in U$, such that $f \cdot N$ is not equal to $N$.

Let $U^\prime$ be a neighborhood of the identity in $S$. Since the action of $S$ on $M$ is locally free by Corollaries~\ref{cor:sl2R_locally_free} and~\ref{cor:twisted_heis_locally_free} and $\mathcal{O}=\mathcal{S}^\perp$, any leaf is locally given by $f^\prime \cdot N$, $f^\prime \in U^\prime$.

If the neighborhood $U$ is small enough, then for any $f \in U$, such that $f \cdot N$ is not equal to $N$, there is $f^\prime \in U^\prime$ such that $f \cdot N=f^\prime \cdot N$.

Since $\widehat{C}^0$ acts transitively on $N$, the map $\widehat{C}^0 \to N$, $\widehat{f} \mapsto \widehat{f} \cdot x_0$, is a locally trivial fiber bundle. So if we fix a neighborhood $\widehat{U}$ of the identity in $\widehat{C}^0$, $\widehat{U} \cdot x_0$ is a neighborhood of $x_0$ in $N$. It follows that if $U^\prime$ and $U$ are sufficiently small, we can choose $f \in U$, $f^\prime \in U^\prime$, such that \[f \cdot N=f^\prime \cdot N\neq N,\] and there is $\widehat{f} \in \widehat{U}$ such that \[\widehat{f}f^{-1}f^\prime \cdot x_0=x_0,\] that is, \[\widehat{f}f^{-1}f^\prime \in H.\] If the neighborhoods $\widehat{U}$, $U^\prime$ and $U$ are sufficiently small, even \[\widehat{f}f^{-1}f^\prime \in H^0,\] where $H^0$ is the identity component in $H$.

By Proposition~\ref{prop:H_compact}, $\mathfrak{h}\subseteq \mathfrak{c}$. Therefore, $H^0 \subseteq C$. Thus, \[\widehat{f}f^{-1}f^\prime \in C.\] Using that $\widehat{f} \in \widehat{C}^0 \subseteq C$ and $f \in U \subset C$, it follows that \[f^\prime \in C.\] If $U$ is small enough, it follows from $\mathfrak{c} \cap \mathfrak{s}=\mathfrak{z(s)}$ that $f^\prime$ is contained in the subgroup generated by $\mathfrak{z(s)}$. But the latter group preserves $N$, contradicting $f^\prime \cdot N\neq N$. Hence, $\widehat{C}^0= C$ and $\widehat{C}^0$ is compact. It follows that $N$ is compact as well.
\end{proof}\vspace{0pt}

Let $\mathfrak{s}\cong\mathfrak{sl}_2(\mathds{R})$. Since $C$ acts transitively on $N$, it follows that $\sigma$ has to be constant. Therefore, $M$ is covered isometrically by the metric product $N \times \widetilde{\textnormal{SL}_2(\mathds{R})}$, where $\widetilde{\textnormal{SL}_2(\mathds{R})}$ is furnished with the metric defined by a positive multiple of the Killing form of $\mathfrak{sl}_2(\mathds{R})$.

Assume $\mathfrak{s}\cong\mathfrak{he}_d^\lambda$. In the same way as above, we see that the map $m:N \to \mathcal{M}$ is a constant, which we denote now by $m$. Due to Proposition~\ref{prop:m_invariant}, $m$ is ad-invariant.

The proof of the last proposition shows even more:

\begin{corollary}\label{cor:isometry_group}
In the situation of this section, the following is true:
\begin{compactenum}
\item Let $s\cong\mathfrak{sl}_2(\mathds{R})$. Then $\textnormal{Isom}^0(M)$ is isomorphic to a central quotient of the group $C \times \widetilde{\textnormal{SL}_2(\mathds{R})}$. Here $C$ is the identity component of the centralizer of the projection of $\Gamma$ to $\textnormal{Isom}(N)$. $C$ acts transitively on $N$.

\item Let $\mathfrak{s}\cong\mathfrak{he}_d^\lambda$. Then $\textnormal{Isom}^0(M)$ is isomorphic to a central quotient of the group $S \times_{Z(S)} C$. Here $C$ is the identity component of the centralizer of the projection of $\Gamma$ to $Z(S)\cdot\textnormal{Isom}_{Z(S)}(N)$. $C$ acts transitively on $N$.
\end{compactenum}
\end{corollary}
\vspace{-1.4em}\begin{proof}
We have seen in the proof of Proposition~\ref{prop:N_compact}, that the compact subgroup $C$ generated by $\mathfrak{c}:=\mathfrak{k}\oplus\mathfrak{a}$ or $\mathfrak{c}:=\mathfrak{k}\oplus\mathfrak{a}\oplus\mathfrak{z(s)}$, respectively, preserves the leaf $N$ and acts transitively on $N$. Therefore, we can see $C$ as a subgroup of $\textnormal{Isom}(N)$. In case~(ii), $C$ is even a subgroup of $\textnormal{Isom}_{Z(S)}(N)$, since $\mathfrak{z(s)}$ is central in $\mathfrak{c}$, so any isometry in $C$ commutes with the action of $Z(S)$.

Because of $M \cong \Gamma \backslash \mkern-5mu \left(N \times \widetilde{\textnormal{SL}_2(\mathds{R})}\right)$ or $M \cong \Gamma \backslash \mkern-5mu \left(S \times_{Z(S)} N\right)$, respectively, an element of $\textnormal{Isom}(N)$ or $\textnormal{Isom}_{Z(S)}(N)$, respectively, induces an isometry of $M$ if and only if it commutes with the projection of $\Gamma$ to $\textnormal{Isom}(N)$ or $Z(S)\cdot\textnormal{Isom}_{Z(S)}(N)$, respectively. Note that by Theorem~\ref{th:geometric_characterization}, $\Gamma$ is a discrete subgroup of the group $\text{Isom}(N) \times \widetilde{\textnormal{SL}_2(\mathds{R})}$ or $S \times_{Z(S)} \text{Isom}_{Z(S)}(N)$, respectively. In the latter case, the projection to $Z(S)\cdot\textnormal{Isom}_{Z(S)}(N)$ is given by \[[f,\psi] \mapsto Z(S) \cdot \psi \subseteq \text{Isom}_{Z(S)}(N),\] $[f,\psi]$ is the equivalence class of $(f,\psi)\in S \times \text{Isom}_{Z(S)}(N)$ in $S \times_{Z(S)} \text{Isom}_{Z(S)}(N)$.

In the case~(ii), $Z(S)$ is contained in $S$ and $C$ and acts as the same group of isometries on $M$. Thus, $S \times_{Z(S)} C$ is defined.

The Lie algebra of $C \times S$ or $S \times_{Z(S)} C$, respectively, is isomorphic to the Lie algebra of the isometry group,  $\mathfrak{isom}(M)=\mathfrak{k}\oplus\mathfrak{a}\oplus\mathfrak{s}$.

The action of $C \times S$ or $S \times_{Z(S)} C$, respectively, on $M$ is clearly locally effective. Additionally, $S$ and $C$ are Lie subgroups of $\textnormal{Isom}^0(M)$. Therefore, the canonical Lie group homomorphism $C \times S \to \textnormal{Isom}^0(M)$ or $S \times_{Z(S)} C \to \textnormal{Isom}^0(M)$, respectively, is surjective and the kernel is a discrete central subgroup of $C \times S$ or $S \times_{Z(S)} C$, respectively.

Finally, the covering $\widetilde{\textnormal{SL}_2(\mathds{R})} \to S$ is central in the case~(i).
\end{proof}\vspace{0pt}

Consider the canonical projection $P:\textnormal{Isom}(N) \times \widetilde{\textnormal{SL}_2(\mathds{R})} \to \widetilde{\textnormal{SL}_2(\mathds{R})}$ (projection to the second component) or $P:S \times_{Z(S)} \textnormal{Isom}_{Z(S)}(N) \to S/Z(S)$ (projection to the first component), respectively. In both situations, $P$ is a continuous homomorphism, which is surjective and has kernel $\textnormal{Isom}(N)$ or kernel isomorphic to $\textnormal{Isom}_{Z(S)}(N)$.

Since $N$ is a compact Riemannian manifold by Proposition~\ref{prop:N_compact}, $\textnormal{Isom}(N)$ is compact. $\textnormal{Isom}_{Z(S)}(N)$ as a closed subgroup of $\textnormal{Isom}(N)$ is compact as well. Thus, in any case, $P$ has compact kernel. In the following lemma, we show that $\Gamma_0:=P(\Gamma)$ is discrete.

\begin{lemma}\label{lem:compact_kernel_discrete_subgroup}
Let $G$ and $G^\prime$ be connected Lie groups and $P:G \to G^\prime$ be a surjective continuous homomorphism with compact kernel. If $\Gamma \subset G$ is discrete, $P(\Gamma) \subset G^\prime$ is discrete as well.
\end{lemma}
\vspace{-1.4em}\begin{proof}
Assume the contrary. Then there is a sequence $\left\{f_k\right\}_{k=0}^\infty$ of distinct elements in $\Gamma$, such that $P(f_k) \to P(f)$ converges in $G^\prime$ as $k \to \infty$ for some $f \in G$.

Since $G \to G^\prime$ is a locally trivial fiber bundle with fiber $\textnormal{ker} (P)$, there is a sequence $\left\{p_k\right\}_{k=0}^\infty$ in $\textnormal{ker} (P)$, such that $f_k p_k \to f$ converges in $G$ as $k \to \infty$. But $\textnormal{ker} (P)$ is compact, hence, we may assume already that $p_k \to p$ converges. It follows that $f_k \to f p^{-1}$ as $k \to \infty$, contradicting the fact that $\Gamma$ is discrete.
\end{proof}\vspace{0pt}

By Lemma~\ref{lem:heis_gamma_0_cocompact}, we know that $\Gamma_0=\overline{\Gamma_0}$ is cocompact in $S/Z(S)$ in the case that $S$ is a twisted Heisenberg group. Thus, $\Gamma_0$ is a uniform lattice in $S/Z(S)$. We will now show that the corresponding result is true in the case that $S$ is locally isomorphic to $\widetilde{\textnormal{SL}_2(\mathds{R})}$. The proof of the following lemma is done in the same way as the proof of Lemma~\ref{lem:heis_gamma_0_cocompact}.

\begin{lemma}\label{lem:sl2R_gamma_0_cocompact}
Let $\mathfrak{s}\cong\mathfrak{sl}_2(\mathds{R})$. Then $\Gamma_0$ is cocompact in $\widetilde{\textnormal{SL}_2(\mathds{R})}$. Thus, $\Gamma_0$ is a uniform lattice.
\end{lemma}
\vspace{-1.4em}\begin{proof}
Consider the covering map $p: N \times \widetilde{\textnormal{SL}_2(\mathds{R})} \to M$ and cover $M$ by evenly covered precompact open sets. Since $M$ is compact, we may choose finitely many of them, such that they still cover $M$. For any such precompact open set, choose one sheet in $N \times \widetilde{\textnormal{SL}_2(\mathds{R})}$. The collection of these sheets is precompact and we denote the closure of them by $A$. Because $M \cong \Gamma \backslash \mkern-5mu \left(N \times \widetilde{\textnormal{SL}_2(\mathds{R})}\right)$, the deck transformations are acting transitively on the fibers. Hence, $\Gamma \cdot A=N \times \widetilde{\textnormal{SL}_2(\mathds{R})}$.

Let $B$ denote the projection of $A$ under the projection $N \times \widetilde{\textnormal{SL}_2(\mathds{R})} \to \widetilde{\textnormal{SL}_2(\mathds{R})}$. $B$ is compact and $B \cdot\Gamma_0 = \widetilde{\textnormal{SL}_2(\mathds{R})}$.

We want to show that $\widetilde{\textnormal{SL}_2(\mathds{R})}/\Gamma_0$ is compact. Equivalently, we have to show that any sequence $\left\{f_k\Gamma_0\right\}_{k=0}^\infty$, $f_k \in \widetilde{\textnormal{SL}_2(\mathds{R})}$, has a convergent subsequence. For all $k$, choose $b_k \in B$ and $\gamma_k \in \Gamma_0$ such that $b_k \gamma_k=f_k$. $B$ is compact, so we may choose a convergent subsequence $\left\{b_{k_j}\Gamma_0\right\}_{j=0}^\infty$. Therefore, $f_{k_j}\Gamma_0=b_{k_j}\Gamma_0$ converges as $j \to \infty$.
\end{proof}\vspace{0pt}

The kernel of the projection $\Gamma \to \Gamma_0$ is (isomorphic to) a discrete subgroup of $\textnormal{Isom}(N)$ or $\textnormal{Isom}_{Z(S)}(N)$, respectively. Note that the kernel is finite, since the latter groups are compact. Any element of the kernel can be written as $(\psi,e) \in \textnormal{Isom}(N) \times \widetilde{\textnormal{SL}_2(\mathds{R})}$ or $[e,\psi] \in S \times_{Z(S)} \textnormal{Isom}_{Z(S)}(N)$, respectively, $e$ being the identity element of the corresponding group.

Consider now the covering maps $p: N \times \widetilde{\textnormal{SL}_2(\mathds{R})} \to M$ and $\pi: S \times_{Z(S)} N \to M$, respectively, defined as in Proposition~\ref{prop:covering}. By definition, \[\psi(x_0)=p(\psi(x_0),e)=p((\psi,e)(x_0,e))=p(x_0,e)=x_0\] in the first case. Analogously, \[\psi(x_0)=\pi([e,\psi(x_0)])=\pi([e,\psi]([e,x_0]))=\pi([e,x_0])=x_0\] in the second case.
It follows that $(\psi,e)$ and $[e,\psi]$, respectively, have the fixed point $(x_0,e)$ on $N \times \widetilde{\textnormal{SL}_2(\mathds{R})}$ or $[e,x_0]$ on $S \times_{Z(S)} N$, respectively. But $\Gamma$ acts freely, so it has to be the identity element. Thus, the kernel of the projection $\Gamma \to \Gamma_0$ is trivial.

It follows that $\Gamma$ projects isomorphically to $\Gamma_0$. In the case~(i) of the theorem, we obtain that $\Gamma$ is the graph of a homomorphism $\varrho:\Gamma_0 \to \textnormal{Isom}(N)$.

$\Gamma$ corresponds to the group of deck transformations. Since $M$ is homogeneous, the centralizer of $\Gamma$ in the isometry group of the covering manifold acts transitively (cf.~\cite{Wo61}, Theorem~2.5). This completes the proof of Theorem~\ref{th:homogeneous_characterization}.

\section{General reductive representation}\label{sec:reductive}

Let $M$ be a compact homogeneous Lorentzian manifold and $G:=\textnormal{Isom}^0(M)$. $G$ acts transitively on the Lorentzian manifold $M$. We consider the isotropy group $H \subseteq G$ of some point $x \in M$. Denote by $\mathfrak{h}\subseteq\mathfrak{g}$ the corresponding Lie algebras.

If $G$ is compact, $H$ is compact as well, since $H$ is a closed subgroup. By Proposition~\ref{prop:compact_algebra}, $\mathfrak{g}$ possesses an ad-invariant symmetric bilinear form $b$, which is positive definite. $G$ is connected, therefore, $b$ is $\textnormal{Ad}(G)$-invariant. Now choose $\mathfrak{m}$ to be the $b$-orthogonal complement to $\mathfrak{h}$ in $\mathfrak{g}$. $\mathfrak{m}$ is $\textnormal{Ad}(H)$-invariant, because $\mathfrak{h}$ and $b$ are.

Suppose in the following that $G$ is not compact. By Proposition~\ref{prop:H_compact}, the connected component of $H$ is compact. To finish the proof of Theorem~\ref{th:homogeneous_reductive}, we have to find a reductive representation of $M$.

According to Theorems~\ref{th:algebraic_classification} and~\ref{th:homogeneous_characterization}, we have a decomposition $\mathfrak{g}=\mathfrak{k}\oplus\mathfrak{a}\oplus\mathfrak{s}$, where $\mathfrak{k}$ is compact semisimple, $\mathfrak{a}$ is abelian and $\mathfrak{s}$ is either isomorphic to $\mathfrak{sl}_2(\mathds{R})$ or to $\mathfrak{he}_d^\lambda$, $\lambda \in \mathds{Z}_+^d$. Moreover, the induced bilinear form $\kappa$ on $\mathfrak{g}$ is Lorentzian by Corollary~\ref{cor:condition_star} and Theorem~\ref{th:algebraic}.

$\kappa$ is ad-invariant and $G$ is connected, so $\kappa$ is $\textnormal{Ad}(G)$-invariant. Since the subgroup generated by $\mathfrak{s}$ acts locally freely on $M$ by Theorem~\ref{th:locally_free}, $\mathfrak{s}\cap\mathfrak{h}=\left\{0\right\}$. It follows from Theorem~\ref{th:algebraic} and Proposition~\ref{prop:H_compact}, that $\kappa$ restricted to $\mathfrak{h} \times \mathfrak{h}$ is positive definite. For this, note that $\mathfrak{z(s)}$ is $\kappa$-isotropic in the case $\mathfrak{s}\cong\mathfrak{he}_d^\lambda$.

Choosing $\mathfrak{m}$ to be the $\kappa$-orthogonal complement of $\mathfrak{h}$ in $\mathfrak{g}$, we are done.

\section{Geometry of homogeneous manifolds}\label{sec:homogeneous_geometry}

In this section, we will describe the geometry (mainly in terms of curvature) of compact homogeneous Lorentzian manifolds $M=(M,g)$, whose isometry groups $\textnormal{Isom}(M)$ have non-compact connected components. According to Theorem~\ref{th:homogeneous_characterization}, $\textnormal{Isom}^0(M)$ contains either a subgroup isomorphic to a central quotient of $\widetilde{\textnormal{SL}_2(\mathds{R})}$ or a subgroup isomorphic to a twisted Heisenberg group $\textnormal{He}_d^\lambda$ or $\overline{\textnormal{He}_d^\lambda}$.

In Paragraph~\ref{sec:homogeneous_general}, we give an overview of the theorems concerning reductive homogeneous semi-Riemannian manifolds, which we will use in the sequel. This presentation is standard, for the most part, we refer the reader to \cite{Ar03}, Chapters~4 and~5. Note that our sign of the Riemannian curvature tensor is different than in \cite{Ar03}.

In the Paragraphs~\ref{sec:homogeneous_sl2R} and~\ref{sec:homogeneous_twisted_heisenberg} we describe the geometry of the homogeneous manifold, if the Lie algebra of $\textnormal{Isom}(M)$ contains $\mathfrak{sl}_2(\mathds{R})$ or $\mathfrak{he}_d^\lambda$, respectively. Especially, we investigate the curvature of Lie groups with Lie algebra $\mathfrak{sl}_2(\mathds{R})$ or $\mathfrak{he}_d^\lambda$, respectively, provided with the bi-invariant metric defined by an ad-invariant Lorentzian scalar product on $\mathfrak{sl}_2(\mathds{R})$ or $\mathfrak{he}_d^\lambda$, respectively.

In Paragraph~\ref{sec:corollaries}, we prove Theorems~\ref{th:isotropy_representation} and~\ref{th:homogeneous_not_Ricci_flat}. These results follow directly from our investigation in Sections~\ref{sec:homogeneous_sl2R} and~\ref{sec:homogeneous_twisted_heisenberg}.

\subsection{Curvature and holonomy of homogeneous semi-Rie\-mannian manifolds}\label{sec:homogeneous_general}

Let $M=(M,g)$ be a homogeneous semi-Riemannian manifold and $x_0 \in M$. By definition, the action of $\textnormal{Isom}(M)$ on $M$ is transitive. Let $G$ be a Lie group and $\rho: G \to \textnormal{Isom}(M)$ be a Lie group homomorphism, such that $\rho(G)$ still acts transitively. Let $H \subseteq G$ be the isotropy group in $G$ of $x_0$, that is, $H$ consists of all elements of $G$ having $x_0$ as a fixed point. Then $H$ is a closed subgroup and $M\cong G/H$. We will assume that $G/H$ is reductive.

Denote by $\mathfrak{g}$ and $\mathfrak{h}$ the Lie algebras of $G$ and $H$, respectively, $\mathfrak{h} \subseteq \mathfrak{g}$. Since $G/H$ is reductive, there is an $\textnormal{Ad}(H)$-invariant vector space $\mathfrak{m}$  complementary to $\mathfrak{h}$ in $\mathfrak{g}$. Note that $\mathfrak{m}$ does not have to be a subalgebra. For $X \in \mathfrak{g}$, we will write $X_{\mathfrak{m}}$ for the $\mathfrak{m}$-component of $X$ with respect to the decomposition $\mathfrak{g}=\mathfrak{m}\oplus\mathfrak{h}$.

Remember the identification of $\mathfrak{isom}(M)$ with $\mathfrak{kill}_c(M)$ in Proposition~\ref{prop:isometry_Killing}, where $X \in \mathfrak{isom}(M)$ was identified with the complete Killing vector field $\widetilde{X}$ defined by $\widetilde{X}(x):=\frac{\partial}{\partial t}(\exp(tX) \cdot x) |_{t=0}$. The differential $d\rho: \mathfrak{g} \to \mathfrak{isom}(M)$ allows us then to obtain a canonical map $\mathfrak{g} \to \mathfrak{kill}_c(M)$, $X \mapsto \widetilde{X}$.

\begin{proposition}\label{prop:tangent_space}
$\mathfrak{m} \cong T_{x_0} M$, where the identification is given by $X \mapsto \widetilde{X}(x_0)$.
\end{proposition}
\vspace{-1.4em}\begin{proof}
Since $M\cong G/H$, the map $\mathfrak{g} \to T_{x_0} M$ defined in the same way as above is surjective and has kernel exactly $\mathfrak{h}$. Thus, $\mathfrak{m} \cong T_{x_0} M$.
\end{proof}\vspace{0pt}

\begin{definition}
The map $H \to \textnormal{GL}(T_{x_0} M)$, $h \to d(\varrho(h))_{x_0}$, is called the \textit{isotropy representation} of $G/H$. $G/H$ is called \textit{isotropy irreducible}, if the isotropy representation is irreducible, that is, there are no other $H$-invariant subspaces in $T_{x_0} M$ than $\left\{0\right\}$ and $T_{x_0} M$. $G/H$ is called \textit{weakly isotropy irreducible}, if any $H$-invariant subspace in $T_{x_0} M$ is either trivial or degenerate (with respect to the metric $g_{x_0}$).
\end{definition}
\begin{remark}
Note that isotropy irreducibility implies weak isotropy irreducibility.
\end{remark}

\begin{proposition}\label{prop:isotropy_representation}
In the situation above, the following is true:
\begin{compactenum}
\item Using the isomorphism $T_{x_0} M \cong \mathfrak{m}$, the isotropy representation of $G/H$ corresponds to the adjoint action of $H$ on $\mathfrak{m}$.

\item The metric $g$ on $M$ corresponds to an $\textnormal{Ad}(H)$-invariant scalar product $\langle \cdot, \cdot \rangle$ on $\mathfrak{m}$. Conversely, any $\textnormal{Ad}(H)$-invariant scalar product $\langle \cdot, \cdot \rangle$ on $\mathfrak{m}$ defines a semi-Riemannian metric $g$ on $M$.
\end{compactenum}
\end{proposition}
\vspace{-1.4em}\begin{proof}
(i) This follows from \cite{Ar03}, Proposition~4.5.

(ii) The proof of Proposition~5.1 in \cite{Ar03} works in the same way also in the case of semi-Riemannian metrics.
\end{proof}\vspace{0pt}

Let $\langle \cdot, \cdot \rangle$ denote the scalar product on $\mathfrak{m}$ associated to $g$.

\begin{definition}
The map $U: \mathfrak{m} \times \mathfrak{m} \to \mathfrak{m}$ is determined by \[2\langle U(X,Y),Z \rangle = \langle [Z,X]_{\mathfrak{m}} , Y \rangle + \langle X , [Z,Y]_{\mathfrak{m}} \rangle\] for all $X,Y,Z \in \mathfrak{m}$.
\end{definition}

Clearly, $U$ is symmetric and bilinear. $U$ vanishes in the case of an ad-invariant scalar product $\langle \cdot,\cdot \rangle$.

Let $\nabla$ denote the Levi-Civita connection of $M$.

\begin{definition}
For $u,v,w,z \in T_xM$, the \textit{Riemannian curvature tensor} $R$ is determined by \[R(u,v,w,z):=g_x (\nabla_{\overline{U}} \nabla_{\overline{V}} \overline{W} - \nabla_{\overline{V}} \nabla_{\overline{U}} \overline{W} - \nabla_{[\overline{U},\overline{V}]} \overline{W},z)\] for any vector fields $\overline{U},\overline{V},\overline{W}$ on $M$ extending $u,v,w$.

If $E$ is a non-degenerate subspace of $(T_x M,g_x)$ and $\left\{v,w\right\}$ a basis of $(E,g_x)$, the \textit{sectional curvature} $K_E$ is defined by \[K_E:=\frac{R(v,w,w,v)}{Q(v,w)},\] where $Q(v,w):= g_x(v,v) g_x(w,w)-g_x(v,w)^2$.

Let $\left\{v_1, \ldots, v_n \right\}$ be an orthonormal basis of $(T_x M,g_x)$ and $\varepsilon_j:= g_x(v_j,v_j)$ for $j=1,\ldots,n$. For $u,w \in T_x M$, the \textit{real-valued Ricci tensor} $\textnormal{Ric}$ is given by \[\textnormal{Ric}(u,w):=\sum\limits_{j=1}^n \varepsilon_j R(u,v_j,v_j,w).\] For $u \in T_x M$, the \textit{vector-valued Ricci tensor} $\textnormal{Ric}$ is determined by \[g_x(\textnormal{Ric}(u),w)=\textnormal{Ric}(u,w)\] for all $w \in T_x M$.

Finally, the scalar curvature $\textnormal{scal}$ is given by \[\textnormal{scal}:=\sum\limits_{j=1}^n \varepsilon_j \textnormal{Ric}(v_j,v_j).\]
\end{definition}

Note that isometries of $M$ respect the curvatures defined above. Since $M$ is homogeneous, it suffices to know the curvature at the single point $x_0$. We now characterize the Levi-Civita connection and the sectional curvature in terms of the decomposition $\mathfrak{g}=\mathfrak{m}\oplus\mathfrak{h}$ of the Lie algebra. We use the identification of $T_{x_0} M$ with $\mathfrak{m}$ given in Proposition~\ref{prop:tangent_space}.

\begin{proposition}\label{prop:curvature}
Let $X,Y \in \mathfrak{m}$. Then
\begin{align*}
\left(\nabla_{\widetilde{X}} \widetilde{Y}\right) \left(x_0\right)=&-\frac{1}{2} [X,Y]_{\mathfrak{m}}+U(X,Y),\\
R(X,Y,Y,X)=&-\frac{3}{4}\langle [X,Y]_{\mathfrak{m}},[X,Y]_{\mathfrak{m}} \rangle-\frac{1}{2}\langle [X,[X,Y]_{\mathfrak{m}}]_{\mathfrak{m}},Y \rangle \\
&-\frac{1}{2}\langle [Y,[Y,X]_{\mathfrak{m}}]_{\mathfrak{m}},X \rangle + \langle Y, [[X,Y]_{\mathfrak{h}},X] \rangle\\
&+\langle U(X,Y), U(X,Y) \rangle - \langle U(X,X), U(Y,Y) \rangle.
\end{align*}
\end{proposition}
\vspace{-1.4em}\begin{proof}
See \cite{Ar03}, Proposition~5.2 and Theorem~5.3.
\end{proof}\vspace{0pt}

For the definition of holonomy and basic results, we follow \cite{Ba09}, Chapter~5.

\begin{definition}
Let $x \in M$ and denote by $\Omega(x)$ the set of closed piecewise smooth paths $\gamma$ starting and ending in $x$. If we restrict to null-homotopic paths $\gamma$, we denote the corresponding set by $\Omega_0(x)$. Additionally, let $P_{\gamma}:T_x M \to T_xM$ be the parallel transport along $\gamma$ with respect to $\nabla$.

Then $\textnormal{Hol}_x(M):=\left\{P_{\gamma} | \gamma \in \Omega(x)\right\}$ is the \textit{holonomy group} of $M=(M,g)$ with respect to the base point $x$. $\textnormal{Hol}_x^0(M):=\left\{P_{\gamma} | \gamma \in \Omega_0(x)\right\}$ is the \textit{restricted holonomy group} of $M=(M,g)$ with respect to the base point $x$.
\end{definition}

\begin{remark}
Since parallel transport is an orthogonal transformation, \[\textnormal{Hol}_x^0(M)\subseteq \textnormal{Hol}_x(M) \subseteq O(T_x M, g_x).\]

If $\delta:[0,1] \to M$ is a piecewise smooth path connecting $\delta(0)=x_0$ and $\delta(1)=x$, the holonomy groups with respect to the base points $x$ and $x_0$ are conjugate to each other: \[\textnormal{Hol}_x(M)=P_\delta \circ \textnormal{Hol}_{x_0}(M) \circ P_\delta^{-1}.\]

Therefore, we know the holonomy of $M$ if we calculate only the holonomy group of $M$ with respect to the base point $x_0$.
\end{remark}

For the following, see \cite{Ba09}, Proposition~5.1:

\begin{proposition}\label{prop:holonomy}
$\textnormal{Hol}_x(M)$ is a Lie subgroup of $O(T_x M, g_x)$. The restricted holonomy group $\textnormal{Hol}_x^0(M)$ is the connected component of the identity in $\textnormal{Hol}_x(M)$.
\end{proposition}

\begin{definition}
The \textit{holonomy algebra} $\mathfrak{hol}_x(M)\subseteq \mathfrak{so}(T_xM,g_x)$ is the Lie algebra of $\textnormal{Hol}_x(M)$.
\end{definition}

If $\widetilde{X}$ is a Killing vector field on $M$, the mapping $(\nabla \widetilde{X})(x_0) :T_{x_0}M \to T_{x_0}M$, $v \mapsto \nabla_v \widetilde{X}$, is skew-symmetric with respect to $g_{x_0}$.

\begin{definition}
Using the identification of $(T_{x_0}M,g_{x_0})$ with $(\mathfrak{m}, \langle \cdot, \cdot \rangle)$ given in Proposition~\ref{prop:tangent_space}, $\mathfrak{so}(T_{x_0}M,g_{x_0}) \cong \mathfrak{so}(\mathfrak{m}, \langle \cdot, \cdot \rangle)$ and we can consider the homomorphism $\Lambda_{\mathfrak{m}}:\mathfrak{m} \to \mathfrak{so}(\mathfrak{m}, \langle \cdot, \cdot \rangle)$ defined by $\Lambda_{\mathfrak{m}}(X)Y:=\nabla_{\widetilde{Y}} \widetilde{X} (x_0)$ for all $X,Y \in \mathfrak{m}$.
\end{definition}

By Proposition~\ref{prop:curvature}, $\Lambda_{\mathfrak{m}}(X)Y=\frac{1}{2}[X,Y]_\mathfrak{m}+U(X,Y)$ for all $X,Y \in \mathfrak{m}$. Since $\mathfrak{m}$ and $\langle \cdot, \cdot \rangle$ are $\textnormal{Ad}(H)$-invariant, $\textnormal{ad}_Z|_{\mathfrak{m}} \in \mathfrak{so}(\mathfrak{m}, \langle \cdot, \cdot \rangle)$ for all $Z \in \mathfrak{h}$.

\begin{proposition}\label{prop:holonomy_algebra}
$\mathfrak{hol}_{x_0}(M)\subseteq \mathfrak{so}(T_{x_0}M,g_{x_0})\cong\mathfrak{so}(\mathfrak{m}, \langle \cdot, \cdot \rangle)$ is given by \[\mathfrak{hol}_{x_0}(M)=\mathfrak{m}_0+[\Lambda_{\mathfrak{m}}(\mathfrak{m}),\mathfrak{m}_0]+[\Lambda_{\mathfrak{m}}(\mathfrak{m}),[\Lambda_{\mathfrak{m}}(\mathfrak{m}),\mathfrak{m}_0]]+\ldots,\]
where $\mathfrak{m}_0$ is the subspace in $\mathfrak{so}(\mathfrak{m}, \langle \cdot, \cdot \rangle)$ spanned by \[\left\{[\Lambda_{\mathfrak{m}}(X),\Lambda_{\mathfrak{m}}(Y)]-\Lambda_{\mathfrak{m}}([X,Y]_{\mathfrak{m}})-\textnormal{ad}_{[X,Y]_{\mathfrak{h}}} \ | \ X,Y \in \mathfrak{m}\right\}.\]
\end{proposition}

\vspace{-1.4em}\begin{proof}
See \cite{KoNo69}, Chapter~X, Corollary~4.2.
\end{proof}\vspace{0pt}

\begin{corollary}\label{cor:holonomy_Lie_group}
Let $G$ be a Lie group furnished with a bi-invariant semi-Rie\-mannian metric $g$ and $x \in G$. Denote by $\mathfrak{g}$ the Lie algebra of $G$ and by $\langle \cdot, \cdot \rangle$ the ad-invariant scalar product on $\mathfrak{g}$ corresponding to $g$. Then \[\mathfrak{hol}_{x}(G)=\textnormal{ad}([\mathfrak{g},\mathfrak{g}])\subseteq \mathfrak{so}(\mathfrak{g}, \langle \cdot, \cdot \rangle).\]
\end{corollary}

\vspace{-1.4em}\begin{proof}
Following Proposition~\ref{prop:holonomy_algebra}, we first calculate $\mathfrak{g}_0$, which is the subspace in $\mathfrak{so}(\mathfrak{g}, \langle \cdot, \cdot \rangle)$ spanned by $\left\{[\Lambda_{\mathfrak{g}}(X),\Lambda_{\mathfrak{g}}(Y)]-\Lambda_{\mathfrak{g}}([X,Y]) | X,Y \in \mathfrak{g}\right\}.$ For any $Z \in \mathfrak{g}$, \[[\Lambda_{\mathfrak{g}}(X),\Lambda_{\mathfrak{g}}(Y)](Z)-\Lambda_{\mathfrak{g}}([X,Y])(Z)=\frac{1}{4}([X,[Y,Z]]-[Y,[X,Z]])-\frac{1}{2}[[X,Y],Z]\] because $U \equiv 0$. By the Jacobi identity, $[X,[Y,Z]]-[Y,[X,Z]]=[[X,Y],Z]$. It follows that $\mathfrak{g}_0$ is equal to $\textnormal{ad}([\mathfrak{g},\mathfrak{g}])$.

Now $\Lambda_{\mathfrak{g}}(\mathfrak{g})=\textnormal{ad}(\mathfrak{g})$. Using that $\textnormal{ad}$ is a Lie algebra homomorphism, \begin{equation*}\mathfrak{hol}_x(G)=\textnormal{ad}\left([\mathfrak{g},\mathfrak{g}]+[\mathfrak{g},[\mathfrak{g},\mathfrak{g}]]+[\mathfrak{g},[\mathfrak{g},[\mathfrak{g},\mathfrak{g}]]]+\ldots \right)=\textnormal{ad}([\mathfrak{g},\mathfrak{g}]). \qedhere\end{equation*}
\end{proof}\vspace{0pt}

\subsection{Isometry group contains a cover of the projective special linear group}\label{sec:homogeneous_sl2R}

In this section, let $M=(M,g)$ be a compact homogeneous Lorentzian manifold and assume that the Lie algebra of $G:=\textnormal{Isom}^0(M)$ contains a direct summand isomorphic to $\mathfrak{sl}_2(\mathds{R})$. Let $H$ be the isotropy group in $G$ of some $x \in M$. Clearly, $M \cong G/H$.

First, we will investigate the isotropy representation of $G/H$ in Section~\ref{sec:isotropy_sl2R}. Later on, in Section~\ref{sec:local_geometry_sl2R}, we describe the local geometry of the manifold $M$ in terms of the holonomy and curvature of $\widetilde{\textnormal{SL}_2(\mathds{R})}$, furnished with the metric given by a multiple of its Killing form, and of a compact homogeneous Riemannian manifold $N$. Finally, we investigate the local geometry of central quotients of $\widetilde{\textnormal{SL}_2(\mathds{R})}$ with the metric defined by a positive multiple of the Killing form of $\mathfrak{sl}_2(\mathds{R})$ in Section~\ref{sec:sl2R}.

\subsubsection{Isotropy representation}\label{sec:isotropy_sl2R}

According to Corollary~\ref{cor:condition_star} and Theorem~\ref{th:algebraic}, the Lie algebra of $G$ decomposes as a $\kappa$-orthogonal direct sum $\mathfrak{g}=\mathfrak{k}\oplus\mathfrak{a}\oplus\mathfrak{s}$, $\kappa$ being the induced bilinear form on $\mathfrak{g}$. Here $\mathfrak{a}$ is abelian, $\mathfrak{k}$ is compact semisimple and $\mathfrak{s}\cong\mathfrak{sl}_2(\mathds{R})$. Moreover, $\kappa$ restricted to $\mathfrak{s} \times \mathfrak{s}$ is a positive multiple of the Killing form and the restriction of $\kappa$ to $(\mathfrak{k}\oplus\mathfrak{a})\times(\mathfrak{k}\oplus\mathfrak{a})$ is positive definite.

Due to Proposition~\ref{prop:H_compact}, it holds for the Lie algebra $\mathfrak{h}$ of $H$, that $\mathfrak{h}\subseteq \mathfrak{k}\oplus\mathfrak{a}$. Choosing $\mathfrak{m}$ to be the $\kappa$-orthogonal complement of $\mathfrak{h}$ in $\mathfrak{g}$ as in the proof of Theorem~\ref{th:homogeneous_reductive}, we have the reductive decomposition $\mathfrak{g}=\mathfrak{m}\oplus\mathfrak{h}$. Note that $\mathfrak{s} \subseteq \mathfrak{m}$ by definition. Therefore, $\mathfrak{m}=\mathfrak{p}\oplus\mathfrak{s}$, where $\mathfrak{p}$ is the $\kappa$-orthogonal complement of $\mathfrak{h}$ in $\mathfrak{k}\oplus\mathfrak{a}$.

By Theorem~\ref{th:homogeneous_characterization}, $M$ is isometric to $\Gamma\backslash\mkern-5mu\left(N \times \widetilde{\textnormal{SL}_2(\mathds{R})}\right)$, where $N$ is a compact homogeneous Riemannian manifold and $S:=\widetilde{\textnormal{SL}_2(\mathds{R})}$ is provided with the metric defined by a positive multiple of the Killing form of $\mathfrak{s}=\mathfrak{sl}_2(\mathds{R})$. Also, $G$ is a central quotient of $C \times S$, where $C$ is a subgroup of the isometry group of $N$ acting transitively on $N$. It follows from the proof of Corollary~\ref{cor:isometry_group}, that $C$ can be identified with the subgroup in $G$ generated by $\mathfrak{c}:=\mathfrak{k}\oplus\mathfrak{a}$.

Let $\langle \cdot,\cdot \rangle$ denote the Lorentzian scalar product on $\mathfrak{m}$ corresponding to the metric on $M$. It follows from the construction, that $\mathfrak{m}=\mathfrak{p}\oplus\mathfrak{s}$ is $\langle \cdot,\cdot \rangle$-orthogonal and $\langle \cdot,\cdot \rangle$ restricted to $\mathfrak{s} \times \mathfrak{s}$ is a multiple of the Killing form. The restriction to $\mathfrak{p} \times \mathfrak{p}$ is Riemannian and corresponds to the metric on $N$. For this, we know that $C$ acts transitively on $N$ and since $C \subset G$, the isotropy group $H_C$ of the point $x$ in $N$ (remember that we can consider $N$ as a leaf in $M$; see Section~\ref{sec:structure}) in the group $C$ is contained in $H$. Moreover, since $\mathfrak{h}\subseteq \mathfrak{c}$, it follows that the Lie algebra of $H_C$ is $\mathfrak{h}$ as well. Thus, $N \cong C/H_C$ and $\mathfrak{c}=\mathfrak{p}\oplus\mathfrak{h}$ is a reductive decomposition. Note that the $\textnormal{Ad}(H_C)$-invariance of $\mathfrak{p}$ follows from the fact that the latter space is the $\kappa$-orthogonal complement of $\mathfrak{h}$.

\begin{proposition}\label{prop:sl2R_isotropy}
In the situation described above, a decomposition of $\mathfrak{m}$ into a direct sum of irreducible $\textnormal{Ad}(H)$-invariant subspaces is given by the $\langle \cdot,\cdot \rangle$-orthogon\-al sum $\mathfrak{p}_1\oplus\ldots\oplus\mathfrak{p}_k\oplus\mathfrak{s}$, where the $\mathfrak{p}_j$ are irreducible subspaces of $\mathfrak{p}$.

Especially, the space $G/H$ is weakly isotropy irreducible if and only if $M$ is isometric to $\widetilde{\textnormal{SL}_2(\mathds{R})}/\Gamma$. In this case, it is also isotropy irreducible.
\end{proposition}
\vspace{-1.4em}\begin{proof}
By Theorem~\ref{th:homogeneous_characterization}, there is a uniform lattice $\Gamma_0$ in $S$ and a homomorphism $\varrho:\Gamma_0 \to \textnormal{Isom}(N)$, such that $\Gamma$ is the graph of $\varrho$. Remember that for an element $(\gamma,\varrho(\gamma)) \in \Gamma$, $\gamma$ is acting on the $S$-part by multiplication from the right (see Section~\ref{sec:structure}).

Denote $p: N \times S \to M$ the covering given by the theorem. Let $e$ be the identity element of $S$ and $\gamma \in \Gamma_0$. Since $C$ acts isometrically on $N$, there is $\psi \in C$, such that $\varrho(\gamma)(x)=\psi(x)$. Consider now the mapping $\varphi_{\gamma}: M \to M$ defined by \[\varphi_{\gamma}(p(y,f)):=p(\psi(y),f\gamma).\] Due to Theorem~\ref{th:homogeneous_characterization}, $C$ centralizes $\Gamma$, therefore, $\varphi_{\gamma}$ is correctly defined. By construction, $\varphi_{\gamma}$ is an isometry and \[x=p(x,e)=p(\varrho(\gamma)(x),\gamma)=p(\psi(x),\gamma),\] that is, $\varphi_{\gamma} \in H$.

$\mathfrak{s}$ is an ideal in $\mathfrak{g}$, especially $\textnormal{Ad}(H)$-invariant. $\mathfrak{p}$ is the $\langle \cdot,\cdot \rangle$-orthogonal complement of $\mathfrak{s}$ in $\mathfrak{m}$, so it is $\textnormal{Ad}(H)$-invariant as well. Thus, $\mathfrak{m}=\mathfrak{p}\oplus\mathfrak{s}$ is a decomposition into invariant subspaces. Since the restriction of $\langle \cdot,\cdot \rangle$ to $\mathfrak{p}$ is positive definite, we obtain by induction an orthogonal decomposition $\mathfrak{p}=\mathfrak{p}_1\oplus\ldots\oplus\mathfrak{p}_k$ into invariant irreducible subspaces. It remains to show that $\mathfrak{s}$ is irreducible.

Because $C$ centralizes $\mathfrak{s}$, the adjoint action of $C$ on $\mathfrak{s}$ is trivial. It follows that $\textnormal{Ad}_{\varphi_{\gamma}}|_{\mathfrak{s}}=\textnormal{Ad}^S_{\gamma}$, $\textnormal{Ad}^S$ being the adjoint action of $S$. It follows that any $\textnormal{Ad}(H)$-invariant subspace of $\mathfrak{s}$ is also $\textnormal{Ad}^S(\Gamma_0)$-invariant. But $\textnormal{Ad}^S(\Gamma_0)$ is Zariski-dense in $\textnormal{Ad}^S(S)$ (cf.~\cite{OnVin00}, Part~I, Chapter~3, Theorem~1.2) and $\mathfrak{s}$ contains no non-trivial proper ideals. Thus, $\mathfrak{s}$ is irreducible.
\end{proof}\vspace{0pt}
\begin{remark}
Note that the subspaces $\mathfrak{p_j}$ in Proposition~\ref{prop:sl2R_isotropy} are $\textnormal{Ad}(H_C)$-invariant, since $H_C \subseteq H$.
\end{remark}

Proposition~\ref{prop:sl2R_isotropy} shows the first part of Theorem~\ref{th:isotropy_representation}.

\subsubsection{Local geometry of the manifold}\label{sec:local_geometry_sl2R}

Due to Theorem~\ref{th:homogeneous_characterization}, $M$ is covered by the metric product $\widetilde{M}:=N \times \widetilde{\textnormal{SL}_2(\mathds{R})}$, where $N$ is a compact homogeneous Riemannian manifold and $S:=\widetilde{\textnormal{SL}_2(\mathds{R})}$ is provided with the metric defined by a positive multiple of the Killing form of $\mathfrak{s}=\mathfrak{sl}_2(\mathds{R})$.

Since the local geometry of $M$ and $\widetilde{M}$ coincide, it suffices to investigate the homogeneous space $\widetilde{M}$. Let $(x,f) \in N \times S$ be arbitrary. Because $\widetilde{M}$ is a metric product of the two homogeneous spaces $S$ and $N$, we can decompose any $v$ in $T_{(x,f)} \widetilde{M}$ uniquely into $v=v_N+v_S$, where $v_N \in T_x N$ and $v_S \in T_f S$. Also, for any $u,v,w,z \in T_{(x,f)} \widetilde{M}$,
\begin{align*}
R(u,v,w,z)&=R^N(u_N,v_N,w_N,z_N)+R^S(u_S,v_S,w_S,z_S),\\
\textnormal{Ric}(u,v)&=\textnormal{Ric}^N(u_N,v_N)+\textnormal{Ric}^S(u_S,v_S),\\
\textnormal{scal}&=\textnormal{scal}^N+\textnormal{scal}^S.
\end{align*}
Here $R^N,\textnormal{Ric}^N,\textnormal{scal}^N$ and $R^S,\textnormal{Ric}^S,\textnormal{scal}^S$ correspond to the homogeneous spaces $N$ and $S$, respectively.

According to the second part of Proposition~5.4 in \cite{Ba09}, \[\textnormal{Hol}_{(x,f)}(\widetilde{M})\cong\textnormal{Hol}_{x}(N) \times \textnormal{Hol}_{f}(S).\] Hence, $\mathfrak{hol}_{(x,f)}(\widetilde{M})=\mathfrak{hol}_{x}(N) \oplus \mathfrak{hol}_{f}(S)$, where $\mathfrak{hol}_{x}(N)$  and $\mathfrak{hol}_{f}(S)$ are embedded canonically in $\mathfrak{hol}_{(x,f)}(\widetilde{M})$.

Thus, it remains to determine the curvature and holonomy algebra of $S$.

\subsubsection{Curvature and holonomy of the two-dimensional special linear group}\label{sec:sl2R}

As above, let $S:=\widetilde{\textnormal{SL}_2(\mathds{R})}$ be furnished with the bi-invariant metric given by a positive multiple $\langle \cdot,\cdot \rangle=\lambda k$ of the Killing form $k$ of $\mathfrak{s}=\mathfrak{sl}_2(\mathds{R})$. Note that $S$ is a symmetric space. Choose $x \in S$ arbitrary.

Let $\left\{e,f,h\right\}$ be an $\mathfrak{sl}_2$-triple of $\mathfrak{s}$, that is, $[h,e]=2e$, $[h,f]=-2f$ and $[e,f]=h$. The Killing form $k$ with respect to the ordered basis $(e,f,h)$ is determined by the matrix \[\begin{pmatrix} 0 & 4 & 0 \\ 4 & 0 & 0 \\ 0 & 0 & 8\end{pmatrix}.\]

Using the ad-invariance of $\langle \cdot,\cdot \rangle$, we obtain $U \equiv 0$, and by Proposition~\ref{prop:curvature}, we have for $X,Y \in \mathfrak{s}$:
\begin{align*}
R(X,Y,Y,X)&=-\frac{3}{4}\langle [X,Y],[X,Y] \rangle-\frac{1}{2}\langle [X,[X,Y]],Y \rangle -\frac{1}{2}\langle [Y,[Y,X]],X \rangle \\
&=-\frac{3}{4}\langle [X,Y],[X,Y] \rangle+\frac{1}{2}\langle [X,Y],[X,Y] \rangle +\frac{1}{2}\langle [Y,X],[Y,X] \rangle\\
&=\frac{1}{4}\langle [X,Y],[X,Y] \rangle.
\end{align*}

If $\left\{X_1,X_2,X_3\right\} \subset \mathfrak{s}$ is a $(\lambda k)$-orthonormal basis, $\varepsilon_j:=\lambda k(X_j,X_j)$ for all $j$, then for any $X \in \mathfrak{s}$,
\begin{align*}
\begin{array}{rcl}
\textnormal{Ric}(X,X)=&\sum\limits_{j=1}^3 \varepsilon_j R(X,X_j,X_j,X)
&=\sum\limits_{j=1}^3 \frac{1}{4}\varepsilon_j\langle [X,X_j],[X,X_j] \rangle\\
=&\sum\limits_{j=1}^3 -\frac{1}{4}\varepsilon_j\langle [X,[X,X_j]],X_j \rangle
&=-\frac{1}{4} k (X,X).
\end{array}
\end{align*}
It follows that $S$ is an Einstein manifold. We conclude \[\textnormal{scal}=\sum\limits_{j=1}^3 \varepsilon_j \textnormal{Ric}(X_j,X_j)=-\frac{3}{4\lambda}.\]

Note that $S$ has constant sectional curvature $K$, since any connected Einstein manifold of dimension 3 has. If we consider $E=\textnormal{span}\left\{e,f\right\}$, \[K=K_E=\frac{\langle [e,f],[e,f] \rangle}{4 (\langle e,e \rangle \langle f,f \rangle - \langle e,f \rangle^2)}=\frac{8 \lambda}{-4 \cdot (4 \lambda)^2}=-\frac{1}{8\lambda}.\]

Due to Corollary~\ref{cor:holonomy_Lie_group}, the holonomy algebra of $S$ is given by \[\mathfrak{hol}_x(S)=\textnormal{ad}([\mathfrak{s},\mathfrak{s}])\subseteq \mathfrak{so}(\mathfrak{s}, \langle \cdot, \cdot \rangle).\] But $\mathfrak{s}$ is semisimple, so $\textnormal{ad}([\mathfrak{s},\mathfrak{s}])=\textnormal{ad}(\mathfrak{s})\cong \mathfrak{s}$ is three-dimensional. $\mathfrak{so}(\mathfrak{s}, \langle \cdot, \cdot \rangle)$ is is three-dimensional as well, therefore, \[\mathfrak{hol}_x(S)=\mathfrak{so}(\mathfrak{s}, \langle \cdot, \cdot \rangle).\]

\subsection{Isometry group contains a twisted Heisenberg group}\label{sec:homogeneous_twisted_heisenberg}

Throughout this paragraph, we assume that the Lie algebra of $G:=\textnormal{Isom}^0(M)$, $M=(M,g)$ a compact homogeneous Lorentzian manifold, contains a direct summand $\mathfrak{s}$ isomorphic to $\mathfrak{he}_d^\lambda$, $\lambda \in \mathds{Z}_+^d$. Let $S$ be the subgroup generated by $\mathfrak{s}$ and $H$ be the isotropy group in $G$ of some $x \in M$.

First, we will choose in Section~\ref{sec:reductive_representation} a slightly different reductive decomposition than in Section~\ref{sec:reductive}, which is adapted to the covering space $S \times_{Z(S)} N$. In Section~\ref{sec:isotropy_heis}, we will investigate the isotropy representation of $M\cong G/H$. Before finally describing the local geometry of the manifold $M$ in Section~\ref{sec:local_geometry_heis}, we describe the curvature and holonomy of twisted Heisenberg algebras furnished with the bi-invariant metric defined by an ad-invariant Lorentzian scalar product on its Lie algebra in Section~\ref{sec:curvature_twisted_Heisenberg}.

\subsubsection{Reductive representation}\label{sec:reductive_representation}

According to Corollary~\ref{cor:condition_star} and Theorem~\ref{th:algebraic}, the Lie algebra of $G$ decomposes as a $\kappa$-orthogonal direct sum $\mathfrak{g}=\mathfrak{s}\oplus\mathfrak{k}\oplus\mathfrak{a}$, $\kappa$ being the induced bilinear form on $\mathfrak{g}$. $\mathfrak{a}$ is abelian and $\mathfrak{k}$ is semisimple. Additionally, $\kappa$ restricted to $\mathfrak{s} \times \mathfrak{s}$ is an ad-invariant Lorentzian scalar product on $\mathfrak{s}$ and the restriction to $(\mathfrak{k}\oplus\mathfrak{a})\times(\mathfrak{k}\oplus\mathfrak{a})$ is positive definite.

Due to Proposition~\ref{prop:H_compact}, it holds for the Lie algebra $\mathfrak{h}$ of $H$, that $\mathfrak{h}\subseteq \mathfrak{z(s)}\oplus\mathfrak{k}\oplus\mathfrak{a}$. By Theorem~\ref{th:locally_free}, $\mathfrak{s}\cap\mathfrak{h}=\left\{0\right\}$. Thus, we can consider the $\kappa$-orthogonal complement $\mathfrak{m}^\prime$  of $\mathfrak{h}$ in $\mathfrak{z(s)}\oplus\mathfrak{k}\oplus\mathfrak{a}$. Then $\mathfrak{m}^\prime$ is $\textnormal{Ad}(H)$-invariant. Note that $\mathfrak{z(s)}\subseteq \mathfrak{m}^\prime$.

$\mathfrak{s}$ is an ideal in $\mathfrak{g}$ and hence $\textnormal{Ad}(G)$-invariant. It follows that $\mathfrak{m}:=\mathfrak{s}+\mathfrak{m}^\prime$ is $\textnormal{Ad}(H)$-invariant and complementary to $\mathfrak{h}$. Thus, $\mathfrak{g}=\mathfrak{m}\oplus \mathfrak{h}$ is a reductive decomposition.

By Theorem~\ref{th:homogeneous_characterization}, $M$ is isometric to $\Gamma\backslash\mkern-5mu\left(S \times_{Z(S)} N\right)$, where $N$ is a compact homogeneous Riemannian manifold and $S$ is provided with the metric defined by an ad-invariant Lorentzian scalar product on $\mathfrak{s}$. Also, $G$ is a central quotient of $S \times_{Z(S)} C$, where $C$ is a subgroup of the isometry group of $N$ acting transitively on $N$. It follows from the proof of Corollary~\ref{cor:isometry_group}, that $C$ can be identified with the subgroup in $G$ generated by $\mathfrak{c}:=\mathfrak{z(s)}\oplus\mathfrak{k}\oplus\mathfrak{a}$.

Since $C \subset G$, the isotropy group $H_C$ of the point $x$ in $N$ (remember that we can consider $N$ as a leaf in $M$; see Section~\ref{sec:structure}) in the group $C$ is contained in $H$. Moreover, since $\mathfrak{h}\subseteq \mathfrak{c}$, it follows that the Lie algebra of $H_C$ is $\mathfrak{h}$ as well. Thus, $N \cong C/H_C$ and $\mathfrak{c}=\mathfrak{m}^\prime\oplus\mathfrak{h}$ is a reductive decomposition. Note that the $\textnormal{Ad}(H_C)$-invariance of $\mathfrak{m}^\prime$ follows from the fact that the latter space is the $\kappa$-orthogonal complement of $\mathfrak{h}$.

The metric on $N$ corresponds to an $\textnormal{Ad}(H_C)$-invariant Riemannian scalar product $(\cdot,\cdot)$ on $\mathfrak{m}^\prime$. Let $\mathfrak{p}$ be the $(\cdot,\cdot)$-orthogonal complement of $\mathfrak{z(s)}$ in $\mathfrak{m}^\prime$. We obtain $\mathfrak{m}^\prime=\mathfrak{z(s)}\oplus\mathfrak{p}$ and $\mathfrak{m}=\mathfrak{s}\oplus\mathfrak{p}$.

Let $\langle \cdot,\cdot \rangle$ denote the Lorentzian scalar product on $\mathfrak{m}$ corresponding to the metric on $M$. It follows from the construction, that the direct sum $\mathfrak{m}=\mathfrak{s}\oplus\mathfrak{p}$ is $\langle \cdot,\cdot \rangle$-orthogonal and $\langle \cdot,\cdot \rangle$ restricted to $\mathfrak{s} \times \mathfrak{s}$ is an ad-invariant Lorentzian scalar product. The restriction to $\mathfrak{p} \times \mathfrak{p}$ is Riemannian and equals $(\cdot,\cdot)$.

Note that $S$ furnished with the metric defined by the restriction of $\langle \cdot, \cdot \rangle$ to $\mathfrak{s} \times \mathfrak{s}$, is a symmetric space.

\subsubsection{Isotropy representation}\label{sec:isotropy_heis}

\begin{proposition}\label{prop:heis_isotropy}
In the situation of Section~\ref{sec:reductive_representation}, a decomposition of $\mathfrak{m}$ into a direct sum of weakly irreducible $\textnormal{Ad}(H)$-invariant subspaces is given by the $\langle \cdot,\cdot \rangle$-orthogonal sum $\mathfrak{s}\oplus\mathfrak{p}_1\oplus\ldots\oplus\mathfrak{p}_k$, where the $\mathfrak{p}_j$ are irreducible subspaces of $\mathfrak{p}$. Furthermore, $\mathfrak{s}$ is not irreducible, but cannot be decomposed into irreducible invariant subspaces.

Especially, the space $G/H$ is not isotropy irreducible. It is weakly isotropy irreducible if and only if $M\cong S/\Gamma$.
\end{proposition}
\vspace{-1.4em}\begin{proof}
By Theorem~\ref{th:homogeneous_characterization}, there is a uniform lattice $\Gamma_0$ in $S/Z(S)$, such that $\Gamma$ projects isomorphically to $\Gamma_0$. Remember that for $[\gamma,\psi_\gamma] \in \Gamma \subset S \times_{Z(S)} \textnormal{Isom}_{Z(S)}(N)$, $[\gamma,\psi_\gamma]$ is acting on $S \times_{Z(S)} N$ by \[[\gamma,\psi_\gamma]([f,x])\mapsto [f\gamma,\psi_\gamma(x)] \textnormal{ (see Section~\ref{sec:structure}).}\]
Denote by $\pi: S \times_{Z(S)} N \to M$ the covering given by the theorem. Let $e$ be the identity element of $S$. Choose $\gamma \in \Gamma_0 \cdot Z(S) \subset S$ arbitrary and $\psi_\gamma$ in $\textnormal{Isom}_{Z(S)}(N)$, such that $[\gamma,\psi_\gamma] \in \Gamma$.

Since $C$ acts transitively on $N$, there is $\psi \in C$, such that $\psi_\gamma(x)=\psi(x)$. Consider now the map $\varphi_{\gamma}: M \to M$ defined by \[\varphi_{\gamma}(\pi([f,y])):=\pi[(f\gamma,\psi(y))].\] Due to Theorem~\ref{th:homogeneous_characterization}, $C\subseteq \textnormal{Isom}_{Z(S)}(N)$ centralizes $\Gamma$, therefore, $\varphi_{\gamma}$ is correctly defined. By construction, $\varphi_{\gamma}$ is an isometry and \[x=\pi([e,x])=\pi([\gamma,\psi_\gamma(x)])=\pi([\gamma,\psi(x)]),\] that is, $\varphi_{\gamma} \in H$.

$\mathfrak{s}$ is an ideal in $\mathfrak{g}$, especially $\textnormal{Ad}(H)$-invariant. $\mathfrak{p}$ is the $\langle \cdot,\cdot \rangle$-orthogonal complement of $\mathfrak{s}$ in $\mathfrak{m}$, so it is $\textnormal{Ad}(H)$-invariant as well. Thus, $\mathfrak{m}=\mathfrak{s}\oplus\mathfrak{p}$ is a decomposition into invariant subspaces. Since the restriction of $\langle \cdot,\cdot \rangle$ to $\mathfrak{p}$ is positive definite, we obtain by induction an orthogonal decomposition $\mathfrak{p}=\mathfrak{p}_1\oplus\ldots\oplus\mathfrak{p}_k$ into invariant irreducible subspaces. It remains to show that $\mathfrak{s}$ is weakly irreducible.

Because $C$ centralizes $\mathfrak{s}$, the adjoint action of $C$ on $\mathfrak{s}$ is trivial. It follows that $\textnormal{Ad}_{\varphi_{\gamma}}|_{\mathfrak{s}}=\textnormal{Ad}^S_{\gamma}$, $\textnormal{Ad}^S$ being the adjoint action of $S$. Since the adjoint action of the center is trivial, $S/Z(S)$ acts on $\mathfrak{s}$ through the adjoint action as well. We denote this action by $\textnormal{Ad}^{S/Z(S)}$.

Thus, any $\textnormal{Ad}(H)$-invariant subspace of $\mathfrak{s}$ is also $\textnormal{Ad}^{S/Z(S)}(\Gamma_0)$-invariant.

Let $N$ denote the nilradical of $S/Z(S)$. $N$ is isomorphic to $\widetilde{\textnormal{He}_d}/Z(\widetilde{\textnormal{He}_d})\cong \mathds{R}^{2d}$. Moreover, $\Gamma_0 \cap N$ is a lattice in $N$ (cf.~\cite{OnVin00}, Part~I, Chapter~2, Theorem~3.6), and $\Gamma_0 \cap N$ is Zariski-dense in $N$ since $N$ is simply-connected and nilpotent (cf.~\cite{OnVin00}, Part~I, Chapter~2, Theorem~2.4). It follows that any subspace of $\mathfrak{s}$ invariant under $\textnormal{Ad}^{S/Z(S)}(\Gamma_0)$ is $\textnormal{ad}(\mathfrak{he}_d)$-invariant.

Let $\mathfrak{i} \subseteq \mathfrak{s}$ be a non-trivial, $\textnormal{ad}(\mathfrak{he}_d)$-invariant subspace. If $\mathfrak{i} \subseteq \mathfrak{he}_d$, it follows that $Z \in \mathfrak{i}$ and $\mathfrak{i}$ is degenerate. Otherwise, $T+X \in \mathfrak{i}$, $X \in \mathfrak{he}_d$. But $[\mathfrak{he}_d,T]$ is equal to the subspace generated by $X_1,Y_1,\ldots,X_d,Y_d$. Thus, all $X_k$ and $Y_k$, $k=1,\ldots, d$, are contained in $\mathfrak{i} + \mathds{R}Z$. Applying $\textnormal{ad}(\mathfrak{he}_d)$-invariance another time, we obtain $Z \in \mathfrak{i}$, so $\mathfrak{i}=\mathfrak{s}$. It follows that $\mathfrak{s}$ is weakly irreducible.

Since $\mathfrak{he}_d$ is an ideal in $\mathfrak{g}$, it is $\textnormal{Ad}(H)$-invariant. Thus, $\mathfrak{s}$ is not irreducible. We have seen above that any non-trivial $\textnormal{Ad}(H)$-invariant subspace $\mathfrak{i}\subseteq \mathfrak{s}$ contains $\mathfrak{z(s)}$. Therefore, $\mathfrak{s}$ cannot be decomposed into irreducible invariant subspaces.
\end{proof}\vspace{0pt}
\begin{remark}
Note that the subspaces $\mathfrak{p_j}$ in Proposition~\ref{prop:heis_isotropy} are $\textnormal{Ad}(H_C)$-invariant, since $H_C \subseteq H$.
\end{remark}

\subsubsection{Curvature and holonomy of a twisted Heisenberg group}\label{sec:curvature_twisted_Heisenberg}

In the following, we consider a twisted Heisenberg group $S$ furnished with the bi-invariant metric defined by an ad-invariant Lorentzian scalar product $\langle \cdot,\cdot \rangle$ on $\mathfrak{s}\cong\mathfrak{he}_d^\lambda$. In this case, $S$ is a symmetric space.

Let $\left\{T, Z, X_1, Y_1, \ldots, X_d, Y_d\right\}$ be a canonical basis of $\mathfrak{s}$. By Proposition~\ref{prop:Lorentz_s_only_one}, all ad-invariant Lorentzian scalar products on $\mathfrak{s}$ are equivalent. So without loss of generality, we may assume that $\left\{X_1, Y_1, \ldots, X_d, Y_d\right\}$ is $\langle \cdot, \cdot \rangle$-orthonormal, \[\langle T,T \rangle=0 \text{ and }\langle T,Z \rangle=1.\] Note that $\langle Z,Z\rangle =0$.

Using the ad-invariance of $\langle \cdot,\cdot \rangle$, we obtain $U \equiv 0$, and by Proposition~\ref{prop:curvature}, we have for $X,Y \in \mathfrak{s}$:
\begin{align*}
R(X,Y,Y,X)&=-\frac{3}{4}\langle [X,Y],[X,Y] \rangle-\frac{1}{2}\langle [X,[X,Y]],Y \rangle -\frac{1}{2}\langle [Y,[Y,X]],X \rangle \\
&=-\frac{3}{4}\langle [X,Y],[X,Y] \rangle+\frac{1}{2}\langle [X,Y],[X,Y] \rangle +\frac{1}{2}\langle [Y,X],[Y,X] \rangle\\
&=\frac{1}{4}\langle [X,Y],[X,Y] \rangle.
\end{align*}
Especially, the sectional curvature vanishes for non-degenerate two-dimensional subspaces of $\mathfrak{he}_d$.

Since \[R(X,Y,Y,X)=\frac{1}{4}\langle [X,Y],[X,Y] \rangle=-\frac{1}{4} \langle [X,[X,Y]],Y \rangle,\] it holds for any $X \in \mathfrak{s}$, that \[\textnormal{Ric}(X,X)=-\frac{1}{4} \textnormal{Tr}(\textnormal{ad}_X^2)=-\frac{1}{4} k(X,X),\] $k$ being the Killing form of $\mathfrak{s}$. It follows that \[\textnormal{Ric}(X,Y)=-\frac{1}{4} k(X,Y)\] for all $X,Y \in \mathfrak{s}$. Especially, \[\textnormal{Ric}(X,Y)=0 \text{ and } \textnormal{Ric}(T+X,T+Y)=\frac{1}{2} \sum\limits_{k=1}^d \lambda_k^2\] for all $X,Y \in \mathfrak{he}_d$. Therefore, \[\textnormal{Ric}(X)=0\] if $X \in \mathfrak{he}_d$ and \[\textnormal{Ric}(T)=\frac{1}{2} \sum\limits_{k=1}^d \lambda_k^2 Z.\] It follows that the vector-valued Ricci tensor is totally isotropic.

Furthermore, $\left\{\frac{1}{\sqrt{2}}(T-Z), \frac{1}{\sqrt{2}}(T+Z), X_1, Y_1, \ldots, X_d, Y_d\right\}$ is an orthonormal basis with respect to $\langle \cdot,\cdot \rangle$. Thus, the scalar curvature $\textnormal{scal}$ is given by \[\frac{1}{2} (\textnormal{Ric}(T+Z,T+Z)- \textnormal{Ric}(T-Z,T-Z))+\sum\limits_{k=1}^d (\textnormal{Ric}(X_k,X_k)+\textnormal{Ric}(Y_k,Y_k)),\] which is identically $0$.

Due to Corollary~\ref{cor:holonomy_Lie_group}, the holonomy algebra of $S$ is given by \[\mathfrak{hol}_x(S)=\textnormal{ad}([\mathfrak{s},\mathfrak{s}])=\textnormal{ad}(\mathfrak{he}_d)\subseteq \mathfrak{so}(\mathfrak{s}, \langle \cdot, \cdot \rangle).\]
The kernel of ad is exactly equal to the center of $\mathfrak{s}$, therefore, $\textnormal{ad}(\mathfrak{he}_d)$ is a $2d$-dimensional abelian subalgebra of $\mathfrak{so}(\mathfrak{s}, \langle \cdot, \cdot \rangle)$.

\subsubsection{Local geometry of the manifold}\label{sec:local_geometry_heis}

We use the reductive representation of $M$ given in Section~\ref{sec:reductive_representation} and use the same notations. For $X \in \mathfrak{m}=\mathfrak{s}\oplus \mathfrak{p}$, we will write $X_{\mathfrak{s}}$ for the $\mathfrak{s}$-component and $X_{\mathfrak{p}}$ for the $\mathfrak{p}$-component.

Let $\left\{T, Z, X_1, Y_1, \ldots, X_d, Y_d\right\}$ be a canonical basis of $\mathfrak{s}$. By Proposition~\ref{prop:Lorentz_s_only_one}, all ad-invariant Lorentzian scalar products on $\mathfrak{s}$ are equivalent. Since the scalar product $\langle \cdot, \cdot \rangle$ restricted to $\mathfrak{s} \times \mathfrak{s}$ is ad-invariant, we may assume without loss of generality, that $\left\{X_1, Y_1, \ldots, X_d, Y_d\right\}$ is $\langle \cdot, \cdot \rangle$-orthonormal, $\langle T,T \rangle=0$ and $\langle T,Z \rangle=1$. Note that also $\langle Z,Z \rangle=0$.

\begin{lemma}\label{lem:decomposition}
Let $X,Y \in \mathfrak{m}$. Then $[X,Y]_{\mathfrak{m}}=[X_{\mathfrak{s}},Y_{\mathfrak{s}}]+[X_{\mathfrak{p}},Y_{\mathfrak{p}}]_{\mathfrak{m}^\prime}$ is a sum of commuting elements.
\end{lemma}
\vspace{-1.4em}\begin{proof}
The claim follows from $[X,Y]=[X_{\mathfrak{s}},Y_{\mathfrak{s}}]+[X_{\mathfrak{p}},Y_{\mathfrak{p}}]$ and $[X_{\mathfrak{p}},Y_{\mathfrak{p}}]\in \mathfrak{c}$.
\end{proof}\vspace{0pt}

For $X \in \mathfrak{c}=\mathds{R}Z \oplus \mathfrak{p} \oplus \mathfrak{h}$, we will write $X_Z$ for the $Z$-component and for an element $X \in \mathfrak{g}=(\mathds{R}T \inplus \mathfrak{he}_d) \oplus \mathfrak{k}\oplus\mathfrak{a}$, we will write $X_T$ for the $T$-component.

\begin{definition}
The map $V: \mathfrak{m} \times \mathfrak{m} \to \mathfrak{m}$ is determined by \[2 \langle V(X,Y) , W\rangle=\langle [W_{\mathfrak{p}},X_{\mathfrak{p}}]_Z, Y_T \rangle+\langle X_T,[W_{\mathfrak{p}},Y_{\mathfrak{p}}]_Z \rangle\] for all $W \in \mathfrak{m}$ and for any $X,Y \in \mathfrak{m}$.
\end{definition}
\begin{remark}
$V$ is symmetric and bilinear. For all $X,Y \in \mathfrak{m}$, $V(X,Y) \in \mathfrak{p}$.
\end{remark}

\begin{definition}
In the situation above, we call the homogeneous space $M$ \textit{special}, if $V \equiv 0$.
\end{definition}

\begin{lemma}\label{lem:special}
In the situation above, the following is true:
\begin{compactenum}
\item $M$ is special if and only if $[\mathfrak{p},\mathfrak{p}]_Z=\left\{0\right\}$.

\item If $M$ is special, $N$ is covered isometrically by the metric product $Z(S) \times N^\prime$ and $M$ is covered isometrically by the metric product $S \times N^\prime$. Here $Z(S)$ is furnished with a homogeneous metric (which is unique up to multiplication with a constant), $S$ is provided with the bi-invariant metric defined in the end of Section~\ref{sec:reductive_representation}, and $N^\prime$ is a homogeneous Riemannian space.
\end{compactenum}
\end{lemma}
\vspace{-1.4em}\begin{proof}
(i) If $[\mathfrak{p},\mathfrak{p}]_Z=\left\{0\right\}$, the statement is obvious.

Conversely, if $V \equiv 0$, \[\langle [W,X]_Z, T \rangle=0\] for all $W,X \in \mathfrak{p}$. Because of $\langle T,Z \rangle =1$, $[\mathfrak{p},\mathfrak{p}]_Z=\left\{0\right\}$ follows.

(ii) As in Section~\ref{sec:orthogonal_distribution}, denote by $\mathcal{O}$ the distribution on $M$ orthogonal to the orbits of $S$ and let $\mathcal{Z}$ be the distribution defined by the action of $Z(S)$. We consider $N$ as a leaf of the foliation defined by the involutive distribution $\mathcal{O}+\mathcal{Z}$ (see Proposition~\ref{prop:integrable_distribution}) through the point $x$, as we did in Section~\ref{sec:structure}.

Let $v,w \in \mathcal{O}_x$. Since $M$ is a homogeneous space, there are $X,Y \in \mathfrak{m}$, such that $\widetilde{X}(x)=v, \widetilde{Y}(x)=w$. Remember that $\widetilde{X},\widetilde{Y}$ are complete Killing vector fields defined by $\widetilde{X}(y):=\frac{\partial}{\partial t} (\exp (t X) \cdot y)|_{t=0}$ for all $y \in M$, $\widetilde{Y}$ is defined in the same way. Because $\mathcal{O}$ is orthogonal to the orbits of $S$ and $\mathfrak{p}$ is orthogonal to $\mathfrak{s}$, $X,Y \in \mathfrak{p}$. We have seen in the proof of Proposition~\ref{prop:N_compact}, that any element of $\textnormal{Isom}^0(M)$ preserves the orbits $\mathcal{S}$ and $\mathcal{O}$. It follows that $\widetilde{X}$ and $\widetilde{Y}$ lie in $\mathcal{O}$.

Let us consider the vector-valued symmetric bilinear form $\omega: \mathcal{O} \times \mathcal{O} \to \mathfrak{s}$ defined in the proof of Proposition~\ref{prop:integrable_distribution}. For any $y \in M$, $\omega_y: \mathcal{O}_y \times \mathcal{O}_y \to \mathfrak{s}$ is defined as \[\omega_x(v,w)=\textnormal{proj}_{\mathcal{S}_y}\left([\widetilde{V},\widetilde{W}]\left(y\right)\right),\] where $\textnormal{proj}_{\mathcal{S}_y}$ denotes the orthogonal projection to $\mathcal{S}_y \cong \mathfrak{s}$ and $\widetilde{V}$ and $\widetilde{W}$ are vector fields in $\mathcal{O}$ extending $v$ and $w$, respectively. We have shown in the proposition that $\omega$ takes only values in $\mathfrak{z(s)}$. But by (i), $[\mathfrak{p},\mathfrak{p}]_Z=\left\{0\right\}$, so $\omega_x$ is trivial. Using that $M$ is homogeneous, it follows that $\omega$ is trivial, hence, $\mathcal{O}$ is involutive.

Now choose $N^\prime$ to be a leaf of the foliation defined by $\mathcal{O}$ through the point $x$. We furnish $N^\prime$ with the induced metric of $N$. Due to Corollary~\ref{cor:isometry_group}~(ii), $\textnormal{Isom}_{Z(S)}N$ acts transitively on $N$. Clearly, $\textnormal{Isom}_{Z(S)}N$ preserves $\mathcal{Z}_N$ by definition, so it also preserves $\mathcal{O}_N$. Thus, if an isometry in $\textnormal{Isom}_{Z(S)}N$ maps $x$ to $y \in N^\prime$, it preserves $N^\prime$. It follows that $N^\prime$ is homogeneous.

Since the action of $S$ is locally free by Corollary~\ref{cor:sl2R_locally_free} and $\mathcal{O}=\mathcal{S}^\perp$, any leaf of $\mathcal{O}$ in a neighborhood of $y \in M$ is given by $f \cdot N_y^\prime$ for some $f \in S$ and $N_y^\prime$ being the leaf of the foliation defined by $\mathcal{O}$ through $y$. Thus, the image of the map $p: S \times N^\prime \to M$, $(f,y) \mapsto f \cdot y$, is open as well as the complement of the image of $p$ in $M$. But $M$ is connected and the image is non-empty, hence, $p$ is surjective.

By construction of the metric on $S \times_{Z(S)} N$, $S \to S \cdot x \subseteq M$ is an isometric covering. It follows that $p$ is a local isometry. Because $S$ is a symmetric space and $N^\prime$ a homogeneous Riemannian space, both manifolds are geodesically complete. Thus, $S \times N^\prime$ is geodesically complete as well. Using that $p$ is surjective and locally isometric, $p$ is a Lorentzian covering map.

If we restrict $p$ to $Z(S) \times N^\prime$, its image is contained in $N$. But the action of $Z(S)$ is locally free, so any leaf of $\mathcal{O}_N$ in a neighborhood of $y \in N$ is given by $f \cdot N_y^\prime$ for some $f \in Z(S)$ and $N_y^\prime$ being the leaf through $y$. Thus, the image of the restriction of $p$ is open as well as the complement of the image is open in $N$. But $N$ is connected and the image is non-empty, hence, the restriction $p|_{Z(S) \times N^\prime}: Z(S) \times N^\prime \to N$ is surjective. The metric on $Z(S)$ is determined by pulling back the metric of an $Z(S)$-orbit of $N$. Because the action of $Z(S)$ is isometric, this metric makes $Z(S)$ into a homogeneous space.
\end{proof}\vspace{0pt}
\begin{remark}
As a converse of Lemma~\ref{lem:special}, if $N$ is covered isometrically by the metric product $Z(S) \times N^\prime$ and the action of $Z(S)$ is defined as the action on the $Z(S)$-component, it is easy to see that $M$ will be special.

Note that $N^\prime$ in Lemma~\ref{lem:special} does not have to be compact. For example, a leaf of the orthogonal distribution of a $\mathds{S}^1$-action on a two-dimensional torus lies dense if the homogeneous metric is chosen appropriately.
\end{remark}

\begin{proposition}\label{prop:U}
Let $X,Y \in \mathfrak{m}$. Then
\[U(X,Y)=U^N(X_{\mathfrak{p}},Y_{\mathfrak{p}})+V(X,Y).\] Here $U^N$ corresponds to the homogeneous space $N$.
\end{proposition}
\vspace{-1.4em}\begin{proof}
Let $W \in \mathfrak{m}$. Since $\mathfrak{he}_d$ is $\langle \cdot,\cdot \rangle$-orthogonal to $\mathfrak{p}$, \[\langle [W_{\mathfrak{s}},X_{\mathfrak{s}}], Y_{\mathfrak{p}} \rangle =0.\] Also, \[\langle [W_{\mathfrak{p}},X_{\mathfrak{p}}]_{\mathfrak{m}^\prime}, Y_{\mathfrak{s}} \rangle = \langle [W_{\mathfrak{p}},X_{\mathfrak{p}}]_Z, Y_T \rangle.\] Because of $(Y_{\mathfrak{p}})_Z=0$, \[\langle [W_{\mathfrak{p}},X_{\mathfrak{p}}]_{\mathfrak{m}^\prime}, Y_{\mathfrak{p}} \rangle=([W_{\mathfrak{p}},X_{\mathfrak{p}}]_{\mathfrak{m}^\prime}, Y_{\mathfrak{p}}).\]
Since $\langle \cdot,\cdot \rangle$ is $\textnormal{ad}(\mathfrak{s})$-invariant, $\langle [W_{\mathfrak{s}},X_{\mathfrak{s}}],Y_{\mathfrak{s}} \rangle+\langle X_{\mathfrak{s}},[W_{\mathfrak{s}},Y_{\mathfrak{s}}] \rangle=0$. Using the definition of $U$, it follows with Lemma~\ref{lem:decomposition} that
\begin{align*}
2 \langle U(X,Y) ,W\rangle&= \langle [W,X]_{\mathfrak{m}},Y \rangle+\langle X,[W,Y]_{\mathfrak{m}} \rangle\\
&=\langle [W_{\mathfrak{s}},X_{\mathfrak{s}}]+[W_{\mathfrak{p}},X_{\mathfrak{p}}]_{\mathfrak{m}^\prime},Y_{\mathfrak{s}}+Y_{\mathfrak{p}} \rangle\\
&+\langle X_{\mathfrak{s}}+X_{\mathfrak{p}},[W_{\mathfrak{s}},Y_{\mathfrak{s}}]+[W_{\mathfrak{p}},Y_{\mathfrak{p}}]_{\mathfrak{m}^\prime} \rangle\\
&=\langle [W_{\mathfrak{s}},X_{\mathfrak{s}}],Y_{\mathfrak{s}} \rangle+\langle X_{\mathfrak{s}},[W_{\mathfrak{s}},Y_{\mathfrak{s}}] \rangle\\
&+\langle [W_{\mathfrak{p}},X_{\mathfrak{p}}]_{\mathfrak{m}^\prime},Y_{\mathfrak{p}} \rangle+\langle X_{\mathfrak{p}}, [W_{\mathfrak{p}},Y_{\mathfrak{p}}]_{\mathfrak{m}^\prime} \rangle\\
&+\langle [W_{\mathfrak{p}},X_{\mathfrak{p}}]_Z, Y_T \rangle+\langle X_T,[W_{\mathfrak{p}},Y_{\mathfrak{p}}]_Z\rangle\\
&=2\langle V(X,Y), W \rangle+( [W_{\mathfrak{p}},X_{\mathfrak{p}}]_{\mathfrak{m}^\prime},Y_{\mathfrak{p}} )+ (X_{\mathfrak{p}}, [W_{\mathfrak{p}},Y_{\mathfrak{p}}]_{\mathfrak{m}^\prime})\\
&=2 ( U^N(X_{\mathfrak{p}},Y_{\mathfrak{p}}),W_{\mathfrak{p}})+2\langle V(X,Y), W \rangle.\end{align*}
Since $(U^N(X_{\mathfrak{p}},Y_{\mathfrak{p}}), Z)=0$, it follows that $U^N(X_{\mathfrak{p}},Y_{\mathfrak{p}})\in \mathfrak{p}$. Using that $V(X,Y)$ in $\mathfrak{p}$ as well, we see that \[2 \langle U(X,Y) ,W\rangle=0\] if $W\in \mathfrak{s}$ and \[\langle U(X,Y) ,W\rangle= ( U^N(X_{\mathfrak{p}},Y_{\mathfrak{p}}),W)+( V(X,Y), W )\] if $W \in \mathfrak{p}$. The result follows.
\end{proof}\vspace{0pt}

\begin{proposition}\label{prop:R}
Let $X,Y \in \mathfrak{m}$. Then
\begin{align*}
R(X,Y,Y,X)&=R^S(X_{\mathfrak{s}},Y_{\mathfrak{s}},Y_{\mathfrak{s}},X_{\mathfrak{s}})+R^N(X_{\mathfrak{p}},Y_{\mathfrak{p}},Y_{\mathfrak{p}},X_{\mathfrak{p}})\\
&-\frac{1}{2}\langle [X_{\mathfrak{p}},[X_{\mathfrak{p}},Y_{\mathfrak{p}}]_{\mathfrak{m}^\prime}]_Z,Y_T \rangle-\frac{1}{2}\langle X_T, [Y_{\mathfrak{p}},[Y_{\mathfrak{p}},X_{\mathfrak{p}}]_{\mathfrak{m}^\prime}]_Z \rangle\\
&+\frac{3}{4}([X_{\mathfrak{p}},Y_{\mathfrak{p}}]_Z,[X_{\mathfrak{p}},Y_{\mathfrak{p}}]_Z)
+2 (V(X,Y),U^N(X_{\mathfrak{p}},Y_{\mathfrak{p}}))\\&+(V(X,Y),V(X,Y))-(V(X,X),V(Y,Y))\\
&-(V(X,X),U^N(Y_{\mathfrak{p}},Y_{\mathfrak{p}}))-(U^N(X_{\mathfrak{p}},X_{\mathfrak{p}}),V(Y,Y)).
\end{align*}
Here $R^S$ and $R^N$ correspond to the homogeneous spaces $S$ and $N$.
\end{proposition}
\vspace{-1.4em}\begin{proof}
We investigate the terms appearing in Proposition~\ref{prop:curvature} separately.

Using Lemma~\ref{lem:decomposition} and the $\langle \cdot,\cdot \rangle$-orthogonality of $\mathfrak{he}_d$ and $\mathfrak{p}$, we have
\begin{align*}
&\langle [X,Y]_{\mathfrak{m}},[X,Y]_{\mathfrak{m}} \rangle\\
=&\langle [X_{\mathfrak{s}},Y_{\mathfrak{s}}],[X_{\mathfrak{s}},Y_{\mathfrak{s}}] \rangle+\langle [X_{\mathfrak{p}},Y_{\mathfrak{p}}]_{\mathfrak{m}^\prime},[X_{\mathfrak{p}},Y_{\mathfrak{p}}]_{\mathfrak{m}^\prime} \rangle\\
=&\langle [X_{\mathfrak{s}},Y_{\mathfrak{s}}],[X_{\mathfrak{s}},Y_{\mathfrak{s}}] \rangle+( [X_{\mathfrak{p}},Y_{\mathfrak{p}}]_{\mathfrak{m}^\prime},[X_{\mathfrak{p}},Y_{\mathfrak{p}}]_{\mathfrak{m}^\prime})-([X_{\mathfrak{p}},Y_{\mathfrak{p}}]_Z,[X_{\mathfrak{p}},Y_{\mathfrak{p}}]_Z),
\end{align*}
where in the end we used that $Z$ is orthogonal to $\mathfrak{p}$ in both metrics, but $\langle \cdot,\cdot \rangle$-isotropic.

Because of $(Y_{\mathfrak{p}})_Z=0$, an easy calculation similar to the ones before yields
\begin{align*}
&\langle [X,[X,Y]_{\mathfrak{m}}]_{\mathfrak{m}},Y \rangle\\
=&
\langle
[X_{\mathfrak{s}},[X_{\mathfrak{s}},Y_{\mathfrak{s}}]],Y_{\mathfrak{s}} \rangle
+
\langle[X_{\mathfrak{p}},[X_{\mathfrak{p}},Y_{\mathfrak{p}}]_{\mathfrak{m}^\prime}]_{\mathfrak{m}^\prime},Y_{\mathfrak{p}}\rangle
+\langle [X_{\mathfrak{p}},[X_{\mathfrak{p}},Y_{\mathfrak{p}}]_{\mathfrak{m}^\prime}]_{\mathfrak{m}^\prime},Y_T \rangle\\
=&
\langle
[X_{\mathfrak{s}},[X_{\mathfrak{s}},Y_{\mathfrak{s}}]],Y_{\mathfrak{s}} \rangle
+
([X_{\mathfrak{p}},[X_{\mathfrak{p}},Y_{\mathfrak{p}}]_{\mathfrak{m}^\prime}]_{\mathfrak{m}^\prime},Y_{\mathfrak{p}})
+\langle [X_{\mathfrak{p}},[X_{\mathfrak{p}},Y_{\mathfrak{p}}]_{\mathfrak{m}^\prime}]_Z,Y_T \rangle.
\end{align*}

Because $\langle \cdot, \cdot \rangle$ is $\textnormal{ad}(\mathfrak{h})$-invariant, \[\langle Y_T, [[X_{\mathfrak{p}},Y_{\mathfrak{p}}]_{\mathfrak{h}},X_{\mathfrak{p}}] \rangle=\langle [Y_T,[X_{\mathfrak{p}},Y_{\mathfrak{p}}]_{\mathfrak{h}}], X_{\mathfrak{p}}\rangle=0.\] Therefore, $\langle Y, [[X,Y]_{\mathfrak{h}},X] \rangle
=( Y_{\mathfrak{m}}, [[X_{\mathfrak{p}},Y_{\mathfrak{p}}]_{\mathfrak{h}},X_{\mathfrak{p}}]_{\mathfrak{m}} ).$

Finally, for the $U$-terms in Proposition~\ref{prop:curvature}, we use Proposition~\ref{prop:U} and the fact, that the image of $V$ lies in $\mathfrak{p}$. The result now follows by summing up the corresponding terms.
\end{proof}\vspace{0pt}

\begin{corollary}\label{cor:he_d_flat}
For $X,Y \in \mathfrak{he}_d \subset \mathfrak{g}$, $R(X,Y,Y,X)=0$.
\end{corollary}
\vspace{-1.4em}\begin{proof}
By Proposition~\ref{prop:R}, \[R(X,Y,Y,X)=R^S(X,Y,Y,X)+(V(X,Y),V(X,Y))-(V(X,X),V(Y,Y)).\] We have seen in Section~\ref{sec:curvature_twisted_Heisenberg}, that $R^S(X,Y,Y,X)=0$. Also, it is obvious that $V(X,Y)=0$ for any $X,Y \in \mathfrak{he}_d \subset \mathfrak{m}$. Thus, $R(X,Y,Y,X)=0$.
\end{proof}\vspace{0pt}

\begin{corollary}\label{cor:curvature_special}
Suppose $M$ is special. Then for $X,Y \in \mathfrak{m}$,
\begin{align*}
R(X,Y,Y,X)&=R^S(X_{\mathfrak{s}},Y_{\mathfrak{s}},Y_{\mathfrak{s}},X_{\mathfrak{s}})+R^N(X_{\mathfrak{p}},Y_{\mathfrak{p}},Y_{\mathfrak{p}},X_{\mathfrak{p}}).
\end{align*}
As above, $R^S$ and $R^N$ correspond to the homogeneous spaces $S$ and $N$.
\end{corollary}
\vspace{-1.4em}\begin{proof}
The claim follows directly from Lemma~\ref{lem:special} and Proposition~\ref{prop:R}.
\end{proof}\vspace{0pt}

\begin{proposition}\label{prop:Ric}
Let $\left\{W_1,\ldots,W_m\right\}$ be any orthonormal basis of $\mathfrak{p}$ and $X\in \mathfrak{m}$. Then
\begin{align*}
\textnormal{Ric}(X,X)&=\textnormal{Ric}^S(X_{\mathfrak{s}},X_{\mathfrak{s}})+\textnormal{Ric}^N(X_{\mathfrak{p}},X_{\mathfrak{p}})-\frac{1}{(Z,Z)} (U^N(X_{\mathfrak{p}},Z),U^N(X_{\mathfrak{p}},Z))\\
&-\frac{1}{2}\sum\limits_{j=1}^m\langle X_T, [W_j,[W_j,X_{\mathfrak{p}}]_{\mathfrak{m}^\prime}]_Z \rangle
+\frac{3}{4}\sum\limits_{j=1}^m([X_{\mathfrak{p}},W_j]_Z,[X_{\mathfrak{p}},W_j]_Z)\\
&+2 \sum\limits_{j=1}^m \langle X_T,[U^N(X_{\mathfrak{p}},W_j),W_j]_Z \rangle
+\frac{1}{4}\sum\limits_{j,k=1}^m (\langle X_T, [W_k,W_j]_Z \rangle)^2\\
&-\sum\limits_{j=1}^m \langle X_T,[U^N(W_j,W_j),X_{\mathfrak{p}}]_Z\rangle
\end{align*}
Here $\textnormal{Ric}^S$ and $\textnormal{Ric}^N$ correspond to the homogeneous spaces $S$ and $N$.
\end{proposition}
\vspace{-1.4em}\begin{proof}
$\left\{\frac{1}{\sqrt{2}}(T-Z),\frac{1}{\sqrt{2}}(T+Z),X_1,Y_1,\ldots,X_d,Y_d, W_1,\ldots, W_m\right\}$ is a basis of $\mathfrak{m}$ being $\langle \cdot,\cdot \rangle$-orthonormal. By definition, \begin{align*}
\textnormal{Ric}(X,X)=&-\frac{1}{2}R(X,T-Z,T-Z,X)+\frac{1}{2}R(X,T+Z,T+Z,X)\\
&+\sum\limits_{j=1}^d(R(X,X_j,X_j,X)+R(X,Y_j,Y_j,X))+\sum\limits_{j=1}^m R(X,W_j,W_j,X).
\end{align*}
We will treat the terms involved in Proposition~\ref{prop:R} separately. For any basis element $Y$, $Y_T=0$ or $Y_{\mathfrak{p}}=0$. Furthermore, $Y_{\mathfrak{p}}\neq 0$ if and only if $Y=W_j$ for some $j$. Also, \[\left\{\frac{Z}{\sqrt{(Z,Z)}},W_1,\ldots,W_m\right\}\] is an orthonormal basis for $\mathfrak{m}^\prime$ with respect to $(\cdot,\cdot)$. So by definition, \[\textnormal{Ric}^N(X_{\mathfrak{p}},X_{\mathfrak{p}})=\frac{1}{(Z,Z)} R^N(X_{\mathfrak{p}},Z,Z,X_{\mathfrak{p}})+\sum\limits_{j=1}^m R^N(X_{\mathfrak{p}},W_j,W_j,X_{\mathfrak{p}}).\] Since $Z$ lies in the center, $U^N(Z,Z)=0$. So using Proposition~\ref{prop:curvature}, it follows that \[R^N(X_{\mathfrak{p}},Z,Z,X_{\mathfrak{p}})=(U^N(X_{\mathfrak{p}},Z),U^N(X_{\mathfrak{p}},Z)).\] This allows us to calculate the corresponding sum of the terms in Proposition~\ref{prop:R} not involving a $V$-expression.

Because of $Y_T=0$ or $Y_{\mathfrak{p}}=0$ for any basis element $Y$, $V(Y,Y)=0$. By definition, \[(V(X,X),W)=\langle [W,X_{\mathfrak{p}}]_{\mathfrak{m}^\prime},X_T \rangle\] for any $W \in \mathfrak{m}^\prime$ and for all $1\leq j \leq m$, \[(V(X,W_j),U^N(X_{\mathfrak{p}},W_j))=\langle X_T,[U^N(X_{\mathfrak{p}},W_j),W_j]_Z \rangle.\]

If $Y$ is a basis element in $\mathfrak{s}$, \[2(V(X,Y),W)=\langle [W_{\mathfrak{p}},X_{\mathfrak{p}}]_Z,Y_T \rangle\] for all $W \in \mathfrak{m}$. Thus, \[V(X,Y)=0 \text{ if }Y \in \mathfrak{he}_d \text{ and }V(X,T-Z)=V(X,T+Z).\]

It is easy to see that \[V(X,W_j)=\frac{1}{2}\sum\limits_{k=1}^d \langle X_T, [W_k,W_j]_Z \rangle W_k.\] Therefore, it holds \[(V(X,W_j),V(X,W_j))=\frac{1}{4}\sum\limits_{k=1}^d (\langle X_T, [W_k,W_j]_Z \rangle)^2\] for all $j$. This finishes the consideration of the terms in Proposition~\ref{prop:R} involving a $V$-expression.
\end{proof}\vspace{0pt}

\begin{corollary}\label{cor:Ric}
Let $\left\{W_1,\ldots,W_m\right\}$ be any orthonormal basis of $\mathfrak{p}$.
\begin{compactenum}
\item If $X,Y \in \mathfrak{he}_d$, $\textnormal{Ric}(X,Y)=0$  and \[\textnormal{Ric}(T+X,T+Y)=\frac{1}{2} \sum\limits_{k=1}^d \lambda_k^2+\frac{1}{4}\sum\limits_{j,k=1}^m (\langle T, [W_k,W_j]_Z \rangle)^2>0.\]

\item The vector-valued Ricci tensor is totally isotropic if and only if \[\textnormal{Ric}^N(X,X)=\frac{1}{(Z,Z)} (U^N(X,Z),U^N(X,Z))
-\frac{3}{4}\sum\limits_{j=1}^m([X,W_j]_Z,[X,W_j]_Z)\] \[ \textnormal{and }
\sum\limits_{j=1}^m \left\{
 [W_j,[W_j,X]_{\mathfrak{m}^\prime}]
-4 [U^N(X,W_j),W_j]
+2[U^N(W_j,W_j),X] \right\}_Z=0
\] for all $X \in \mathfrak{p}$. In this case, we have for all $X \in \mathfrak{m}$, that \[\textnormal{Ric}(X,X)=\textnormal{Ric}^S(X_{\mathfrak{s}},X_{\mathfrak{s}})+\frac{1}{4}\sum\limits_{j,k=1}^m (\langle X_T, [W_k,W_j]_Z \rangle)^2.\]
\end{compactenum}
\end{corollary}
\vspace{-1.4em}\begin{proof}
(i) This directly follows from Proposition~\ref{prop:Ric} and the results of Section~\ref{sec:curvature_twisted_Heisenberg}.

(ii) Assume that the image of the vector-valued Ricci tensor is a totally isotropic subspace of $\mathfrak{m}$. Since $\langle \cdot,\cdot \rangle$ is a Lorentzian scalar product, it follows that the image of $\textnormal{Ric}$ has dimension at most one. By~(i), $\textnormal{Ric}(X,Y)=0$ if $X,Y \in \mathfrak{he}_d$. Therefore, if $X \in \mathfrak{he}_d$, $\textnormal{Ric}(X)$ is $\langle \cdot,\cdot \rangle$-orthogonal to $\mathfrak{he}_d$, that is, \[\textnormal{Ric}(X) \in \mathds{R}Z \oplus \mathfrak{p}.\] Because of $\textnormal{Ric}(T,T)>0$ and due to the totally isotropic image of the Ricci tensor, the image of $\textnormal{Ric}$ is equal to $\mathds{R}Z$. Hence, $\textnormal{Ric}(X,X)=0$ for all $X \in \mathfrak{p}$ (which proves the first equation) and $\textnormal{Ric}(T+X,T+X)=\textnormal{Ric}(X,X)$ for all $X \in \mathfrak{p}$ (which proves the second equation).

Conversely, if both equations are fulfilled, it holds \[\textnormal{Ric}(X,X)=\textnormal{Ric}^S(X_{\mathfrak{s}},X_{\mathfrak{s}})+\frac{1}{4}\sum\limits_{j,k=1}^m (\langle X_T, [W_k,W_j]_Z \rangle)^2\] for all $X \in \mathfrak{m}$ by Proposition~\ref{prop:Ric}. Using the results of Section~\ref{sec:curvature_twisted_Heisenberg}, we see that for any $X \in \mathfrak{m}$, $\text{Ric}(X)=0$ unless $X_T\neq 0$. Thus, the image of the vector-valued Ricci tensor is contained in $\mathds{R}Z$, especially totally isotropic.
\end{proof}\vspace{0pt}

\begin{corollary}\label{cor:Ric_special}
In the situation of this section, the following is true:
\begin{compactenum}
\item $M$ is not Ricci-flat.

\item If $M$ is special, \[\textnormal{Ric}(X,Y)=\textnormal{Ric}^S(X_{\mathfrak{s}},Y_{\mathfrak{s}})+\textnormal{Ric}^N(X_{\mathfrak{p}},Y_{\mathfrak{p}})\] for all $X,Y \in \mathfrak{m}$. In this case, the vector-valued Ricci tensor is totally isotropic if and only if $N$ is Ricci-flat, that is, flat.
\end{compactenum}
\end{corollary}
\vspace{-1.4em}\begin{proof}
The first part follows from Corollary~\ref{cor:Ric}~(i). For the second part, note that $[\mathfrak{p},\mathfrak{p}]_Z=\left\{0\right\}$ due to Lemma~\ref{lem:special}~(i). Hence, by Proposition~\ref{prop:Ric}, \[\textnormal{Ric}(X,X)=\textnormal{Ric}^S(X_{\mathfrak{s}},X_{\mathfrak{s}})+\textnormal{Ric}^N(X_{\mathfrak{p}},X_{\mathfrak{p}})-\frac{1}{(Z,Z)} (U^N(X_{\mathfrak{p}},Z),U^N(X_{\mathfrak{p}},Z)).\]

For all $W \in \mathfrak{m}^\prime$, \[2 (U^N(X_{\mathfrak{p}},Z),W)=([W,X_{\mathfrak{p}}]_{\mathfrak{m}^\prime},Z)+(X_{\mathfrak{p}},[W,Z]_{\mathfrak{m}^\prime})=0,\] because $Z$ lies in the center and $[\mathfrak{p},\mathfrak{p}]_{\mathfrak{m}^\prime}\subseteq\mathfrak{p}$ is orthogonal to $Z$. Thus, \[\textnormal{Ric}(X,Y)=\textnormal{Ric}^S(X_{\mathfrak{s}},X_{\mathfrak{s}})+\textnormal{Ric}^M(X_{\mathfrak{p}},X_{\mathfrak{p}})\] for all $X \in \mathfrak{m}$.

By Corollary~\ref{cor:Ric}~(ii), $\textnormal{Ric}$ is totally isotropic if and only if the two equations given in the corollary are fulfilled. The second one is true because of $[\mathfrak{p},\mathfrak{p}]_Z=\left\{0\right\}$, the first one is now equivalent to $\textnormal{Ric}^N(X,X)=0$ for all $X \in \mathfrak{p}$.

Using $U^N(Z,Z)=0$ and Proposition~\ref{prop:curvature}, it follows that \[R^N(X,Z,Z,X)=(U^N(X,Z),U^N(X,Z))=0\] for all $X \in \mathfrak{m}^\prime$. Thus, $\textnormal{Ric}^N(Z,Z)=0$, which finishes the proof. Note that any homogeneous Riemannian manifold that is Ricci-flat, is also flat.
\end{proof}\vspace{0pt}

\begin{proposition}\label{prop:scal}
Let $\left\{W_1,\ldots,W_m\right\}$ be any orthonormal basis of $\mathfrak{p}$. Then the scalar curvature is given by \[\textnormal{scal}=\textnormal{scal}^N-\frac{2}{(Z,Z)}\sum\limits_{j=1}^m (U^N(W_j,Z),U^N(W_j,Z))+\frac{3}{4}\sum\limits_{j,k=1}^m ([W_j,W_k]_Z,[W_j,W_k]_Z).\]
Here $\textnormal{scal}^N$ denotes the scalar curvature of the homogeneous space $N$.
\end{proposition}
\vspace{-1.4em}\begin{proof}
For the calculation of the scalar curvature, we choose the $\langle \cdot,\cdot \rangle$-orthonormal basis $\left\{\frac{1}{\sqrt{2}}(T-Z),\frac{1}{\sqrt{2}}(T+Z),X_1,Y_1,\ldots,X_d,Y_d, W_1,\ldots, W_m\right\}$.

Because of Corollary~\ref{cor:Ric}~(i), \[\textnormal{Ric}(X,X)=0 \text{ if }X \in \mathfrak{he}_d \text{ and }\textnormal{Ric}(T+Z,T+Z)=\textnormal{Ric}(T-Z,T-Z).\] Using Proposition~\ref{prop:Ric}, \begin{align*}
\textnormal{scal}=&\sum\limits_{j=1}^m \textnormal{Ric}^N(W_j,W_j)-\frac{1}{(Z,Z)}\sum\limits_{j=1}^m (U^N(W_j,Z),U^N(W_j,Z))\\
&+\frac{3}{4}\sum\limits_{j,k=1}^m([W_j,W_k]_Z,[W_j,W_k]_Z)\\
=&\textnormal{scal}^N-\frac{1}{(Z,Z)}\textnormal{Ric}^N(Z,Z)-\frac{1}{(Z,Z)}\sum\limits_{j=1}^m (U^N(W_j,Z),U^N(W_j,Z))\\
&+\frac{3}{4}\sum\limits_{j,k=1}^m([W_j,W_k]_Z,[W_j,W_k]_Z).
\end{align*}
Now \[\textnormal{Ric}^N(Z,Z)=\sum\limits_{j=1}^m R^N(Z,W_j,W_j,Z)=\sum\limits_{j=1}^m R^N(W_j,Z,Z,W_j).\] Using the fact that $Z$ lies in the center, $U(Z,Z)=0$, and by Proposition~\ref{prop:curvature}, \[R^N(X,Z,Z,X)=(U^N(X,Z),U^N(X,Z))\] for $X \in \mathfrak{p}$. The claim follows.
\end{proof}\vspace{0pt}

\begin{corollary}\label{cor:scal_special}
If $M$ is special, $\textnormal{scal}=\textnormal{scal}^N$.
\end{corollary}
\vspace{-1.4em}\begin{proof}
Due to Lemma~\ref{lem:special}~(i), $[\mathfrak{p},\mathfrak{p}]_Z=\left\{0\right\}$. So for all $W \in \mathfrak{m}^\prime$, \[2 (U^N(X_{\mathfrak{p}},Z),W)=([W,X_{\mathfrak{p}}]_{\mathfrak{m}^\prime},Z)+(X_{\mathfrak{p}},[W,Z]_{\mathfrak{m}^\prime})=0,\] because $Z$ lies in the center and $[\mathfrak{p},\mathfrak{p}]_{\mathfrak{m}^\prime}\subseteq \mathfrak{p}$ is orthogonal to $Z$. The claim now follows immediately from Proposition~\ref{prop:scal}.
\end{proof}\vspace{0pt}

We conclude with determining the holonomy algebra of $M$ in the special case, that $M$ is special.

\begin{proposition}\label{prop:holonomy_special}
Let $M$ be special and $x \in N \subset M$ as in Section~\ref{sec:reductive_representation}. Then $\mathfrak{hol}_x(N)\subseteq\mathfrak{so}(\mathfrak{p}, ( \cdot,\cdot))$, so there is a natural injection $\mathfrak{hol}_x(N) \to \mathfrak{so}(\mathfrak{m}, \langle \cdot,\cdot \rangle)$. Also, $\textnormal{ad}({\mathfrak{he}_d})\subseteq \mathfrak{so}(\mathfrak{m}, \langle \cdot,\cdot \rangle)$. The holonomy algebra of $M$ is then determined by
\[\mathfrak{hol}_x(M)=\textnormal{ad}({\mathfrak{he}_d})\oplus\mathfrak{hol}_x(N)\subseteq \mathfrak{so}(\mathfrak{m}, \langle \cdot,\cdot \rangle).\]
\end{proposition}
\vspace{-1.4em}\begin{proof}
We know already, that $\langle \cdot,\cdot \rangle$ restricted to $\mathfrak{s} \times \mathfrak{s}$ is ad($\mathfrak{s}$)-invariant. But $\mathfrak{p}$ commutes with $\mathfrak{s}$ and is orthogonal to $\mathfrak{s}$, thus, $\langle \cdot,\cdot \rangle$ is ad($\mathfrak{s}$)-invariant. Especially, $\textnormal{ad}({\mathfrak{he}_d})\subseteq \mathfrak{so}(\mathfrak{m}, \langle \cdot,\cdot \rangle)$.

By Proposition~\ref{prop:holonomy_algebra}, \[\mathfrak{hol}_{x}(N)=\mathfrak{m}^\prime_0+[\Lambda_{\mathfrak{m}^\prime}(\mathfrak{m}^\prime),\mathfrak{m}^\prime_0]+[\Lambda_{\mathfrak{m}^\prime}(\mathfrak{m}^\prime),[\Lambda_{\mathfrak{m}^\prime}(\mathfrak{m}^\prime),\mathfrak{m}^\prime_0]]+\ldots,\]
where $\mathfrak{m}^\prime_0$ is the subspace in $\mathfrak{so}(\mathfrak{m}^\prime, ( \cdot, \cdot ))$ spanned by \[\left\{[\Lambda_{\mathfrak{m}^\prime}(X),\Lambda_{\mathfrak{m}^\prime}(Y)]-\Lambda_{\mathfrak{m}^\prime}([X,Y]_{\mathfrak{m}^\prime})-\textnormal{ad}_{[X,Y]_{\mathfrak{h}}} \ | \ X,Y \in \mathfrak{m}^\prime\right\}.\]

Now \[\Lambda_{\mathfrak{m}^\prime}(X)Y=\frac{1}{2}[X,Y]_{\mathfrak{m}^\prime}+U^N(X,Y)\] for all $X,Y \in \mathfrak{m}^\prime$. Using \[2(U^N(X,Y),Z)=([Z,X]_{\mathfrak{m}^\prime},Y)+(X,[Z,Y]_{\mathfrak{m}^\prime})=0,\] we have $U^N(X,Y) \in \mathfrak{p}$, and it follows from the first part of Lemma~\ref{lem:special}, that $[\mathfrak{m}^\prime,\mathfrak{m}^\prime]_{\mathfrak{m}^\prime} \subseteq \mathfrak{p}$. Thus, $\Lambda_{\mathfrak{m}^\prime}(X)Y \in \mathfrak{p}$ for all $X,Y \in \mathfrak{m}^\prime$.

Since $Z$ lies in the center and $U^N(X,Z)=0$ for all $X \in \mathfrak{m}^\prime$ as shown in the proof of Corollary~\ref{cor:scal_special}, $\Lambda_{\mathfrak{m}^\prime}(X)Z=0$. Hence, $\Lambda_{\mathfrak{m}^\prime}(\mathfrak{m}^\prime)\subseteq\mathfrak{so}(\mathfrak{p}, ( \cdot,\cdot))$.

For any $X,Y,W \in \mathfrak{m}^\prime$, \[([[X,Y]_{\mathfrak{h}},W],Z)=-(W,[[X,Y]_{\mathfrak{h}},Z])=0\] because $(\cdot,\cdot)$ is ad($\mathfrak{h}$)-invariant. Therefore, $\textnormal{ad}_{[X,Y]_{\mathfrak{h}}}$, which annihilates $Z$, is in $\mathfrak{so}(\mathfrak{p}, ( \cdot,\cdot))$ as well. We conclude that $\mathfrak{hol}_x(N)\subseteq\mathfrak{so}(\mathfrak{p}, ( \cdot,\cdot))$, so there is an injection $\mathfrak{hol}_x(N) \to \mathfrak{so}(\mathfrak{m}, \langle \cdot,\cdot \rangle)$ using the canonical injection of $\mathfrak{so}(\mathfrak{p}, ( \cdot,\cdot))$ into $\mathfrak{so}(\mathfrak{m}, \langle \cdot,\cdot \rangle)$.

For all $X,Y \in \mathfrak{m}$, \[\Lambda_{\mathfrak{m}}(X)Y=\frac{1}{2}[X,Y]_{\mathfrak{m}}+U(X,Y)=\Lambda_{\mathfrak{s}}(X_{\mathfrak{s}})Y_{\mathfrak{s}}+\Lambda_{\mathfrak{m}^\prime}(X_{\mathfrak{p}})Y_{\mathfrak{p}}\] by Lemma~\ref{lem:decomposition} and Proposition~\ref{prop:U} and because of $V \equiv 0$. Additionally, we have $[X,Y]_{\mathfrak{h}}=[X_{\mathfrak{p}},Y_{\mathfrak{p}}]_{\mathfrak{h}}$.

We have seen in Section~\ref{sec:curvature_twisted_Heisenberg}, that the holonomy algebra of $S$ is given by ad($\mathfrak{he}_d$). The required result now follows from Proposition~\ref{prop:holonomy_algebra}.
\end{proof}\vspace{0pt}

\subsection{Isotropy representation and Ricci-flat manifolds}\label{sec:corollaries}

Propositions~\ref{prop:sl2R_isotropy} and~\ref{prop:heis_isotropy} show Theorem~\ref{th:isotropy_representation}.

Also the proof of Theorem~\ref{th:homogeneous_not_Ricci_flat} is now immediate: Due to the results of Section~\ref{sec:homogeneous_sl2R}, a compact homogeneous Lorentzian manifold $M$ is not Ricci-flat, if the Lie algebra of its isometry group contains a direct summand isomorphic to $\mathfrak{sl}_2(\mathds{R})$. By Corollary~\ref{cor:Ric}, the same is true in the case of a direct summand isomorphic to $\mathfrak{he}_d^\lambda$, $\lambda \in \mathds{Z}_+^d$. So by Theorem~\ref{th:homogeneous_characterization}, $\textnormal{Isom}^0(M)$ is compact.

\newpage
%%%%%%%%%%%%%%%%%%%%%%%%%%%%%%%%%%%%%%%%%%%%%%%%%%%%%%%%%%%%%%%%%%%%%%%%%%%%%%%%%%%%%%%%%%%%

\begin{appendix}

%%%%%%%%%%%%%%%%%%%%%%%%%%%%%%%%%%%%%%%%%%%%%%%%%%%%%%%%%%%%%%%%%%%%%%%%%%%%%%%%%%%%%%%%%%%%

\pagestyle{plain}
\fancyhead[ER]{Appendix}
\fancyhead[OL]{Complete Jordan decomposition}
\chapter{Complete Jordan decomposition}

Let $V$ be a real vector space and $\mathfrak{g} \subset \mathfrak{gl}(V)$ a semisimple subalgebra. In the following, we will introduce a Jordan-type decomposition of elements of $\mathfrak{g}$:

\begin{proposition}\label{prop:Jordan}
For any $X \in \mathfrak{g}$, there are commuting elements $E,H,N \in \mathfrak{g}$, such that $X=E+H+N$, $E$ and $H$ are semisimple with only purely imaginary and real eigenvalues, respectively, and $N$ is nilpotent. We call this the \textit{complete additive Jordan decomposition}.

Moreover, $X$ is nilpotent or semisimple with only purely imaginary or real eigenvalues, respectively, if and only if $\textnormal{ad}_X$ is nilpotent or semisimple with only purely imaginary or real eigenvalues, respectively.
\end{proposition}

\begin{proof}
We follow \cite{Hel92}, Chapter~IX, Paragraph~7. According to Lemma~7.1, there is a unique decomposition \[\exp(X)=ehu\] with commuting elements $e,h,u$ in $\textnormal{GL}(V)$ and $e$ is semisimple with all its (complex) eigenvalues having modulus one, $h$ is semisimple with only positive eigenvalues and $(u-\textnormal{id}_V)$ is nilpotent. This decomposition is called the \textit{complete multiplicative Jordan decomposition}. The Further Result~A.6 at the end of Chapter~IX states that $e,h,u$ are contained in the adjoint group $\textnormal{Int}(\mathfrak{g})$.

Using the usual additive Jordan decomposition of $X$, we can write $X$ uniquely as a sum \[X=E^\prime+H^\prime+N^\prime\] of commuting elements $E^\prime,H^\prime,N^\prime \in \mathfrak{gl}(V)$, where $N^\prime$ is nilpotent, $E^\prime$ and $H^\prime$ are semisimple and have purely imaginary or real eigenvalues only, respectively. (If $J \in \textnormal{GL}(V)$ such that $JXJ^{-1}$ is in Jordan canonical form, take $JN^\prime J^{-1}$ to be the strictly upper triangular matrix, $JE^\prime J^{-1}$ and $JH^\prime J^{-1}$ to be the imaginary and real part of the diagonal, respectively.)

Then $e^\prime:=\exp(E^\prime),h^\prime:=\exp(H^\prime),u^\prime:=\exp(N^\prime)$ commute and \[\exp(X)=e^\prime h^\prime u^\prime.\] Clearly, $e^\prime$ and $h^\prime$ are semisimple and all its eigenvalues have modulus one or are positive, respectively, and $(u^\prime-\textnormal{id}_V)$ is nilpotent. By uniqueness of the multiplicative Jordan decomposition, \[e=e^\prime, \ h=h^\prime, \ u=u^\prime.\]

Lemma~7.3 and the subsequent remark assert that $e,h,u$ are contained in one-parameter groups of $\textnormal{Int}(\mathfrak{g})$, so \[e=\exp(E),\ h=\exp(H),\ u=\exp(N)\] with $E,H,N \in \mathfrak{g}$. From the proof of Lemma~7.3 it is clear that \[H=H^\prime.\] Since $u=u^\prime$, $\exp(N-N^\prime)=\textnormal{id}_V$. But $N-N^\prime$ is nilpotent, therefore \[N=N^\prime.\] Thus, also \[E=E^\prime.\]

In summary, we have shown that for any $X \in \mathfrak{g}$, the elements of its complete additive Jordan decomposition are in $\mathfrak{g}$ as well.

It is clear that $\textnormal{ad}_X$ is nilpotent if $X \in \mathfrak{g}$ is nilpotent. The following is an adaptation of the proof of Lemma~2 in \cite{Bo07}, Paragraph~5.4. Let $X \in \mathfrak{g}$ be semisimple, such that its eigenvalues $\lambda_1, \ldots, \lambda_j$ are all purely imaginary or real, respectively. Let $v_1,\ldots,v_j$ be corresponding eigenvalues in $V^\mathds{C}$, the complexification of $V$. Consider the basis $\left\{M_{kl}\right\}_{k,l=1,\ldots,j}$ of $\mathfrak{gl}(V)$ defined by \[M_{kl}v_m=\delta_{lm}v_k\] for all $k,l,m$. Then \[\textnormal{ad}_X(M_{kl})v_m=XM_{kl}v_m-M_{kl}Xv_m=\delta_{lm}(\lambda_k-\lambda_m)v_k.\] Therefore, \[\textnormal{ad}_X(M_{kl})=(\lambda_k-\lambda_l)M_{kl},\] so $\textnormal{ad}_X$ is semisimple with only purely imaginary or real eigenvalues, respectively. So if \[X=E+H+N\] is the complete additive Jordan decomposition of $X \in \mathfrak{g}$, \[\textnormal{ad}_X=\textnormal{ad}_E+\textnormal{ad}_H+\textnormal{ad}_N\] is the one of $\textnormal{ad}_X$. Since $\mathfrak{g}$ is semisimple, the adjoint representation of $\mathfrak{g}$ is injective. Thus, $X \in \mathfrak{g}$ is nilpotent or semisimple with only purely imaginary or real eigenvalues, respectively, if and only if $\textnormal{ad}_X$ is nilpotent or semisimple with only purely imaginary or real eigenvalues, respectively.
\end{proof}
\newpage
%%%%%%%%%%%%%%%%%%%%%%%%%%%%%%%%%%%%%%%%%%%%%%%%%%%%%%%%%%%%%%%%%%%%%%%%%%%%%%%%%%%%%%%%%%%%

\pagestyle{plain}
\fancyhead[ER]{References}
\fancyhead[OL]{References}
\renewcommand{\bibname}{References}

\newpage

\end{appendix}

\end{document}